\documentclass[a4paper,11pt,reqno]{amsart}
\usepackage[margin=0.8in]{geometry}
\usepackage[colorlinks]{hyperref}

\usepackage{amsthm,amsmath,amssymb,dsfont}
\usepackage{mathrsfs,amsfonts,functan,extarrows,mathtools}

\newtheorem{theorem}{Theorem}[section]
\newtheorem{proposition}{Proposition}[section]
\newtheorem{definition}{Definition}[section]

\newtheorem{lemma}{Lemma}[section]

\numberwithin{equation}{section}
\allowdisplaybreaks
\def\p{\partial}
\def\d{\mathop{}\!\mathrm{d}}
\def\no{\nonumber}
\def\eps{\varepsilon}
\def\div{\mathrm{div}}

\def\A{\mathcal{A}}
\def\B{\mathcal{B}}
\def\C{\mathcal{C}}
\def\P{\mathcal{P}}

\def\R{\mathbb{R}}

\def\w{\mathfrak{w}}
\def\L{\mathcal{L}}
\def\Q{\mathcal{Q}}
\def\T{{\mathbb{T}^3}}
\def\ker{\mathrm{ker}}
\def\IN{\mathrm{in}}

\begin{document}
\title[Hydrodynamic Limit of Vlasov-Maxwell-Boltzmann] 
			{Hydrodynamic Limit of the Incompressible Navier-Stokes-Fourier-Maxwell System with Ohm's Law \\from the Vlasov-Maxwell-Boltzmann System: \\Hilbert Expansion Approach}

\author[N. Jiang]{Ning Jiang}
\address[Ning Jiang]
    {\newline\indent School of Mathematics and Statistics, Wuhan University, Wuhan, 430072, P. R. China}
\email{njiang@whu.edu.cn}

\author[Y.-L. Luo]{Yi-Long Luo}
\address[Yi-Long Luo]
    {\newline\indent School of Mathematics, South China University of Technology, Guangzhou, 510641, P. R. China}
\email{luoylmath@scut.edu.cn}

\author[T.-F. Zhang]{Teng-Fei Zhang$^\dag$}
\address[Teng-Fei Zhang]
    {\newline\indent School of Mathematics and Physics, China University of Geosciences, Wuhan, 430074, P. R. China}
\email{zhangtf@cug.edu.cn}

\thanks{$^\dag$ Corresponding author.}
\date{} 

\begin{abstract}
We prove a global-in-time limit from the two-species Vlasov-Maxwell-Boltzmann system to the two-fluid incompressible Navier-Stokes-Fourier-Maxwell system with Ohm's law. Besides the techniques developed for the classical solutions to the Vlasov-Maxwell-Boltzmann equations in the past years, such as the nonlinear energy method and micro-macro decomposition are employed, key roles are played by the decay properties of both the electric field and the wave equation with linear damping of the divergence free magnetic field. This is a companion paper of [N. Jiang and Y.-L. Luo, \emph{Ann. PDE} 8 (2022), no. 1, Paper No. 4, 126 pp] in which Hilbert expansion is not employed.

\smallskip\noindent\textit{\textsc{Keywords}}. Hydrodynamic limits; Hilbert expansion; Two-species Vlasov-Maxwell-Boltzmann system; Navier-Stokes-Fourier-Maxwell system; Ohm's law

\smallskip\noindent\textit{\textsc{2020 Mathematics Subject Classification}}. 35Q20, 35Q30, 35Q61, 76P05, 76W05
\end{abstract}

\maketitle





\section{Introduction}
\subsection{The two-species Vlasov-Maxwell-Boltzmann system and its hydrodynamic limits}

The two-species Vlasov-Maxwell-Boltzmann (VMB) system describes the evolution of a gas of two species of oppositely charged particles (cations of charge $q^+ > 0$ and mass $m^+>0$, and anions of charge $-q^- <0$ and mass $m^->0$), subject to auto-induced electromagnetic forces.
Such a gas of charged particles equipped with a global neutrality condition is usually called a plasma. The unknowns $F^+(t,x,v) \geq 0 $ and $F^-(t,x,v) \geq 0 $ represent respectively the particle number densities of the positive charged ions (i.e., cations), and the negative charged ions (i.e., anions), which are at position $x\in \mathbb{T}^3$ with velocity $v\in \mathbb{R}^3$, at time $t\geq 0$. The VMB system reads as follows:
	\begin{align}\label{VMB-0}\tag{VMB}
	\left\{\begin{array}{rcl}
	  \partial_t F^+ + v\cdot \nabla_x F^+ + \tfrac{q^+}{m^+}(E + v \times B)\cdot\nabla_v F^+ & = &  Q(F^+, F^+) + Q(F^+, F^-)\,,
	  \\[5pt]
	  	\partial_t F^- + v\cdot \nabla_x F^- - \tfrac{q^-}{m^-}(E + v \times B)\cdot\nabla_v F^- & = & Q(F^-, F^-) + Q(F^-, F^+)\,,
	  \\[5pt] \displaystyle
	  	\mu_0\eps_0\partial_t E - \nabla_x \times B & = & -\mu_0\int_{\mathbb{R}^3}(q^+F^+- q^- F^-)v\,\mathrm{d}v\,,
	  \\[5pt]
	    \partial_t B + \nabla_{\!x} \times E & = & 0\,,
	  \\[5pt] \displaystyle
	  	\div_x E & = & \tfrac{1}{\eps_0}\int_{\mathbb{R}^3}(q^+F^+- q^- F^-)\,\mathrm{d}v\,,
	  \\[5pt]
	  	\div_x B & = & 0\,.
	\end{array}\right.
	\end{align}
Here the evolution of the densities $F^\pm$ are governed by the Vlasov-Boltzmann equations in the first two equations of \eqref{VMB-0}. This means that the variations of densities along the particle trajectories are subject to the influence of an auto-induced Lorentz force and inter-particle collisions in the gas. The electric field $E(t,x)$ and the magnetic field $B(t,x)$, which are generated by the motion of the particles in the plasma itself, are governed by the Maxwell equation. It consists of the Amp\`ere equation, Faraday's equation and Gauss' laws, representing in the third, fourth, and the last two equations, respectively.
The vacuum permeability and permittivity (or say, the magnetic and electric constants) are denoted, respectively, by the physical coefficients $\mu_0, \eps_0 >0$. Both species of particles are assumed here to have the same mass $m^\pm=m>0$ and charge $q^\pm = q >0$.

In the VMB system \eqref{VMB-0}, the collision between particles is given by the standard Boltzmann collision operator $Q(f,h)$, with respect to two species of density functions $f(v)$, $h(v)$. The bilinear operator $Q(f,h)(v)$ is defined as
\begin{equation}\label{Clisn-Oprt}
 Q(f,h) = \int_{\R^3} \int_{\mathbb{S}^2} ( f' h'_* - f h_* ) b ( v - v_* , \sigma ) \d \sigma \d v_* \,,
\end{equation}
where the parameter $\sigma \in \mathbb{S}^2$, and the cross-section function is taken the form of
\begin{equation}
  b( v - v_* , \sigma ) = |v - v_*|^\gamma \hat{b} (\cos \theta).
\end{equation}
For convenience, we consider the hard potential case $\gamma \in [0, 1]$, and take the cutoff factor $\hat{b} (\cos \theta)$ such that $ \int_{\mathbb{S}^2} \hat{b} (\cos \theta) \d \sigma = 1 $.
Note that we have used here the standard abbreviations
$$ f = f(v)\,, \ f' = f (v') \,, \ h_* = h (v_*) \,, \ h'_* = h (v'_*) $$
with the relation between pre- and post-collisional velocities that
\begin{equation}\label{v-v-collision}
  \begin{aligned}
    v' = v + [ ( v - v_* ) \cdot \sigma ] \sigma \,, \quad
    v'_* = v_* - [ ( v - v_* ) \cdot \sigma ] \sigma \,,
  \end{aligned}
\end{equation}
which is derived from the conservation of momentum and energy during the collision process.

In this paper, our goal is to investigate the hydrodynamic limit of the two-species VMB system, in the regime where Knudsen number is small, is exactly the two-fluid incompressible Navier-Stokes-Fourier-Maxwell (in brief, NSFM) system with Ohm's law.
More precisely, by considering the non-dimensional variables in the Vlasov-Boltzmann equations in \eqref{VMB-0}, and by a change of units in electric and magnetic fields $E$ and $B$, we are led to the scaled two-species Vlasov-Maxwell-Boltzmann system, as follows:
\begin{equation}\label{VMB-f} \tag*{(VMB)$_\eps$}
	\left\{
	\begin{array}{rcl}
		\eps \partial_t F^{\pm}_\eps + v \cdot \nabla_x F^{\pm}_\eps \pm ( \eps E_\eps + v \times B_\eps ) \cdot \nabla_v F^{\pm}_\eps & = & \frac{1}{\eps} Q ( F^{\pm}_\eps , F^{\pm}_\eps  ) + \frac{1}{\eps} Q( F^{\pm}_\eps ,F^{\mp}_\eps) \,,
	\\[3pt]
		F^{\pm}_\eps & = & M ( 1 + \eps g^{\pm}_\eps ) \,,
	\\[5pt] \displaystyle
		\partial_t E_\eps - \nabla_x \times B_\eps & = & - \frac{1}{\eps} \int_{\R^3} (g^+_\eps - g^-_\eps) v M \d v \,,
		\\[5pt]
				\div_x E_\eps & = & \int_{\R^3} (g^+_\eps - g^-_\eps) M \d v \,,
	\\[5pt]
		\partial_t B_\eps + \nabla_x \times E_\eps & = & 0 \,,
		\\[3pt]
				\div_x B_\eps & = & 0 \,.
	\end{array}
	\right.
\end{equation}
Note that, the dimensionless Knudsen number and Mach number are taken the same order $\eps>0$.
The second equation represents the fluctuations $g^{\pm}_\eps$ is around the global Maxwellian distribution $M (v) = \tfrac{1}{\sqrt{2 \pi}^3} \exp{\left(-\tfrac{|v|^2}{2}\right)}$.

We will prove the two-fluid incompressible NSFM equations with Ohm's law can be obtained, in the above hydrodynamic regime, from the scaled VMB system \ref{VMB-f}. The limit system are rewritten as:
	\begin{equation}\label{NSFM} \tag{NSFM}
 		\left\{
		\begin{array}{rclc}
			\partial_t u + u \cdot \nabla_x u - \mu \Delta_x u + \nabla_x p & = & \tfrac{1}{2} n E + \tfrac{1}{2} j \times B \,, & \quad \div_x u = 0 \,,  \\[2pt]
			\partial_t \theta + u \cdot \nabla_x \theta - \kappa \Delta_x \theta & = & 0 \,, & \\[2pt]
			\partial_t E - \nabla_x \times B & = & - j \,, & \quad \div_x E = n \,, \\[2pt]
			\partial_t B + \nabla_x \times E & = & 0 \,, & \quad \div_x B = 0 \,, \\[2pt]
			j = n u +  \sigma ( - \tfrac{1}{2} & \nabla_x n & + E + u \times B ) \,, &
		\end{array}
		\right.
	\end{equation}
where $\mu$, $\kappa$ and $\sigma$ are positive constants to be determined, see Section \ref{Sec:Formal-Analysis}.

\subsection{Historical remarks}

There have been extensive researches on the well-posedness of the VMB system. DiPerna-Lions developed a theory of global-in-time renormalized solutions with large initial data to the Boltzmann equation \cite{DL-Annals1989}, Vlasov-Maxwell equations \cite{DL-CPAM1989} and Vlasov-Poisson-Boltzmann equation \cite{Lions-Kyoto1994, Lions-Kyoto1994-2}. However, there are still severe difficulties in establishing the global well-posedness theory for VMB. The major one is that the a priori bounds provided by physical laws are not enough to prove the existence of global solutions, even in the renormalized solutions sense. Recently, Ars\`enio and Saint-Raymond \cite{Arsenio-SaintRaymond, Arsenio-SRM} established global-in-time renormalized solutions with large initial data for VMB, with both cut-off and non-cutoff collision kernels. We emphasize that, as far as we know, the renormalized solutions are still the only existing theory for solutions without any smallness requirements on initial data. On the other line, in the context of classical solutions, through a so-called nonlinear energy method, Guo \cite{Guo-Inventions2003} constructed a classical solution of VMB near the global Maxwellian. Guo's work inspired many results on VMB with more general collision kernels, among which we only mention here the results related to the most general collision kernels with or without angular cutoff assumptions, see \cite{DLYZ-KRM2013, DLYZ-CMP2017, FLLZ-2018} for example.

Hydrodynamic limits from kinetic equations have been an active research field for decades. Among many research results in this field, the most successful program is the so-called BGL program (named after Bardos-Golse-Levermore \cite{BGL-CPAM1993}), which aimed at establishing the asymptotic limit between DiPerna-Lions' renormalized solutions of the Boltzmann equation and the Leray solutions of incompressible Navier-Stokes equations. The BGL program was completed by Golse and Saint-Raymond \cite{G-SRM-Invent2004, G-SRM2009}, and the readers are referred to see \cite{Masmoudi-SRM-CPAM2003, JM-CPAM2017} for the case in a domain with boundary.

Concerning the VMB system, the corresponding hydrodynamic limits are much harder, even at the formal level, since it is coupled with Maxwell equations which are essentially hyperbolic. In a recent remarkable breakthrough \cite{Arsenio-SaintRaymond}, Ars\'enio and Saint-Raymond proved, in one hand, the existence of renormalized solutions of VMB, as mentioned above. More importantly, they also justified, in the other hand, various limits (depending on different scalings) towards incompressible viscous electro-magneto-hydrodynamics. Among these limits, the most singular one is from renormalized solutions of two-species VMB to dissipative solutions of the two-fluid incompressible NSFM system \eqref{NSFM} with Ohm's law. The proofs in \cite{Arsenio-SaintRaymond}, justifying the weak limit from a sequence of solutions of VMB \ref{VMB-f} to a dissipative solution of incompressible NSFM \eqref{NSFM}, are extremely hard. Part of the reasons are, besides many difficulties of the existence of renormalized solutions of VMB itself, our current understanding for the incompressible NSFM with Ohm's law is far from complete. From the viewpoint of mathematical analysis, NSFM have a behavior which is more similar to the much less understood incompressible Euler equations than to the Navier-Stokes equations. That is the reason the authors in \cite{Arsenio-SaintRaymond} considered the so-called dissipative solutions of NSFM rather than the usual weak solutions. The dissipative solutions were introduced by Lions for 3-dimensional incompressible Euler equations (see \cite[\S 4.4]{Lions-1996}). In the content of classical solutions, the first two authors of considered in a recent work \cite{JL-2022annPDE} the hydrodynamic limit for the two-species VMB to two-fluid incompressible NSFM with Ohm’s law. They proved the uniform estimates with respect to Knudsen number $\eps$ for the fluctuations. Another previous hydrodynamic limit result for the VMB in the framework of classical solutions belongs to Jang \cite{Jang-ARMA2008}. In \cite{Jang-ARMA2008}, it was taken a special scaling that the magnetic effect appears only as a higher order. As a consequence, it vanishes in the limit as the Knudsen number $\eps\rightarrow 0$. So in the limiting equations derived in \cite{Jang-ARMA2008}, there is no equation for the magnetic field.

The studies of incompressible NSFM just started in recent years (for the introduction of physical background, see \cite{Biskamp, Davidson}).  For weak solutions, the existence of global in time Leray type weak solutions are completely open, even in 2-dimensional case.  A first breakthrough comes from Masmoudi \cite{Masmoudi-JMPA2010}, who proved the existence and uniqueness of global strong solutions of incompressible NSFM in 2-dimensional case, for the initial data $(u^{\IN}, E^{\IN}, B^{\IN})\in L^2(\mathbb{R}^2)\times (H^s(\mathbb{R}^2))^2$ with $s>0$. In fact, the system in \cite{Masmoudi-JMPA2010} is a little different from the NSFM in this paper, but the spirit of analytic analysis still remains same. Note that in \cite{Masmoudi-JMPA2010}, the divergence-free condition of the magnetic field $B$ or the decay property of the linear part arising from Maxwell's equations is \emph{not} used.  Ibrahim and Keraani \cite{Ibrahim-Keraani-2011-SIMA} considered the data $(u^{\IN}, E^{\IN}, B^{\IN}) \in \dot B^{1/2}_{2,1}(\mathbb{R}^3)\times (\dot H^{1/2}(\mathbb{R}^3))^2$ for 3-dimension, and $(u^{\IN}, E^{\IN}, B^{\IN})\in \dot B^0_{2,1}(\mathbb{R}^2)\times (L^2_{\log}(\mathbb{R}^2))^2$ for 2-dimensional case. Later on, German, Ibrahim and Masmoudi \cite{GIM2014} refined the previous results by applying a fixed-point argument to obtain mild solutions, but with the initial velocity field in the natural Navier-Stokes space $H^{1/2}$. In their results the regularity of the initial velocity and electromagnetic fields is lowered. Furthermore, they employed an $L^2L^\infty$-estimate on the velocity field, which significantly simplifies the fixed-point arguments used in \cite{Ibrahim-Keraani-2011-SIMA}. For some other related asymptotic problems, say, the derivation of the MHD from the Navier-Stokes-Maxwell system in the context of weak solutions, see Ars\'enio-Ibrahim-Masmoudi \cite{AIM-ARMA-2015}. Recently, in \cite{JL-CMS-2018} the first two authors of the current paper proved the global classical solutions of the incompressible NSFM with small initial data, by employing the decay properties of both the electric field and  the wave equation with linear damping of the divergence free magnetic field. This key idea was already used in \cite{GIM2014}. Thus, it is natural to construct a series of classical solutions to VMB around the solutions of the incompressible NSFM with Ohm's law established in \cite{JL-CMS-2018}. This is the main concern of the current paper. We note that for the Boltzmann equation, Guo \cite{Guo-CPAM-2006} proved the incompressible Navier-Stokes limit from the Boltzmann equation in classical solutions framework. Our result in this paper is in the same spirit as that in \cite{Guo-CPAM-2006}, but for much more involved incompressible NSFM limit from two-species VMB.

\subsection{Hilbert expansion and well-prepared initial data}

In this paper, we will investigate the issues in the classical regimes. Working in the classical solutions has some mathematical convenience, comparing to renormalized solutions. We have relatively much better understandings for both VMB and incompressible NSFM, at least for classical solutions near equilibriums. Furthermore, we can employ some properties of the incompressible NSFM specific for classical solutions, for example, the decay properties of both the electric field, and the wave equation with linear damping of the divergence free magnetic field. In fact, the analog of these properties in VMB will play an essential role in the current work.

First, the initial data of the scaled VMB system \ref{VMB-f} are imposed as
\begin{equation}\label{IC-1}
	( F^{\pm}_\eps (0, x, v)\,, \ E_\eps (0,x) \,, \ B_\eps (0,x)  ) = ( F^{\pm, \IN}_\eps(x,v)\,, \ E^{\IN}_\eps(x) \,, \ B^{\IN}_\eps(x)  ) 
\end{equation}
with the compatibility condition $\div_x B^{\IN}_\eps = 0$.

It is well-known that the classical solutions to the system \ref{VMB-f} obey the following global conservation laws of mass, momentum and  energy:
	\begin{align}\label{Totl-Consv-Law}
    \begin{cases}
      \frac{\d}{\d t} \iint_{\mathbb{T}^3 \times \mathbb{R}^3} F^{\pm}_\eps \d v \d x = 0 \,,
    \\[5pt] 
      \frac{\d }{ \d t } \left( \iint_{\mathbb{T}^3 \times \mathbb{R}^3} v ( F^+_\eps + F^-_\eps ) \d v \d x + \eps \int_{\mathbb{T}^3} E_\eps \times B_\eps \d x \right) = 0 \,,
    \\[5pt] 
      \frac{\d }{ \d t } \left( \iint_{\mathbb{T}^3 \times \mathbb{R}^3} |v|^2 ( F^+_\eps + F^-_\eps ) \d v \d x + \eps^2 \int_{\mathbb{T}^3} | E_\eps |^2 + | B_\eps |^2 \d x \right) = 0 \,.
    \end{cases}
	\end{align}
Moreover, there holds
	\begin{equation}\label{Conservtn-B}
	  \frac{\d}{\d t} \int_{\T} B_\eps (t,x) \d x = 0 \,.
	\end{equation}
Assuming the initial data \eqref{IC-1}  has the same mass, total momentum and energy as the global equilibrium $(M(v)\,, 0 \,, 0)$, the conservation laws \eqref{Totl-Consv-Law}-\eqref{Conservtn-B} can be rewritten as
\begin{equation}\label{Totl-Consv-Law-1}
    \left\{\begin{array}{l} \displaystyle
       \iint_{\mathbb{T}^3 \times \mathbb{R}^3} F^{\pm}_\eps \d v \d x = \iint_{\mathbb{T}^3 \times \mathbb{R}^3} M \d v \d x \,,
    \\[7pt] \displaystyle
       \iint_{\mathbb{T}^3 \times \mathbb{R}^3} v ( F^+_\eps + F^-_\eps ) \d v \d x + \eps \int_{\mathbb{T}^3} E_\eps \times B_\eps \d x = 0 \,,
    \\[7pt] \displaystyle
      \iint_{\mathbb{T}^3 \times \mathbb{R}^3} |v|^2 ( F^+_\eps + F^-_\eps ) \d v \d x + \eps^2 \int_{\mathbb{T}^3} | E_\eps |^2 + | B_\eps |^2 \d x  = 2 \iint_{\mathbb{T}^3 \times \mathbb{R}^3} |v|^2 M \d v \d x \,,
    \\[7pt] \displaystyle
      \int_{\T} B_\eps(t,x) \d x = 0 \,.
    \end{array}\right.
\end{equation}
In other words, the relation \eqref{Totl-Consv-Law-1} holds for $( F^{\pm}_\eps, E_\eps, B_\eps )$ at any positive time provided that it holds initially.

The goal of this paper is to justify the hydrodynamic limits from two-species VMB system \ref{VMB-f} to NSFM system \eqref{NSFM} with Ohm's law, by employing the Hilbert expansion approach. As shown in Section \ref{Sec:Formal-Analysis}, we consider the solutions to \ref{VMB-f} of the form:
\begin{equation}\label{Soltn-Fom}
    \left\{\begin{array}{rcl} \displaystyle
      F^{\pm}_\eps (t,x,v) &=&  M (v) \left\{ 1  + \eps \left[ g_0^{\pm} (t,x,v) + \eps {g}_1^{\pm} (t,x,v) + \eps^2 {g}_2^{\pm} (t,x,v) + \eps g^{\pm}_{R, \eps} (t,x,v) \right] \right\} \,,
    \\[3pt]
      E_\eps (t,x) &=& E_0 (t,x) + \eps {E}_1 (t,x) + \eps E_{R, \eps} (t,x) \,,
    \\[3pt]
      B_\eps (t,x) &=& B_0 (t,x) + \eps {B}_1 (t,x) + \eps B_{R, \eps} (t,x) \,.
    \end{array}\right.
\end{equation}
Here the leading order terms $g_0^{\pm}$ are the infinitesimal Maxwellians, given by
\begin{equation}\label{g0-0}
  g_0^{\pm} (t,x,v) = \rho_0^{\pm} (t,x) + u_0 (t,x) \cdot v + \theta_0 (t,x) ( \tfrac{|v|^2}{2} - \tfrac{3}{2} ) \,.
\end{equation}
Denote $\rho_0(t,x) = \tfrac{1}{2} ( \rho_0^+ (t,x) + \rho_0^- (t,x) ) $ and $n_0 (t,x) = \rho_0^+ (t,x) - \rho_0^- (t,x) = \div_x E_0 $ with $\rho_0 + \theta_0 = 0$, then the functions $(u_0 , \theta_0, E_0, B_0 )$ 
 solve the NSFM system \eqref{NSFM}. The initial data
	\begin{equation}\label{IN-NSFM}
	  ( u_0 , \theta_0 , E_0 , B_0 ) (0,x)  = ( u_0^{\IN} , \theta_0^{\IN} , E_0^{\IN} , B_0^{\IN} ) (x)
	\end{equation}
is assumed to satisfy the compatibility conditions $\div_x (u_0^{\IN},B_0^{\IN}) = 0$. We further assume
\begin{equation}\label{IC-av-B0}
	\begin{aligned}
		\int_{\T} B_0^{\IN} \d x = 0 \,.
	\end{aligned}
\end{equation}
Lemma \ref{lemm:bnd-NSMF} below ensures that the solutions $( u_0 , \theta_0 , E_0 , B_0 )$ and $n_0$ to system \eqref{IN-NSFM} are all globally well-defined, if the initial data are taken in appropriate Sobolev spaces.

The following order terms ${g}_i^{\pm} (t,x,v)$ ($i=1,2$), as expressed in \eqref{g1-g2}, depends on some macroscopic quantities $(\rho_1, u_1, \theta_1, n_1, E_1, B_1)$, which can be determined by \eqref{Od:eps(10)} and \eqref{Spec-Chos*} in Section \ref{Sec:Formal-Analysis} below. In fact, we assume that the initial data
\begin{equation}\label{IC-LM}
  ( {E}_1 (0, x) , {B}_1 (0, x) ) = ( {E}^{\IN}_1 (x) , {B}^{\IN}_1 (x) ) 
\end{equation}
satisfies the compatibility condition $\div_x {B}^{\IN}_1 = 0$.
In addition, the magnetic field obeys $\tfrac{\d}{\d t} \int_{\T} {B}_1 \d x = 0$, due to the formulation of its equation in \eqref{Od:eps(10)}. This ensures $ \int_{\T} {B}_1 \d x = 0 $, provided that it holds initially:
	\begin{equation}\label{IC-4}
	  \int_{\T} {B}^{\IN}_1 \d x = 0 \,.
	\end{equation}
Lemma \ref{lemm:bnd-linearMaxwl}, equipped with the above initial conditions, 
will show that the functions $ {E}_1 (t,x)$, $ {B}_1 (t,x)$ and $ {n}_1(t,x)$ are globally and uniquely determined, once some functions $(\rho_1, u_1, \theta_1)$ are well-chosen to satisfy \eqref{Spec-Chos*}. As a result, the quantities $ {g}_i^{\pm} (t,x,v)$ for $i = 1 , 2$ in \eqref{g1-g2} are already known.

By the formal analysis in Section \ref{Sec:Formal-Analysis}, the remainder terms $(g_{R,\eps}^\pm, E_{R,\eps} , B_{R,\eps})$ are determined through the system
\begin{equation}\label{Remd-Equ-Bref}
  \left\{
    \begin{array}{rll}
      \eps \partial_t G_{R,\eps} + v \cdot \nabla_x G_{R,\eps} + \mathcal{T} (v \times B_0) \cdot \nabla_v G_{R,\eps} + & \mathcal{T} (v \times B_{R,\eps}) \cdot \nabla_v G_0 &
    \\[2mm]
      - E_{R,\eps} \cdot v \mathcal{T}_1 + \frac{1}{\eps} \mathbb{L} G_{R,\eps} = & \eps H_{R,\eps} \,,
    \\[2mm]
      \partial_t E_{R,\eps} - \nabla_x \times B_{R,\eps} = & - \frac{1}{\eps} \left\langle G_{R,\eps} \cdot \mathcal{T}_1 v \right\rangle \,, \
      	& \div_x E_{R,\eps} = \left\langle G_{R,\eps} \cdot \mathcal{T}_1 \right\rangle \,,
    \\[2mm]
      \partial_t B_{R,\eps} + \nabla_x \times E_{R,\eps} = & 0 \,, \ & \div_x B_{R,\eps} = 0 \,,
    \end{array}
  \right.
\end{equation}
with $\mathcal{T}_1 = (1, -1)^\top$ and $\mathcal{T} = diag(1,-1)$. Here $G_{R,\eps} = (g_{R,\eps}^+, g_{R,\eps}^-)^\top$, $ G_0 = (g_0^+ ,g_0^-)^\top $. The term $H_{R,\eps}$ will be specified later, see \eqref{Remnd-H} in Section \ref{Sec:Formal-Analysis}.

The initial data of the remainder system \eqref{Remd-Equ-Bref} is imposed on
	\begin{equation}\label{IC-Remd}
	  ( G_{R,\eps} (0 , x,v) , E_{R, \eps} (0,x) , B_{R,\eps} (0,x) ) = ( G_{R, \eps}^{\IN}(x,v) , E_{R,\eps}^{\IN}(x) , B_{R,\eps}^{\IN} (x) ),
	\end{equation}
satisfying the compatibility $\div_x B_{R, \eps}^{\IN} = 0$, where $ G_{R , \eps}^{\IN} (x,v) = (g_{R,\eps}^{+, \IN} (x,v),g_{R,\eps}^{-, \IN} (x,v))^\top$.

Therefore, the initial data \eqref{IC-1} of \ref{VMB-f} are given through the following well-prepared formulations:
	\begin{equation}\label{IC-3}
	  \left\{
	  \begin{array}{rl}
	    F_\eps^{\pm, \IN} (x,v) =& M (v) \left\{ 1 + \eps \left[ g_0^{\pm , \IN} (x,v) + \eps {g}_1^{\pm, \IN} (x,v) + \eps^2 {g}_2^{\pm , \IN} (x,v) + \eps g^{\pm , \IN}_{R,\eps} (x,v) \right] \right\} \,, \\[2mm]
	    E_\eps^{\IN} (x) =& E_0^{\IN} (x) + \eps {E}_1^{\IN} (x) + \eps E_{R,\eps}^{\IN} (x) \,, \\[2mm]
	    B_\eps^{\IN} (x) =& B_0^{\IN} (x) + \eps {B}_1^{\IN} (x) + \eps B_{R,\eps}^{\IN} (x) \,.
	  \end{array}\right.
	\end{equation}
Here
\begin{equation*}
  g_0^{\pm , \IN} (x,v) = \pm \tfrac{1}{2} \div_x E_0^{\IN} (x) - \theta_0^{\IN} (x) + u_0^{\IN} (x) \cdot v + \theta_0^{\IN} (x) ( \tfrac{|v|^2}{2} - \tfrac{3}{2} ) \,,
\end{equation*}
and ${g}_k^{\pm, \IN} (x,v)$ ($k = 1,2$) have the same forms as ${g}_k^{\pm} (t,x,v)$ defined in \eqref{g1-g2}, just replacing the vectors $ (u_0 , \theta_0 , E_0 , B_0) $ and $({E}_1 , {B}_1)$ by $( u_0^{\IN} , \theta_0^{\IN} , E_0^{\IN} , B_0^{\IN} )$ and $({E}^{\IN}_1 , {B}^{\IN}_1)$, respectively.
As is mentioned before, the initial data $(F_\eps^{\pm, \IN}, E_\eps^{\IN}, B_\eps^{\IN})$ in \eqref{IC-3} are supposed to obey the relations \eqref{Totl-Consv-Law-1}, and $B_0^{\IN} (x), B_1^{\IN} (x)$ satisfy, respectively, \eqref{IC-av-B0} and \eqref{IC-4}.

\subsection{Notations and main results}

In order to state clearly our results, some notations should be introduced firstly. $A \lesssim B$ means $A \leq C B$ for some harmless constant $C> 0$. $A \sim B$ stands for $C_1 B \leq A \leq C_2 B$ for some constants $C_1, C_2 > 0$. We denote by the symbol $\langle \, \cdot \, \rangle$ the integral of a function $f(v)$ in $\mathbb{R}^3$ with respect to the probability measure $M\!\d v$, i.e.,
$$ \langle f \rangle : = \int_{\mathbb{R}^3} f (v) M \d v \,. $$
Let the multi-indices $m$ and $\beta$ be $ m = (m_1 , m_2, m_3 ) \,, \ \beta = (\beta_1 , \beta_2 , \beta_3 )$, and let
$$ \partial^m_\beta : = \partial^m \partial_\beta = \partial^{m_1}_{x_1} \partial^{m_2}_{x_2} \partial^{m_3}_{x_3} \partial^{\beta_1}_{v_1} \partial^{\beta_2}_{v_2} \partial^{\beta_3}_{v_3} \,,$$
where $m$ is related to the space derivatives, and $\beta$ corresponds to the velocity derivatives. If each component of $m$ is not greater than that of $\tilde{m}$'s, we denote by $ m \leq \tilde{m} $. The symbol $ m < \tilde{m} $ means $ m \leq \tilde{m} $ and $ |m| < | \tilde{m} | $, where $|m| = m_1 + m_2 + m_3$.

In what follows, the capital symbol $G = (g^+, g^-)^\top$ denotes a column vector in $\mathbb{\R}^2$. For example, we will use the notations $G_0 = \left( \begin{array}{c}
		g^+_0 \\
		g^-_0
	\end{array} \right) $, and $G_{R,\eps} = \left( \begin{array}{c}
		g_{R,\eps}^+ \\
		g_{R,\eps}^-
	\end{array} \right)$.

Now we introduce some basic properties of the linearized operators, which can also be found in \cite{Arsenio-SaintRaymond}. For the Boltzmann collision operator \eqref{Clisn-Oprt}, the collision frequency is defined as
\begin{equation}\label{Clisn-Frqc}
  \nu (v) := \int_{\mathbb{R}^3} \int_{\mathbb{S}^2} |v - v_*|^\gamma \hat{b} ( \cos \theta ) \d \sigma \d v_* = \int_{\mathbb{R}^3} |v - v_*|^\gamma  \d v_* \,.
\end{equation}
We define the weight $\w(v)$ by $ \w(v) =  (1+|v|^2)^{1/2} $. Moreover, the linearized collision operator $\mathbb{L}$ acting on a two-species vector $G = (g^+, g^-)^\top$ is defined through
\begin{equation}\label{Def-L-bb}
  \mathbb{L} G =
  \left(
    \begin{array}{c}
      \mathcal{L} g^+ + \mathcal{L}^\flat (g^+ , g^-) \\
      \mathcal{L} g^- + \mathcal{L}^\flat (g^+ , g^-)
    \end{array}
  \right) \,,
\end{equation}
where we have used the linearized Boltzmann collision operator
\begin{equation}\label{Def-L-cal}
	\begin{aligned}
		\mathcal{L} g = - \tfrac{1}{M} \big[ Q ( M g, M ) + Q ( M , M g ) \big] \,,
	\end{aligned}
\end{equation}
and the linearized inter-species collision operator
\begin{equation}\label{Def-L-flat}
	\begin{aligned}
		\mathcal{L}^\flat (g , h) = - \tfrac{1}{M} \big[ Q ( M g, M ) + Q ( M , M h ) \big]  \,.
	\end{aligned}
\end{equation}
For simplicity, we denote
\begin{equation}\label{Q-cal}
	\begin{aligned}
		\mathcal{Q} (g , h) = \tfrac{1}{M} Q ( M g , M h ) \,.
	\end{aligned}
\end{equation}

The kernel of $\mathbb{L}$ (denoted by $\ker \, \mathbb{L}$) is spanned by, see \cite{Arsenio-SaintRaymond},
 \begin{equation}\label{L-ker}
   \left( \begin{array}{c} 1 \\ 0 \end{array} \right) \,, \ \left( \begin{array}{c} 0 \\ 1 \end{array} \right) \,, \ \left( \begin{array}{c} v \\ v \end{array} \right) \,, \left( \begin{array}{c} \tfrac{|v|^2}{2} - \tfrac{3}{2} \\ \tfrac{|v|^2}{2} - \tfrac{3}{2} \end{array} \right) \,.
 \end{equation}
The projection operator $\mathbb{P}$ from $L^2_v$ to  $\ker \mathbb{L}$ is defined as
\begin{equation}\label{Projc-P}
  \mathbb{P} G = \rho^+ \left( \begin{array}{c} 1 \\ 0 \end{array} \right) + \rho^- \left( \begin{array}{c} 0 \\ 1 \end{array} \right) + u \cdot \left( \begin{array}{c} v \\ v \end{array} \right) + \theta \left( \begin{array}{c} \tfrac{|v|^2}{2} - \tfrac{3}{2} \\ \tfrac{|v|^2}{2} - \tfrac{3}{2} \end{array} \right) \,,
\end{equation}
where $\rho^\pm = \langle g^\pm \rangle$, $u = \left\langle v \tfrac{g^+ + g^-}{2} \right\rangle$ and $\theta = \left\langle ( \tfrac{|v|^2}{3} - 1 ) \tfrac{g^+ + g^-}{2} \right\rangle $. It is easy to infer, from \cite{BGL-CPAM1993}, that
\begin{equation}
	\begin{aligned}
		\ker \L = \mathrm{span} \{ 1, v, \tfrac{|v|^2}{2} - \tfrac{3}{2} \},
	\end{aligned}
\end{equation}
and that, the projection $\P_{\L} : L^2(M\!\d v) \rightarrow \ker \L$ is
\begin{equation}\label{Projc-PL}
	\begin{aligned}
		\P_{\L} g = \langle g \rangle + \langle g v \rangle \cdot v + \langle g ( \tfrac{|v|^2}{3} - 1 ) \rangle ( \tfrac{|v|^2}{2} - \tfrac{3}{2} ) \,.
	\end{aligned}
\end{equation}
Furthermore, the symbol $\mathbb{I}$ denotes the identity operator, and $\mathbb{P}^\perp = \mathbb{I} - \mathbb{P}$, $\P_{\L}^\perp = Id - \P_{\L}$. The linear operator $\mathfrak{L}$ is further defined as
\begin{equation}\label{Def-L-frk}
	\begin{aligned}
		\mathfrak{L} g = - \Q ( g, 1 ) - \Q ( 1 , - g ) \,,
	\end{aligned}
\end{equation}
which satisfies $\ker \mathfrak{L} = \mathrm{span} \{ 1 \}$, see \cite{Arsenio-SaintRaymond}.

Let $L^p_x$ and $H^s_x$ be the standard $L^p$ and Sobolev spaces on $x \in \T$, respectively. The norms of $\| \cdot \|_{L^p_x}$ and $\| \cdot \|_{H^s_x}$ are
\begin{equation*}
	\begin{aligned}
		\| f \|_{L^p_x} := \Big( \int_{\mathbb{T}^3} |f (x)|^p \d x \Big)^{\tfrac{1}{p}} \ (1 \leq p < \infty) \,, \ \| f \|_{L^\infty_x} = \sup_{x \in \T} |f (x)| \,, \ \| f \|^2_{H^s_x} = \sum_{|m| \leq s} \| \partial^m f \|^2_{L^2_x} \,.
	\end{aligned}
\end{equation*}
For a function $a(v) > 0$, the weighted Hilbert space $L^2_v(a)$ is endowed with the norm
$$\| g \|_{L^2_v(a)} : = \Big( \int_{\mathbb{R}^3} |g(v)|^2 a(v) M \d v \Big)^{\tfrac{1}{2}} \,.$$
For the function $h(x,v)$, one can define a norm
\begin{equation*}
  \| h \|_{L^2_{x,v} (a)} := \Big( \int_{\mathbb{T}^3} \int_{\mathbb{R}^3} |h(x,v)|^2 a(v) \d v \d x \Big)^{\tfrac{1}{2}} \,.
\end{equation*}
Furthermore, the spaces $H^s_x L^2_v (a)$ and $H^s_{x,v}(a)$ are defined by
$$ \| h \|^2_{H^s_x L^2_v (a)} : = \sum_{|m| \leq s} \| \partial^m h \|^2_{L^2_{x,v} (a)} \,, \ \| h \|^2_{H^s_{x,v}(a)} : = \sum_{|m|+|\beta| \leq s} \| \partial^m_\beta h \|^2_{L^2_{x,v} (a)} \,, $$
respectively. Moreover, the space $\widetilde{\mathbf{H}}^s_{x,v} (a)$ is endowed with the norm
\begin{equation*}
	\begin{aligned}
		\| h \|^2_{\widetilde{\mathbf{H}}^s_{x,v} (a)} : = \sum_{\substack{|m| + |\beta| \leq s , \beta \neq 0}} \| \p^m_\beta h \|^2_{L^2_{x,v} (a)} \,.
	\end{aligned}
\end{equation*}
In the case $a (v) \equiv 1$, we denote briefly $X = X (1)$, where the space $X$ could be $L^2_v$, $L^2_{x,v}$, $H^s_x L^2_v$, $ H^s_{x,v}$ or $\widetilde{\mathbf{H}}^s_{x,v}$.

Next, we define the following energy and the dissipation rate functionals.
\begin{definition}[Energy and dissipation rate]
  For $N \geq 4$, $l \geq 0$ and $\w (v) = \sqrt{1 + |v|^2}$, we define the energy functional $ \mathbb{E}_{N,l} (G , E, B)$ as
  \begin{equation}\label{Energ-1}
    \begin{aligned}
      \mathbb{E}_{N,l} (G , E, B) = & \| E \|^2_{H^{N+1}_x} + \| B \|^2_{H^{N+1}_x} + \| G \|^2_{H^{N+1}_x L^2_v} + \| \mathbb{P}^\perp G \|^2_{H^N_{x,v} (\w^{2l})} \\
      & + \| \mathbb{P}^\perp G \|^2_{\widetilde{\mathbf{H}}^{N+1}_{x,v} (\w^{2l})} + \mathcal{E}_{0, N+5} (t) + \mathcal{E}_{1, N+3} (t) \,,
    \end{aligned}
  \end{equation}
and the dissipation rate functional $\mathbb{D}_{N,l} (G,E,B)$ as
  \begin{equation}\label{Dspt-Rate}
    \begin{aligned}
      \mathbb{D}_{N,l} (G,E,B) = & \| E \|^2_{H^{N-1}_x} + \| \nabla_x B \|^2_{H^{N-1}_x} + \| \partial_t B \|^2_{H^{N-1}_x} + \mathcal{D}_{0,N+5} (t) + \mathcal{D}_{1,N+3}(t) \\
      & + \| G \|^2_{H^{N+1}_x L^2_v} + \tfrac{1}{\eps^2} \| \mathbb{P}^\perp G \|^2_{\widetilde{\mathbf{H}}^{N+1}_{x,v} (\w^{2l} \nu)} \\
      & + \tfrac{1}{\eps^2} \| \mathbb{P}^\perp G \|^2_{H^{N+1}_x L^2_v(\nu)} + \tfrac{1}{\eps^2} \| \mathbb{P}^\perp G \|^2_{H^N_{x,v}(\w^{2l} \nu)}
    \end{aligned}
  \end{equation}
where the functionals $\mathcal{E}_{0,N+5}(t)$, $\mathcal{D}_{0,N+5} (t)$, $\mathcal{E}_{1,N+3} (t)$ and $\mathcal{D}_{1,N+3} (t)$ are given in \eqref{E0s}-\eqref{D0s}, and \eqref{E1M}-\eqref{D1M}, respectively.
\end{definition}


As a correspondence, we employ the initial energy functional for the remainder term:
\begin{equation*}
  \begin{aligned}
   {\mathbb{E}}_{N,l} ( G_{R,\eps}^{\IN} , E_{R,\eps}^{\IN} , B_{R,\eps}^{\IN} ) = & \| E_{R,\eps}^{\IN} \|^2_{H^{N+1}_x} + \| B_{R,\eps}^{\IN} \|^2_{H^{N+1}_x} + \| G_{R,\eps}^{\IN} \|^2_{H^{N+1}_x L^2_v} + \| \mathbb{P}^\perp G_{R,\eps}^{\IN} \|^2_{H^N_{x,v} (\w^{2l})} \\
   & + \| \mathbb{P}^\perp G_{R,\eps}^{\IN} \|^2_{\widetilde{\mathbf{H}}^{N+1}_{x,v} (\w^{2l} )} + \mathcal{E}_{0, N+5}^{\IN} + \mathcal{E}_{1,N+3}^{\IN} \,,
  \end{aligned}
\end{equation*}
where the quantities $\mathcal{E}_{0, N+5}^{\IN}$ and $\mathcal{E}_{1,N+3}^{\IN}$ are given in \eqref{E0s-in} and \eqref{E1M-in}, respectively, by taking $s=N+5$, $M=N+3$.

We now are in a position to state our main result, as follows:
\begin{theorem}\label{Main-Thm}
  Let $N \geq 4$ and $l \geq 2 \gamma + 1$. Assume that
	  \begin{enumerate}
	  	\item $u_0^{\IN}(x)$, $\theta_0^{\IN}(x)$, $E_0^{\IN}(x)$ and $B_0^{\IN}(x)$ given in \eqref{IN-NSFM}-\eqref{IC-av-B0} satisfy the initial hypotheses in Lemma \ref{lemm:bnd-NSMF} with $s = N+5$;

	  	\item ${E}_1^{\IN}(x)$ and ${B}_1^{\IN}(x)$ given in \eqref{IC-LM}-\eqref{IC-4} satisfy the initial assumptions in Lemma \ref{lemm:bnd-linearMaxwl} with $M = N+3$;

	  	\item $E_{R,\eps}^{\IN}(x)$, $B_{R,\eps}^{\IN}(x)$ and $G_{R,\eps}^{\IN}(x)$ given in \eqref{IC-Remd} satisfy ${\mathbb{E}}_{N,l} ( G_{R,\eps}^{\IN} , E_{R,\eps}^{\IN} , B_{R,\eps}^{\IN} ) < \infty$;

	  	\item The initial data $( F^{\pm, \IN}_\eps(x,v)\,, \ E^{\IN}_\eps(x) \,, \ B^{\IN}_\eps(x)  )$ of the scaled VMB system \ref{VMB-f} are imposed on the well-prepared forms \eqref{IC-3}.
	  \end{enumerate}

  There exist two sufficiently small constants $\eps_0$, $\eta_0 > 0$, depending only on $\mu$, $\sigma$, $\kappa$, $l$ and $N$, such that if
  \begin{equation}\label{Thm-IC-Small}
    {\mathbb{E}}_{N,l} ( G_{R,\eps}^{\IN} , E_{R,\eps}^{\IN}, B_{R,\eps}^{\IN} ) \leq \eta_0 \,, \quad \forall \eps \in (0 , \eps_0) \,,
  \end{equation}
  then the system \ref{VMB-f} with the initial data \eqref{IC-3} admits a global-in-time classical solution $ ( F_\eps^{\pm} (t,x,v)$, $E_\eps(t,x) , B_\eps(t,x) ) $, belonging to $L^\infty( \mathbb{R}^+ ; H^{N+1}_{x,v} \times H^{N+1}_x \times H^{N+1}_x )$, and being of the form \eqref{Soltn-Fom}, i.e.,
  \begin{equation*}
  	\left\{\begin{array}{rl}
  		F^{\pm}_\eps (t,x,v) =&  M (v) \left\{ 1  + \eps \left[ g_0^{\pm} (t,x,v) + \eps {g}_1^{\pm} (t,x,v) + \eps^2 {g}_2^{\pm} (t,x,v) + \eps g^{\pm}_{R, \eps} (t,x,v) \right] \right\} \,,\\ [2mm]
  		E_\eps (t,x) =& E_0 (t,x) + \eps {E}_1 (t,x) + \eps E_{R, \eps} (t,x) \,, \\[2mm]
  		B_\eps (t,x) =& B_0 (t,x) + \eps {B}_1 (t,x) + \eps B_{R, \eps} (t,x) \,,
  	\end{array}\right.
  \end{equation*}
  where $g_0^\pm$ in \eqref{g0-0} depend on $(\rho_0, n_0, u_0, \theta_0)$, and $g_k^\pm$ $(k = 1,2)$ in \eqref{g1-g2} depend on $(u_0, \theta_0, E_0)$ and $(u_1, \theta_1, E_1)$. Here $(u_0, \theta_0, E_0, B_0)$ subjecting to the two fluid incompressible NSFM equations \eqref{NSFM} with initial data \eqref{IN-NSFM}-\eqref{IC-av-B0} are solved in Lemma \ref{lemm:bnd-NSMF}. $(u_1, \theta_1, E_1, B_1)$ obeying the linear Maxwell-type equations \eqref{Od:eps(10)}-\eqref{Spec-Chos*} with initial data \eqref{IC-LM}-\eqref{IC-4} are solved in Lemma \ref{lemm:bnd-linearMaxwl}.

	Moreover, there exists an instant energy functional $\mathbb{E}_{N,l} (R_{R,\eps}, E_{R, \eps}, B_{R,\eps})$, such that
\begin{equation}
	\begin{aligned}
      \sup_{t \geq 0} \mathbb{E}_{N,l} (G_{R,\eps}, E_{R, \eps}, B_{R,\eps}) (t) + \int_0^\infty \mathbb{D}_{N,l} (G_{R,\eps}, & E_{R, \eps}, B_{R,\eps})(t) \d t \\
      & \leq C {\mathbb{E}}_{N,l} (G^{\IN}_{R,\eps}, E^{\IN}_{R, \eps}, B^{\IN}_{R,\eps})
    \end{aligned}
\end{equation}
holds for some constant $C > 0$, depending only on $\mu$, $\sigma$, $\kappa$, $l$ and $N$.
\end{theorem}

\subsection{Main ideas and novelties}

The hydrodynamic limits from the scaled two-species VMB system \ref{VMB-f} to the incompressible two-fluid NSFM equations \eqref{NSFM} with Ohm's law is a singular limit problem. Usually, there are two frameworks to deal with the limit: Moment method and Hilbert expansion approach. The former one is to obtain a uniform bound (with respect to the small parameter $\eps$) on the solutions to the original scaled singular equations, and then to extract a convergent subsequence converging to the solutions of the target (limit) equations as $\eps \to 0$. The later one is, based on the obtained solutions for the limiting equations, to construct a sequence of special solutions of the original scaled singular equations for small parameter $\eps$ around the limit equations. We remark that the first method is usually much harder than the second approach, since the uniform bounds are very difficult to be established in most situations. The readers are referred to \cite{JL-2022annPDE}, in which the uniform convergence result belongs to a result of the first method, as a companion of this paper.

In the current paper, the Hilbert expansion approach will be adopted to help us justify the limit. The most advantage of this way lies in the observation the remainder equations have weaker singularity than the original system, although the remainder equations remains still nonlinear and singular. To be more precise, after utilizing the Hilbert expansion, the nonlinearity and singularity are separated. So, it will be relatively easier to get an energy bound, uniformly in small $\eps > 0$, of the remainder system.

Our first goal is to derive the limit NSFM equations \eqref{NSFM} with Ohm's law, by a formal derivations. The point is to apply the order analysis to the local conservation formulas, which is derived from projecting the VMB equations \ref{VMB-f} into the fluid part $\ker \mathbb{L}$. Deriving the remainder equations is the next goal of the formal derivations. In the Hilbert expansion approach, while the limit equations can be totally determined, we also hope to get an expansion formulation with less expanded terms, and its corresponding remainder system admitting certain uniform bounds. Thus, we hope to find an optimal truncation on the expansion formulation. Our truncated formulas are exactly \eqref{Soltn-Fom}, in which $(g_0^\pm, E_0, B_0 )$ is completely determined by the limit NSFM equations \eqref{NSFM}, and $g_k^\pm (k=1,2)$, $E_1, B_1$ are the auxiliary functions to formulate the NSFM equations \eqref{NSFM}.
More precisely, as expressed in \eqref{g1-g2} in \S \ref{sub:truncation_expansion}, the first order term $\tfrac{{g}_1^+ + {g}_1^-}{2}$ contains two parts: one part belongs to the $\ker^\perp \L$, which can be uniquely determined, while the other one belonging to $\ker \L$ is with undetermined functions $(\rho_1, u_1, \theta_1)$. The term $\tfrac{{g}_1^+ - {g}_1^-}{2}$ contains a similar structures with completely determined $\ker^\perp(\L + \mathfrak{L})$ part and undetermined $\ker(\L + \mathfrak{L})$ part (indicating by the function $n_1$). We choose here some special functions $(\rho_1, u_1, \theta_1)$ and $(E_1, B_1, n_1)$ to satisfy the relations \eqref{Od:eps(10)}-\eqref{Spec-Chos*}. On the other hand, the expression \eqref{g1-g2} also indicates that, the second order terms $\frac{g_2^+ + g_2^-}{2}$ and $\frac{g_2^+ - g_2^-}{2}$ are with undetermined $(\rho_2, u_2, \theta_2, n_2)$ (and $(\rho_1, u_1, \theta_1, n_1)$). The undetermined functions $(\rho_2, u_2, \theta_2, n_2)$ require us to construct the next order expansion to determine these functions. In order to obtain an expansion formulation with as less terms as possible, we take these functions as zero values for simplicity. Therefore, it is reasonable to guess some missing information, arising from choosing such special auxiliary functions $g_k^\pm (k=1,2)$, $E_1$ and $B_1$, would be contained in the remainder terms $(g_{R,\eps}^\pm, E_{R,\eps}, B_{R,\eps})$ with order of $\eps$. Note that we prove in our main result Theorem \ref{Main-Thm} that the remainders themselves are uniformly bounded in appropriate spaces. In this sense, our truncation on the expansion formulation is optimal. The corresponding remainder system \eqref{Remd-Equ-Bref} is derived in the formal analysis, see Section \ref{Sec:Formal-Analysis} below.

The main novelty of this paper is to obtain the uniformly global energy bounds of the remainder system \eqref{Remd-Equ-Bref} with small initial data and small Knudsen number $\eps$. The key points have two aspects: (1) to find enough dissipative or decay structures required in deriving the global energy bounds with uniformly small initial data. (2) to seek some structures to deal with the singularity. About finding dissipative or decay structures in the first point, it will be divided into two main ingredients: the pure spatial derivative estimates and the mixed $(x,v)$-derivative estimates.

Considering the pure spatial derivative estimates, the first step is to find the fluid dissipation of $G_{R, \eps}$ by employing the so-called micro-macro decomposition scheme. As inspired by \cite{Guo-Inventions2003}, we first introduce a \emph{seventeen moments} $\mathfrak{B}$ (see \eqref{Basis-B} later) and define a projection $\mathcal{P}_{\mathfrak{B}}$ from $L^2_v$ to $\mathrm{span} \{ \mathfrak{B} \} \subseteq L^2_v$ (see \eqref{S-3} below). We project the remainder system ($G_{R, \eps}$-equation) in \eqref{Remd-Equ-Bref} into $\ker \mathbb{L}$ by the projection operator $\mathbb{P}$ in \eqref{Projc-P}. The fluid variables $(\rho^\pm_{R,\eps}, u_{R,\eps}, \theta_{R,\eps}, E_{R,\eps}, B_{R,\eps})$ of $G_{R,\eps}$ are subject to the macroscopic equations \eqref{MM-10}. The $G_{R, \eps}$-equation in \eqref{Remd-Equ-Bref} is also projected into $\mathrm{span} \{ \mathfrak{B} \}$ by the projection operator $\mathcal{P}_{\mathfrak{B}}$ in \eqref{S-3}. The fluid variables obey the equations \eqref{MM-8} too. Based on the two groups of equations \eqref{MM-8} and \eqref{MM-10}, the fluid dissipation $\| \mathbb{P} G_{R,\eps} \|^2_{H^{N+1}_x L^2_v} \sim \| (\rho^+_{R,\eps}, \rho^-_{R, \eps}, u_{R, \eps}, \theta_{R, \eps}) \|^2_{H^{N+1}_x}$ is constructed in Lemma \ref{Lm-Mic-Mac-decomp}. Note that there are several average quantities on $(\rho^+_{R,\eps}, \rho^-_{R, \eps}, u_{R, \eps}, \theta_{R, \eps})$ in the right-hand side of \eqref{Mic-Mac-Inq} in Lemma \ref{Lm-Mic-Mac-decomp}. They will be controlled by employing the conservation laws \eqref{Totl-Consv-Law-1}, see Lemma \ref{Lm-Integral-Bnd}.

The second step on estimating the pure spatial derivative bounds is to seek some decay structures of the electric field $E_{R,\eps}$ and the magnetic field $B_{R,\eps}$. Observe that Faraday's law $\partial_t B_{R,\eps} + \nabla_x \times E_{R,\eps} = 0$ in \eqref{Remd-Equ-Bref} does not have explicit dissipative term. Together with the Amp\`ere equation $\partial_t E_{R,\eps} - \nabla_x \times B_{R,\eps} = - \tfrac{1}{\eps} \langle G_{R,\eps} \cdot \mathcal{T}_1 v \rangle$ and the Gaussian law $\div_x B_{R,\eps} = 0$ in \eqref{Remd-Equ-Bref}, $B_{R,\eps}$ obeys the formulation \eqref{Mxw-1}, i.e.,
\begin{equation*}
	\partial_{tt} B_{R,\eps} - \Delta_x B_{R,\eps} = \tfrac{1}{\eps} \nabla_x \times \langle \mathbb{P}^\perp G_{R,\eps} \cdot \mathcal{T}_1 v \rangle \,,
\end{equation*}
which seems to be a hyperbolic system with source term $\tfrac{1}{\eps} \nabla_x \times \langle \mathbb{P}^\perp G_{R,\eps} \cdot \mathcal{T}_1 v \rangle$. The decay or dissipation of \eqref{Mxw-1} remains not enough. The \emph{key point} is the discovery of the kinetic version of Ohm's law \eqref{Mxw-4}, namely,
\begin{equation*}
	\begin{aligned}
		\tfrac{1}{\eps} \langle G_{R,\eps} \cdot \mathcal{T}_1 v \rangle =  \sigma E_{R,\eps} - \tfrac{1}{2} \sigma \nabla_x ( \rho_{R,\eps}^+ - \rho_{R,\eps}^- ) + \mathcal{K} ( \mathbb{P}^\perp G_{R,\eps} ) \,.
	\end{aligned}
\end{equation*}
Note that $\langle G_{R,\eps} \cdot \mathcal{T}_1 v \rangle = \langle \mathbb{P}^\perp G_{R,\eps} \cdot \mathcal{T}_1 v \rangle$. Then $B_{R,\eps}$ satisfies the damped wave system \eqref{Mxw-5}, i.e.,
\begin{equation*}
	\begin{aligned}
		\partial_{tt} B_{R,\eps} - \Delta_x B_{R,\eps} + \sigma \partial_t B_{R,\eps} = \nabla_x \times \mathcal{K} ( \mathbb{P}^\perp G_{R,\eps} ) \,.
	\end{aligned}
\end{equation*}
Furthermore, by $E_{R,\eps}$-equation in \eqref{Remd-Equ-Bref} and the (kinetic) Ohm's law \eqref{Mxw-4}, $E_{R,\eps}$ satisfies \eqref{Mxw-6}, i.e.,
\begin{equation*}
	\begin{aligned}
		\partial_t E_{R,\eps} + \sigma E_{R,\eps} - \tfrac{1}{2} \sigma \nabla_x \div_x E_{R,\eps} - \nabla_x \times B_{R,\eps} + \mathcal{K} ( \mathbb{P}^\perp G_{R,\eps} ) = 0 \,.
	\end{aligned}
\end{equation*}
By the above two equations (or say, \eqref{Mxw-5}-\eqref{Mxw-6}), the decay structures $\| \partial_t B_{R,\eps} \|^2_{H^N_x} + \| \nabla_x B_{R,\eps} \|^2_{H^{N-1}_x}$ and $\| E_{R,\eps} \|^2_{H^{N-1}_x} + \| \div_x E_{R,\eps} \|^2_{H^{N-1}_x}$ are constructed in Lemma \ref{Lm-Mxw-Dec}.

The third step on deriving the pure spatial derivative estimates should be focused on the kinetic dissipation related to $G_{R,\eps}$. Actually, the hypocoercivity property of the linearized collision operator $\mathbb{L}$ \eqref{L-Proty-2} in Lemma \ref{Lm-L-Propty} will reduce to the required kinetic dissipation $\tfrac{1}{\eps^2} \| \mathbb{P}^\perp G_{R, \eps} \|^2_{H^{N+1}_x L^2_v (\nu)}$ as in Lemma \ref{Lm-Unif-Spatial-Bnd}. Thus, the process of the pure spatial derivative estimates is finished.

However, the pure spatial derivative estimate \eqref{Unif-Spatial-Bnd} in Lemma \ref{Lm-Unif-Spatial-Bnd} is not closed. The first reason is that the quantity $ \| \nabla_v \mathbb{P}^\perp G_{R,\eps} \|^2_{H^N_x L^2_v}$ coming from the Lorentz force term $\mathcal{T} (v \times B_0) \cdot \nabla_v G_{R,\eps}$ is not yet controlled, which motivates us to consider the mixed $(x,v)$-derivative estimates of the remainder system \eqref{Remd-Equ-Bref}. This is different from the Boltzmann equation. As in \cite{JXZ-IUM-2018}, the Boltzmann equation has the closed pure spatial derivative estimates, because it does not have the structures like Lorentz force. The second reason is that the quantities $ \| H_{R,\eps} \|^2_{H^N_x L^2_v}$, $ \| \mathcal{P}_{\mathfrak{B}} H_{R,\eps} \|^2_{H^N_x L^2_v}$ and $\sum_{|m|\leq N+1} \int_{\T} \langle \partial^m H_{R,\eps} \cdot \partial^m G_{R,\eps}  \rangle \d x$ in the right-hand of \eqref{Unif-Spatial-Bnd} have not yet been dominated. As in \eqref{Remnd-H}, the quantity $H_{R,\eps}$ contains $\mathcal{Q} (g_{R,\eps}^\pm, g_{R,\eps}^\pm)$. These quantities should be controlled by employing Lemma \ref{Lm-L-Q}, in which the weighted norms with the weight $\w (v) = (1+|v|^2)^{1/2}$ are required. Therefore, we will derive the estimates on the mixed $(x,v)$-derivatives with weighted norms, as in Lemma \ref{lemm:remainder-apriori} and Lemma \ref{lemm:remainder-apriori2}. Thanks to the hypocoercivity of $\mathbb{L}$, there will also generate two kinetic dissipative contributions $\tfrac{1}{\eps^2} \| \mathbb{P}^\perp G_{R,\eps} \|^2_{H^N_{x,v} (\w^{2l} \nu)} $ and $\tfrac{1}{\eps^2} \| \mathbb{P}^\perp G_{R,\eps} \|^2_{\widetilde{\bf H}^{N+1}_{x,v} (\w^{2l} \nu)} $. Finally, Lemma \ref{Lm-HR-Square} and Lemma \ref{Lm-HR-Inner-Product} will dominate the quantities related to $H_{R,\eps}$, so that the closed estimate \eqref{ED-2} in Proposition \ref{Prop-Est} is derived.

On the other hand, one should carefully deal with the singularity appeared in the remainder system of $G_{R,\eps}$ \eqref{Remd-Equ-Bref}. This mainly contains the following three aspects:
\begin{enumerate}
	\item The above constructed kinetic dissipative structures $\tfrac{1}{\eps^2} \| \mathbb{P}^\perp G_{R,\eps} \|^2_{H^{N+1}_x L^2_v (\nu)}$, $\tfrac{1}{\eps^2} \| \mathbb{P}^\perp G_{R,\eps} \|^2_{H^N_{x,v} (\w^{2l} \nu)} $ and $\tfrac{1}{\eps^2} \| \mathbb{P}^\perp G_{R,\eps} \|^2_{\widetilde{\bf H}^{N+1}_{x,v} (\w^{2l} \nu)} $ are all with singular factor $\tfrac{1}{\eps^2}$. They are important to absorb the singularity of kinetic part $\mathbb{P}^\perp G_{R,\eps}$ appeared in $G_{R,\eps}$-equation \eqref{Remd-Equ-Bref}.

	\item There is a singular term $- \tfrac{1}{\eps} \langle G_{R,\eps} \cdot \mathcal{T}_1 v \rangle$ in the Amp\`ere equation in \eqref{Remd-Equ-Bref}, i.e., $\partial_t E_{R,\eps} - \nabla_x \times B_{R,\eps} = - \tfrac{1}{\eps} \langle G_{R,\eps} \cdot \mathcal{T}_1 v \rangle$. Thanks again to the kinetic version of Ohm's law in \eqref{Mxw-4},
	the singularity in Amp\`ere equation can be replaced by the terms without singularity.

	\item In the process of deriving the kinetic dissipation in \S \ref{Sec:Lmm-USB}, there are two singular quantities in \eqref{Sp-10} and \eqref{Sp-12} required to be dominated, namely,
	\begin{equation*}
			- \tfrac{1}{\eps} \int_{\T} \langle \partial^m [\mathcal{T} (v \times B_{R,\eps}) \cdot \nabla_v G_0 ] \cdot \partial^m G_{R,\eps} \rangle \d x
	\end{equation*}
    and
    \begin{equation*}
    		- \tfrac{1}{\eps} \sum_{0 \neq m' \leq m} C_m^{m'} \int_{\T} \langle [ \mathcal{T} ( v \times \partial^{m'} B_0 ) \cdot \nabla_v \partial^{m-m'} G_{R,\eps} ] \cdot \partial^m G_{R,\eps} \rangle \d x \,.
    \end{equation*}
    respectively. Observe that, after splitting $G_{R,\eps} = \mathbb{P} G_{R,\eps} + \mathbb{P}^\perp G_{R,\eps}$, the singularity $\tfrac{1}{\eps}$ in the front of $\mathbb{P}^\perp G_{R,\eps}$ can be absorbed by the kinetic dissipative structures with singular factor $\tfrac{1}{\eps^2}$. However, singularity of the fluid part $\mathbb{P} G_{R,\eps}$ seems hard to be controlled. Fortunately, the key cancellation relations \eqref{Sp-9} and \eqref{Sp-11} will be used to overcome this issue, namely,
    \begin{align*}
    		\langle [ \mathcal{T} (v \times \partial^{m - m'} B_{R,\eps} ) \cdot \nabla_v \partial^{m'} G_0 ] \cdot \partial^m \mathbb{P} G_{R,\eps} \rangle & = 0 \,, \\
    		\langle [ \mathcal{T} (v \times \partial^{m'} B_0) \cdot \nabla_v \partial^{m-m'} \mathbb{P} G_{R,\eps} ] \cdot \partial^m \mathbb{P} G_{R,\eps} \rangle & = 0 \,.
    \end{align*}
\end{enumerate}

\subsection{Organization of this paper}

In the next section, the formal derivation of limit equations will be obtained by Hilbert expansion approach. The corresponding remainder system \eqref{Remd-Equ-Bref} are also derived. In Section \ref{Sec:Prem}, the analysis on collision operators is first given. Moreover, the energy estimates for NSFM system \eqref{NSFM} and the linear Maxwell-type system \eqref{Od:eps(10)}-\eqref{Spec-Chos*} are also constructed. In Section \ref{Sec:UEB}, one first derives the uniform energy bounds for the remainder system. The proof of Theorem \ref{Main-Thm} can thus be completed. Section \ref{Sec:detailed_lemmas} is devoted to the detailed proofs for those lemmas stated in previous section. Precisely, \S \ref{Sec:Lmm-MM} is used to deal with the fluid variables, i.e., to prove Lemma \ref{Lm-Mic-Mac-decomp}. In \S \ref{Sec:Lmm-MD}, one focuses on deriving the estimates of $(E_{R,\eps}, B_{R,\eps})$, namely, proving Lemma \ref{Lm-Mxw-Dec}. \S \ref{Sec:Lmm-USB} is devoted to justifying Lemma \ref{Lm-Unif-Spatial-Bnd}, i.e., deriving the kinetic dissipation and completing the pure spatial derivative estimates. In \S \ref{Sec:Lmm-RA12}, we prove Lemmas \ref{lemm:remainder-apriori}-\ref{lemm:remainder-apriori2} associated to the mixed derivative estimates. In \S \ref{Sec:HR}, we justify Lemmas \ref{Lm-HR-Square}-\ref{Lm-HR-Inner-Product} on estimates of $H_{R,\eps}$. In Appendix \ref{Appendix-A}, the proof of Lemma \ref{lemm:bnd-linearMaxwl} about the estimates on $(\rho_1, u_1, \theta_1, n_1, B_1)$ is given.

\section{Formal Analysis: Hilbert Expansion}
\label{Sec:Formal-Analysis}

In this section, the goal is to derive formally the limit from the scaled VMB model \ref{VMB-f} to the NSFM model \eqref{NSFM} by Hilbert expansion approach. The remainder equations will be obtained by truncating the expansion ansatz.

\subsection{Derivation of limit equations}\label{Subsec:Order}

By the relation $ F^{\pm}_\eps = M ( 1 + \eps g^{\pm}_\eps ) $, one can rewrite the system \ref{VMB-f} as
\begin{equation}\label{VMB-g}
  \left\{
    \begin{array}{rl} \displaystyle
      \eps \partial_t \begin{pmatrix} g^+_\eps \\ g^-_\eps \end{pmatrix}
      + v \cdot \nabla_x \begin{pmatrix} g^+_\eps \\ g^-_\eps \end{pmatrix}
      + ( \eps E_\eps + v \times B_\eps ) & \cdot \nabla_v \begin{pmatrix} g^+_\eps \\ -g^-_\eps \end{pmatrix}
      - E_\eps \cdot v \begin{pmatrix} 1 + \eps g^+_\eps \\ -1 - \eps g^-_\eps \end{pmatrix}
    \\[7pt] \displaystyle
      = & - \frac{1}{\eps} \mathbb{L} \begin{pmatrix} g^+_\eps \\ g^-_\eps \end{pmatrix}
      		+ \begin{pmatrix}
              \mathcal{Q} ( g^+_\eps , g^+_\eps ) + \mathcal{Q} ( g^+_\eps , g^-_\eps ) \\
              \mathcal{Q} ( g^-_\eps , g^-_\eps ) + \mathcal{Q} ( g^-_\eps , g^+_\eps )
            \end{pmatrix} \,,
  \\[7pt]
    \partial_t E_\eps - \nabla_x \times B_\eps = & - \frac{1}{\eps} \int_{\R^3} (g^+_\eps - g^-_\eps) v M \d v \,,
  \\[4pt]
  	\div_x E_\eps = & \int_{\R^3} (g^+_\eps - g^-_\eps) M \d v \,,
  \\[3pt]
    \partial_t B_\eps + \nabla_x \times E_\eps = & 0 \,,
  \\[2pt]
  	\div_x B_\eps = & 0 \,,
  \end{array}
  \right.
\end{equation}
where the linearized collision operator $\mathbb{L}$ is defined in \eqref{Def-L-bb}, and $\mathcal{Q} (\cdot, \cdot)$ is given in \eqref{Q-cal}.

We take the following expansion ansatz, with respect to the parameter $\eps \rightarrow 0$, that
\begin{equation}\label{Ansatz-1}
    g^{\pm}_\eps = g^{\pm}_0 + \eps g^{\pm}_1 + \eps^2 g^{\pm}_2 + \cdots \,, \quad E_\eps = E_0 + \eps E_1 + \cdots \,, \quad B_\eps = B_0 + \eps B_1 + \cdots \,.
\end{equation}
Plugging the ansatz \eqref{Ansatz-1} into the system \eqref{VMB-g}, the balance of orders will provide a framework which helps us analyze the hydrodynamic limit.

\subsubsection{Order of \texorpdfstring{$\mathcal{O} (\frac{1}{\eps})$}{O(1/epsilon)}}

The order of $\mathcal{O}(\frac{1}{\eps})$ reads
\begin{equation}\label{Od:eps^-1(1)}
  \mathbb{L} \left(
                \begin{array}{c}
                  g_0^+ \\
                  g_0^-
                \end{array}
              \right) = 0 \,, \quad \int_{\R^3} ( g_0^+ - g_0^- ) v M \d v = 0 \,.
\end{equation}
Recalling the structure of the kernel $\ker \mathbb{L}$ given in \eqref{L-ker} above, one derives
  \begin{equation}\label{g0}
    g_0^{\pm} = \rho_0^{\pm} + u_0 \cdot v + \theta_0 ( \tfrac{|v|^2}{2} - \tfrac{3}{2} ) \,,
  \end{equation}
where $\rho_0^{\pm} = \langle g_0^{\pm} \rangle $, $ u_0 = \langle g_0^{+} v \rangle = \langle g_0^{-} v \rangle $ and $ \theta_0 = \langle g_0^+ (\frac{|v|^2}{3} - 1  ) \rangle = \langle g_0^- (\frac{|v|^2}{3} - 1  ) \rangle $. Remark that the second relation in \eqref{Od:eps^-1(1)} is exactly consistent with the structure of $\ker \mathbb{L}$, namely, $\langle g_0^{+} v \rangle = \langle g_0^{-} v \rangle$.

\subsubsection{Order of \texorpdfstring{$\mathcal{O}(1)$}{O(1)}}
	The order of $\mathcal{O}(1)$ in the system \eqref{VMB-g} reads
    \begin{equation}\label{Od:1(1)}
    \begin{aligned}
      v \cdot \nabla_x
    \left(
     \begin{array}{c}
        g_0^+ \\
        g_0^-
      \end{array}
    \right)
    +
    (v \times B_0) \cdot & \nabla_v
    \left(
     \begin{array}{c}
        g_0^+ \\
        g_0^-
      \end{array}
    \right)
    - E_0 \cdot v
    \left(
     \begin{array}{c}
        1 \\
        -1
      \end{array}
    \right) \\
    = & -
    \mathbb{L}\left(
     \begin{array}{c}
        g_1^+ \\
        g_1^-
      \end{array}
    \right)
    +
    \left(
     \begin{array}{c}
        \Q ( g_0^+ , g_0^+ ) + \Q ( g_0^+ , g_0^- ) \\
        \Q ( g_0^- , g_0^- ) + \Q ( g_0^- , g_0^+ )
      \end{array}
    \right) \,,
    \end{aligned}
    \end{equation}
    and
    \begin{equation}\label{Od:1(2)}
      \left\{
        \begin{array}{l}
          \partial_t E_0 - \nabla_x \times B_0 = - j_0 \,, \ \ \, \qquad \qquad \quad \div_x E_0 = n_0 \,, \\
          \partial_t B_0 + \nabla_x \times E_0 = 0 \,, \qquad \qquad \qquad \quad \div_x B_0 = 0 \,, \\
          j_0 = u_1^+ - u_1^- = \int_{\R^3} ( g_1^+ - g_1^- ) v M \d v \,, \ n_0 = \rho_0^+ - \rho_0^- \,.
        \end{array}
      \right.
    \end{equation}
    The expression of $g_0^{\pm}$ \eqref{g0} implies
    \begin{equation}\label{Od:1-1}
      ( v \times B_0 ) \cdot \nabla_v g_0^{\pm} = - ( u_0 \times B_0 ) \cdot v \,.
    \end{equation}
  Let
    \begin{equation*}
    	\begin{aligned}
    		\A(v) = v \otimes v - \tfrac{|v|^2}{3} Id \,, \  
    		\B(v) = v ( \tfrac{|v|^2}{2} - \tfrac{5}{2} )  \,, \  
    \text{ and } \C (v) = \tfrac{1}{4} |v|^4 - \tfrac{5}{2} |v|^2 + \tfrac{15}{4} 
    	\end{aligned}
    \end{equation*}
  be the Burnett functions, see \cite{BGL-CPAM1993}. They belong to the kernel orthogonal space of the linearized Boltzmann operator, denoting by $\ker^\perp L$, with $L$ being defined in \eqref{Def-L-cal}. A straightforward calculation yields
    \begin{equation}\label{Od:1-2}
      v \cdot \nabla_x g_0^{\pm}
      = \div_x u_0 + \nabla_x ( \rho_0^{\pm} + \theta_0 ) \cdot v + \tfrac{2}{3} \div_x u_0 ( \tfrac{|v|^2}{2} - \tfrac{3}{2} ) + \A(v) : \nabla_x u_0 + \B(v) \cdot \nabla_x \theta_0 \,.
    \end{equation}
    Denote $\rho_0 = \frac{\rho_0^+ + \rho_0^-}{2}$. Note that
    \begin{equation*}
    	\begin{aligned}
    		\L g_1^+ + \L g_1^- + & \L^\flat ( g_1^+ , g_1^- ) + \L^\flat ( g_1^- , g_1^+ )
    		= & \L ( g_1^+ + g_1^- ) - \Q ( g_1^+ + g_1^- , 1 ) - \Q ( 1 , g_1^+ + g_1^- ) \\
        = & 2 \L ( g_1^+ + g_1^- ) \,.
    	\end{aligned}
    \end{equation*}
    We also have, due to the fact $\Q ( g, g ) = \frac{1}{2} \L ( g^2 )$ for $g \in \ker \L$ (see \cite{BGL-CPAM1993}), that
    \begin{equation*}
    	\begin{aligned}
    		\tfrac{1}{2} \left[ \Q ( g_0^+ , g_0^+ ) + \Q ( g_0^+ , g_0^- ) + \Q ( g_0^- , g_0^- ) + \Q ( g_0^- , g_0^+ ) \right]
    		& = \tfrac{1}{2} \Q ( g_0^+ + g_0^- , g_0^+ + g_0^- )
          = \tfrac{1}{4} \L \left[ ( g_0^+ + g_0^- )^2 \right]  \\
        &  = \L ( \A(v) : u_0 \otimes u_0 + 2 \theta_0 u_0 \cdot \B (v) + \theta_0^2 \C (v) ) \,,
    	\end{aligned}
    \end{equation*}
    and
    \begin{align*}
    		\tfrac{1}{4} ( g_0^+ + g_0^- )^2
    		= & \rho_0^2 + 2 ( \rho_0 + \theta_0 ) u_0 \cdot v + 2 \rho_0 \theta_0 ( \tfrac{|v|^2}{2} - \tfrac{3}{2} ) + \tfrac{|v|^2}{3} |u_0|^2 \\
    		& + \A(v) : u_0 \otimes u_0 + 2 \theta_0 u_0 \cdot \B (v) + \theta_0^2 \C (v) \,.
    \end{align*}

Summing up the two equations in \eqref{Od:1(1)} and multiplying by $\frac{1}{2}$, together with the equalities \eqref{Od:1-1}-\eqref{Od:1-2}, ensure us to derive
    \begin{multline}\label{Od:1-3}
       \div_x u_0 + v \cdot \nabla_x (\rho_0 + \theta_0) + \tfrac{2}{3} \div_x u_0 ( \tfrac{|v|^2}{2} - \tfrac{3}{2} ) + \A(v) : \nabla_x u_0 + \B(v) \cdot \nabla_x \theta_0 \\
       = - 2 \L ( \tfrac{g_1^+ + g_1^-}{2} ) + \L ( \A(v) : u_0 \otimes u_0 + 2 \theta_0 u_0 \cdot \B (v) + \theta_0^2 \C (v) ) \,.
    \end{multline}
Performing projection on the equation \eqref{Od:1-3} into $\ker \L $, by recalling the definition $\P_{\L} : L^2(M\!\d v) \rightarrow \ker \L$ in \eqref{Projc-PL}, reduces to $\div_x u_0 = 0$ and $\nabla_x (\rho_0 + \theta_0) = 0$. This yields,
    \begin{equation}\label{Od:1(3)}
      \div_x u_0 = 0 \,, \ \rho_0 + \theta_0 = 0 \,.
    \end{equation}
Besides, the kernel orthogonal part of \eqref{Od:1-3} is
    \begin{equation}\label{Od:1-4}
      \L ( \tfrac{g_1^+ + g_1^-}{2} )= \L ( \tfrac{1}{2} u_0 \otimes u_0 : \A(v) + \theta_0 u_0 \cdot \B(v) + \tfrac{1}{2} \theta_0^2 \C(v) - \tfrac{1}{2} \nabla_x u_0 : \widehat{\A} (v) - \tfrac{1}{2} \nabla_x \theta_0 \cdot \widehat{\B} (v) ) \,,
    \end{equation}
  where the terms $\widehat{\A}(v)$, $ \widehat{\B}(v) \in \ker^\perp \L $ is uniquely determined in $ \ker^\perp \L $, satisfying
    \begin{equation*}
      \L \widehat{\A}(v) = \A(v) \,, \ \textrm{and} \quad \L \widehat{\B}(v) = \B(v) \,,
    \end{equation*}
  Moreover, there exits functions $ \varphi\,, \ \psi : \R^+ \rightarrow \R^+$ such that $ \widehat{\A}(v) = \varphi(|v|) \A (v) $ and $ \widehat{\B}(v) = \psi(|v|) \B (v)$, see \cite{BGL-CPAM1993} for instance.

As a result, the equation \eqref{Od:1-4} yields, for some undetermined functions $\rho_1(t,x)$, $u_1(t,x)$ and $\theta_1(t,x)$, that
    \begin{equation}\label{Od:1(4)}
      \begin{aligned}
        \tfrac{g_1^+ + g_1^-}{2} = & \rho_1 + u_1 \cdot v + \theta_1 ( \tfrac{|v|^2}{2} - \tfrac{3}{2} ) \\
         + & \tfrac{1}{2} u_0 \otimes u_0 : \A(v) + \theta_0 u_0 \cdot \B(v) + \tfrac{1}{2} \theta_0^2 \C(v)  - \tfrac{1}{2} \nabla_x u_0 : \widehat{\A} (v) - \tfrac{1}{2} \nabla_x \theta_0 \cdot \widehat{\B} (v) \\
         : = & l_1 + l_1^\perp \in \ker \L + \ker^\perp \L \,.
      \end{aligned}
    \end{equation}

By definitions of $\mathcal{L}$, $\mathcal{L}^\flat$ and $\mathfrak{L}$ in \eqref{Def-L-cal}, \eqref{Def-L-flat} and \eqref{Def-L-frk}, respectively, we can write
\begin{equation*}
	\begin{aligned}
		\L g_1^+ + & \L^\flat ( g_1^+ , g_1^- ) - \L g_1^- - \L^\flat ( g_1^- , g_1^+ ) \\
		= & \L ( g_1^+ - g_1^- ) - \Q (g_1^+ , 1) - \Q(1, g_1^-) + \Q ( g_1^- , 1 ) + \Q ( 1 , g_1^+ ) = ( \L + \mathfrak{L} ) ( g_1^+ - g_1^- ) \,.
	\end{aligned}
\end{equation*}
By \eqref{g0}, we have
\begin{equation*}
	\begin{aligned}
		\Q ( g_0^+ , g_0^+ ) + \Q ( g_0^+ , g_0^- ) - \Q ( g_0^- , g_0^- ) - \Q ( g_0^- , g_0^+ ) = \Q ( g_0^+ - g_0^- , g_0^+ + g_0^- ) \\
		= 2 n_0 \Q ( 1 , u_0 \cdot v + \theta_0 ( \tfrac{|v|^2}{2} - \tfrac{3}{2} ) ) = n_0 ( \L + \mathfrak{L} ) ( u_0 \cdot v + \theta_0 ( \tfrac{|v|^2}{2} - \tfrac{3}{2} ) ) \,.
	\end{aligned}
\end{equation*}
Thus, it follows from subtracting the second equality from the first equality in \eqref{Od:1(1)} and multiplying by $\frac{1}{2}$ that
    \begin{equation}\label{Od:1-5}
      \begin{aligned}
        - ( - \tfrac{1}{2} \nabla_x n_0 + E_0 + & u_0 \times B_0 ) \cdot v \\
        = & - ( \L + \mathfrak{L} ) ( \tfrac{g_1^+ - g_1^-}{2} ) + \tfrac{n_0}{2} ( \L + \mathfrak{L} ) ( u_0 \cdot v + \theta_0 ( \tfrac{|v|^2}{2} - \tfrac{3}{2} ) ) \,.
      \end{aligned}
    \end{equation}
As shown in \cite{Arsenio-SaintRaymond}, for $\Phi (v) = v \in L^2(M\!\d v)$ and $\Psi (v) = \tfrac{|v|^2}{2} - \tfrac{3}{2} \in L^2(M\!\d v)$, there exist two functions $ \widetilde{\Phi} \,, \ \widetilde{\Psi} \in \ker^\perp(\L + \mathfrak{L}) $, such that
    \begin{equation}\label{L+L-Dual}
      ( \L + \mathfrak{L} ) \widetilde{\Phi} = \Phi \quad \textrm{and} \quad ( \L + \mathfrak{L} ) \widetilde{\Psi} = \Psi \,,
    \end{equation}
    which can be uniquely determined in $\ker^\perp(\L + \mathfrak{L})$. Furthermore, there exist two scalar functions $\alpha \,, \ \beta \,: \R^+ \rightarrow \R$, such that
    $$ \widetilde{\Phi} (v) = \alpha(|v|) \Phi (v) \quad \textrm{and} \quad \widetilde{\Psi} (v) = \beta(|v|) \Psi (v)  \,, $$
    which implies that
    $$ \int_{\R^3} \Phi_i (v) \widetilde{\Phi}_j (v) M \d v = \tfrac{1}{2} \sigma \delta_{ij} \,.  $$
This defines the electrical conductivity $ \sigma = \frac{2}{3} \int_{\R^3} \Phi \cdot \widetilde{\Phi} M \d v > 0$. One can also obtain the energy conductivity $ \lambda = \int_{\R^3} \Psi \cdot \widetilde{\Psi} M \d v  > 0$.

As a consequence, the equation \eqref{Od:1-5} implies that, for some undetermined scalar function $n_1(t,x)$,
    \begin{equation}\label{Od:1(5)}
      \tfrac{g_1^+ - g_1^-}{2} = \tfrac{1}{2} n_1 + \tfrac{1}{2} n_0 u_0 \cdot v + \tfrac{1}{2} n_0 \theta_0 ( \tfrac{|v|^2}{2} - \tfrac{3}{2} ) + ( - \tfrac{1}{2} \nabla_x n_0 + E_0 + u_0 \times B_0 ) \cdot \widetilde{\Phi} (v) \,.
    \end{equation}
	Therefore, projecting \eqref{Od:1(5)} to $\ker \L$ reduces to
    \begin{equation}\label{Od:1(6)}
      \left\{
        \begin{array}{l}
          n_1 = \langle g_1^+ - g_1^- \rangle = \rho_1^+ - \rho_1^- \,, \qquad \theta_1^+ - \theta_1^- = n_0 \theta_0 \,, \\
          j_0 = u_1^+ - u_1^- = n_0 u_0 + \sigma ( - \tfrac{1}{2} \nabla_x n_0 + E_0 + u_0 \times B_0 ) \,,
        \end{array}
      \right.
    \end{equation}
  where $ u_1^{\pm} = \langle g_1^{\pm} v \rangle \,, \ \theta_1^{\pm} = \langle g_1^{\pm} ( \tfrac{|v|^2}{3} - 1 ) \rangle$. The last equation is exactly the Ohm's law.

\subsubsection{Order of \texorpdfstring{$\mathcal{O}(\eps)$}{O(epsilon)}}

  The order of $\mathcal{O}(\eps)$ in the system \eqref{VMB-g} reads
    \begin{multline}\label{Od:eps(1)}
        \partial_t
        \left(
          \begin{array}{c}
            g_0^+ \\
            g_0^-
          \end{array}
        \right)
        +
        v \cdot \nabla_x
        \left(
          \begin{array}{c}
            g_1^+ \\
            g_1^-
          \end{array}
        \right)
        +
        E_0 \cdot \nabla_v
        \left(
          \begin{array}{c}
            g_0^+ \\
            - g_0^-
          \end{array}
        \right)
        +
        ( v \times B_0 ) \cdot \nabla_v
        \left(
          \begin{array}{c}
            g_1^+ \\
            - g_1^-
          \end{array}
        \right) \\
      + ( v \times B_1 ) \cdot \nabla_v
        \left(
          \begin{array}{c}
            g_0^+ \\
            - g_0^-
          \end{array}
        \right)
        -
        E_0 \cdot v
        \left(
          \begin{array}{c}
            g_0^+ \\
            - g_0^-
          \end{array}
        \right)
        -
        E_1 \cdot v
        \left(
          \begin{array}{c}
            1 \\
            - 1
          \end{array}
        \right) \\
      = - \mathbb{L}
        \left(
          \begin{array}{c}
            g_2^+ \\
            g_2^-
          \end{array}
        \right)
        +
        \left(
          \begin{array}{c}
            \Q(g_0^+ , g_1^+ ) + \Q(g_1^+ , g_0^+ ) + \Q(g_0^+ , g_1^- ) + \Q(g_1^+ , g_0^- ) \\
            \Q(g_0^- , g_1^- ) + \Q(g_1^- , g_0^- ) + \Q(g_0^- , g_1^+ ) + \Q(g_1^- , g_0^+ )
          \end{array}
        \right) \,,
    \end{multline}
    and
    \begin{equation}\label{Od:eps(2)}
      \left\{
        \begin{array}{ll}
          \partial_t E_1 - \nabla_x \times B_1 = - ( u_2^+ - u_2^- ) \,, \quad  & \div_x E_1 = n_1 \,, \\
          \partial_t B_1 + \nabla_x \times E_1 = 0 \,, & \div_x B_1 = 0 \,.
        \end{array}
      \right.
    \end{equation}
The equation \eqref{Od:eps(1)} can be rewritten equivalently as
    \begin{equation}\label{Od:eps-1}
      \begin{aligned}
        & \partial_t
        \left(
          \begin{array}{c}
            \tfrac{g_0^+ + g_0^-}{2} \\[3pt]
            \tfrac{g_0^+ - g_0^-}{2}
          \end{array}
        \right)
        + v \cdot \nabla_x
        \left(
          \begin{array}{c}
            \tfrac{g_1^+ + g_1^-}{2} \\[3pt]
            \tfrac{g_1^+ - g_1^-}{2}
          \end{array}
        \right)
        + E_0 \cdot \nabla_v
        \left(
          \begin{array}{c}
            \tfrac{g_0^+ - g_0^-}{2} \\
            \tfrac{g_0^+ + g_0^-}{2}
          \end{array}
        \right)
        + ( v \times B_0 ) \cdot \nabla_v
        \left(
          \begin{array}{c}
            \tfrac{g_1^+ - g_1^-}{2} \\
            \tfrac{g_1^+ + g_1^-}{2}
          \end{array}
        \right) \\
        & \hspace*{1cm} + ( v \times B_1 ) \cdot \nabla_v
        \left(
          \begin{array}{c}
            \tfrac{g_0^+ - g_0^-}{2} \\
            \tfrac{g_0^+ + g_0^-}{2}
          \end{array}
        \right)
        - E_0 \cdot v
        \left(
          \begin{array}{c}
            \tfrac{g_0^+ - g_0^-}{2} \\
            \tfrac{g_0^+ + g_0^-}{2}
          \end{array}
        \right)
        - E_1 \cdot v
        \left(
          \begin{array}{c}
            0 \\
            1
          \end{array}
        \right) \\
    & = -
        \left(
          \begin{array}{c}
             2 \L ( \tfrac{g_2^+ + g_2^-}{2} ) \\
            ( \L + \mathfrak{L} ) ( \tfrac{g_2^+ - g_2^-}{2} )
          \end{array}
        \right)
        +
        \left(
          \begin{array}{c}
            \Q( \tfrac{g_0^+ + g_0^-}{2} , g_1^+ + g_1^- ) + \Q( g_1^+ + g_1^- ,\tfrac{g_0^+ + g_0^-}{2} ) \\
           - n_0 \L ( \tfrac{g_1^+ + g_1^-}{2} ) + \tfrac{1}{2} n_0 ( \L + \mathfrak{L} )  ( \tfrac{g_1^+ + g_1^-}{2} ) + \Q ( g_1^+ - g_1^- , \tfrac{g_0^+ + g_0^-}{2} )
          \end{array}
        \right) \,.
      \end{aligned}
    \end{equation}
Recalling the expression of $g_0^{\pm}$ \eqref{g0} and the relations \eqref{Od:1(3)}, the following terms on the left-hand side of \eqref{Od:eps-1} reduces to
      \begin{align}\label{Od:eps-2}
        \no & \partial_t
        \left(
          \begin{array}{c}
            \tfrac{g_0^+ + g_0^-}{2} \\
            \tfrac{g_0^+ - g_0^-}{2}
          \end{array}
        \right)
        + E_0 \cdot \nabla_v
        \left(
          \begin{array}{c}
            \tfrac{g_0^+ - g_0^-}{2} \\
            \tfrac{g_0^+ + g_0^-}{2}
          \end{array}
        \right)
        + ( v \times B_1 ) \cdot \nabla_v
        \left(
          \begin{array}{c}
            \tfrac{g_0^+ - g_0^-}{2} \\
            \tfrac{g_0^+ + g_0^-}{2}
          \end{array}
        \right)
      \\ \no & \quad
        - E_0 \cdot v
        \left(
          \begin{array}{c}
            \tfrac{g_0^+ - g_0^-}{2} \\
            \tfrac{g_0^+ + g_0^-}{2}
          \end{array}
        \right)
        - E_1 \cdot v
        \left(
          \begin{array}{c}
            0 \\
            1
          \end{array}
        \right)
        \\ \no &
      =
        \left(
          \begin{array}{c}
            \partial_t \rho_0 + ( \partial_t u_0 - \tfrac{1}{2} n_0 E_0 ) \cdot v + \partial_t \theta_0 ( \tfrac{|v|^2}{2} - \tfrac{3}{2} )  \\
            \tfrac{1}{2} \partial_t n_0 + ( \theta_0 E_0 - u_0 \times B_1 - E_1 ) \cdot v - \tfrac{2}{3} E_0 \cdot u_0 ( \tfrac{|v|^2}{2} - \tfrac{3}{2} )
          \end{array}
        \right)
      \\ & \quad
        -
        \left(
          \begin{array}{c}
            0 \\
            E_0 \otimes u_0 : \A (v) + \theta_0 E_0 \cdot \B(v)
          \end{array}
        \right) \,.
      \end{align}

We next deal with the second term 
on the left-hand side of \eqref{Od:eps-1}. By the expression of $\frac{g_1^+ + g_1^-}{2}$ in \eqref{Od:1(4)}, we can compute
    \begin{equation*}
    	\begin{aligned}
    		v \cdot \nabla_x ( \tfrac{g_1^+ + g_1^-}{2} ) = v \cdot \nabla_x l_1 + v \cdot \nabla_x l_1^\perp \,,
    	\end{aligned}
    \end{equation*}
where
    \begin{equation}\label{Od:eps-3}
      \begin{aligned}
        v \cdot \nabla_x l_1 = \div_x u_1 + \nabla_x ( \rho_1 + \theta_1 ) \cdot v + \tfrac{2}{3} \div_x u_1 ( \tfrac{|v|^2}{2} - \tfrac{3}{2} ) + \nabla_x u_1 : \A (v) + \nabla_x \theta_1 \cdot \B(v)
      \end{aligned}
    \end{equation}
with both $\A (v)$ and $\B (v)$ belonging to $\ker^\perp \L$. It is derived from projecting the term $v \cdot \nabla_x l_1^\perp$ into $\ker \L$ that
    \begin{equation}\label{Od:eps-4}
      \left\langle
        v \cdot \nabla_x l_1^\perp
         \left(
          \begin{array}{c}
            1 \\
            v \\
            \tfrac{|v|^2}{3} - 1
          \end{array}
         \right)
      \right\rangle =
      \left(
          \begin{array}{c}
            0 \\
            \div_x ( u_0 \otimes u_0 ) - \mu \Delta_x u_0 - \tfrac{1}{3} \nabla_x ( |u_0|^2 ) \\
            \tfrac{5}{3} \div_x ( \theta_0 u_0 ) - \tfrac{5}{3} \kappa \Delta_x \theta_0
          \end{array}
         \right) \,,
    \end{equation}
where $\mu = \frac{1}{30} \int_{\R^3} \varphi (|v|) |v|^4 M \d v > 0$ and $\kappa = \tfrac{1}{15} \int_{\R^3} \widehat{\B}(v) \cdot \B(v) M \d v > 0$, see \cite{BGL-CPAM1993}.

On the other hand, the expression of $ \frac{g_1^+ - g_1^- }{2} $ in \eqref{Od:1(5)} implies
    \begin{align}\label{Od:eps-5}
        v \cdot \nabla_x ( \tfrac{g_1^+ - g_1^-}{2} )
      = & \tfrac{1}{2} \div_x ( n_0 u_0 ) + \tfrac{1}{2} ( \nabla_x n_1 + \tfrac{5}{2} \nabla_x ( n_0 \theta_0 ) ) \cdot v + \tfrac{1}{3} \div_x ( n_0 u_0 ) ( \tfrac{|v|^2}{2} - \tfrac{3}{2} ) \no \\
        & + \nabla_x ( - \tfrac{1}{2} \nabla_x n_0 + E_0 + u_0 \times B_0 ) : v \otimes \widetilde{\Phi} (v) \\ \no
        & + \tfrac{1}{2} \nabla_x ( n_0 u_0  ) : \A(v) + \tfrac{1}{2} \nabla_x ( n_0 \theta_0 ) \cdot \B(v) \,.
    \end{align}
Moreover, the expressions \eqref{Od:1(4)} and \eqref{Od:1(5)} yield that
    \begin{equation}\label{Od:eps-6}
      \begin{aligned}
        ( v \times B_0 ) \cdot \nabla_v ( \tfrac{g_1^+ - g_1^-}{2} ) = & - ( \tfrac{1}{2} n_0 u_0 \times B_0 ) \cdot v - B_{0i} ( - \tfrac{1}{2} \nabla_x n_0 + E_0 + u_0 \times B_0 )_l \epsilon_{ijk} v_j \partial_{v_k} \widetilde{\Phi}_l (v) \,,
      \end{aligned}
    \end{equation}
    and
    \begin{equation}\label{Od:eps-7}
      \begin{aligned}
        ( v \times B_0 ) \cdot \nabla_v ( \tfrac{g_1^+ + g_1^-}{2} ) = & - ( u_1 \times B_0 ) \cdot v - ( B_0 \times u_0 ) \otimes u_0 : \A(v) - ( n_0 u_0 \times B_0 ) \cdot \B(v) \\
        & + \tfrac{1}{2} \partial_l u_{0m} \epsilon_{ijk} v_j \partial_{v_k} \widehat{\A}_{lm} (v) + \tfrac{1}{2} B_{0i} \partial_l \theta_0 \epsilon_{ijk} v_j \partial_{v_k} \widehat{\B}_l (v) \,,
      \end{aligned}
    \end{equation}
    where $\epsilon_{ijk} \ (1 \leq i,j,k \leq 3)$ is defined as
    \begin{equation}\label{eps-ijk}
    	\epsilon_{ijk} =
    	\left\{
    	  \begin{aligned}
    	  	1 \,, & \quad \textrm{ if the order of } \{ i,j,k \} \ \textrm{is even} \,, \\
    	  	- 1 \,, & \quad \textrm{ if the order of } \{ i,j,k \} \ \textrm{is odd} \,, \\
    	  	0 \,, & \quad \ \textrm{if there are two numbers identical among } i,j,k.
    	  \end{aligned}
    	\right.
    \end{equation}

It is an easy matter to derive that
    \begin{equation}\label{Od:eps-8}
        \P_{\L} \Big( B_{0i} ( - \tfrac{1}{2} \nabla_x n_0 + E_0 + u_0 \times B_0 )_l \epsilon_{ijk} v_j \partial_{v_k} \widetilde{\Phi}_l (v) \Big)
        = \big[ \tfrac{1}{2} \sigma ( - \tfrac{1}{2} \nabla_x n_0 + E_0 + u_0 \times B_0 ) \times B_0 \big] \cdot v \,,
    \end{equation}
    and
    \begin{equation}\label{Od:eps-9}
      \begin{aligned}
        \P_{\L + \mathfrak{L}} \Big( \tfrac{1}{2} \partial_l u_{0m} \epsilon_{ijk} v_j \partial_{v_k} \widehat{\A}_{lm} (v) \Big) = \P_{\L + \mathfrak{L}} \Big( \tfrac{1}{2} B_{0i} \partial_l \theta_0 \epsilon_{ijk} v_j \partial_{v_k} \widehat{\B}_l (v) \Big) = 0 \,,
      \end{aligned}
    \end{equation}
    where the projection $\P_{\L + \mathfrak{L}} : L^2(M\!\d v) \rightarrow \ker (\L + \mathfrak{L}) \subset L^2(M\!\d v) $ is given by $ \P_{\L + \mathfrak{L}} ( g ) = \langle g \rangle \,. $

Observing that, for all $ g , h \in L^2(M\!\d v) $,
  $$
    \Q(g,h) + \Q(h,g) \in \ker^\perp \L \,, \ \textrm{and} \quad \Q(g,h) \in \ker^\perp ( \L + \mathfrak{L} ) \,,
  $$
we perform projection on the first equation of \eqref{Od:eps-1} into $\ker \L$, and combine together the relations \eqref{Od:eps-2}$_1$, \eqref{Od:eps-3}, \eqref{Od:eps-4}, \eqref{Od:eps-6} and \eqref{Od:eps-8}, then it follows that
    \begin{equation*}
      \left\{
        \begin{array}{rl}
          \partial_t \rho_0 + \div_x u_1 = & 0 \,, \\[2pt]
          \partial_t u_0 + \div_x ( u_0 \otimes u_0 ) - & \mu \Delta_x u_0 + \nabla_x ( \rho_1 + \theta_1 - \tfrac{1}{3} |u_0|^2 ) \\[4pt]
          = & \tfrac{1}{2} n_0 E_0 + \tfrac{1}{2} [ n_0 u_0 + \sigma ( - \tfrac{1}{2} \nabla_x n_0 + E_0 + u_0 \times B_0 ) ] \times B_0  \,, \\[4pt]
          \partial_t \theta_0 + \tfrac{2}{3} \div_x u_1 + & \tfrac{5}{3} \div_x (\theta_0 u_0) - \tfrac{5}{3} \kappa \Delta_x \theta_0 = 0\,.
        \end{array}
      \right.
    \end{equation*}
This, together with \eqref{Od:1(2)}, \eqref{Od:1(3)} and \eqref{Od:1(6)}, indicates the limit NSFM system, that
    \begin{equation}\label{Limit-Equ}
      \left\{
        \begin{array}{rll}
          \partial_t u_0 + u_0 \cdot \nabla_x u_0 - \mu \Delta_x u_0 + \nabla_x p_0 = & \tfrac{1}{2} n_0 E_0 + \tfrac{1}{2} j_0 \times B_0 \,, & \ \div_x u_0 = 0 \,,  \\
          \partial_t \theta_0 + u_0 \cdot \nabla_x \theta_0 - \kappa \Delta_x \theta_0 = & 0 \,, & \ \rho_0 + \theta_0 = 0 \,, \\
          \partial_t E_0 - \nabla_x \times B_0 = & - j_0 \,, & \ \div_x E_0 = n_0 \,, \\
          \partial_t B_0 + \nabla_x \times E_0 = & 0 \,, & \ \div_x B_0 = 0 \,, \\
          j_0 = n_0 u_0 +  \sigma ( - \tfrac{1}{2} & \nabla_x n_0 + E_0 + u_0 \times B_0 ) \,, &
        \end{array}
      \right.
    \end{equation}
    and
    \begin{equation}\label{Od:eps(3)}
    	\div_x u_1 = \partial_t \theta_0 \,, \quad
      \nabla_x ( \rho_1 + \theta_1 ) = - \partial_t u_0 - u_0 \cdot \nabla_x u_0 + \mu \Delta_x u_0 + \tfrac{1}{3} \nabla_x ( |u_0|^2 ) + \tfrac{1}{2} j_0 \times B_0 \,.
    \end{equation}

    On the other hand, from projecting the first equation of \eqref{Od:eps-1} into $\ker^\perp \L$,
    \begin{equation}\label{Od:eps-10}
      \begin{aligned}
        \L ( \tfrac{g_2^+ + g_2^-}{2} ) = & \Q ( \tfrac{g_0^+ + g_0^- }{2} , \tfrac{g_1^+ + g_1^- }{2} ) + \Q ( \tfrac{g_1^+ + g_1^- }{2} , \tfrac{g_0^+ + g_0^- }{2} ) \\
        & - \tfrac{1}{2} \nabla_x u_1 : \A(v) - \tfrac{1}{2} \nabla_x \theta_1 \cdot \B(v) - \tfrac{1}{2} \P_{\L}^\perp ( v \cdot \nabla_x l_1^\perp ) \\
        & + \tfrac{1}{2 \sigma} B_{0i} ( j_{0l} - n_0 u_{0l} ) \epsilon_{ijk} \P_{\L}^\perp ( v_j \partial_{v_k} \widetilde{\Phi}_l (v) ) \,,
      \end{aligned}
    \end{equation}
where $l_1^\perp$ is defined in \eqref{Od:1(4)}.

Since we have $ \Q(g,h) + \Q(h,g) = \L (gh) $ for all $g , h \in \ker \L$ (see \cite{BGL-CPAM1993}, for instance), and $ \tfrac{g_1^+ + g_1^-}{2} = l_1 + l_1^\perp $ with $l_1 \in \ker \L \,, \ l_1^\perp \in \ker^\perp \L$, there holds
    \begin{equation}\label{Od:eps-11}
      \begin{aligned}
        & \Q ( \tfrac{g_0^+ + g_0^- }{2} , \tfrac{g_1^+ + g_1^- }{2} ) + \Q ( \tfrac{g_1^+ + g_1^- }{2} , \tfrac{g_0^+ + g_0^- }{2} ) \\
        = & \Q ( \tfrac{g_0^+ + g_0^- }{2} , l_1 ) + \Q ( l_1 , \tfrac{g_0^+ + g_0^- }{2} ) + \Q ( \tfrac{g_0^+ + g_0^- }{2} , l_1^\perp ) + \Q ( l_1^\perp , \tfrac{g_0^+ + g_0^- }{2} ) \\
        = & \L ( l_1 \tfrac{g_0^+ + g_0^- }{2} ) + \Q ( \tfrac{g_0^+ + g_0^- }{2} , l_1^\perp ) + \Q ( l_1^\perp , \tfrac{g_0^+ + g_0^- }{2} ) \,.
      \end{aligned}
    \end{equation}
By direct calculation, one has
    \begin{equation*}
      \begin{aligned}
        l_1 \tfrac{g_0^+ + g_0^- }{2} = & ( \rho_0 \rho_1 + u_0 \cdot u_1 + \tfrac{3}{2} \theta_0 \theta_1 ) + ( \rho_0 u_1 + \rho_1 u_0 + \theta_0 u_1 + \theta_1 u_0 ) \cdot v \\
        & + ( \tfrac{2}{3} u_0 \cdot u_1 + 2 \theta_0 \theta_1 + \rho_0 \theta_1 + \rho_1 \theta_0 ) ( \tfrac{|v|^2}{2} - \tfrac{3}{2} ) \\
        & + ( u_0 \otimes u_1 + u_1 \otimes u_0 ) : \A(v) + ( \theta_0 u_1 + \theta_1 u_0 ) \cdot \B(v) + \theta_0 \theta_1 \C(v) \,,
      \end{aligned}
    \end{equation*}
which implies that
    \begin{equation}\label{Od:eps-12}
      \begin{aligned}
        \L ( l_1 \tfrac{g_0^+ + g_0^- }{2} ) = \L ( ( u_0 \otimes u_1 + u_1 \otimes u_0 ) : \A(v) + ( \theta_0 u_1 + \theta_1 u_0 ) \cdot \B(v) + \theta_0 \theta_1 \C(v) ) \,.
      \end{aligned}
    \end{equation}
Thus, by using the expression of $l_1^\perp$ and by plugging \eqref{g0}, \eqref{Od:eps-11}-\eqref{Od:eps-12} into \eqref{Od:eps-10} above, one sees that, for some undetermined functions $\rho_2(t,x)$, $u_2(t,x)$ and $\theta_2(t,x)$,
      \begin{align}\label{Od:eps(4)}
         \tfrac{g_2^+ + g_2^-}{2}  = & \rho_2 + u_2 \cdot v + \theta_2 ( \tfrac{|v|^2}{2} - \tfrac{3}{2} ) \no \\
           & + ( u_0 \otimes u_1 + u_1 \otimes u_0 ) : \A(v) +  ( \theta_0 u_1 + \theta_1 u_0  ) \cdot \B(v) + \theta_0 \theta_1 \C(v) \\ \no
           & - \tfrac{1}{2} \nabla_x u_1 : \widehat{\A} (v) - \tfrac{1}{2} \nabla_x \theta_1 \cdot \widehat{\B} (v) + \sum \Gamma_0^+ \Upsilon^+ (v) \,.
      \end{align}
    Here we have used the summation symbol $ \sum$, $ \Gamma_0^+  $ explicitly depends on the $u_0, n_0, \theta_0 , B_0 , E_0$ only, and $ \Upsilon^+ (v) \in \ker^\perp \L $. More precisely,
      \begin{align}\label{Od:eps(5)}
        \no \sum \Gamma_0^+ & \Upsilon^+ (v) =  - \tfrac{5}{4} \theta_0 u_0 \otimes u_0 : \widehat{\A} (v) - \tfrac{5}{4} \theta_0 \nabla_x \theta_0 \cdot \widehat{\B} (v) + \tfrac{5}{4} \theta_0 u_0 \otimes u_0 : \A (v) + \tfrac{5}{2} \theta_0^2 u_0 \cdot \B(v) \\
        \no & + \tfrac{5}{4} \theta_0^3 \C (v) + \tfrac{1}{2 \sigma} B_{0i} ( j_{0l} - n_0 u_{0l} ) \epsilon_{ijk} \L^{-1} [ \P^\perp_{\L} ( v_j \partial_{v_k} \widetilde{\Phi}_l (v) ) ] \\
        \no & + \tfrac{1}{2} u_0 \otimes u_0 \otimes u_0 : \L^{-1} ( \Q (v, \A(v)) + \Q(  \A(v) , v ) ) \\
        \no & + \theta_0 u_0 \otimes u_0 : \L^{-1} ( \Q (v, \B(v)) + \Q(  \B(v) , v ) ) - \tfrac{1}{2} \nabla_x u_0 \otimes u_0 : \L^{-1} \P^\perp_{\L} \Q( \widehat{\A} (v) , v ) \\
        \no & + \tfrac{1}{2} \theta_0^2 u_0 \cdot \L^{-1} ( \Q (v, \C(v)) + \Q(  \C(v) , v ) ) - \tfrac{1}{2} u_0 \otimes \nabla_x u_0 : \L^{-1} \P^\perp_{\L} \Q(v , \widehat{\A} (v)) \\
        \no &  - \tfrac{1}{2} \nabla_x \theta_0 \otimes u_0 : \L^{-1} \P^\perp_{\L} \Q( \widehat{\B} (v) , v ) -  \tfrac{1}{2} u_0 \otimes \nabla_x \theta_0 : \L^{-1} \P^\perp_{\L} \Q( v , \widehat{\B} (v) ) \\
         & + \tfrac{1}{4} \theta_0 u_0 \otimes u_0 : \L^{-1} ( \Q( |v|^2 , \A(v)) + \Q ( \A(v) , |v|^2 ) ) \\
        \no & + \tfrac{1}{2} \theta_0^2 u_0 \cdot \L^{-1} ( \Q( |v|^2 , \B(v)) + \Q ( \B(v) , |v|^2 ) ) \\
        \no & + \tfrac{1}{4} \theta_0^3 \L^{-1} ( \Q( |v|^2 , \C(v)) + \Q ( \C(v) , |v|^2 ) ) - \tfrac{1}{4} \nabla_x \theta_0^2 \cdot \L^{-1} \P^\perp_{\L} ( v \C(v) ) \\
        \no & - \tfrac{1}{4} \theta_0 \nabla_x u_0 : \L^{-1} ( \Q( |v|^2 , \widehat{\A}(v)) + \Q ( \widehat{\A}(v) , |v|^2 ) ) \\
        \no & - \tfrac{1}{4} \theta_0 \nabla_x \theta_0 \cdot \L^{-1} ( \Q( |v|^2 , \widehat{\B}(v)) + \Q ( \widehat{\B}(v) , |v|^2 ) )  \\
        \no & - \tfrac{1}{4} \nabla_x ( u_0 \otimes u_0 ) : \L^{-1} \P^\perp_{\L} ( v \otimes \A(v) ) - \tfrac{1}{2} \nabla_x ( \theta_0 u_0 ) : \L^{-1} \P^\perp_{\L} ( v \otimes \B(v) )  \\
        \no & + \tfrac{1}{4} \nabla_x^2 u_0 : \L^{-1} \P^\perp_{\L} ( v \otimes \widehat{\A} (v) ) + \tfrac{1}{4} \nabla_x^2 \theta_0 : \L^{-1} \P^\perp_{\L} ( v \otimes \widehat{\B} (v) ) \,,
      \end{align}
where the symbol $\L^{-1}$ represents $( \L |_{\ker^\perp \L} )^{-1} : \ker^\perp \L \rightarrow \ker^\perp \L$ (a one to one and onto map).

Furthermore, one derives from projecting \eqref{Od:eps-1}$_2$ into $\ker ( \L + \mathfrak{L} ) $, and from combining with \eqref{Od:eps-2}$_2$, \eqref{Od:eps-5}, \eqref{Od:eps-7} and \eqref{Od:eps-9}, that
    \begin{equation}\label{Limit-Equ-n0}
      \partial_t n_0 + \div_x j_0 = 0 \,.
    \end{equation}
Note that the equation \eqref{Limit-Equ-n0} is consistent with the limit system \eqref{Limit-Equ}. Indeed, it can also be derived from taking the divergence operator $\div_x$ on the $E_0$-equation in \eqref{Limit-Equ}, due to the relations $ \div_x ( \nabla_x \times B_0 ) = 0$ and $ \div_x E_0 = n_0 $.

Correspondingly, projecting equation \eqref{Od:eps-1} into $\ker^\perp ( \L + \mathfrak{L} )$ reduces to
      \begin{align}\label{Od:eps-13}
        \no ( \L + \mathfrak{L} ) ( \tfrac{g_2^+  - g_2^-}{2} )
      = & - ( \tfrac{1}{2} \nabla_x n_1 - E_1 - u_0 \times B_1 - u_1 \times B_0 + \tfrac{5}{4} \nabla_x ( n_0 \theta_0 ) + \theta_0 E_0 ) \cdot v \\
        \no & - \tfrac{1}{3} ( \div_x ( n_0 u_0 ) - 2 E_0 \cdot u_0 ) ( \tfrac{|v|^2}{2} - \tfrac{3}{2} )           \\
        \no & - (  \tfrac{1}{2} \nabla_x ( n_0 u_0 ) - ( E_0 + B_0 \times u_0 ) \otimes u_0 ) : \A(v) \\
        \no & - ( \tfrac{1}{2} \nabla_x ( n_0 \theta_0 ) + \theta_0 E_0 - n_0 u_0 \times B_0 ) \cdot \B(v) - \tfrac{1}{2} B_{0i} \partial_l \theta_0 \epsilon_{ijk} v_j \partial_{v_k} \widehat{\B}_l (v) \\
        \no & - \tfrac{1}{\sigma} \nabla_x ( j_0 - n_0 u_0 ) : \P^\perp_{\L + \mathfrak{L}} (  v \otimes \widetilde{\Phi} (v) ) - \tfrac{1}{2} B_{0i} \partial_l u_{0m} \epsilon_{ijk} v_j \partial_{v_k} \widehat{\A}_{lm} (v) \\
        & - n_0 \L ( \tfrac{g_1^+ + g_1^-}{2} ) + \tfrac{1}{2} ( \L + \mathfrak{L} ) ( \tfrac{g_1^+ + g_1^-}{2} ) + \Q( g_1^+ - g_1^- , \tfrac{g_0^+ + g_0^-}{2} ) \,,
      \end{align}
where the operator $ \P^\perp_{\L + \mathfrak{L}} = Id - \P_{\L + \mathfrak{L}} $. Recalling \eqref{g0} and \eqref{Od:1(5)} again, a simple calculation implies
    \begin{equation}\label{Od:eps-14}
      \begin{aligned}
        \Q ( g_1^+ - g_1^- , \tfrac{g_0^+ + g_0^-}{2} )
      = & ( n_1 + n_0 \theta_0 ) u_0 \cdot \Q ( 1 , v ) + \tfrac{1}{2} ( n_1 + n_0 \theta_0 ) \theta_0 \Q( 1 , |v|^2 ) \\
        & + \tfrac{n_0}{2} \L (  u_0 \otimes u_0 : \A(v) + 2 \theta_0 u_0 \cdot \B(v) + \theta_0^2 \C(v) ) \\
        & + \tfrac{1}{\sigma} \theta_0 ( j_0 - n_0 u_0 ) \cdot \Q( \widetilde{\Phi} (v) , |v|^2 - 5 ) + \tfrac{2}{\sigma} ( j_0 - n_0 u_0  ) \otimes u_0 : \Q ( \widetilde{\Phi} (v) , v ) \,.
      \end{aligned}
    \end{equation}
As a result, combined with \eqref{Od:1-4}, this yields, for some undetermined function $n_2(t,x)$,
    \begin{equation}\label{Od:eps(6)}
      \begin{aligned}
        \tfrac{g_2^+ - g_2^-}{2} = & \tfrac{1}{2} n_2 - ( \tfrac{1}{2} \nabla_x n_1 - E_1 - u_0 \times B_1 - u_1 \times B_0 ) \cdot \widetilde{\Phi}(v) + \tfrac{1}{2} u_1 \cdot v + \tfrac{1}{2} \theta_1 ( \tfrac{|v|^2}{2} - \tfrac{3}{2} ) \\
        & + n_1 u_0 \cdot ( \L + \mathfrak{L} )^{-1} \Q(1 , v) + \tfrac{1}{2} n_1 \theta_0 ( \L + \mathfrak{L} )^{-1} \Q(1 , |v|^2) + \sum \Gamma^-_0 \Upsilon^-(v) \,,
      \end{aligned}
    \end{equation}
where $ \Gamma_0^-  $ in the summation term explicitly depends on the $u_0, n_0, \theta_0 , B_0 , E_0$ only, and $ \Upsilon^- (v) \in \ker^\perp ( \L + \mathfrak{L} ) $. More precisely,
      \begin{align}\label{Od:eps(7)}
        \no \sum \Gamma_0^- & \Upsilon^- (v) = - ( \tfrac{5}{4} \nabla_x ( n_0 \theta_0 ) + \theta_0 E_0 ) \cdot \widetilde{\Phi} (v) - \tfrac{1}{3} ( \div_x ( n_0 u_0 ) - 2 E_0 \cdot u_0 ) \widetilde{\Psi} (v) \\
        \no & + \tfrac{1}{4} u_0 \otimes u_0 : \A(v) + \tfrac{1}{2} \theta_0 u_0 \cdot \B(v) + \tfrac{1}{4} \theta_0^2 \C(v) - \tfrac{1}{4} \nabla_x u_0 : \widehat{\A} (v) - \tfrac{1}{4} \nabla_x \theta_0 \cdot \widehat{\B}(v) \\
        \no & + n_0 \theta_0 u_0 \cdot  ( \L + \mathfrak{L} )^{-1} \Q(1 , v) + \tfrac{1}{2} n_0 \theta_0^2 ( \L + \mathfrak{L} )^{-1} \Q(1 , |v|^2 ) \\
        \no & - ( \tfrac{1}{2} \nabla_x ( n_0 u_0 ) - ( E_0 + B_0 \times u_0 ) \otimes u_0 ) : ( \L + \mathfrak{L} )^{-1} \A(v) \\
        \no & - ( \tfrac{1}{2} \nabla_x ( n_0 \theta_0 ) + \theta_0 E_0 - n_0 u_0 \times B_0 ) \cdot ( \L + \mathfrak{L} )^{-1} \B(v) \\
        \no & - \tfrac{1}{2} B_{0i} \partial_l \theta_0 \epsilon_{ijk} ( \L + \mathfrak{L} )^{-1} ( v_j \partial_{v_k} \widehat{\B}_l (v) ) - \tfrac{1}{2} B_{0i} \partial_{0m} \epsilon_{ijk} ( \L + \mathfrak{L} )^{-1} (  v_j \partial_{v_k} \widehat{\A}_{lm} (v) ) \\
        & - \tfrac{1}{\sigma} \nabla_x ( j_0 - n_0 u_0  ) : ( \L + \mathfrak{L} )^{-1} \P^\perp_{\L + \mathfrak{L}} ( v \otimes \widetilde{\Phi} (v) ) \\
        \no & + \tfrac{1}{\sigma} \theta_0 ( j_0 - n_0 u_0 ) \cdot ( \L + \mathfrak{L} )^{-1} \Q ( \widetilde{\Phi} (v) , |v|^2 - 5 ) + \tfrac{n_0}{2} \nabla_x u_0 : ( \L + \mathfrak{L} )^{-1} \A(v) \\
        \no & + \tfrac{2}{\sigma} ( j_0 - n_0 u_0 ) \otimes u_0 : ( \L + \mathfrak{L} )^{-1} \Q ( \widetilde{\Phi} (v) , v ) + \tfrac{n_0}{2} \nabla_x \theta_0 \cdot ( \L + \mathfrak{L} )^{-1} \B(v) \,,
      \end{align}
where the symbol $ ( \L + \mathfrak{L} )^{-1} $ represents $ \Big( ( \L + \mathfrak{L} )|_{\ker^\perp (\L + \mathfrak{L})} \Big)^{-1} : \ker^\perp (\L + \mathfrak{L}) \rightarrow \ker^\perp (\L + \mathfrak{L}) $.

By projecting the equality \eqref{Od:eps(6)} into $\ker \L$, there holds
      \begin{equation}\label{Od:eps(8)}
        \left\{
          \begin{array}{l}
            n_2 = \langle g_2^+ - g_2^- \rangle \,, \qquad \theta_2^+ - \theta_2^- = \theta_1 + n_1 ( \overline{V} \cdot u_0 + C \theta_0 ) + \sum C_{\Upsilon^-} \Gamma_0^- \,, \\
            u_2^+ - u_2^- = n_1 ( u_0 \cdot \mathbf{M} + \theta_0 V ) + u_1 + \sigma ( - \tfrac{1}{2} \nabla_x n_1 + E_1 + u_0 \times B_1 + u_1 \times B_0 ) + \sum \Gamma_0^- U_{\Upsilon^-} \,,
          \end{array}
        \right.
      \end{equation}
with constant value quantities $ \mathbf{M} \in \R^{3 \times 3} $, $ V , \overline{V} , U_{\Upsilon^-} \in \R^3 $ and $ C , C_{\Upsilon^-} \in \R $, given as
      \begin{equation*}
      \begin{array}{rll}
        \mathbf{M} = & 2 \int_{\R^3} ( \L + \mathfrak{L} )^{-1} \Q(1 ,v) \otimes v M \d v \,, \quad \quad \quad \ & V =  \int_{\R^3} ( \L + \mathfrak{L} )^{-1} \Q(1 , |v|^2) v M \d v \,,
      \\[2mm]
        \overline{V} = & 2 \int_{\R^3} ( \L + \mathfrak{L} )^{-1} \Q(1 , v) ( \tfrac{|v|^2}{3} - 1 ) M \d v \,, \quad \ & U_{\Upsilon^-} =  2 \int_{\R^3} \Upsilon^- (v) v M \d v \,,
      \\[2mm]
        C = & \int_{\R^3} ( \L + \mathfrak{L} )^{-1} \Q( 1 , |v|^2 ) ( \tfrac{|v|^2}{3} - 1 ) M \d v \,, \quad & C_{\Upsilon^-} =  2 \int_{\R^3} \Upsilon^- (v) ( \tfrac{|v|^2}{3} - 1 ) M \d v \,.
      \end{array}
      \end{equation*}

In summary, the two-fluid incompressible NSFM system with Ohm's law is first derived in \eqref{Limit-Equ}, which completely determines the leading order terms $g_0^\pm$ in \eqref{g0} and $(E_0, B_0)$. The expressions of $\frac{g_1^+ + g_1^-}{2}$ and $\frac{g_1^+ - g_1^-}{2}$ are given, respectively, in \eqref{Od:1(4)} and \eqref{Od:1(5)}, with undetermined functions $(\rho_1, u_1, \theta_1, n_1) (t,x)$.
The relations \eqref{Od:eps(3)} cannot completely determine $(\rho_1, u_1, \theta_1 )$. However, once a pair $(\rho_1, u_1, \theta_1)$ is fixed, the functions $(E_1, B_1, n_1)$ can be determined by \eqref{Od:eps(2)} and \eqref{Od:eps(8)}, namely,
      \begin{align}\label{Od:eps(10)}
      	\begin{cases}
      		\partial_t E_1 - \nabla_x \times B_1 = - ( u_2^+ - u_2^- ) \,, \quad \div_x E_1 = n_1 \,, \\
      		\partial_t B_1 + \nabla_x \times E_1 = 0 \,, \hspace*{2.2cm} \div_x B_1 = 0 \,, \\
      		u_2^+ - u_2^- = n_1 ( u_0 \cdot \mathbf{M} + \theta_0 V ) + u_1  \\
      		\hspace*{2cm} + \sigma ( - \tfrac{1}{2} \nabla_x n_1 + E_1 + u_0 \times B_1 + u_1 \times B_0 ) + \sum \Gamma_0^- U_{\Upsilon^-} \,,
      	\end{cases}
      \end{align}
Finally, $\frac{g_2^+ + g_2^-}{2}$ and $\frac{g_2^+ - g_2^-}{2}$ are, respectively, formally expressed in \eqref{Od:eps(4)} and \eqref{Od:eps(6)} with undetermined $(\rho_2, u_2, \theta_2, n_2) (t,x)$ and $(\rho_1, u_1, \theta_1, n_1) (t,x)$.

\subsection{Truncation of Hilbert expansion} 
\label{sub:truncation_expansion}


When using the Hilbert expansion approach, one hopes to expand the number of terms as small as possible, if the limit equation had been derived. It is possible now to truncate the expansion ansatz \eqref{Ansatz-1} from the first order terms $(g_1^\pm, E_1, B_1)$ with order $\eps$.

At first, noticing \eqref{Od:eps(3)}, one chooses a special function pair $(\rho_1, u_1, \theta_1)$ solving
      \begin{equation}\label{Spec-Chos*}
        \left\{
          \begin{array}{rl}
            \Delta_x \rho_1 (t,x) = & \tfrac{1}{6} \Delta_x |u_0|^2 - \tfrac{1}{2} \div_x ( u_0 \cdot \nabla_x u_0 - \tfrac{1}{2} j_0 \times B_0) \,,  \\
            \Delta_x \phi (t,x) = & \partial_t \theta_0 (t,x) \,, \qquad \int_{\mathbb{T}^3} \phi \d x = \int_{\T} \rho_1 \d x = 0 \,,\\
            u_1 (t , x) = & \nabla_x \phi (t , x) \,, \qquad \rho_1 = \theta_1 \,.
          \end{array}
        \right.
      \end{equation}
This indicates, $(E_1, B_1, n_1)$ can be uniquely determined by \eqref{Od:eps(10)}.

Furthermore, the undetermined $(\rho_2, u_2, \theta_2, n_2)$ in \eqref{Od:eps(4)} and \eqref{Od:eps(6)} are taken as
      \begin{equation*}
      	\begin{aligned}
      		\rho_2 = \theta_2 = n_2 = 0 \,, \ u_2 = 0 \,.
      	\end{aligned}
      \end{equation*}
Consequently, the terms $g_1^\pm$ and $g_2^\pm$ are determined by
      \begin{align}\label{g1-g2}
          \begin{cases}
            \tfrac{ {g}_1^+ + {g}_1^-}{2} = {\rho}_1 + {u}_1 \cdot v + {\theta}_1 ( \tfrac{|v|^2}{2} - \tfrac{3}{2} ) +  \tfrac{1}{2} u_0 \otimes u_0 : \A(v) + \theta_0 u_0 \cdot \B(v) + \tfrac{1}{2} \theta_0^2 \C(v) \\
          \qquad\qquad  - \tfrac{1}{2} \nabla_x u_0 : \widehat{\A} (v) - \tfrac{1}{2} \nabla_x \theta_0 \cdot \widehat{\B} (v) \,, \\
          \tfrac{ {g}_1^+ - {g}_1^-}{2} = \tfrac{1}{2} {n}_1 + \tfrac{1}{2} n_0 u_0 \cdot v + \tfrac{1}{2} n_0 \theta_0 ( \tfrac{|v|^2}{2} - \tfrac{3}{2} ) + ( - \tfrac{1}{2} \nabla_x n_0 + E_0 + u_0 \times B_0 ) \cdot \widetilde{\Phi} (v) \,, \\
          \tfrac{ {g}_2^+ + {g}_2^-}{2}  =  ( u_0 \otimes {u}_1 + {u}_1 \otimes u_0 ) : \A(v) +  ( \theta_0 {u}_1 + {\theta}_1 u_0  ) \cdot \B(v) + \theta_0 {\theta}_1 \C(v) \\
           \qquad \qquad - \tfrac{1}{2} \nabla_x {u}_1 : \widehat{\A} (v) - \tfrac{1}{2} \nabla_x {\theta}_1 \cdot \widehat{\B} (v) + \sum \Gamma_0^+ \Upsilon^+ (v) \,, \\
           \tfrac{ {g}_2^+ - {g}_2^-}{2} = - ( \tfrac{1}{2} \nabla_x {n}_1 - {E}_1 - u_0 \times {B}_1 - {u}_1 \times B_0 ) \cdot \widetilde{\Phi}(v) + \tfrac{1}{2} {u}_1 \cdot v + \tfrac{1}{2} {\theta}_1 ( \tfrac{|v|^2}{2} - \tfrac{3}{2} ) \\
        \qquad \qquad + {n}_1 u_0 \cdot ( \L + \mathfrak{L} )^{-1} \Q(1 , v) + \tfrac{1}{2} {n}_1 \theta_0 ( \L + \mathfrak{L} )^{-1} \Q(1 , |v|^2) + \sum \Gamma^-_0 \Upsilon^-(v) \,,
          \end{cases}
      \end{align}
with the functions $ ({\rho}_1$, $ {u}_1$, $ {\theta}_1)$ and $(E_1, B_1, n_1)$ given in \eqref{Od:eps(10)}-\eqref{Spec-Chos*}.

As a result, we write the Hilbert expansion ansatz as the following truncated form:
      \begin{equation}\label{Ansatz}
      	g_{\eps}^{\pm} = g_0^{\pm} + \eps {g}_1^{\pm} + \eps^2 {g}_2^{\pm} + \eps g_{R , \eps}^{\pm} \,, \quad
      	E_\eps = E_0 + \eps {E}_1 + \eps E_{R, \eps} \,, \quad
      	B_\eps = B_0 + \eps {B}_1 + \eps B_{R, \eps} \,.
      \end{equation}
Plugging this ansatz \eqref{Ansatz} into \eqref{VMB-g}, and combining the order analysis in \S \ref{Subsec:Order} above, the remainder system \eqref{Remd-Equ-Bref} obeyed by the remainders $(g_{R,\eps}^\pm, E_{R,\eps}, B_{R,\eps})$ can be obtained, as follows:
\begin{equation} 
  \left\{
    \begin{array}{rll}
      \eps \partial_t G_{R,\eps} + v \cdot \nabla_x G_{R,\eps} + \mathcal{T} (v \times B_0) \cdot \nabla_v G_{R,\eps} + & \mathcal{T} (v \times B_{R,\eps}) \cdot \nabla_v G_0 &
    \\[2mm]
      - E_{R,\eps} \cdot v \mathcal{T}_1 + \frac{1}{\eps} \mathbb{L} G_{R,\eps} = & \eps H_{R,\eps} \,,
    \\[2mm]
      \partial_t E_{R,\eps} - \nabla_x \times B_{R,\eps} = & - \frac{1}{\eps} \left\langle G_{R,\eps} \cdot \mathcal{T}_1 v \right\rangle \,, \
      	& \div_x E_{R,\eps} = \left\langle G_{R,\eps} \cdot \mathcal{T}_1 \right\rangle \,,
    \\[2mm]
      \partial_t B_{R,\eps} + \nabla_x \times E_{R,\eps} = & 0 \,, \ & \div_x B_{R,\eps} = 0 \,,
    \end{array}
  \right.
\end{equation}
with $\mathcal{T}_1 = (1, -1)^\top$ and $\mathcal{T} = diag(1,-1)$. Here
\begin{align}\label{Remnd-H}
    H_{R,\eps} = & - \eps \mathcal{T} {E}_1 \cdot \nabla_v G_{R,\eps} + \eps \mathcal{T} {E}_1 \cdot v G_{R,\eps} - \eps \mathcal{T} E_{R,\eps} \cdot \nabla_v G_{R,\eps} + \eps \mathcal{T} E_{R,\eps} \cdot v G_{R,\eps} \no \\
    & -  \mathcal{T} E_0 \cdot \nabla_v G_{R,\eps} +  \mathcal{T} E_0 \cdot v G_{R,\eps} - \mathcal{T} (v \times {B}_1) \cdot \nabla_v G_{R,\eps} - \mathcal{T} (v \times B_{R,\eps}) \cdot \nabla_v G_{R,\eps} \no \\
    & - \mathcal{T} E_{R,\eps} \cdot \nabla_v ( G_0 + \eps G_1 + \eps^2 G_2 ) + \mathcal{T} E_{R,\eps} \cdot v ( G_0 +  \eps G_1 + \eps^2 G_2 ) \no \\
    & + \frac{1}{\eps} \Gamma_0 G_{R,\eps} - \mathcal{T} (v \times B_{R,\eps}) \cdot \nabla_v ( G_1 + \eps G_2 ) + \eps
        \left(
          \begin{array}{c}
            \Q ( {g}_2^+ , g_{R, \eps}^+ + g_{R , \eps}^- ) \\
            \Q ( {g}_2^- , g_{R, \eps}^+ + g_{R, \eps}^-  )
          \end{array}
        \right) \no \\
        & +
        \left(
          \begin{array}{c}
            \Q ( {g}_1^+ , g_{R, \eps}^+ + g_{R , \eps}^- ) + \Q ( g_{R , \eps}^+ , {g}_1^+ + {g}_1^- ) \\
            \Q ( {g}_1^- , g_{R, \eps}^+ + g_{R, \eps}^- ) + \Q ( g_{R , \eps}^- , {g}_1^+ + {g}_1^- )
          \end{array}
        \right)
        + \eps
        \left(
          \begin{array}{c}
            \Q (g_{R, \eps}^+, {g}_2^+ + {g}_2^- ) \\
            \Q (g_{R, \eps}^-, {g}_2^+ + {g}_2^- )
          \end{array}
        \right) \\ \no
        & + \left(
          \begin{array}{c}
            \Q ( g_{R, \eps}^+ , g_{R , \eps}^+ + g_{R , \eps}^- ) \\
            \Q ( g_{R, \eps}^- , g_{R , \eps}^+ + g_{g , \eps}^- )
          \end{array}
        \right)
         +
        \left(
          \begin{array}{c}
            \mathcal{R}^+ \\
            \mathcal{R}^-
          \end{array}
        \right) \,,
\end{align}
where
\begin{equation}\label{Remd-Gamma0}
  \Gamma_0 G_{R,\eps} = \left(
          \begin{array}{c}
            \Q ( g_0^+ , g_{R, \eps}^+ + g_{R, \eps}^- ) + \Q ( g_{R , \eps}^+ , g_0^+ + g_0^- ) \\
            \Q ( g_0^- , g_{R, \eps}^- + g_{R , \eps}^- ) + \Q ( g_{R , \eps}^- , g_0^+ + g_0^- )
          \end{array}
        \right) \,,
\end{equation}
and the symbols $\mathcal{R}^\pm$ are defined as
	\begin{align}\label{Remd-Terms}
		\no \left(
		\begin{array}{c}
			\mathcal{R}^+ \\
			\mathcal{R}^-
		\end{array}
		\right)
		= & - \partial_t
		\left(
		\begin{array}{c}
			{g}_1^+ + \eps {g}_2^+ \\
			{g}_1^- + \eps {g}_2^-
		\end{array}
		\right)
		- v \cdot \nabla_x
		\left(
		\begin{array}{c}
			{g}_2^+ \\
			{g}_2^-
		\end{array}
		\right)
		- {E}_1 \cdot \nabla_v
		\left(
		\begin{array}{c}
			g_0^+ + \eps {g}_1^+ + \eps^2 {g}_2^+ \\
			- ( g_0^- + \eps {g}_1^- + \eps^2 {g}_2^-
		\end{array}
		\right) \\
		\no & - E_0 \cdot \nabla_v
		\left(
		\begin{array}{c}
			{g}_1^+ + \eps {g}_2^+ \\
			- ( {g}_1^- + \eps {g}_2^- )
		\end{array}
		\right)
		- ( v \times B_0 ) \cdot \nabla_v
		\left(
		\begin{array}{c}
			{g}_2^+ \\
			- {g}_2^-
		\end{array}
		\right)\\
		\no & + ( v \times {B}_1 ) \cdot \nabla_v
		\left(
		\begin{array}{c}
			{g}_1^+ + \eps {g}_2^+ \\
			- ( {g}_1^- + \eps {g}_2^- )
		\end{array}
		\right)
		+ {E}_1 \cdot v
		\left(
		\begin{array}{c}
			g_0^+ + \eps {g}_1^+ + \eps^2 {g}_2^+ \\
			- ( g_0^- + \eps {g}_1^- + \eps^2 {g}_2^- )
		\end{array}
		\right) \\
		\no & + E_0 \cdot v
		\left(
		\begin{array}{c}
			{g}_1^+ + \eps {g}_2^+ \\
			- ( {g}_1^- + \eps {g}_2^- )
		\end{array}
		\right)
		+
		\left(
		\begin{array}{c}
			\Q ( g_0^+ , {g}_2^+ + {g}_2^- ) + \Q ( {g}_2^+ , g_0^+ + g_0^- ) \\
			\Q ( g_0^- , {g}_2^+ + {g}_2^- ) + \Q ( {g}_2^- , g_0^+ + g_0^- )
		\end{array}
		\right) \\
		& +
		\left(
		\begin{array}{c}
			\Q ( {g}_1^+ + \eps {g}_2^+ , {g}_1^+ + \eps {g}_2^+ ) + \Q ( {g}_1^+ + \eps {g}_2^+ , {g}_1^- + \eps {g}_2^- ) \\
			\Q ( {g}_1^- + \eps {g}_2^- , {g}_1^- + \eps {g}_2^- ) + \Q ( {g}_1^- + \eps {g}_2^-, {g}_1^+ + \eps {g}_2^+ )
		\end{array}
		\right) \,.
	\end{align}

\section{Preliminaries}\label{Sec:Prem}

This section is devoted to collect some useful estimates for deriving the uniform bounds of the remainder equations.

\subsection{Estimates on collision operators}

\begin{lemma}\label{Lm-L-Propty}
  For any $G, H \in L^2_v \cap L^2_v(\nu)$, we have
  \begin{equation}\label{L-Proty-1}
    \langle \mathbb{L} G \cdot H \rangle = \langle G \cdot \mathbb{L} H \rangle
  \end{equation}
  and $\mathbb{L} G = 0$ if and only if $G = \mathbb{P} G$. Moreover, there is $\Lambda > 0$ such that
  \begin{equation}\label{L-Proty-2}
    \langle \mathbb{L} G \cdot G \rangle \geq \Lambda \| \mathbb{P}^\perp G \|^2_{L^2_v(\nu)} \,.
  \end{equation}
\end{lemma}


\begin{lemma}\label{Lm-nu-Propty}
  The collision frequency $\nu(v)$ defined in \eqref{Clisn-Frqc} is smooth, and for every $v \in \mathbb{R}^3$
  \begin{align}\label{nu-Propty-1}
    \nu(v) \sim ( 1 + |v| )^\gamma.
  \end{align}
  For any multi-index $\beta \neq 0$, we have
  \begin{align}\label{nu-Propty-2}
    \sup_{v \in \mathbb{R}^3} | \partial_\beta \nu(v) | < \infty \,.
  \end{align}
  Furthermore, if the velocities $v, v_*, v' , v_*' \in \mathbb{R}^3$ satisfy
  \begin{align}\label{Clisn-Law}
    v + v_* = v' + v_*' \,, \quad |v|^2 + |v_*|^2 = |v'|^2 + |v_*'|^2 \,,
  \end{align}
  then it holds, for some generic positive constant $C > 0$, that
  \begin{align}\label{nu-Propty-3}
    \nu(v) + \nu(v_*) \leq C \left( \nu(v') + \nu(v_*') \right).
  \end{align}
\end{lemma}


The detailed proof of Lemma \ref{Lm-L-Propty} is referred to \cite{Arsenio-SaintRaymond}. The proof of \eqref{nu-Propty-1} in Lemma \ref{Lm-nu-Propty} can be found in \cite{LS-2010-KRM}. The inequalities \eqref{nu-Propty-2} and \eqref{nu-Propty-3} have been justified in \cite{Guo-Inventions2003}, for instance. So, we omit the proofs here.

Next, some mixed higher-order derivative estimates on the collision operators $\mathbb{L}$ and $\mathcal{Q}$ are also gathered here. The proof can also be found in \cite{Guo-Inventions2003,Guo-CPAM-2006}, for instance.

\begin{lemma}\label{Lm-L-Q}
  For the hard potential with $\gamma \in [0,1]$ and weighed function $\w (v) = \sqrt{1 + |v|^2}$, there exist $C_{|\beta|}, C > 0$, such that
  \begin{equation}\label{Coercive-L}
    \left\langle \w^{2 l} \partial^m_\beta \mathbb{L} G \cdot \partial^m_\beta G \right\rangle \geq \tfrac{1}{2} \| \w^l \partial^m_\beta G \|^2_{L^2_v(\nu)} - C_{|\beta|} \| \partial^m G \|^2_{L^2_v(\nu)} \,,
  \end{equation}
  \begin{equation}\label{Bnd-Q-1}
    \begin{aligned}
      \left\langle \w^{2 l} \partial^m_\beta \mathcal{Q} (g_1, g_2) \cdot \partial^m_\beta g_3 \right\rangle \leq & C \Big( \| \w^l \partial^{m_1}_{\beta_1} g_1 \|_{L^2_v(\nu)} \| \w^l \partial^{m_2}_{\beta_2} g_2 \|_{L^2_v} \\
      & + \| \w^l \partial^{m_1}_{\beta_1} g_1 \|_{L^2_v} \| \w^l \partial^{m_2}_{\beta_2} g_2 \|_{L^2_v(\nu)} \Big) \| \w^l \partial^m_\beta g_3 \|_{L^2_v(\nu)} \,,
    \end{aligned}
  \end{equation}
  and
  \begin{equation}\label{Bnd-Q-2}
    \begin{aligned}
      \| \partial^m_\beta \mathcal{Q} (g_1 , g_2) \|_{L^2_v} \leq & C \Big( \| \partial^{m_1}_{\beta_1} g_1 \|_{L^2_v} \| \w^{2 \gamma} \partial^{m_2}_{\beta_2} g_2 \|_{L^2_v}  + \| \w^{2 \gamma} \partial^{m_1}_{\beta_1} g_1 \|_{L^2_v} \| \partial^{m_2}_{\beta_2} g_2 \|_{L^2_v} \Big) \,,
    \end{aligned}
  \end{equation}
  where $l \geq 0$ and the summation is for $|m|+|\beta| \leq N$ with $\beta_1 + \beta_2 \leq \beta$ and $m_1 + m_2 \leq m$ componentwise.
\end{lemma}


\subsection{Energy bounds for NSFM system and linear Maxwell-type system}
\label{sec:estimates_for_nsmf_and_linear_maxwell}

In the remainder system \eqref{Remd-Equ-Bref}, the terms $g_i^\pm$ ($i = 0, 1, 2$), playing actually a role of parameter functions, are determined by fluid variables $(u_0, \theta_0, E_0, B_0)$ and $(\rho_1, u_1, \theta_1, n_1, E_1, B_1)$. These fluid variables are estimated here, such that all of $g_i^\pm$ ($i = 0, 1, 2$) can be bounded, which is helpful in deriving the uniform bounds for the remainder system with respect to the small Knudsen number $\eps > 0$. Note that $(u_0, \theta_0, E_0, B_0)$ solves the NSFM equations \eqref{NSFM} (or say, \eqref{Limit-Equ}) with initial data \eqref{IN-NSFM}, and $(\rho_1, u_1, \theta_1, n_1, E_1, B_1)$ satisfies the coupled system of linear Maxwell equations \eqref{Od:eps(10)} and a elliptic system \eqref{Spec-Chos*}, equipped with initial data \eqref{IC-LM}.

First, we consider the NSFM system \eqref{Limit-Equ} with initial data \eqref{IN-NSFM}. For any fixed integer $s\ge 2$, we define the energy and dissipative rate functionals $\mathcal{E}_{0,s}(t)$ and $\mathcal{D}_{0,s}(t)$, respectively, as follows,
	\begin{align}\label{E0s}
		\mathcal{E}_{0,s}(t) =\ & \|u_0\|^2_{H^s_x} + \|\theta_0\|^2_{H^s_x}
			+ \tfrac{3}{2}\|E_0\|^2_{H^s_x} + \tfrac{5}{4} \|n_0\|^2_{H^s_x}
			+ (\tfrac{3}{2} {- \delta + \delta \sigma })\|B_0\|^2_{H^s_x}
			\no\\
			& + (1- \delta)\|\p_t B_0\|^2_{H^s_x} + \|\nabla_x B_0\|^2_{H^s_x} + \delta \|\p_t B_0 + B_0\|^2_{H^s_x}, \\
		\label{D0s} \mathcal{D}_{0,s}(t) =\ & \mu\|\nabla_x u_0\|^2_{H^s_x} + \tfrac{\kappa}{2} \|\nabla_x \theta_0 \|^2_{H^s_x}
			+ \sigma \|E_0\|^2_{H^s_x}
			+ \tfrac{3}{2}\sigma \|n_0\|^2_{H^s_x} + \tfrac{1}{2}\sigma \|\nabla_x n_0\|^2_{H^s_x}
			+ \delta \|\nabla_x B_0\|^2_{H^s_x} \no\\
			& + (\sigma- \delta)\|\p_t B_0\|^2_{H^s_x}
			+ \tfrac{1}{2}\sigma \sum_{m \le s} \| - \tfrac{1}{2} \nabla_x \p^m_x n_0 + \p^m_x E_0 + (\p^m_x u_0) \times B_0\|^2_{L^2_x} \,,
	\end{align}
where $\delta > 0$ is a sufficiently small constant. Moreover, the initial energy is denoted by
	\begin{align}\label{E0s-in}
	  \mathcal{E}_{0,s}^{\IN} = \|u_0^{\IN}\|^2_{H^s_x} + \|\theta_0^{\IN}\|^2_{H^s_x}
			+ \|E_0^{\IN}\|^2_{H^s_x} + \|\div_x E_0^{\IN}\|^2_{H^s_x} + \|\nabla_x \times E_0^{\IN}\|^2_{H^s_x} + \|B_0^{\IN}\|^2_{H^s_x}.
	\end{align}
It is an obvious fact that $\mathcal{E}_{0,s}(0) \le C \mathcal{E}_{0,s}^{\IN}$ for some positive constant $C=C(\mu,\kappa,\sigma) >0$. Besides, we point out that for $2 \le s_1 \le s_2$, there hold the relations $\mathcal{E}_{0,s_1} \le \mathcal{E}_{0,s_2}$, $\mathcal{D}_{0,s_1} \le \mathcal{D}_{0,s_2}$, and $\mathcal{E}_{0,s_1}^{\IN} \le \mathcal{E}_{0,s_2}^{\IN}$.

\begin{lemma}[\cite{JL-CMS-2018}]\label{lemm:bnd-NSMF}
	Assume that the initial data $(u_0^{\IN}, \theta_0^{\IN}, E_0^{\IN}, B_0^{\IN}) \in H^s_x \times H^s_x \times H^{s+1}_x \times H^{s+1}_x$ with $s\ge 2$. There is a small positive constant $\lambda_0(s)= \lambda(s, \mu, \kappa, \sigma)$ depending only on $s$, the viscous coefficient $\mu$, the heat conductivity coefficient $\kappa$, and the electric resistivity $\sigma$, such that if the initial energy satisfies $\mathcal{E}_{0,s}^{\IN} \le \lambda_0(s)$. Then the NSFM system \eqref{Limit-Equ}-\eqref{IN-NSFM} admits a unique global-in-time solution $(u_0, \theta_0, E_0, B_0)$ satisfying
		\begin{align}
		  & u_0, \theta_0 \in L^\infty(\mathbb{R}^+, H^s_x) \cap L^2(\mathbb{R}^+, \dot H^{s+1}_x) \\\no
		  & E_0 \in L^\infty(\mathbb{R}^+, H^s_x), \quad n_0 (= \div_x E_0) \in L^\infty(\mathbb{R}^+, H^s_x) \cap L^2(\mathbb{R}^+, \dot H^{s+1}_x) \\\no
		  & B_0 \in L^\infty(\mathbb{R}^+, H^{s+1}_x), \quad \p_t B_0 (= - \nabla_x \times B_0) \in L^\infty(\mathbb{R}^+, H^s_x).
		\end{align}
	Moreover, there holds
		\begin{align}\label{esm:E-0s}
		  \tfrac{\d}{\d t} \mathcal{E}_{0,s}(t) + \mathcal{D}_{0,s}(t) \le 0 \quad \forall \ t \geq 0 \,.
		\end{align}
	Consequently, there exists a constant $C=C(\mu,\kappa,\sigma)>0$ such that
		\begin{align}\label{Bnd-NSFM}
			\sup_{t\ge 0} (\|u_0\|^2_{H^s_x} + & \|\theta_0\|^2_{H^s_x} + \|E_0\|^2_{H^s_x} + \|n_0\|^2_{H^s_x} + \|B_0\|^2_{H^s_x} + \|\p_t B_0 \|^2_{H^s_x}) \\ \no
			& + \int_0^\infty (\mu \|\nabla_x u_0\|^2_{H^s_x} + \kappa \|\nabla_x \theta_0\|^2_{H^s_x} + \sigma \|\nabla_x n_0 \|^2_{H^s_x}) \d t
			\le C \mathcal{E}_{0,s}^{\IN}.
		\end{align}
\end{lemma}

Next, one turns to consider the linear Maxwell-type system \eqref{Od:eps(10)}-\eqref{Spec-Chos*} with the initial conditions \eqref{IC-LM} and \eqref{IC-4}. The condition \eqref{IC-4} implies that $\int_{\T} B_1 \d x = 0$. As a result, the following Poincar\'e inequality holds:
	\begin{align}
	  \| B_1\|_{L^2_x} \le C \|\nabla_x B_1 \|_{L^2_x} \,,
	\end{align}
which will be used in deriving the energy bounds of the linear Maxwell-type equations \eqref{Od:eps(10)}-\eqref{Spec-Chos*}.

For the system \eqref{Od:eps(10)}-\eqref{Spec-Chos*}, introduce the following energy functional $\mathcal{E}_{1,M} (t)$ and dissipative rate $\mathcal{D}_{1,M} (t)$, for any integer $M\ge 1$, that
	\begin{align}\label{E1M}
		\mathcal{E}_{1,M}(t) = & \| E_1 \|^2_{H^M_x} + \| n_1 \|^2_{H^M_x}
			+ (1 - \delta + \delta \sigma) \| B_1 \|^2_{H^M_x} + \|\nabla_x B_1 \|^2_{H^M_x} \no \\
			& + (1- \delta)\| \p_t B_1 \|^2_{H^M_x} + \delta \| \p_t B_1 + B_1 \|^2_{H^M_x} \,, \\
		\label{D1M}\mathcal{D}_{1,M}(t) = & \tfrac{\sigma}{2} \| E_1 \|^2_{H^M_x}
			+ \tfrac{3}{4}\sigma \| n_1 \|^2_{H^M_x} + \tfrac{1}{4} \sigma \| \nabla_x n_1 \|^2_{H^M_x}
			+ \tfrac{\delta}{2} \|\nabla_x B_1 \|^2_{H^M_x} \no \\
			& + \tfrac{(\sigma- \delta)}{2} \| \p_t B_1 \|^2_{H^M_x} + \tfrac{1}{2} \| \nabla_x u_1 \|^2_{H^M_x} \,,
	\end{align}
with the constant $\delta = \tfrac{1}{2} \min\{1, \sigma\} \in (0, \tfrac{1}{2}]$ such that all coefficients are positive. Furthermore, the initial energy is defined as
	\begin{align}\label{E1M-in}
	  \mathcal{E}_{1,M}^{\IN} = \| E_1^{\IN} \|^2_{H^M_x} + \| \div_x E_1^{\IN} \|^2_{H^M_x} + \| \nabla_x \times E_1^{\IN} \|^2_{H^M_x} + \| B_1^{\IN} \|^2_{H^M_x} \,.
	\end{align}
It is easy to see $\mathcal{E}_{1,M}(0) \le C \mathcal{E}_{1,M}^{\IN}$.

\begin{lemma}\label{lemm:bnd-linearMaxwl}
	Let $ E_1^{\IN}, B_1^{\IN} \in H_x^{M+1} $ with $M\ge 1$. Assume that there exists some small constant $\lambda_1(M+2) = \lambda_1(M, \mu,\kappa,\sigma) \in (0, \lambda_0(M+2)]$, such that $\mathcal{E}_{0,M+1}^{\IN} \le \lambda_1(M+2)$.
	Then smooth solution $( \rho_1, u_1, \theta_1, E_1, B_1)$ to the linear Maxwell-type equation \eqref{Od:eps(10)}-\eqref{Spec-Chos*} obey the following bounds
		\begin{align}
		  & \label{bnd:Spec-u} \| u_1 \|^2_{H^{M+1}_x} (t) \le C ( 1 + \mathcal{E}_{0, M+1} (t) ) \mathcal{D}_{0,M+1} (t) \,, \\
		  & \label{bnd:Spec-rho} \| \rho_1 \|^2_{H^{M+1}_x} (t) = \| \theta_1 \|^2_{H^{M+1}_x} (t) \le C \mathcal{E}_{0, M+1} (t) \big( 1 + \mathcal{E}_{0, M+1} (t) \big) \mathcal{D}_{0,M+1} (t) \,,
		\end{align}
    and
		\begin{align}
			& \label{bnd:Spec-energy-1} \tfrac{\d}{\d t} [ \mathcal{E}_{1, M} (t) + \widetilde C_M \mathcal{E}_{0, M+2} (t) ] + \left[ \mathcal{D}_{1, M} (t) + \mathcal{D}_{0, M+2} (t) \right] \le 0 \,, \ \forall \, t \geq 0 \,,
		\end{align}
	where the constants $C=C(\mu, \kappa, \sigma)>0$ and $\widetilde C_M = \widetilde C_M(M, \mu, \kappa, \sigma) >1$. Moreover, for any $t \ge 0$, the following energy bounds holds:
		\begin{align}\label{bnd:Spec-energy-2}
			\big( \| E_1 \|^2_{H^M_x} + \| \div_x E_1 \|^2_{H^M_x} + \| \nabla_x \times E_1 \|^2_{H^M_x} + \| B_1 \|^2_{H^M_x} \big) (t) \le C \big(\mathcal{E}_{0,M+2}^{\IN} + \mathcal{E}_{1,M}^{\IN} \big).
		\end{align}
\end{lemma}

The proof of Lemma \ref{lemm:bnd-linearMaxwl} will be given in Appendix \ref{Appendix-A}.


\section{Uniform Bounds for Remainder System: Proof of Main Theorem} \label{Sec:UEB}

In this section, one aims at establishing a uniform energy bound on the remainder system \eqref{Remd-Equ-Bref}, as a crucial part in proving our main result (Theorem \ref{Main-Thm}). For sake of exposition, we rewrite \eqref{Remd-Equ-Bref} by dropping the subscript $\eps$, namely,
\begin{equation}\label{Remd-Equ-GR}
  \left\{
    \begin{array}{rll}
      \eps \partial_t G_{R} + v \cdot \nabla_x G_{R} + \mathcal{T} (v \times B_0) \cdot \nabla_v G_{R} + & \mathcal{T} (v \times B_{R}) \cdot \nabla_v G_0  &
    \\[2mm]
       - E_{R} \cdot v \mathcal{T}_1 + \frac{1}{\eps} \mathbb{L} G_{R} = & \eps H_{R} \,, &
    \\[2mm]
      \partial_t E_{R} - \nabla_x \times B_{R} = & - \frac{1}{\eps} \left\langle G_{R} \cdot \mathcal{T}_1 v \right\rangle \,,
      	& \div_x E_{R} = \left\langle G_{R} \cdot \mathcal{T}_1 \right\rangle \,,
    \\[2mm]
      \partial_t B_{R} + \nabla_x \times E_{R} = & 0 \,,
      	& \div_x B_{R} = 0 \,.
    \end{array}
  \right.
\end{equation}
The proof on the uniform bounds of the remainder system \eqref{Remd-Equ-GR} will mainly be divided into two steps: 1) The pure spatial derivative estimates; 2) The mixed $(x,v)$-derivative estimates. By inspired by \cite{Guo-CPAM-2006} related to the Boltzmann equation, one first employs the so-called micro-macro decomposition approach to establish the pure spatial derivative estimates of the kinetic function $G_R (t,x,v)$. Then one tries to seek a structure of Ohm's law from the kinetic equation, which can be used to study the decay mechanism of the Maxwell system on $(E_R, B_R) (t,x)$. The corresponding expected energy inequalities thereby are constructed. At the end, the mixed $(x,v)$-derivative estimates shall be derived in order to close the energy estimates.

\subsection{Pure spatial derivative estimates}

This subsection is to estimate the pure spatial derivatives of $(G_R, E_R, B_R)$ in the remainder system \eqref{Remd-Equ-GR} (or say equivalently, \eqref{Remd-Equ-Bref}) by the micro-macro decomposition method.

\subsubsection{Estimates on fluid variables}

By the definition of the projection operator $\mathbb{P}$ in \eqref{Projc-P}, the fluid part $\mathbb{P} G_R$ of $G_R$ reads
\begin{equation}\label{Projc-GR}
  \mathbb{P} G_R = \rho_R^+ \left( \begin{array}{c} 1 \\ 0 \end{array} \right) + \rho_R^- \left( \begin{array}{c} 0 \\ 1 \end{array} \right) + u_R \cdot \left( \begin{array}{c} v \\ v \end{array} \right) + \theta_R \left( \begin{array}{c} \tfrac{|v|^2}{2} - \tfrac{3}{2} \\ \tfrac{|v|^2}{2} - \tfrac{3}{2} \end{array} \right) \,,
\end{equation}
where $\rho_R^\pm = \langle g_R^\pm \rangle$, $u_R = \langle \tfrac{g_R^+ + g_R^-}{2} v \rangle$ and $\theta_R = \langle \tfrac{g_R^+ + g_R^-}{2} (\tfrac{|v|^2}{3} - 1) \rangle $. Our goal is to estimate the fluid variables $\rho_R^\pm (t,x)$, $u_R(t,x)$ and $\theta_R (t,x)$ in terms of $\mathbb{P}^\perp G_R$.

We first introduce the following so-called \emph{seventeen moments} set associated to $L^2_v$ \cite{Guo-Inventions2003}:
\begin{equation}\label{Basis-B}
  \mathfrak{B} = \bigg\{ \left( \begin{array}{c} 1 \\ 0 \end{array} \right) , \left( \begin{array}{c} 0 \\ 1 \end{array} \right) , \left( \begin{array}{c} v_i \\ 0 \end{array} \right) , \left( \begin{array}{c} 0 \\ v_i \end{array} \right) , \left( \begin{array}{c} v_i^2 \\ v_i^2 \end{array} \right) , \left( \begin{array}{c} v_i |v|^2 \\ v_i |v|^2 \end{array} \right) , \left( \begin{array}{c} v_j v_k \\ v_j v_k \end{array} \right) ; 1 \leq i \leq 3, 1 \leq j < k \leq 3 \bigg\} \,,
\end{equation}
which is one of the key structures in the process of micro-macro decomposition. One claims that $\mathfrak{B}$ is linearly independent in $L^2_v$. Indeed, suppose that
\begin{equation}\label{S-1}
  \begin{aligned}
    k_+ \left( \begin{array}{c} 1 \\ 0 \end{array} \right) & + k_- \left( \begin{array}{c} 0 \\ 1 \end{array} \right) + \sum_{i=1}^3 k_{i+} \left( \begin{array}{c} v_i \\ 0 \end{array} \right)  + \sum_{i=1}^3 k_{i-} \left( \begin{array}{c} 0 \\ v_i \end{array} \right) \\
    & + \sum_{i=1}^3 k_{i} \left( \begin{array}{c} v_i^2 \\ v_i^2 \end{array} \right) + \sum_{i=1}^3 \bar{k}_{i} \left( \begin{array}{c} v_i |v|^2 \\ v_i |v|^2 \end{array} \right) + \sum_{1 \leq i < j \leq 3} k_{ij} \left( \begin{array}{c} v_i v_j \\ v_i v_j \end{array} \right) = 0 \,,
  \end{aligned}
\end{equation}
then it is derived, from taking $L^2_v$-inner product in \eqref{S-1} by multiplying each element in the set $\mathfrak{B}$, that
\begin{equation}\label{S-2}
	\begin{aligned}
		& k_{i \pm} + 5 \bar{k}_i = 0 \,, \ k_+ + k_- + 6 k_i = 0 \,, \ k_{i+} + k_{i-} + 2 \bar{k}_i = 0 \ (1 \leq i \leq 3) \,, \\
		& k_\pm + \sum\limits_{i=1}^3 k_i = 0 \,, \ k_{ij} = 0 \ (1 \leq i < j \leq 3) \,,
	\end{aligned}
\end{equation}
where the relations $\int_{\R^3} M \d v = 1$, $\int_{\R^3} |v|^2 M \d v = 3$ and $\int_{\R^3} |v|^4 M \d v = 15$ have been used. A simple calculation implies that the linear system \eqref{S-2} admits only zero solution, namely, $k_\pm = k_{i \pm} = k_i = \bar{k}_i = k_{ij} = 0$.

Furthermore, we define a projection $\mathcal{P}_{\mathfrak{B}} : L^2_v \rightarrow \mathrm{span}\{\mathfrak{B}\} \subset L^2_v$ by
\begin{equation}\label{S-3}
  \mathcal{P}_{\mathfrak{B}} f : = \sum_{\zeta \in \mathfrak{B}} \langle f \cdot \zeta \rangle \zeta \,, \quad \forall \ f \in L^2_v \,.
\end{equation}

In order to estimate the fluid part $\mathbb{P} G_R$, or precisely, $\rho_R^\pm (t,x)$, $u_R(t,x)$ and $\theta_R(t,x)$, the treatment contains two main ingredients: (1) projecting the $G_R$-equation in \eqref{Remd-Equ-GR} into $\mathrm{span} \{ \mathfrak{B} \}$ yields the fluid relations \eqref{MM-8} below; (2) projecting the $G_R$-equation into $\ker \mathbb{L}$ leads to the macroscopic equations \eqref{MM-10} below. Based on the relations \eqref{MM-8} and \eqref{MM-10}, the following lemma can be derived:

\begin{lemma}\label{Lm-Mic-Mac-decomp}
  Assume that $(G_R, E_R, B_R)$ is a sufficiently smooth solution to the remainder system \eqref{Remd-Equ-GR}. Then there is a constant $C_1 > 0$ such that
  \begin{equation}\label{Mic-Mac-Inq}
    \begin{aligned}
      & \| u_R \|^2_{H^{N+1}_x} + \| \theta_R \|^2_{H^{N+1}_x} + \| \rho_R^+ \|^2_{H^{N+1}_x} + \| \rho_R^- \|^2_{H^{N+1}_x} + \| \div_x E_R \|^2_{H^N_x}
    \\
      & \leq C_1 \eps \tfrac{\d}{\d t} \widetilde{\mathcal{A}}_N (t) + \tfrac{C_1}{\eps^2} \| \mathbb{P}^\perp G_R \|^2_{H^{N+1}_x L^2_v} + C_1 \eps^2 \| \mathcal{P}_{\mathfrak{B}} H_R \|^2_{H^N_x L^2_v} + C_1 \mathcal{E}^{\IN}_{0,N+2} \| B_R \|^2_{H^N_x} \\
      & \quad + C_1 \Big[ \Big( \int_{\T} u_R \d x \Big)^2 + \Big( \int_{\T} \theta_R \d x \Big)^2 + \Big( \int_{\T} \rho^+_R \d x \Big)^2 + \Big( \int_{\T} \rho^-_R \d x \Big)^2 \Big]
    \end{aligned}
  \end{equation}
  for $\eps$ sufficiently small, where the quantity $ \widetilde{\mathcal{A}}_N (t) $ is defined as
    \begin{align}\label{Tilde-AN}
      \no \widetilde{\mathcal{A}}_N (t)
    = & \sum_{|m| \leq N} \Big\{ \sum_{i,j=1}^3 \int_{\T} \langle \partial^m \mathbb{P}^\perp G_R \cdot \zeta_{ij} \rangle \partial_j \partial^m u_R^i \d x - \tfrac{1}{4} \int_{\T} \partial^m u_R \cdot \nabla_x \partial^m ( \rho_R^+ + \rho_R^- ) \d x
    \\
      & \qquad \qquad + \sum_{i=1}^3 \int_{\T} \langle \partial^m \mathbb{P}^\perp G_R \cdot \zeta_i \rangle \partial_i \partial^m \theta_R \d x  \\
      \no & \qquad \qquad + \tfrac{1}{4} \sum_{i=1} \int_{\T} \big( \langle \partial^m \mathbb{P}^\perp G_R \cdot \zeta_{i+} \rangle \partial_i \partial^m \rho_R^+ + \langle \partial^m \mathbb{P}^\perp G_R \cdot \zeta_{i-} \rangle \partial_i \partial^m \rho_R^- \big) \d x \Big\}\,.
    \end{align}
  Here $\zeta_{ij}(v)$, $\zeta_i(v)$ and $\zeta_{i\pm} (v)$ are some fixed linear combinations of the basis $\mathfrak{B}$.
\end{lemma}

The proof of Lemma \ref{Lm-Mic-Mac-decomp} will be carried out in \S \ref{Sec:Lmm-MM} later. Remark that the last three terms of \eqref{Mic-Mac-Inq} require to be estimated. The term $\eps^2 \sum\limits_{|m| \leq N} \| \mathcal{P}_{\mathfrak{B}} \partial^m H_R \|^2_{L^2_{x,v}}$ will be controlled by estimating the mixed derivatives. The term $\mathcal{E}^{\IN}_{0,N+2} \| B_R \|^2_{H^N_x}$ can be dominated by finding enough decay structures of the Maxwell equations on $E_R$ and $B_R$ in the remainder system \eqref{Remd-Equ-GR}. By analyzing the conservation laws of mass, momentum and energy of \eqref{Remd-Equ-GR}, one can give an estimation on the last term involving integrals in \eqref{Mic-Mac-Inq}.

\subsubsection{Estimates on the integrals of fluid variables in \eqref{Mic-Mac-Inq}}

Recalling the truncated expansion \eqref{Soltn-Fom}, the well-prepared initial data \eqref{IC-3} and the conservation laws \eqref{Totl-Consv-Law-1}, it is derived that the leading term $(g_0^\pm , E_0, B_0)$ satisfies
\begin{equation}\label{Conservtn-g0}
  \begin{array}{l}
    \iint_{\T \times \R^3} g_0^\pm M \d v \d x = 0 \,, \quad \iint_{\T \times \R^3} |v|^2 (g_0^+ + g_0^-) M \d v \d x = 0 \,, \\[2mm]
    \iint_{\T \times \R^3} v ( g_0^+ + g_0^- ) M \d v \d x + \int_{\T} E_0 \times B_0 \d x = 0\,, \quad \int_{\T} B_0 \d x = 0 \,,
  \end{array}
\end{equation}
and that the remainder term $(g_R^\pm , E_R, B_R)$ obeys
\begin{equation}\label{Conservtn-gR}
  \begin{array}{l}
    \iint_{\T \times \R^3} ( g_R^\pm + \overline{g}_1^\pm ) M \d v \d x = 0 \,, \quad \iint_{\T} B_R \d x = 0 \,, \\ [2mm]
    \int_{\T \times \R^3} [ v ( g_R^+ + g_R^- ) + v ( {g}_1^+ + {g}_1^- ) ] M \d v \d x \\ [2mm]
    \qquad + \int_{\T} [ E_0 \times ( {B}_1 + B_R ) + ( {E}_1 + E_R ) \times B_0 + \eps E_R \times B_R ] \d x = 0 \,, \\[2mm]
    \iint_{\T \times \R^3} [ |v|^2 ( g_R^+ + g_R^- ) + |v|^2 ( {g}_1^+ + {g}_1^- ) ] M \d v \d x \\ [2mm]
    \qquad + \int_{\T} [ |E_0 + \eps {E}_1 + \eps E_R|^2 + |B_0 + \eps {B}_1 + B_R|^2 ] \d x = 0 \,,
  \end{array}
\end{equation}
where the facts $\int_{\R^3} {g}_2^\pm M \d v = \int_{\R^3} |v|^2 ( {g}_2^+ + {g}_2^-) M \d v = 0$ and $\int_{\R^3} v ( {g}_2^+ + {g}_2^-) M \d v = 0$ have been utilized. Combined with the relations
\begin{equation*}
    \begin{array}{lll}
      \int_{\R^3} g_R^\pm M \d v = \rho_R^\pm \,, & \int_{\R^3} v \tfrac{g_R^+ + g_R^- }{2} M \d v = u_R \,, & \int_{\R^3} ( \tfrac{|v|^2}{3} - 1 ) \tfrac{g_R^+ + g_R^- }{2} M \d v = \theta_R \,, \\ [2mm]
      \int_{\R^3} {g}_1^\pm M \d v = {\theta}_1 \pm \tfrac{1}{2} {n}_1 \,, \ & \int_{\R^3} v ( {g}_1^+ + {g}_1^- ) M \d v = 2 {u}_1 \,, \ & \int_{\R^3} |v|^2 ( {g}_1^+ + {g}_1^- ) M \d v = 12 {\theta}_1 \,,
    \end{array}
\end{equation*}
the conservation laws \eqref{Conservtn-gR} implies
\begin{equation}\label{Conservtn-Mac}
  \begin{array}{l}
    \int_{\T} ( \rho_R^\pm + {\theta}_1 \pm \tfrac{1}{2} {n}_1 ) \d x = 0 \,, \quad \int_{\T} B_R \d x = 0 \,, \\[2mm]
    \int_{\T} (u_R + {u}_1) \d x + \tfrac{1}{2} \int_{\T} [ E_0 \times ( {B}_1 + B_R ) + ( {E}_1 + E_R ) \times B_0 + \eps E_R \times B_R ] \d x = 0 \,, \\[2mm]
    \int_{\T} ( \theta_R + 3 {\theta}_1 ) \d x + \int_{\T} |E_0 + \eps {E}_1 + \eps E_R|^2 + |B_0 + \eps {B}_1 + \eps B_R|^2 \d x = 0 \,.
  \end{array}
\end{equation}

The above conservation laws \eqref{Conservtn-Mac} imply the following lemma to control the terms involving integrals of \eqref{Mic-Mac-Inq} in Lemma \ref{Lm-Mic-Mac-decomp}.

\begin{lemma}\label{Lm-Integral-Bnd}
  Assume that $(G_R, E_R , B_R)$ is a sufficiently smooth solution to the remainder system \eqref{Remd-Equ-GR} with the well-prepared initial data \eqref{IC-Remd}-\eqref{IC-3}. Then there is a constant $C_2 > 0$ such that
  \begin{equation}\label{Poincare-BR}
    \| B_R \|_{L^2_x} \leq C_2 \| \nabla_x B_R \|_{L^2_x} \,,
  \end{equation}
  and
  \begin{equation}\label{Integral-Inqu}
    \begin{aligned}
      & \Big( \int_{\T} u_R \d x \Big)^2 + \Big( \int_{\T} \theta_R \d x \Big)^2 + \Big( \int_{\T} \rho_R^+ \d x \Big)^2 + \Big( \int_{\T} \rho_R^- \d x \Big)^2 \\
      & \leq C_2 \Big[ \mathcal{D}_{0,2} (t) + \mathcal{D}_{1,2} (t) + \eps^2 \big( \| E_R \|^4_{L^2_x} + \| B_R \|^2_{L^2_x} \| \nabla_x B_R \|^2_{L^2_x} \big) \Big]
    \end{aligned}
  \end{equation}
  for $\eps$ small enough, where the quantities $ \mathcal{D}_{0,2} (t) $ and $ \mathcal{D}_{1,2} (t) $ are mentioned in Lemma \ref{lemm:bnd-NSMF} and \ref{lemm:bnd-linearMaxwl} respectively.
\end{lemma}

\begin{proof}
The inequality \eqref{Poincare-BR} is derived from the Poincar\'e inequality, by noticing the zero mean value property of $B_R$ in \eqref{Conservtn-Mac}.

Next one focuses on the inequality \eqref{Integral-Inqu}. For the first equation of \eqref{Conservtn-Mac}, there holds
  \begin{align}\label{I-1}
    \no \Big( \int_{\T} \rho_R^\pm \d x \Big)^2 \leq & \Big( \int_{\T} ( {\theta}_1 \pm \tfrac{1}{2} {n}_1 ) \d x \Big)^2 \leq C ( \| {\theta}_1 \|^2_{L^2_x} + \| {n}_1 \|^2_{L^2_x} ) \\
    \leq & C \mathcal{E}^{\IN}_{0,2} ( 1 + \mathcal{E}^{\IN}_{0,2} ) \mathcal{D}_{0,2} (t) + C \mathcal{D}_{1,2} (t) \,,
  \end{align}
  where the last inequality is implied by Lemmas \ref{lemm:bnd-NSMF}-\ref{lemm:bnd-linearMaxwl}.

Combined with Lemmas \ref{lemm:bnd-NSMF} and \ref{lemm:bnd-linearMaxwl}, the second equality of \eqref{Conservtn-Mac} reduces to
  \begin{align}\label{I-2}
    \no \Big( \int_{\T} u_R \d x \Big)^2 \leq & C \Big( \int_{\T} {u}_1 \d x \Big)^2 + C \Big( \int_{\T} E_0 \times ( {B}_1 + B_R ) \d x \Big)^2 \\
    \no & + C \Big( \int_{\T} ( {E}_1 + E_R ) \times B_0 \d  x \Big)^2 + C \eps^2 \Big( \int_{\T} E_R \times B_R \d x \Big)^2 \\
    \no \leq & C \| {u}_1 \|^2_{L^2_x} + C \| E_0 \|^2_{L^2_x} ( \| {B}_1 \|^2_{L^2_x} + \| B_R \|^2_{L^2_x} ) \\
    & + C \| B_0 \|^2_{L^2_x} ( \| {E}_1 \|^2_{L^2_x} + \| E_R \|^2_{L^2_x} ) + C \eps^2 \| E_R \|^2_{L^2_x} \| B_R \|^2_{L^2_x} \\
    \no \leq & C \mathcal{D}_{1,2}(t) + C \eps^2 \| E_R \|^2_{L^2_x} \| B_R \|^2_{L^2_x} \\
    \no & + C \mathcal{E}^{\IN}_{0,2} ( \mathcal{D}_{1,2} (t) + \| E_R \|^2_{L^2_x} + \| \nabla_x B_R \|^2_{L^2_x} ) \\
    \no \leq & C ( 1 + \mathcal{E}^{\IN}_{0,2} ) \mathcal{D}_{1,2} (t) + C \mathcal{E}^{\IN}_{0,2} ( \| E_R \|^2_{L^2_x} + \| \nabla_x B_R \|^2_{L^2_x} ) + C \eps^2 \| E_R \|^2_{L^2_x} \| B_R \|^2_{L^2_x} \,,
  \end{align}
where the inequality \eqref{Poincare-BR} has also been used.

Furthermore, combined with the last conservation law in \eqref{Conservtn-Mac}, the inequality \eqref{Poincare-BR}, Lemmas \ref{lemm:bnd-NSMF}-\ref{lemm:bnd-linearMaxwl} and the inequality $ \| B_0 \|^4_{L^2_x} \leq C \| \nabla_x B_0 \|^2_{L^2_x} \| B_0 \|^2_{H^1_x} $ together, it is derived that
  \begin{align}\label{I-3}
    \no \Big( \int_{\T} \theta_R \d x \Big)^2 \leq & C \Big( \int_{\T} {\theta}_1 \d x \Big)^2 + C ( \| E_0 \|^4_{L^2_x} + \| B_0 \|^4_{L^2_x} ) \\
    & + C \eps^2 ( \| {E}_1 \|^4_{L^2_x} + \| {B}_1 \|^4_{L^2_x} + \| E_R \|^4_{L^2_x} + \| B_R \|^4_{L^2_x} ) \\
    \no \leq & C \mathcal{E}^{\IN}_{0,2} ( 1 + \mathcal{E}^{\IN}_{0,2} ) \mathcal{D}_{0,2}(t) + C \eps^2 ( \mathcal{E}^{\IN}_{0,2} + \mathcal{E}^{\IN}_{1,2} ) \mathcal{D}_{1,2} (t) \\
    \no & + C \eps^2 ( \| E_R \|^4_{L^2_x} + \| B_R \|^2_{L^2_x} \| \nabla_x B_R \|^2_{L^2_x} ) \,.
  \end{align}
Consequently, the above inequalities \eqref{I-1}, \eqref{I-2} and \eqref{I-3} conclude the inequality \eqref{Integral-Inqu}, and the proof of Lemma \ref{Lm-Integral-Bnd} is hence completed.
\end{proof}

\subsubsection{Estimates on the electric-magnetic filed \texorpdfstring{$(E_R, B_R)$}{(E\_R, B\_R)}}

This subsection is devoted to estimating the term $C_1 \mathcal{E}^{\IN}_{0,N+2} \| B_R \|^2_{H^N_x}$ in \eqref{Mic-Mac-Inq}. The \emph{key point} is to seek enough dissipation or decay properties on $B_R$ by making use of the Maxwell equations, i.e., the last four equations in \eqref{Remd-Equ-GR}
\begin{equation}\label{Mxw-GR}
  \left\{
    \begin{array}{l}
      \partial_t E_R - \nabla_x \times B_R = - \frac{1}{\eps} \langle G_R \cdot \mathcal{T}_1 v \rangle = - \frac{1}{\eps} \langle \mathbb{P}^\perp G_R \cdot \mathcal{T}_1 v \rangle \,, \\[2mm]
      \partial_t B_R + \nabla_x \times E_R = 0 \,, \quad \div_x E_R = \rho_R^+ - \rho_R^- \,, \quad \div_x B_R = 0 \,.
    \end{array}
  \right.
\end{equation}
Notice that the second equation in \eqref{Mxw-GR}, Faraday's law, does not have explicit dissipative term. Besides that, a singular term $- \frac{1}{\eps} \langle \mathbb{P}^\perp G_R \cdot \mathcal{T}_1 v \rangle$ arises in the right-hand side \eqref{Mxw-GR}. If one takes $\partial_t$ on the evolution equation of $B_R$, then by combining together with the equation of $E_R$, there holds $$\partial_{tt} B_R + \nabla_x \times ( \nabla_x \times B_R ) = \tfrac{1}{\eps} \nabla \times \langle \mathbb{P}^\perp G_R \cdot \mathcal{T}_1 v \rangle.$$
This, together with the fact $ \nabla_x \times ( \nabla_x \times B_R ) = - \Delta_x B_R $ under the divergence-free property $\div_x B_R = 0$, implies
\begin{equation}\label{Mxw-1}
  \partial_{tt} B_R - \Delta_x B_R = \tfrac{1}{\eps} \nabla_x \times \langle \mathbb{P}^\perp G_R \cdot \mathcal{T}_1 v \rangle.
\end{equation}

However, the dissipation of \eqref{Mxw-1} still remains not enough to reach our aim. We try to derive the Ohm's law from the microscopic equation of $G_R$ in \eqref{Remd-Equ-GR}, which will supply a decay term $\partial_t B_R$. More precisely, by performing inner product with $\mathcal{T}_1 = (1, -1)^\top$ in the first equation of \eqref{Remd-Equ-GR}, we get
\begin{multline}\label{Mxw-2}
  \eps \partial_t ( G_R \cdot \mathcal{T}_1 ) + v \cdot \nabla_x ( G_R \cdot \mathcal{T}_1 ) + (v \times B_0) \cdot \nabla_v ( G_R \cdot \mathcal{T}_2 ) - 2 E_R \cdot v \\
  + (v \times B_R) \cdot \nabla_v ( G_0 \cdot \mathcal{T}_2 ) + \tfrac{1}{\eps} ( \mathcal{L} + \mathfrak{L} ) ( G_R \cdot \mathcal{T}_1 ) - \eps H_R \cdot \mathcal{T}_1 = 0 \,,
\end{multline}
where $\mathcal{T}_2 = (1, 1)^\top$, and we have used the following relation:
\begin{equation*}
  \mathbb{L} G_R \cdot \mathcal{T}_1 = ( \mathcal{L} + \mathfrak{L} ) ( G_R \cdot \mathcal{T}_1 ) \,.
\end{equation*}
Recalling the property \eqref{L+L-Dual} of the operator $\mathcal{L} + \mathfrak{L}$ in Section \ref{Sec:Formal-Analysis}, one knows that, for $\Phi(v) = v \in L^2_v$, there exists a unique function $\widetilde{\Phi} (v) \in \ker^\perp (\mathcal{L} + \mathfrak{L})$ such that
$$ (\mathcal{L} + \mathfrak{L}) \widetilde{\Phi} = \Phi \,. $$
From multiplying by $\widetilde{\Phi}(v) M$ in \eqref{Mxw-2} and integrating over $v \in \R^3$, it follows:
\begin{align}\label{Mxw-3}
  \no \tfrac{1}{\eps} \langle G_R \cdot \mathcal{T}_1 v \rangle = & - \eps \partial_t \langle G_R \cdot \mathcal{T}_1 \widetilde{\Phi}(v) \rangle - \langle v \cdot \nabla_x (G_R \cdot \mathcal{T}_1 \widetilde{\Phi}(v) ) \rangle \\
  & + \sigma E_R - \langle (v \times B_0) \cdot \nabla_v (G_R \cdot \mathcal{T}_2 ) \widetilde{\Phi}(v) \rangle \\
  \no & - \langle (v \times B_R) \cdot \nabla_v (G_0 \cdot \mathcal{T}_2) \widetilde{\Phi}(v) \rangle + \eps \langle H_R \cdot \mathcal{T}_1 \widetilde{\Phi}(v) \rangle \,.
\end{align}
A simple calculation reduces to
\begin{align}\label{Mxw-4}
  \no \tfrac{1}{\eps} \langle G_R \cdot \mathcal{T}_1 v \rangle = & - \eps \partial_t \langle \mathbb{P}^\perp G_R \cdot \mathcal{T}_1 \widetilde{\Phi}(v) \rangle - \langle v \cdot \nabla_x ( \mathbb{P}^\perp G_R \cdot \mathcal{T}_1 ) \widetilde{\Phi}(v) \rangle \\
  \no & + \langle (v\times B_0) ( \mathbb{P}^\perp G_R \cdot \mathcal{T}_2 ) \cdot \nabla_v \widetilde{\Phi}(v) \rangle + \eps \langle H_R \cdot \mathcal{T}_1 \widetilde{\Phi}(v) \rangle \\
  & + \sigma E_R - \tfrac{1}{2} \sigma \nabla_x ( \rho_R^+ - \rho_R^- ) + \sigma ( u_R \times B_0 + u_0 \times B_R) \\
  \no := & \sigma E_R - \tfrac{1}{2} \sigma \nabla_x ( \rho_R^+ - \rho_R^- ) + \mathcal{K} ( \mathbb{P}^\perp G_R ) \,.
\end{align}
We emphasize that Ohm's law \eqref{Mxw-4} cancels exactly the singularity in the right-hand side of \eqref{Mxw-GR}. Thus, substituting \eqref{Mxw-4} into \eqref{Mxw-1} implies
\begin{equation}\label{Mxw-5}
  \partial_{tt} B_R - \Delta_x B_R + \sigma \partial_t B_R = \nabla_x \times \mathcal{K} ( \mathbb{P}^\perp G_R ) \,,
\end{equation}
where we have used Faraday's law $\partial_t B_R + \nabla_x \times E_R = 0$ and $\nabla_x \times (\nabla_x f) = 0$ for any function $f(x)$. Therefore, the decay term $\partial_t B_R$ of $B_R$-equation has been found.

Moreover, combining together the first equation (the Amp\`ere equation) in \eqref{Mxw-GR} and Ohm's law \eqref{Mxw-4}, one can infer that $E_{R,\eps}$ satisfies
\begin{equation}\label{Mxw-6}
	\begin{aligned}
		\partial_t E_R + \sigma E_R - \tfrac{1}{2} \sigma \nabla_x \div_x E_R - \nabla_x \times B_R + \mathcal{K} ( \mathbb{P}^\perp G_R ) = 0 \,,
	\end{aligned}
\end{equation}
which contains the decay $\sigma E_R$.

Based on \eqref{Mxw-5}-\eqref{Mxw-6}, the following lemma is obtained, with its proof being provided in \S \ref{Sec:Lmm-MD} later.
\begin{lemma}\label{Lm-Mxw-Dec}
  Assume that $(G_R, E_R, B_R)$ is a sufficiently smooth solution to \eqref{Remd-Equ-GR} with the well-prepared initial \eqref{IC-Remd}-\eqref{IC-3}. There is a small constant $\lambda_R (N+2) \in ( 0 , \lambda_1 (N+2) ]$, depending only on $N$, $\sigma$, $\mu$, $\kappa$, such that if
  \begin{equation*}
    \mathcal{E}^{\IN}_{0,N+2} \leq \lambda_R (N+2) \,,
  \end{equation*}
  then the following inequality holds
  \begin{align}\label{Decay-Inq}
    \no & \tfrac{\d}{\d t} \Big( \| E_R \|^2_{H^{N-1}_x} + ( 1 - \delta + \sigma \delta ) \| B_R \|^2_{H^{N-1}_x} + \| \nabla_x B_R \|^2_{H^{N-1}_x} + (1 - \delta) \| \partial_t B_R \|^2_{H^{N-1}_x} \\
    \no & \quad + \delta \| \partial_t B_R + B_R \|^2_{H^{N-1}_x} + \eps^2 \| \nabla_x \times \langle \mathbb{P}^\perp G_R \cdot \mathcal{T}_1 \widetilde{\Phi}(v) \rangle \|^2_{H^{N-1}_x} \Big) \\
    \no & \quad + \sigma \| E_R \|^2_{H^{N-1}_x} + \delta \| \nabla_x B_R \|^2_{H^{N-1}_x} + ( \sigma - \delta ) \| \partial_t B_R \|^2_{H^{N-1}_x} + \sigma \| \div_x E_R \|^2_{H^{N-1}_x} \\
  & \leq C_3 \eps \tfrac{\d}{\d t} \sum_{|m| \leq N-1} \int_{\T} \langle \partial^m G_R \cdot \mathcal{T}_1 \widetilde{\Phi} (v) \rangle \cdot ( \partial^m E_R - \delta \nabla_x \times \partial^m B_R ) \d x \\
    \no & \quad + \tfrac{C_3}{\eps^2} \| \mathbb{P}^\perp G_R \|^2_{H^{N+1}_x L^2_v (\nu)} + C_3 \eps^2 \| H_R \|^2_{H^N_x L^2_v} + C_3 \sqrt{\mathcal{E}^{\IN}_{0,N+2}} \| u_R \|^2_{H^N_x},
  \end{align}
  for sufficiently small $\eps$ and $C_3 > 0$, depending only on $N$, $\sigma$, $\mu$ and $\kappa$.
\end{lemma}


\subsubsection{Kinetic dissipation of \texorpdfstring{$G_R$}{G\_R}}

In this subsection, the goal is to seek the kinetic dissipative structure $\tfrac{1}{\eps^2} \| \mathbb{P}^\perp G_R \|^2_{H^{N+1}_x L^2_v (\nu)}$ with singularity of order $\frac{1}{\eps^2}$ in the pure spatial derivative estimates. This kinetic dissipation is also one of important structures to deal with the singularity appearing in the remainder system \eqref{Remd-Equ-GR}.

Recall the definition of the projection operator $\mathbb{P}$ in \eqref{Projc-P}, i.e.,
\begin{equation*}
  \mathbb{P} G_R = \rho_R^+ \left( \begin{array}{c} 1 \\ 0 \end{array} \right) + \rho_R^- \left( \begin{array}{c} 0 \\ 1 \end{array} \right) + u_R \cdot \left( \begin{array}{c} v \\ v \end{array} \right) + \theta_R \left( \begin{array}{c} \tfrac{|v|^2}{2} - \tfrac{3}{2} \\ \tfrac{|v|^2}{2} - \tfrac{3}{2} \end{array} \right) \,,
\end{equation*}
then there exists a generic constant $C_0 > 0$, such that
\begin{equation}\label{Sp-1}
  \begin{aligned}
    C_0 \| \mathbb{P} G_R \|^2_{H^{N+1}_x L^2_v} \leq \| ( \rho_R^+, \rho_R^-, u_R, \theta_R ) \|^2_{H^{N+1}_x} \leq \tfrac{1}{C_0} \| \mathbb{P} G_R \|^2_{H^{N+1}_x L^2_v} \,.
  \end{aligned}
\end{equation}
Moreover, the inequality \eqref{Poincare-BR} in Lemma \ref{Lm-Integral-Bnd} implies
\begin{equation}\label{Sp-2}
  \| B_R \|^2_{H^N_x} \leq C \| \nabla_x B_R \|^2_{H^{N-1}_x} \,.
\end{equation}

Based on Lemmas \ref{Lm-Mic-Mac-decomp}, \ref{Lm-Integral-Bnd} and \ref{Lm-Mxw-Dec}, the coercivity \eqref{L-Proty-2} in Lemma \ref{Lm-L-Propty}, the above inequalities \eqref{Sp-1}-\eqref{Sp-2} yield the following lemma.
\begin{lemma}\label{Lm-Unif-Spatial-Bnd}
  Let $N \geq 2$ and $(G_R , E_R, B_R)$ be a sufficiently smooth solution to the remainder system \eqref{Remd-Equ-GR} with the well-prepared initial data \eqref{IC-Remd}-\eqref{IC-3}. There is a small constant $\lambda_R(N+2) \in ( 0 , \lambda_1 (N+2) ]$, depending only on $N$, $\sigma$, $\mu$, $\kappa$, such that if
  \begin{equation}\label{Sp-3}
    \mathcal{E}^{\IN}_{0,N+5} \leq \lambda_R(N+5) \,,
  \end{equation}
  where the quantity $\mathcal{E}^{\IN}_{0,N+2}$ and $\lambda_1 (N+2)$ are mentioned in Lemma \ref{lemm:bnd-NSMF} and \ref{lemm:bnd-linearMaxwl} respectively, then there exist some constants  $C_4 > 1$, $C_5\,, C_6 \,, C_7 > 0$, depending only on $N$, $\sigma$, $\mu$, $\kappa$, such that, for sufficiently small $\eps$,
  \begin{align}\label{Unif-Spatial-Bnd}
    \no & \tfrac{\d}{\d t} \Big[ \| E_R \|^2_{H^{N-1}_x} + (1 - \delta + \sigma \delta) \| B_R \|^2_{H^{N-1}_x} + ( 1 - \delta ) \| \partial_t B_R \|^2_{H^{N-1}_x} + \delta \| \partial_t B_R +  B_R \|^2_{H^{N-1}_x} \\
    \no & \qquad - C_3 \eps \mathcal{A}_N (t) + \eps^2 \| \nabla_x \times \langle \mathbb{P}^\perp G_R \cdot \mathcal{T}_1 \widetilde{\Phi} (v) \rangle \|^2_{H^{N-1}_x} \\
    \no & \qquad + C_4 \big( \| E_R \|^2_{H^{N+1}_x} + \| B_R \|^2_{H^{N+1}_x} + \| G_R \|^2_{H^{N+1}_x L^2_v} \big) + C_5 \big( \mathcal{E}_{1,N} (t) + \widetilde{C}_N \mathcal{E}_{0,N+2} (t) \big) \Big] \\
    \no & \quad + C_6 \Big( \| E_R \|^2_{H^{N-1}_x} + \| \nabla_x B_R \|^2_{H^{N-1}_x} + \| \partial_t B_R \|^2_{H^{N-1}_x} + \mathcal{D}_{0,N+1} (t) \\
    & \qquad + \| \div_x E_R \|^2_{H^{N}_x} + \| \mathbb{P} G_R \|^2_{H^{N+1}_x L^2_v} + \tfrac{1}{\eps^2} \| \mathbb{P}^\perp G_R \|^2_{H^{N+1}_x L^2_v (\nu)} \Big) \\
    \no & \leq C_7 \Big[ \eps^2 \| H_R \|^2_{H^N_x L^2_v} + \eps^2 \| \mathcal{P}_{\mathfrak{B}} H_R \|^2_{H^N_x L^2_v} + \| B_R \|^2_{H^{N+1}_x} \mathcal{D}_{0,N+2} (t) + \mathcal{E}^{\IN}_{0,N+2} \| \nabla_v \mathbb{P}^\perp G_R \|^2_{H^N_x L^2_v} \\
    \no & \quad + \eps^2 ( \| E_R \|^4_{L^2_x} + \| B_R \|^2_{L^2_x} \| \nabla_x B_R \|^2_{L^2_x} ) + \sum_{|m|\leq N+1} \int_{\T} \langle \partial^m H_R \cdot \partial^m G_R  \rangle \d x \Big]
  \end{align}
where the constant $C_3 > 0$ is taken as in Lemma \ref{Lm-Mxw-Dec}, and the quantity $\mathcal{A}_N (t)$ is defined as
  \begin{equation}\label{Sp-4}
    \mathcal{A}_N (t) = \widetilde{\mathcal{A}}_N (t) + \sum_{|m| \leq N-1} \int_{\T} \langle \partial^m G_R \cdot \mathcal{T}_1 \widetilde{\Phi} (v) \rangle \cdot ( \partial^m E_R - \delta \nabla_x \times \partial^m B_R ) \,,
  \end{equation}
  with $ \widetilde{\mathcal{A}}_N (t) $ being given as \eqref{Tilde-AN} in Lemma \ref{Lm-Mic-Mac-decomp}.
\end{lemma}

The proof of Lemma \ref{Lm-Unif-Spatial-Bnd} will be given in \S \ref{Sec:Lmm-USB} later.

\subsection{Estimates on the mixed \texorpdfstring{$(x,v)$}{(x,v)}-derivatives with weighted norms}
\label{sec:mixed-derivatives}

In this subsection, we try to deal with the quantities $ \| H_R \|^2_{H^N_x L^2_v}$ (or $ \| \mathcal{P}_{\mathfrak{B}} H_R \|^2_{H^N_x L^2_v}$), $ \| \nabla_v \mathbb{P}^\perp G_R \|^2_{H^N_x L^2_v}$ and $\sum_{|m|\leq N+1} \int_{\T} \langle \partial^m H_R \cdot \partial^m G_R  \rangle \d x$ in the right-hand of \eqref{Unif-Spatial-Bnd}. As shown in \eqref{Remnd-H}, the quantity $H_R$ contains $\mathcal{Q} (g_R^\pm, g_R^\pm)$. These quantities should be controlled by employing Lemma \ref{Lm-L-Q}, in which the norms with a weight $\w (v) = (1+|v|^2)^{1/2}$ are required. We mention that the quantity $ \| \nabla_v \mathbb{P}^\perp G_R \|^2_{H^N_x L^2_v}$ comes from the Lorentz force term $\mathcal{T} (v \times B_0) \cdot \nabla_v G_R$, which indicates that the mixed $(x,v)$-derivative estimates on the remainder system \eqref{Remd-Equ-GR} is necessary for deriving the closed energy bounds. This is different from the Boltzmann equation. As pointed out in \cite{JXZ-IUM-2018}, the Boltzmann equation possesses the closed pure spatial derivative estimates, since the Lorentz force is not involved in there.

Note that the hydrodynamic part $\mathbb{P}G_R$ can be bounded by
	\begin{align}\label{Weight-Hydro-Control}
	  \|\w^l \p^m_{\beta} \mathbb{P}G_R\|_{L^2_{x,v}} \le C \|\p^m \mathbb{P}G_R\|_{L^2_{x,v}}.
	\end{align}
It suffices to estimate the remaining microscopic part $\w^l \p^m_{\beta} \mathbb{P}^\perp G_R$ with $|m|+|\beta| \le N$.
\begin{lemma}\label{lemm:remainder-apriori}
	Let $N\ge 2$ and $l \ge 0$. Assume that $(G_R, E_R, B_R)$ is a smooth solution of the remainder system \eqref{Remd-Equ-GR}. Then there exists some constant $C>0$ such that for sufficiently small $\eps>0$, it holds that
		\begin{align}\label{esm:remainder-mix}
		  & \tfrac{\d}{\d t} \| \mathbb{P}^\perp G_R\|^2_{H^N_{x,v} (\w^l)} + \tfrac{1}{2 \eps^2} \| \mathbb{P}^\perp G_R\|^2_{H^N_{x,v} (\w^l \nu) } \\
		  \no \le & \tfrac{C}{\eps^2} \| \mathbb{P}^\perp G_R\|^2_{H^{N+1}_x L^2_v (\nu)}
		 	+ C \big( \|E_R\|^2_{H^{N-1}_x} + \|\nabla_x B_R\|^2_{H^{N-1}_x} + \| \mathbb{P} G_R\|^2_{H^{N+1}_x L^2_v} \big)
		 	\\\no &
		 	+ \sum_{|m|+|\beta| \le N} \int_{\mathbb{T}^3} \langle \p^m_{\beta} H_R \cdot \w^{2l} \p^m_{\beta} \mathbb{P}^\perp G_R \rangle \d x.
		\end{align}
\end{lemma}


Observe that in \eqref{esm:remainder-mix}, the highest order of the derivatives of $G_R$ in the right-hand side is $N+1$, as appeared in the terms $ \| \mathbb{P}^\perp G_R \|^2_{H^{N+1}_x L^2_v (\nu)}$ and $ \| \mathbb{P} G_R \|^2_{H^{N+1}_x L^2_v}$.
However, the highest order of the derivatives in the left-hand side of \eqref{esm:remainder-mix} is $N$. It seems to be unreasonable.
One recalls that the dissipation on $E_R$ obtained in \eqref{Unif-Spatial-Bnd} is exactly $\| E_R \|^2_{H^{N-1}_x} + \| \div_x E_R \|^2_{H^N_x}$ (where the last quantity does not play a role). Nevertheless, the dissipation on the pure spatial derivatives of $G_R$ in \eqref{Unif-Spatial-Bnd} is $ \| \mathbb{P} G_R \|^2_{H^{N+1}_x L^2_v} + \tfrac{1}{\eps^2} \| \mathbb{P}^\perp G_R \|^2_{H^{N+1}_x L^2_v (\nu)} $.
This indicates that, one order of spatial derivative on $E_R$ is transformed to $G_R$, such that the highest order of derivatives of $E_R$ in \eqref{esm:remainder-mix} is $N-1$, i.e., $\| E_R \|^2_{H^{N-1}_x}$. For details, the readers could see the derivation of \eqref{eq:term-I1} in \S \ref{Sec:Lmm-RA12} later.

Notice that the term $ \|\nabla_v \mathbb{P}^\perp G_R \|^2_{H^N_x L^2_v}$ in \eqref{Unif-Spatial-Bnd} cannot be dealt with the inequality \eqref{esm:remainder-mix} in Lemma \ref{lemm:remainder-apriori}. The norm $ \| \mathbb{P}^\perp G_R\|^2_{\widetilde{\mathbf{H}}^{N+1}_{x,v} (\w^{2l} )}$ is thereby investigated in the following lemma, where the highest order of the spatial derivatives is at most $N$. As the similar reasons illustrated above, this is the key transforming one order spatial derivative of $E_R$ to $G_R$.

\begin{lemma}\label{lemm:remainder-apriori2}
	Let $N\ge 2$ and $l \ge 0$. Assume that $(G_R, E_R, B_R)$ is a smooth solution of the remainder system \eqref{Remd-Equ-GR}. Then there exists some constant $C>0$ such that for sufficiently small $\eps>0$, it holds that
		\begin{align}\label{esm:remainder-mix2}
		  & \tfrac{\d}{\d t} \| \mathbb{P}^\perp G_R\|^2_{\widetilde{\mathbf{H}}^{N+1}_{x,v} (\w^{2l}} + \tfrac{1}{2 \eps^2} \| \mathbb{P}^\perp G_R\|^2_{\widetilde{\mathbf{H}}^{N+1}_{x,v} (\w^{2l} \nu)} \\
		  \no \le & \tfrac{C}{\eps^2} \| \mathbb{P}^\perp G_R\|^2_{H^{N+1}_x L^2_v (\nu)} + C \big( \|E_R\|^2_{H^{N-1}_x} + \|\nabla_x B_R\|^2_{H^{N-1}_x} + \| \mathbb{P} G_R \|^2_{H^{N+1}_x L^2_v} \big) \\
		 \no & + C \eps^2 \| H_R\|^2_{H^N_x L^2_v} + \sum_{|m|+|\beta| \le N+1, \beta \neq 0} \int_{\mathbb{T}^3} \langle \p^m_{\beta} H_R \cdot \w^{2l} \p^m_{\beta} \mathbb{P}^\perp G_R \rangle \d x.
		\end{align}
\end{lemma}

The proofs of Lemmas \ref{lemm:remainder-apriori}-\ref{lemm:remainder-apriori2} are given in \S \ref{Sec:Lmm-RA12} later.

\subsection{Estimates on the higher order term \texorpdfstring{$H_R$}{H\_R}}

In this subsection, the main goal is to control the quantities related to $H_R$ in \eqref{Unif-Spatial-Bnd}, \eqref{esm:remainder-mix} and \eqref{esm:remainder-mix2}, namely, to control
\begin{align*}
	\begin{array}{cc}
	\displaystyle
  \eps^2 \| \mathcal{P}_{\mathfrak{B}} H_R \|^2_{H^N_x L^2_v}\,, \quad \eps^2 \| H_R \|^2_{H^N_x L^2_v} \,,\quad
  & \sum_{|m|\leq N+1} \int_{\T} \langle \partial^m H_R \cdot \partial^m G_R  \rangle \d x
\\[5pt]
	\displaystyle
  \sum_{|m| + |\beta| \leq N} \int_{\T} \langle \partial^m_\beta H_R \cdot w^{2l} \partial^m_\beta \mathbb{P}^\perp G_R \rangle \d x \,, \qquad
  & \sum_{\substack{|m| + |\beta| \leq N + 1 , \beta \neq 0}} \int_{\T} \langle \partial^m_\beta H_R \cdot \w^{2l} \partial^m_\beta \mathbb{P}^\perp G_R \rangle \d x \,,
  \end{array}
\end{align*}
where $H_R$ is defined in \eqref{Remnd-H} (dropping the index $\eps$).

By the definition of the operator $\mathcal{P}_{\mathfrak{B}}$ in \eqref{S-3}, we have
\begin{equation}\label{HR-1}
  \eps^2 \| \mathcal{P}_{\mathfrak{B}} H_R \|^2_{H^N_x L^2_v} \leq C \eps^2 \| H_R \|^2_{H^N_x L^2_v}.
\end{equation}
So it suffices to estimate the last four terms.
\begin{lemma}\label{Lm-HR-Square}
  Let $N \geq 4$ and $\gamma \in [0,1]$. There exists a constant $C_* > 0$, depending only on $\mathcal{E}^{\IN}_{0,N+5}$, $\mathcal{E}^{\IN}_{1 , N+3}$, $\mu$, $\kappa$, $\sigma$ and $N$, such that
  \begin{align}\label{HR-Square}
    \eps^2 \| H_R \|^2_{H^N_x L^2_v}
    \leq & C_* ( \eps^2 + \mathcal{E}^{\IN}_{0,N+5} ) \mathbb{D}_{N, 2 \gamma + 1} (G_R, E_R, B_R) \\ \no
    & + C (1 + \eps^2 ) \mathbb{E}_{N, 2 \gamma + 1} (G_R, E_R, B_R) \mathbb{D}_{N, 2 \gamma + 1} (G_R, E_R, B_R) \\ \no
    & + C_* \big( \mathcal{D}_{0,N+5} (t) + \mathcal{D}_{1,N+3} (t) \big)
  \end{align}
  holds for some constant $C > 0$, depending only on $\mu$, $\kappa$, $\sigma$ and $N$, where the functionals $ \mathbb{E}_{N,2 \gamma + 1} (G_R$, $E_R, B_R)$ and $\mathbb{D}_{N,2 \gamma + 1} (G_R, E_R, B_R)$ are defined in \eqref{Energ-1} and \eqref{Dspt-Rate}, respectively.
\end{lemma}

\begin{lemma}\label{Lm-HR-Inner-Product}
  Let $N \geq 4$ and $l \geq 0$. Then there is a constant $C > 0$ such that
  \begin{align}\label{HR-Inner-Product}
    \no & \Bigg| \sum_{|m| + |\beta| \leq N} \int_{\T} \langle H_R \cdot \w^{2l} \partial^m_\beta \mathbb{P}^\perp G_R \rangle \d x \Bigg| + \Bigg| \sum_{\substack{|m| + |\beta| \leq N+1 , \beta \neq 0 }} \int_{\T} \langle H_R \cdot \w^{2l} \partial^m_\beta \mathbb{P}^\perp G_R \rangle \d x \Bigg| \\
    \no & + \Bigg| \sum_{|m| \leq N+1} \int_{\T} \langle \partial^m H_R \cdot \partial^m G_R \rangle \d x \Bigg| \\
     \leq & C \big( \eps + \sqrt{\mathcal{E}^{\IN}_{0, N+5}} \big) \mathbb{D}_{N,l} (G_R, E_R, B_R)  + C \eps {\mathbb{E}}^\frac{1}{2}_{N,l} (G_R, E_R, B_R) \mathbb{D}_{N,l} (G_R, E_R, B_R) \,,
  \end{align}
  with ${\mathbb{E}}_{N,l} (G_R, E_R, B_R)$ and $\mathbb{D}_{N,l} (G_R, E_R, B_R)$ being defined in \eqref{Energ-1} and \eqref{Dspt-Rate}, respectively.
\end{lemma}

\subsection{Close the energy bounds: proof of Theorem \ref{Main-Thm}}

Based on above Lemmas \ref{Lm-Unif-Spatial-Bnd}-\ref{lemm:remainder-apriori2}, and Lemma \ref{lemm:bnd-NSMF}-\ref{lemm:bnd-linearMaxwl} in Section \ref{Sec:Prem}, one can derive the closed energy estimates on the remainder system \eqref{Remd-Equ-GR}.

\begin{proposition}\label{Prop-Est}
	Let $N \geq 4$ and $l \geq 2 \gamma + 1$. Let $(G_R, E_R, B_R)$ be the smooth solution over time interval $[0, T]$ to the remainder system \eqref{Remd-Equ-GR}. Assume the condition \eqref{Sp-3} in Lemma \ref{Lm-Unif-Spatial-Bnd} holds. Then there is a sufficiently small $\eps_0 \in (0,1)$, such that if $0 < \eps < \eps_0$, then there exist two energy functionals $\mathfrak{E}_{N,l} (G_R, E_R, B_R)$ and $\mathfrak{D}_{N,l} (G_R, E_R, B_R)$ such that
	\begin{equation}\label{ED-1}
		\begin{aligned}
			\mathfrak{E}_{N,l} (G_R, E_R, B_R) (t) \sim \mathbb{E}_{N,l} (G_R, E_R, B_R) (t) \,, \\ \mathfrak{D}_{N,l} (G_R, E_R, B_R) (t) \sim \mathbb{D}_{N,l} (G_R, E_R, B_R) (t) \,,
		\end{aligned}
	\end{equation}
    and moreover, the following energy inequality holds
    \begin{equation}\label{ED-2}
    	\begin{aligned}
    		\frac{\d}{\d t} \mathfrak{E}_{N,l} (G_R, E_R, B_R) (t) & + \mathfrak{D}_{N,l} (G_R, E_R, B_R) (t) \\
    		& \leq C_\star \mathfrak{E}_{N,l} (G_R, E_R, B_R) (t) \mathfrak{D}_{N,l} (G_R, E_R, B_R) (t),
    	\end{aligned}
    \end{equation}
    for all $t \in [0,T]$ and for some constant $C_\star > 0$.
\end{proposition}

Once this proposition holds, a continuity argument will conclude our main result, Theorem \ref{Main-Thm}. Indeed, one first chooses a small $\zeta_0 > 0$ such that if $\mathfrak{E}_{N,l} (G_R, E_R, B_R) (0) \leq \zeta_0$,
\begin{equation*}
	C_\star \mathfrak{E}_{N,l} (G_R, E_R, B_R) (0) \leq C_\star \zeta_0 < \tfrac{1}{4} \,.
\end{equation*}
Define a number
\begin{equation*}
	\begin{aligned}
		T_\star = \sup \{ \tau \geq 0 ; \sup_{t \in [0, \tau]} C_\star \mathfrak{E}_{N,l} (G_R, E_R, B_R) (t) < \tfrac{1}{2} \} \geq 0 \,.
	\end{aligned}
\end{equation*}
The continuity of $ \mathfrak{E}_{N,l} (G_R , E_R , B_R ) (t) $ implies $T_\star > 0$.

Claim that $T_\star =  \infty$. If $T_\star < \infty$, the inequality \eqref{ED-2} with the definition of $T_\star$ implies that, for all $t \in [0,T_\star]$,
\begin{equation*}
	\tfrac{\d}{\d t} \mathfrak{E}_{N,l} (G_R, E_R, B_R) (t) + \tfrac{1}{2} \mathfrak{D}_{N,l} (G_R, E_R, B_R) (t) \leq 0 \,.
\end{equation*}
From integrating the above inequality on $[0,t]$, there holds
$$\mathfrak{E}_{N,l} (G_R, E_R, B_R) (t) \leq \mathfrak{E}_{N,l} (G_R, E_R, B_R)(0) \leq \zeta_0 \,, $$
and hence, for all $t \in [0,T_\star]$,
$$ C_\star \mathfrak{E}_{N,l} (G_R, E_R, B_R) (t) \leq C_\star \zeta_0 < \tfrac{1}{4} < \tfrac{1}{2} \,.$$

The continuity of $ \mathfrak{E}_{N,l} (G_R , E_R , B_R ) (t) $ yields that, there is a small $\varsigma > 0$, such that
\begin{equation*}
	\begin{aligned}
		\sup_{t \in [0, T_\star + \varsigma]} C_\star \mathfrak{E}_{N,l} (G_R, E_R, B_R) (t) \leq \tfrac{3}{8} < \tfrac{1}{2} \,,
	\end{aligned}
\end{equation*}
which contradicts to the definition of $T_\star$. Thus, $T_\star = \infty$ and the proof of Theorem \ref{Main-Thm} is completed.

\begin{proof}[Proof of Proposition \ref{Prop-Est}]
Firstly, the number $\eta_0$ in \eqref{Thm-IC-Small} is assumed sufficiently small such that the initial conditions of Lemma \ref{lemm:bnd-NSMF}-\ref{lemm:bnd-linearMaxwl} hold for $s = N+5$ and $M = N+3$. This ensures the validity of all previous lemmas in this section.

We multiply by $\rho$ on the inequalities \eqref{esm:remainder-mix} in Lemma \ref{lemm:remainder-apriori} and \eqref{esm:remainder-mix2} in Lemma \ref{lemm:remainder-apriori2}. We then add the obtained results and $\Upsilon$ times of the inequality \eqref{bnd:Spec-energy-1} with the case $M = N+3$ in Lemma \ref{lemm:bnd-linearMaxwl} to the inequality \eqref{Unif-Spatial-Bnd} in Lemma \ref{Lm-Unif-Spatial-Bnd}. Noticing the relation \eqref{HR-1}, we thereby obtain
\begin{align}\label{Thm-inq-1}
  \no & \tfrac{\d}{\d t} \Big( \mathbb{E}_p(t) + \rho \| \mathbb{P}^\perp G_R \|^2_{H^N_{x,v} (\w^{2l})} + \rho \| \mathbb{P}^\perp G_R \|^2_{\widetilde{\mathbf{H}}^{N+1}_{x,v} (\w^{2l} )} \Big) \\
  \no & \quad + C_6 \mathbb{D}_p(t) + \tfrac{\rho}{2 \eps^2} \| \mathbb{P}^\perp G_R \|^2_{H^N_{x,v} (\w^{2l} \nu) } + \tfrac{\rho}{2 \eps^2} \| \mathbb{P}^\perp G_R \|^2_{\widetilde{\mathbf{H}}^{N+1}_{x,v} (\w^{2l} \nu)} \\
  \no & \leq 2 C \rho \mathbb{D}_p (t)  + C(1+\eps^2) \mathbb{E}_{N,l} (G_R  , E_R , B_R) \mathbb{D}_{N,l} (G_R  , E_R , B_R) + C \eps^2 \| H_R \|^2_{H^N_x L^2_v} \\
  & \quad + C \mathcal{E}^{\IN}_{0,N+2} \| \nabla_v \mathbb{P}^\perp G_R \|^2_{H^N_x L^2_v} + C \sum_{|m| + |\beta| \leq N} \int_{\T} \langle \partial^m_\beta H_R \cdot \w^{2l} \partial^m_\beta \mathbb{P}^\perp G_R \rangle \d x \\
  \no & \quad + C \sum_{|m| \leq N+1} \int_{\T} \langle \partial^m H_R \cdot \partial^m G_R \rangle \d x + C \sum_{\substack{|m| + |\beta| \leq N+1 , \beta \neq 0} }  \int_{\T} \langle \partial^m_\beta H_R \cdot \w^{2l} \partial^m_\beta \mathbb{P}^\perp G_R \rangle \d x \,,
\end{align}
where the constant $\rho > 0$ is small, $\Upsilon > 0$ is large enough to be determined. The functionals $\mathbb{E}_{N,l} (G_R  , E_R , B_R) $ and $ \mathbb{D}_{N,l} (G_R  , E_R , B_R) $ are defined in \eqref{Energ-1} and \eqref{Dspt-Rate} respectively, and the energy $\mathbb{E}_p(t)$ and $\mathbb{D}_p(t)$ are defined as
\begin{align}
  \no \mathbb{E}_p(t) = & \| E_R \|^2_{H^{N-1}_x} + (1 - \delta + \sigma \delta) \| B_R \|^2_{H^{N-1}_x} + ( 1 - \delta ) \| \partial_t B_R \|^2_{H^{N-1}_x} + \delta \| \partial_t B_R +  B_R \|^2_{H^{N-1}_x} \\
     \no & - C_3 \eps \mathcal{A}_N (t) + \eps^2 \| \nabla_x \times \langle \mathbb{P}^\perp G_R \cdot \mathcal{T}_1 \widetilde{\Phi} (v) \rangle \|^2_{H^{N-1}_x } + \Upsilon ( \mathcal{E}_{1,N+3} (t) + \widetilde{C}_{N+3} \mathcal{E}_{0,N+5} (t) )  \\
    \no & + C_4 \big( \| E_R \|^2_{H^{N+1}_x} + \| B_R \|^2_{H^{N+1}_x} + \| G_R \|^2_{H^{N+1}_x L^2_v} \big) + C_5 \big( \mathcal{E}_{1,N} (t) + \widetilde{C}_N \mathcal{E}_{0,N+2} (t) \big) \,,
\end{align}
and
\begin{align}
  \no  \mathbb{D}_p(t) = & \| E_R \|^2_{H^{N-1}_x} + \| \nabla_x B_R \|^2_{H^{N-1}_x} + \| \partial_t B_R \|^2_{H^{N-1}_x} + \tfrac{\Upsilon}{C_6} ( \mathcal{D}_{1,N+3} (t) + \mathcal{D}_{0,N+5} (t) ) \\
  \no & + \| \div_x E_R \|^2_{H^{N}_x} + \| \mathbb{P} G_R \|^2_{H^{N+1}_x L^2_v} + \tfrac{1}{\eps^2} \| \mathbb{P}^\perp G_R \|^2_{H^{N+1}_x L^2_v (\nu)} \,.
\end{align}
If we choose a small $\rho > 0$ such that $ 2 C \rho \leq \tfrac{1}{2} C_6 $, then the term $ 2 C \rho \mathbb{D}_p (t) $ in the right-hand side of \eqref{Thm-inq-1} can be absorbed by the left-hand term $ C_6 \mathbb{D}_p(t) $. Furthermore, the bound \eqref{HR-Square} concerning the term $\eps^2 \| H_R \|^2_{H^N_x L^2_v}$ implies that there is a constant $C_0  > 0$, such that
\begin{equation}\label{Thm-inq-2}
  \begin{aligned}
     C \eps^2 \| H_R \|^2_{H^N_x L^2_v} \leq & C_0 \big( \mathcal{D}_{0,N+5} (t) + \mathcal{D}_{1,N+3} (t) \big) + C_0 ( \eps^2 + \mathcal{E}^{\IN}_{0,N+5} ) \mathbb{D}_{N , 2 \gamma + 1} ( G_R, E_R, B_R ) \\
    & + C_0 (1+\eps^2) \mathbb{E}_{N , 2 \gamma + 1} ( G_R, E_R, B_R ) \mathbb{D}_{N , 2 \gamma + 1} ( G_R, E_R, B_R ) \,.
  \end{aligned}
\end{equation}
We take $\Upsilon = \tfrac{C_6 + C_0}{4} > 0$ satisfying $C_0 \big( \mathcal{D}_{0,N+5} (t) + \mathcal{D}_{1,N+3} (t) \big) \le \tfrac{1}{4} C_6 \mathbb{D}_p (t)$.
Remark that the condition $l \geq 2 \gamma + 1$ is required to satisfy
 \begin{equation*}
   \begin{aligned}
     \mathbb{E}_{N , 2 \gamma + 1} ( G_R, E_R, B_R ) \leq \mathbb{E}_{N , l} ( G_R, E_R, B_R ) \,, \ \mathbb{D}_{N , 2 \gamma + 1} ( G_R, E_R, B_R ) \leq \mathbb{D}_{N , l} ( G_R, E_R, B_R ) \,.
   \end{aligned}
 \end{equation*}

Therefore, one has already constructed an energy functional
\begin{equation}\label{Energ-instant-found}
  \begin{aligned}
    \mathfrak{E}_{N,l} (G_R, E_R , B_R) = \mathbb{E}_p(t) + \rho \| \mathbb{P}^\perp G_R \|^2_{H^N_{x,v} (\w^{2l}) } + \rho \| \mathbb{P}^\perp G_R \|^2_{\widetilde{\mathbf{H}}^{N+1}_{x,v} (\w^{2l} )}
  \end{aligned}
\end{equation}
and the dissipation rate functional
\begin{equation}
	\begin{aligned}
		\mathfrak{D}_{N,l} (G_R, E_R , B_R) = \tfrac{1}{8} C_6 \mathbb{D}_p (t) + \tfrac{\rho}{4 \eps^2} \| \mathbb{P}^\perp G_R \|^2_{H^N_{x,v} (\w^{2l} \nu)} + \tfrac{\rho}{4 \eps^2} \| \mathbb{P}^\perp G_R \|^2_{\widetilde{\mathbf{H}}^{N+1}_{x,v} (\w^{2l} \nu)}
	\end{aligned}
\end{equation}
for $l \geq 2 \gamma + 1$. Observe that the quantity $\mathcal{A}_N (t)$ in \eqref{Sp-4} satisfies $|\mathcal{A}_N (t)| \lesssim \mathbb{E}_{N,l} (G_R, E_R , B_R)$.
Therefore, the first relation in \eqref{ED-1} holds for sufficiently small $\eps > 0$. Obviously,
\begin{equation*}
	\begin{aligned}
		\mathfrak{D}_{N,l} (G_R, E_R , B_R) \sim \mathbb{D}_{N,l} (G_R, E_R , B_R) \,.
	\end{aligned}
\end{equation*}

Note that the term $C \mathcal{E}^{\IN}_{0,N+2} \| \nabla_v \mathbb{P}^\perp G_R \|^2_{H^N_x L^2_v}$ in the right-hand side of \eqref{Thm-inq-1} can be dominated by
\begin{equation*}
  C \mathcal{E}^{\IN}_{0,N+2} \| \mathbb{P}^\perp G_R \|^2_{\widetilde{\mathbf{H}}^{N+1}_{x,v} (\w^{2l} \nu)} \leq C \mathcal{E}^{\IN}_{0,N+5} \eps^2 \mathbb{D}_{N,l} (G_R  , E_R , B_R) \,,
\end{equation*}
which means the this term can be absorbed by the term $C' \mathbb{D}_{N,l} (G_R  , E_R , B_R) $, with a constant $C' = \tfrac{1}{4} \min\{ C_6 , \rho \} > 0$.

Recall that the initial condition \eqref{Thm-IC-Small} implies
\begin{equation*}
	\begin{aligned}
		\mathcal{E}_{0, N+5}^{\IN} \leq \mathbb{E}_{N, l} (G_{R,\eps}^{\IN}, E_{R,\eps}^{\IN}, B_{R,\eps}^{\IN}) \leq \eta_0 \,.
	\end{aligned}
\end{equation*}
Combining the initial condition \eqref{Thm-IC-Small} and the inequality \eqref{HR-Inner-Product} in Lemma \ref{Lm-HR-Inner-Product} together, the above arguments conclude that the following inequality holds
\begin{align}\label{Thm-Inq-3}
  & \tfrac{\d}{\d t} \mathfrak{E}_{N,l} (G_R, E_R, B_R) + 2 \mathfrak{D}_{N,l} (G_R, E_R, B_R) \\\no
  & \leq \widetilde{C} ( \eps + \eps^2 + \sqrt{\eta_0} + \eta_0 ) \mathfrak{D}_{N,l} (G_R, E_R, B_R)
  		+ \widetilde{C} \mathfrak{E}_{N,l} (G_R, E_R, B_R)  \mathfrak{D}_{N,l} (G_R, E_R, B_R),
\end{align}
for $l \geq 2 \gamma + 1 $ and some constant $\widetilde{C} > 0$.

Take sufficiently small constants $\eps_0$ and $\eta_0$, depending only on $\mu$, $\sigma$, $\kappa$ and $N$, such that for all $\eps \in (0, \eps_0)$, then there holds $ \widetilde{C} ( \eps + \eps^2 + \sqrt{\eta_0} + \eta_0 ) \leq 1 $. This concludes Proposition \ref{Prop-Est}.
\end{proof}

\section{Detailed Proofs of Lemmas}
\label{Sec:detailed_lemmas}

\subsection{Estimates on Fluid Variables: Proof of Lemma \ref{Lm-Mic-Mac-decomp}}
\label{Sec:Lmm-MM}

In this subsection, the goal is to deal with the fluid part of $G_R$ in terms of the kinetic part of $G_R$, namely, to prove Lemma \ref{Lm-Mic-Mac-decomp}.

\begin{proof}[Proof of Lemma \ref{Lm-Mic-Mac-decomp}]
	First, we project the first equation of $G_R$ in \eqref{Remd-Equ-GR} into $\mathrm{span}\{ \mathfrak{B} \}$. Noticing $G_R =\mathbb{P} G_R + \mathbb{P}^\perp G_R$, the $G_R$-equation can be rewritten as
	\begin{align}\label{MM-1}
			\eps \partial_t \mathbb{P} G_R + v \cdot \nabla_x \mathbb{P} G_R + \mathcal{T} (v \times B_0) \cdot \nabla_v \mathbb{P} G_R & + \mathcal{T} (v \times B_R) \cdot \nabla_v G_0 \\ \no
			& - E_R \cdot v \mathcal{T}_1 = \Theta(\mathbb{P}^\perp G_R) + \eps H_R \,,
	\end{align}
	where
	\begin{equation}\label{MM-2}
		\Theta(\mathbb{P}^\perp G_R) = - ( \eps \partial_t + v \cdot \nabla_x + \tfrac{1}{\eps} \mathbb{L} + (v \times B_0) \cdot \nabla_v ) \mathbb{P}^\perp G_R \,.
	\end{equation}
	Recalling the definition of the projection $\mathbb{P}$ in \eqref{Projc-P}, a simple calculation implies that the left-hand of the equation \eqref{MM-1} is
	\begin{align}\label{MM-5}
		\no & \eps \partial_t \mathbb{P} G_R + v \cdot \nabla_x \mathbb{P} G_R - E_R \cdot v \mathcal{T}_1 + \mathcal{T} (v \times B_0) \cdot \nabla_v \mathbb{P} G_R + \mathcal{T} (v \times B_R) \cdot \nabla_v G_0 \\
		\no & = \eps \partial_t ( \rho_R^+ - \tfrac{3}{2} \theta_R ) \left( \begin{array}{c} 1 \\ 0 \end{array} \right) + \eps \partial_t ( \rho_R^- - \tfrac{3}{2} \theta_R ) \left( \begin{array}{c} 0 \\ 1 \end{array} \right) \\
		\no & \quad + \sum_{i=1}^3 \big[ \eps \partial_t u_R^i + \partial_i (\rho_R^+ - \tfrac{3}{2} \theta_R) - E_R^i - (u_0 \times B_R)^i - (u_R \times B_0)^i \big] \left( \begin{array}{c} v_i \\ 0 \end{array} \right) \\
		\no & \quad + \sum_{i=1}^3 \big[ \eps \partial_t u_R^i + \partial_i (\rho_R^- - \tfrac{3}{2} \theta_R) + E_R^i + (u_0 \times B_R)^i + (u_R \times B_0)^i \big] \left( \begin{array}{c} v_i \\ 0 \end{array} \right) \\
		\no & \quad + \sum_{i=1}^3 \tfrac{1}{2} \partial_i \theta_R \left( \begin{array}{c} v_i |v|^2 \\ v_i |v|^2 \end{array} \right) + \sum_{i=1}^3 ( \tfrac{1}{2} \eps \partial_t \theta_R + \partial_i u_R^i ) \left( \begin{array}{c} v_i^2 \\ v_i^2 \end{array} \right) \\
		& \quad + \sum_{1 \leq i < j \leq 3} ( \partial_i u_R^j + \partial_j u_R^i ) \left( \begin{array}{c} v_i v_j \\ v_i v_j \end{array} \right) \,.
	\end{align}
	By the definition of the projection operator $\mathcal{P}_{\mathfrak{B}}$ in \eqref{S-3}, one knows that
	\begin{align}\label{MM-6}
		\no \mathcal{P}_{\mathfrak{B}} \Theta (\mathbb{P}^\perp G_R) = & \Theta_R^+ (t,x) \left( \begin{array}{c} 1 \\ 0 \end{array} \right)  + \Theta_R^- (t,x) \left( \begin{array}{c} 0 \\ 1 \end{array} \right)  + \sum_{i=1}^3 \Theta_R^{i+} (t,x) \left( \begin{array}{c} v_i \\ 0 \end{array} \right) + \sum_{i=1}^3 \Theta_R^{i-} (t,x) \left( \begin{array}{c} 0 \\ v_i \end{array} \right) \\
		& + \sum_{i=1}^3 \Theta_R^i (t,x) \left( \begin{array}{c} v_i^2 \\ v_i^2 \end{array} \right) + \sum_{i=1}^3 \widetilde{\Theta}_R^i (t,x) \left( \begin{array}{c} v_i |v|^2 \\ v_i |v|^2 \end{array} \right) + \sum_{1 \leq i < j \leq 3} \Theta_R^{ij}(t,x) \left( \begin{array}{c} v_i v_j \\ v_i v_j \end{array} \right) \,,
	\end{align}
	and
	\begin{align}\label{MM-7}
		\no \mathcal{P}_{\mathfrak{B}} H_R = & h_R^+ (t,x) \left( \begin{array}{c} 1 \\ 0 \end{array} \right)  + h_R^- (t,x) \left( \begin{array}{c} 0 \\ 1 \end{array} \right)  + \sum_{i=1}^3 h_R^{i+} (t,x) \left( \begin{array}{c} v_i \\ 0 \end{array} \right) + \sum_{i=1}^3 h_R^{i-} (t,x) \left( \begin{array}{c} 0 \\ v_i \end{array} \right) \\
		& + \sum_{i=1}^3 h_R^i (t,x) \left( \begin{array}{c} v_i^2 \\ v_i^2 \end{array} \right) + \sum_{i=1}^3 \widetilde{h}_R^i (t,x) \left( \begin{array}{c} v_i |v|^2 \\ v_i |v|^2 \end{array} \right) + \sum_{1 \leq i < j \leq 3} h_R^{ij}(t,x) \left( \begin{array}{c} v_i v_j \\ v_i v_j \end{array} \right) \,.
	\end{align}
	Consequently, together with the equations \eqref{MM-5}-\eqref{MM-7}, projecting equation \eqref{MM-1} into $\mathrm{span}\{ \mathfrak{B} \}$ implies
	\begin{equation}\label{MM-8}
		\left\{
		\begin{array}{l}
			\eps \partial_t (\rho_R^\pm - \tfrac{3}{2} \theta_R) = \Theta_R^\pm + \eps h_R^\pm \,, \\[3pt]
			\eps \partial_t u_R^i + \partial_i ( \rho_R^\pm - \tfrac{3}{2} \theta_R ) \mp E_R^i \mp (u_0 \times B_R)^i \mp (u_R \times B_0)^i = \Theta_R^{i \pm} + \eps h_R^{i \pm} \,, \\[3pt]
			\tfrac{1}{2} \partial_i \theta_R = \widetilde{\Theta}_R^i + \eps \widetilde{h}_R^i \,, \qquad \tfrac{1}{2} \eps \partial_t \theta_R + \partial_i u_R^i = \Theta_R^i + \eps h_R^i \,, \\[3pt]
			\partial_i u_R^j + \partial_j u_R^i = \Theta_R^{ij} + \eps h_R^{ij} \ (i \neq j) \,.
		\end{array}
		\right.
	\end{equation}

Second, we project the $G_R$-equation into $\ker \mathbb{L}$. By the splitting $G_R = \mathbb{P} G_R + \mathbb{P}^\perp G_R$, the $G_R$-equation in \eqref{Remd-Equ-GR} can be rewritten as
	\begin{equation}\label{MM-9}
		\begin{aligned}
			& \eps \partial_t G_R + v \cdot \nabla_x \mathbb{P} G_R + \mathcal{T} (v \times B_0) \cdot \nabla_v \mathbb{P} G_R + \mathcal{T} (v \times B_R) \cdot \nabla_v G_0 - E_R \cdot v \mathcal{T}_1 + \tfrac{1}{\eps} \mathbb{L} G_R \\[3pt]
			& = - ( v \cdot \nabla_x \mathbb{P}^\perp G_R + \mathcal{T} (v \times B_0) \cdot \nabla_v \mathbb{P}^\perp G_R ) + \eps H_R \,.
		\end{aligned}
	\end{equation}
	Projecting it into $\ker \mathbb{L}$ reads
	\begin{equation}\label{MM-10}
		\left\{
		\begin{array}{l}
			\eps \partial_t \rho_R^+ + \div_x u_R = \eps \left\langle H_R \cdot \mathcal{T}^+ \right\rangle \,, \quad \eps \partial_t \rho_R^- + \div_x u_R = \eps \left\langle H_R \cdot \mathcal{T}^- \right\rangle \,, \\[3pt]
			\eps \partial_t u_R + \nabla_x \big( \tfrac{\rho_R^+ + \rho_R^-}{2} \big) + \nabla_x \theta_R = \frac{\eps}{2} \left\langle H_R \cdot \mathcal{T}_2 v \right\rangle \\[3pt]
			\qquad\qquad - \frac{1}{2} \left\langle [ v \cdot \nabla_x \mathbb{P}^\perp G_R + \mathcal{T} (v \times B_0) \cdot \nabla_v \mathbb{P}^\perp  G_R ] \cdot \mathcal{T}_2 v \right\rangle \,, \\[3pt]
			\eps \partial_t \theta_R + \frac{2}{3} \div_x u_R = \frac{\eps}{2} \left\langle H_R \cdot \mathcal{T}_2 (\tfrac{|v|^2}{3} - 1) \right\rangle - \frac{1}{2} \left\langle v \cdot \nabla_x \mathbb{P}^\perp G_R \cdot \mathcal{T}_2 (\tfrac{|v|^2}{3} - 1) \right\rangle \,.
		\end{array}
		\right.
	\end{equation}
    where $\mathcal{T}^+ = (1,0)^\top$, $\mathcal{T}^- = (0,1)^\top$, $\mathcal{T}_2 = (1,1)^\top$.

    Third, together with the equations \eqref{MM-8} and \eqref{MM-10}, the inequality \eqref{Mic-Mac-Inq} in Lemma \ref{Lm-Mic-Mac-decomp} will be derived.

\noindent\textbf{Estimates on $u_R$.} From the last two equations of \eqref{MM-8}, we derive
	\begin{align}\label{MM-11}
		\no & - \Delta_x \partial^m u_R^i = - \sum_{j \neq i } \partial_{jj} \partial^m u_R^i - \partial_{ii} \partial^m u_R^i \\
		\no & = - \sum_{j \neq i} \partial_j \partial^m ( - \partial_i u_R^j + \Theta_R^{ij} + \eps h_R^{ij} ) - \partial_i \partial^m ( - \tfrac{1}{2} \eps \partial_t \theta_R + \Theta_R^i + \eps h_R^i ) \\
		\no & = \sum_{j \neq i} \partial_i \partial^m ( - \tfrac{1}{2} \eps \partial_t \theta_R + \Theta_R^j + \eps h_R^j ) + \partial_i \partial^m ( \tfrac{1}{2} \eps \partial_t \theta_R - \Theta_R^i - \eps h_R^i ) - \sum_{j \neq i} \partial_j \partial^m ( \Theta_R^{ij} + \eps h_R^{ij} ) \\
		& = \partial_{ii} \partial^m u_R^i + \sum_{j \neq i } ( \partial_i \partial^m \Theta_R^j  - \partial_j \partial^m \Theta_R^{ij}  )- 2 \partial_i \partial^m \Theta_R^i
		 + \eps \Big[ \sum_{j \neq i} ( \partial_i \partial^m h_R^j - \partial_j \partial^m h_R^{ij} ) - 2 \partial_i \partial^m h_R^i \Big] \,.
	\end{align}
	Noticing that there are certain linear combinations $\zeta_{ij}$ of the basis $\mathfrak{B}$ such that
		\begin{multline}\label{MM-12}
			\no \sum_{j \neq i } ( \partial_i \partial^m \Theta_R^j  - \partial_j \partial^m \Theta_R^{ij}  )- 2 \partial_i \partial^m \Theta_R^i \\
			= \sum_{j=1}^3 \partial_j \partial^m \Big\{ \langle ( - \eps \partial_t \mathbb{P}^\perp G_R - v \cdot \nabla_x \mathbb{P}^\perp G_R
			- \mathcal{T} (v \times B_0) \cdot \nabla_v \mathbb{P}^\perp G_R ) \cdot \zeta_{ij} \rangle - \tfrac{1}{\eps} \langle \mathbb{L} (\mathbb{P}^\perp G_R) \cdot \zeta_{ij} \rangle \Big\} \,,
		\end{multline}
	and
	\begin{equation}\label{MM-13}
		\sum_{j \neq i} ( \partial_i \partial^m h_R^j - \partial_j \partial^m h_R^{ij} ) - 2 \partial_i \partial^m h_R^i = \sum_{j=1}^3 \partial_j \partial^m \langle H_R \cdot \zeta_{ij} \rangle \,,
	\end{equation}
	hence one has
	\begin{align}\label{MM-14}
			- \Delta_x \partial^m u_R^i - \partial_{ii} \partial^m u_R^i
		= & \sum_{j=1}^3 \partial_j \Big\{ - \eps \langle \partial_t \partial^m \mathbb{P}^\perp G_R \cdot \zeta_{ij} \rangle - \langle v \cdot \nabla_x \partial^m \mathbb{P}^\perp G_R \cdot \zeta_{ij} \rangle \\ \no
		& - \langle \partial^m [ \mathcal{T} (v \times B_0) \cdot \nabla_v \mathbb{P}^\perp G_R ] \cdot \zeta_{ij} \rangle - \tfrac{1}{\eps} \langle \mathbb{L} \partial^m \mathbb{P}^\perp G_R \cdot \zeta_{ij} \rangle \Big\} + \eps \sum_{j=1}^3 \partial_j \partial^m \langle H_R \cdot \zeta_{ij} \rangle \,.
	\end{align}

Multiplying by $\partial^m u_R^i$ in \eqref{MM-14}, integrating by parts over $x \in \T$ and summing up for $1 \leq i \leq 3$ imply that
	\begin{align}\label{MM-15}
		\no \| \nabla_x \partial^m u_R \|^2_{L^2_x} + \| \div_x \partial^m u_R \|^2_{L^2_x}
	\leq & \sum_{i,j=1}^3 \int_{\T} - \langle \eps \partial_t \mathbb{P}^\perp \partial_j \partial^m G_R \cdot \zeta_{ij} \rangle \partial^m u_R^i \d x \\
		& + \tfrac{C}{\eps} \big( \| \nabla_x \partial^m \mathbb{P}^\perp G_R \|_{L^2_{x,v}} + \| \partial^m \mathbb{P}^\perp G_R \|_{L^2_{x,v}} \big) \| \nabla_x \partial^m u_R \|_{L^2_x} \\
		\no & + C \eps \| \mathcal{P}_{\mathfrak{B}} \partial^m H_R \|_{L^2_{x,v}} \| \nabla_x \partial^m u_R \|_{L^2_x} \\
		\no & + \sum_{i,j=1}^3 \int_{\T} \langle \partial^m [ \mathcal{T} (v \times B_0) \cdot \nabla_v \mathbb{P}^\perp G_R ] \cdot \zeta_{ij} \rangle \partial_j \partial^m u_R^i \d x \,.
	\end{align}
The last term in \eqref{MM-15} is calculated by
	\begin{align*}
		& \sum_{i,j=1}^3 \int_{\T} \langle \partial^m [ \mathcal{T} (v \times B_0) \cdot \nabla_v \mathbb{P}^\perp G_R ] \cdot \zeta_{ij} \rangle \partial_j \partial^m u_R^i \d x \\
		& = - \sum_{i,j=1}^3 \sum_{m' \leq m} C_m^{m'} \int_{\T} \int_{\R^3} \epsilon_{kpq} \partial^{m-m'} B_0^q \mathcal{T} \mathbb{P}^\perp \partial^{m'} G_R \cdot ( v_q \partial_{v_k} \zeta_{ij} ) M \d v \partial_j \partial^m u_R^i \d x \\
		& \leq C \sum_{m' \leq m} \| \partial^{m-m'} B_0 \|_{L^\infty_x} \| \partial^{m'} \mathbb{P}^\perp G_R \|_{L2_{x,v}} \| \nabla_x \partial^m u_R \|_{L^2_x} \,,
	\end{align*}
where the symbol $\epsilon_{kpq}$ is defined in \eqref{eps-ijk}. Thus, the inequality \eqref{MM-15} reduces to
	\begin{align}\label{MM-16}
		\no & \| \nabla_x \partial^m u_R \|^2_{L^2_x} + \| \div_x \partial^m u_R \|^2_{L^2_x} \\
		\no & \leq \sum_{i,j = 1} \int_{\T} \langle - \eps \partial_t \partial_j \partial^m \mathbb{P}^\perp G_R \cdot \zeta_{ij} \rangle \partial^m u_R^i \d x + C \eps \| \mathcal{P}_{\mathfrak{B}} \partial^m H_R \|_{L^2_{x,v}} \| \nabla_x \partial^m u_R \|_{L^2_x} \\
		& \quad + \tfrac{C}{\eps} \big( \| \nabla_x \partial^m \mathbb{P}^\perp G_R \|_{L^2_{x,v}} + \| \partial^m \mathbb{P}^\perp G_R \|_{L^2_{x,v}} \big) \| \nabla_x \partial^m u_R \|_{L^2_x} \\
		\no & \quad + C \sum_{m' \leq m} \| \partial^{m-m'} B_0 \|_{L^\infty_x} \| \partial^{m'} \mathbb{P}^\perp G_R \|_{L^2_{x,v}} \| \nabla_x \partial^m u_R \|_{L^2_x} \,.
	\end{align}

Next we estimate the first term in the right-hand side of \eqref{MM-16}. We write that
		\begin{align}\label{MM-17}
			& \sum_{i,j=1}^3 \int_{\T} \langle - \eps \partial_t \partial_j \partial^m \mathbb{P}^\perp G_R \cdot \zeta_{ij} \rangle \partial^m u_R^i \d x \\\no
			& = \sum_{i,j=1}^3 \tfrac{\d}{\d t} \int_{\T} \eps \langle - \partial_j \partial^m \mathbb{P}^\perp G_R \cdot \zeta_{ij} \rangle \partial^m u_R^i \d x + \sum_{i,j=1}^3 \int_{\T} \langle \partial_j \partial^m \mathbb{P}^\perp G_R \cdot \zeta_{ij} \rangle \eps \partial^m \partial_t u_R^i \d x \,.
		\end{align}
	Together with the third equation of \eqref{MM-10}, the last term in \eqref{MM-17} can be bounded by
	\begin{align*}
		& \sum_{i,j=1}^3 \int_{\T} \langle \partial_j \partial^m \mathbb{P}^\perp G_R \cdot \zeta_{ij} \rangle  \partial^m \Big[ - \partial_i ( \tfrac{\rho_R^+ + \rho_R^-}{2} ) - \partial_i \theta_R + \tfrac{\eps}{2} \left\langle H_R \cdot \left( \begin{array}{c} v_i \\ v_i \end{array} \right) \right\rangle \\
		& \qquad - \tfrac{1}{2} \left\langle v \cdot \nabla_x \mathbb{P}^\perp G_R \cdot \left( \begin{array}{c} v_i \\ v_i \end{array} \right) \right\rangle - \tfrac{1}{2} \left\langle \mathcal{T} (v \times B_0) \cdot \nabla_v \mathbb{P}^\perp G_R \cdot \left( \begin{array}{c} v_i \\ v_i \end{array} \right) \right\rangle \Big] \d x \\
		& \leq C \| \nabla_x \partial^m \mathbb{P}^\perp G_R \|_{L^2_{x,v}} \Big\{ \| \nabla_x \partial^m \rho_R^+ \|_{L^2_x} + \| \nabla_x \partial^m \rho_R^- \|_{L^2_x} + \| \nabla_x \partial^m \theta_R \|_{L^2_x} + \eps \| \mathcal{P}_{\mathfrak{B}} \partial^m H_R \|_{L^2_{x,v}} \Big\} \\
		& \quad + C \| \nabla_x \partial^m \mathbb{P}^\perp G_R \|^2_{L^2_{x,v}} + \sum_{i,j=1}^3 \int_{\T} \langle \partial_j \partial^m \mathbb{P}^\perp G_R \cdot \zeta_{ij} \rangle \partial^m \langle \eps_{ipq} B_0^q v_p \mathcal{T} \mathbb{P}^\perp G_R \cdot \mathcal{T}_2 \rangle \d x \,,
	\end{align*}
	with $\mathcal{T}_2 = (1,1)^\top$. Note that the last term in the above inequality can be bounded by
	\begin{equation*}
		C \sum_{m'\leq m} \| \nabla_x \partial^m \mathbb{P}^\perp G_R \|_{L^2_{x,v}} \| \partial^{m'} \mathbb{P}^\perp G_R \|_{L^2_{x,v}} \| \partial^{m-m'} B_0 \|_{L^\infty_x} \,.
	\end{equation*}
	Then, there holds
	\begin{align}\label{MM-18}
		\no & \sum_{i,j=1}^3 \int_{\T} \langle \partial_j \partial^m \mathbb{P}^\perp G_R \cdot \zeta_{ij} \rangle \eps \partial^m \partial_t u_R^i \d x \\
		& \leq C \| \nabla_x \partial^m \mathbb{P}^\perp G_R \|_{L^2_{x,v}} \Big( \| \nabla_x \partial^m \rho_R^+ \|_{L^2_x} + \| \nabla_x \partial^m \rho_R^- \|_{L^2_x} + \| \nabla_x \partial^m \theta_R \|_{L^2_x} \\
		\no & \quad + \| \nabla_x \partial^m \mathbb{P}^\perp G_R \|_{L^2_{x,v}} + \eps \| \mathcal{P}_{\mathfrak{B}} \partial^m H_R \|_{L^2_{x,v}} + \sum_{m'\leq m}  \| \partial^{m'} \mathbb{P}^\perp G_R \|_{L^2_{x,v}} \| \partial^{m-m'} B_0 \|_{L^\infty_x} \,.
	\end{align}
	Consequently, by combining the relations \eqref{MM-16}, \eqref{MM-17}, \eqref{MM-18} and the Young's inequality together, we gain that, for $\eps \in (0,1]$,
	\begin{align}\label{MM-19}
		\no \tfrac{1}{2} \| & \nabla_x \partial^m u_R \|^2_{L^2_x} + \| \div_x \partial^m u_R \|^2_{L^2_x} \\
		\no \leq & \tfrac{\d}{\d t} \sum_{i,j=1}^3 \int_{\T} \eps \langle \partial^m \mathbb{P}^\perp G_R \cdot \zeta_{ij} \rangle \partial_j \partial^m u_R^i \d x + C \eps^2 \| \mathcal{P}_{\mathfrak{B}} \partial^m H_R \|^2_{L^2_{x,v}} \\
		\no & + \frac{C}{\eps^2} \big( \| \nabla_x \partial^m \mathbb{P}^\perp G_R \|^2_{L^2_{x,v}} + \| \partial^m \mathbb{P}^\perp G_R \|^2_{L^2_{x,v}} \big) + C \sum_{m' \leq m} \| \partial^{m-m'} B_0 \|^2_{L^\infty_x} \| \partial^{m'} \mathbb{P}^\perp G_R \|^2_{L^2_{x,v}} \\
		& + C \eps^2 \big( \| \nabla_x \rho_R^+ \|^2_{L^2_x} + \| \nabla_x \rho_R^- \|^2_{L^2_x} + \| \nabla_x \theta_R \|^2_{L^2_x} \big) \,.
	\end{align}

    \noindent\textbf{Estimates on $\theta_R$.} From the third equation of \eqref{MM-8}, i.e., $\frac{1}{2} \partial_i \theta_R = \widetilde{\Theta}_R^i + \eps \widetilde{h}_R^i$, we can infer that
	\begin{align*}
		- \tfrac{1}{2} \Delta_x \partial^m \theta_R = & - \tfrac{1}{2} \sum_{i=1}^3 \partial_{ii} \partial^m \theta_R = - \sum_{i=1}^3 \partial_i \partial^m (  \tfrac{1}{2} \partial_i \theta_R ) =  - \sum_{i=1}^3 \partial_i \partial^m ( \widetilde{\Theta}_R^i + \eps \widetilde{h}_R^i ) \,.
	\end{align*}
	Observe that there exist some certain linear combinations $\zeta_i \ (1 \leq i \leq 3)$ of $\mathfrak{B}$, such that
	\begin{align*}
		- \sum_{i=1}^3 \partial_i \partial^m ( \widetilde{\Theta}_R^i + \widetilde{h}_R^i ) = & \sum_{i=1}^3 \partial_i \partial^m \Big\{ \langle - ( \eps \partial_t \mathbb{P}^\perp G_R + v \cdot \nabla_x \mathbb{P}^\perp G_R ) \cdot \zeta_i \rangle \\
		& \qquad - \langle \mathcal{T} (v \times B_0) \cdot \nabla_v \mathbb{P}^\perp G_R \cdot \rangle - \frac{1}{\eps} \langle \mathbb{L} ( \mathbb{P}^\perp G_R) \cdot \zeta_i \rangle + \langle H_R \cdot \zeta_i \rangle \Big\} \,.
	\end{align*}
	There thereby holds
		\begin{align}\label{MM-20}
			\no - \tfrac{1}{2} \Delta_x \partial^m \theta_R = & \sum_{i=1}^3 \partial_i \partial^m \Big\{ \langle - ( \eps \partial_t \mathbb{P}^\perp G_R + v \cdot \nabla_x \mathbb{P}^\perp G_R  ) \cdot \zeta_i \rangle - \tfrac{1}{\eps} \langle \mathbb{L} ( \mathbb{P}^\perp G_R) \cdot \zeta_i \rangle \\
			& - \langle \mathcal{T} (v \times B_0) \cdot \nabla_v \mathbb{P}^\perp G_R \cdot \zeta_i \rangle \Big\} + \eps \sum_{i=1}^3 \partial_i \partial^m \langle H_R \cdot \zeta_i \rangle \,.
		\end{align}
	From taking $L^2_x$-inner product with $\partial^m \theta_R$ in \eqref{MM-20}, we get
	\begin{align}\label{MM-21}
		\no \tfrac{1}{2} \| \nabla_x \partial^m \theta_R \|^2_{L^2_x}
		& \leq \sum_{i=1}^3 \int_{\T} \langle - \eps \partial_t \mathbb{P}^\perp \partial_i \partial^m G_R \cdot \zeta_i \rangle \partial^m \theta_R \d x + \tfrac{C}{\eps} \| \partial^m \mathbb{P}^\perp G_R \|_{L^2_{x,v}} \| \nabla_x \partial^m \theta_R \|_{L^2_x} \\
		\no & + C \| \nabla_x \partial^m \mathbb{P}^\perp G_R \|_{L^2_{x,v}} \| \nabla_x \partial^m \theta_R \|_{L^2_x} + C \eps \| \mathcal{P}_{\mathfrak{B}} \partial^m H_R \|_{L^2_{x,v}} \| \nabla_x \partial^m \theta_R \|_{L^2_x} \\
		&  + C \sum_{m' \leq m} \| \partial^{m-m'} B_0 \|_{L^\infty_x} \| \partial^{m'} \mathbb{P}^\perp G_R \|_{L^2_{x,v}} \| \nabla_x \partial^m \theta_R \|_{L^2_x} \,.
	\end{align}

Recalling the fourth equation in \eqref{MM-10}, the first term in the right-hand side of \eqref{MM-21} reduces to
	\begin{align}\label{MM-22}
		\no \sum_{i=1}^3 & \int_{\T} \langle - \eps \partial_t \mathbb{P}^\perp \partial_i \partial^m G_R \cdot \zeta_i \rangle \partial^m \theta_R \d x \\
		\no = & \tfrac{\d}{\d t} \sum_{i=1}^3 \int_{\T} \langle \eps \partial^m \mathbb{P}^\perp G_R \cdot \zeta_i \rangle \partial_i \partial^m \theta_R \d x - \tfrac{2}{3} \sum_{i=1}^3 \int_{\T} \div_x \partial^m u_R \langle \partial_i \partial^m \mathbb{P}^\perp G_R \cdot \zeta_i \rangle \d x \\
		& - \tfrac{1}{2} \sum_{i=1}^3 \int_{\T} \langle \partial_i \partial^m \mathbb{P}^\perp G_R \cdot \zeta_i \rangle \langle v \cdot \nabla_x \partial^m \mathbb{P}^\perp G_R \cdot \mathcal{T}_2 (\tfrac{|v|^2}{3} - 1) \rangle \d x \\
		\no & + \tfrac{\eps}{2} \sum_{i=1}^3 \int_{\T} \langle \partial_i \partial^m \mathbb{P}^\perp G_R \cdot \zeta_i \rangle \langle \partial^m H_R \cdot \mathcal{T}_2 ( \tfrac{|v|^2}{3} - 1 ) \rangle \d x \\
		\no \leq & \tfrac{\d}{\d t} \sum_{i=1}^3 \int_{\T} \langle \eps \partial^m \mathbb{P}^\perp G_R \cdot \zeta_i \rangle \partial_i \partial^m \theta_R \d x + C \| \nabla_x \partial^m \mathbb{P}^\perp G_R \|^2_{L^2_{x,v}} \\
		\no & + C \| \nabla_x \partial^m \mathbb{P}^\perp G_R \|_{L^2_{x,v}} \| \nabla_x \partial^m u_R \|_{L^2_x} + C \eps \| \nabla_x \partial^m \mathbb{P}^\perp G_R \|_{L^2_{x,v}} \| \mathcal{P}_{\mathfrak{B}} \partial^m H_R \|_{L^2_{x,v}} \,.
	\end{align}
	From plugging \eqref{MM-22} into \eqref{MM-21} and using Young's inequality, it follows that
	\begin{align}\label{MM-23}
		\no \tfrac{1}{4} \| \nabla_x \partial^m \theta_R \|^2_{L^2_x}
		\leq & \tfrac{\d}{\d t} \sum_{i=1}^3 \int_{\T} \langle \eps \partial^m \mathbb{P}^\perp G_R \cdot \zeta_i \rangle \partial_i \partial^m \theta_R \d x + \tfrac{C}{\eps^2} \big( \| \partial^m \mathbb{P}^\perp G_R \|^2_{L^2_{x,v}} + \| \nabla_x \partial^m \mathbb{P}^\perp G_R \|^2_{L^2_{x,v}} \big) \\
		& + C \eps^2 \big( \| \nabla_x \partial^m u_R \|^2_{L^2_x} + \| \mathcal{P}_{\mathfrak{B}} \partial^m H_R \|^2_{L^2_{x,v}} \big) + C \sum_{m' \leq m} \| \partial^{m-m'} B_0 \|^2_{L^\infty_x} \| \partial^{m'} \mathbb{P}^\perp G_R \|^2_{L^2_{x,v}} \,.
	\end{align}

\noindent\textbf{Estimates on $\rho_R^\pm$.} Next we focus on the quantities $\| \nabla_x \partial^m \rho_R^\pm \|^2_{L^2_x}$. It follows from the second equation of \eqref{MM-8},
	\begin{align}\label{MM-24}
		\no - \Delta_x \partial^m \rho_R^\pm \pm \div_x \partial^m E_R \pm \div_x \partial^m ( u_0 & \times B_R ) \pm \div_x \partial^m (u_R \times B_0) \\
		= & \sum_{i=1}^3 \partial_i \partial^m ( \eps \partial_t u_R^i - \frac{3}{2} \partial_i \theta_R - \Theta_R^{i\pm} - \eps h_R^{i\pm} ) \,.
	\end{align}
	Note that there are some fixed linear combinations $\zeta_{i\pm} \ (1\leq i \leq 3)$ of $\mathfrak{B}$, such that
	\begin{align*}
		- \sum_{i=1}^3 \partial_i \partial^m \Theta_R^{i\pm} - \sum_{i=1}^3 \partial_i \partial^m h_R^{i \pm} = \sum_{i=1}^3 \partial_i \partial^m \Big\{ \langle - ( \eps \partial_t \mathbb{P}^\perp G_R + v \cdot \nabla_x \mathbb{P}^\perp G_R ) \cdot \zeta_{i\pm} \rangle \\
		- \frac{1}{\eps} \langle \mathbb{L} ( \mathbb{P}^\perp G_R ) \cdot \zeta_{i\pm} \rangle - \langle \mathcal{T} (v \times B_0) \cdot \nabla_v \mathbb{P}^\perp G_R \cdot \zeta_{i\pm} \rangle \Big\} + \sum_{i=1}^3 \partial_i \partial^m \langle H_R \cdot \zeta_{i \pm} \rangle \,.
	\end{align*}
	Then we have
	\begin{align}\label{MM-25}
		\no - \Delta_x \partial^m \rho_R^\pm \pm \div_x \partial^m E_R = \mp \div_x \partial^m ( u_0 \times B_R + u_R \times B_0 ) + \sum_{i=1}^3 \eps \partial_t \partial_i \partial^m u_R^i \\
		+ \sum_{i=1}^3 \langle - \eps \partial_t \partial^m \partial_i \mathbb{P}^\perp G_R \cdot \zeta_{i\pm} \rangle - \sum_{i=1}^3 \langle v \cdot \nabla_x \partial_i \partial^m \mathbb{P}^\perp G_R \cdot \zeta_{i \pm} \rangle + \eps \sum_{i=1}^3 \partial_i \partial^m \langle H_R \cdot \zeta_{i \pm} \rangle \\
		\no - \frac{1}{\eps} \sum_{i=1}^3 \langle \mathbb{L} \partial_i \partial^m \mathbb{P}^\perp G_R \cdot \zeta_{i \pm} \rangle - \sum_{i=1}^3 \partial_i \partial^m \langle \mathcal{T} (v \times B_0) \cdot \nabla_v \mathbb{P}^\perp G_R \cdot \zeta_{i\pm} \rangle \,.
	\end{align}
Taking $L^2_x$-inner product by multiplying by $\partial^m \rho_R^\pm$, integrating by parts over $x \in \T$ and using Young's inequality, we obtain
	\begin{align}\label{MM-26}
		\no \tfrac{1}{2} \| & \nabla_x \partial^m \rho_R^\pm \|^2_{L^2_x} \pm \int_{\T} \div_x \partial^m E_R \partial^m \rho_R^\pm \d x \\
		\no \leq & \sum_{i=1}^3 \int_{\T} \eps \partial_t \partial_i \partial^m u_R^i \cdot \partial^m \rho_R^\pm \d x + \sum_{i=1}^3 \int_{\T} \langle - \eps \partial_t \partial_i \partial^m \mathbb{P}^\perp G_R \cdot \zeta_{i\pm} \rangle \partial^m \rho_R^\pm \d x \\
		& + \tfrac{C}{\eps^2} \big( \| \partial^m \mathbb{P}^\perp G_R \|^2_{L^2_{x,v}} + \| \nabla_x \partial^m \mathbb{P}^\perp G_R \|^2_{L^2_{x,v}} \big) + C \eps^2 \| \mathcal{P}_{\mathfrak{B}} \partial^m H_R \|^2_{L^2_{x,v}} \\
		\no & + C \sum_{m' \leq m} \big( \| \partial^{m-m'} u_0 \|^2_{L^\infty_x} + \| \partial^{m-m'} B_0 \|^2_{L^\infty_x} \big) ( \| \partial^{m'} \mathbb{P}^\perp G_R \|^2_{L^2_{x,v}} + \| \partial^{m'} u_R \|^2_{L^2_x} + \| \partial^{m'} B_R \|^2_{L^2_x} ) \,.
	\end{align}
	Recalling the equation \eqref{MM-10}, the first two terms in the right-hand side of \eqref{MM-26} can be dealt with as follows:
	\begin{align}\label{MM-27}
		\no \sum_{i=1}^3 \int_{\T} \eps \partial_t \partial_i \partial^m u_R^i \cdot \partial^m \rho_R^\pm \d x \leq & \tfrac{\d}{\d t} \int_{\T} - \eps \partial^m u_R \cdot \nabla_x \partial^m \rho_R^\pm \d x + \| \div_x \partial^m u_R \|^2_{L^2_x} \\
		& + C \eps \| \mathcal{P}_{\mathfrak{B}} \partial^m H_R \|_{L^2_{x,v}} \| \div_x \partial^m u_R \|_{L^2_x} \,,
	\end{align}
	and
	\begin{align}\label{MM-28}
		& \sum_{i=1}^3 \int_{\T} \langle - \eps \partial_t \partial_i \partial^m \mathbb{P}^\perp G_R \cdot \zeta_{i\pm} \rangle \partial^m \rho_R^\pm \d x \\ \no
	& \leq \tfrac{\d}{\d t} \sum_{i=1}^3 \int_{\T} \langle \eps \partial^m \mathbb{P}^\perp G_R \cdot \zeta_{i\pm} \rangle \cdot \partial_i \partial^m \rho_R^\pm \d x
		+ C \| \nabla_x \partial^m \mathbb{P}^\perp G_R \|_{L^2_{x,v}} \| \div_x \partial^m u_R \|_{L^2_x} \\
		\no & \quad + C \eps \| \nabla_x \partial^m \mathbb{P}^\perp G_R \|_{L^2_{x,v}} \| \mathcal{P}_{\mathfrak{B}} \partial^m H_R \|_{L^2_{x,v}} \,.
	\end{align}
Thus, combining the inequalities \eqref{MM-26}-\eqref{MM-28} yields
	\begin{align}\label{MM-29}
		\no \tfrac{1}{2} \| & \nabla_x \partial^m \rho_R^\pm \|^2_{L^2_x} \pm \int_{\T} \div_x \partial^m E_R \partial^m \rho_R^\pm \d x \\
		\leq &  \tfrac{\d}{\d t} \Big\{ \int_{\T} - \eps \partial^m u_R \cdot \nabla_x \partial^m \rho_R^\pm \d x + \sum_{i=1}^3 \int_{\T} \langle \eps \partial^m \mathbb{P}^\perp G_R \cdot \zeta_{i\pm} \rangle \cdot \partial_i \partial^m \rho_R^\pm \d x \Big\} \\
		\no & + 2 \| \div_x \partial^m E_R \|^2_{L^2_x} +  \tfrac{C}{\eps^2} \big( \| \partial^m \mathbb{P}^\perp G_R \|^2_{L^2_{x,v}} + \| \nabla_x \partial^m \mathbb{P}^\perp G_R \|^2_{L^2_{x,v}} \big) + C \eps^2 \| \mathcal{P}_{\mathfrak{B}} \partial^m H_R \|^2_{L^2_{x,v}} \\
		\no & + C \sum_{m' \leq m} \big( \| \partial^{m-m'} u_0 \|^2_{L^\infty_x} + \| \partial^{m-m'} B_0 \|^2_{L^\infty_x} \big) ( \| \partial^{m'} \mathbb{P}^\perp G_R \|^2_{L^2_{x,v}} + \| \partial^{m'} u_R \|^2_{L^2_x} + \| \partial^{m'} B_R \|^2_{L^2_x} ) \,.
	\end{align}

Recalling the last second equation of \eqref{Remd-Equ-GR}, i.e., $\div_x E_R = \langle G_R \cdot \mathcal{T}_1 \rangle = \rho_R^+ - \rho_R^-$, it follows from summing up the two cases of $``\pm"$ in \eqref{MM-29}, that
	\begin{align}\label{MM-30}
		\no \tfrac{1}{2} \| & \nabla_x \partial^m \rho_R^+ \|^2_{L^2_x} + \tfrac{1}{2} \| \nabla_x \partial^m \rho_R^- \|^2_{L^2_x} + \| \div_x \partial^m E_R \|^2_{L^2_x} \\
		\leq & \tfrac{\d}{\d t} \Big\{ \int_{\T} - \eps \partial^m u_R \cdot \nabla_x \partial^m (\rho_R^+ + \rho_R^-) \d x \\
		\no & + \sum_{i=1}^3 \Big[ \int_{\T} \langle \eps \partial^m \mathbb{P}^\perp G_R \cdot \zeta_{i+} \rangle \partial_i \partial^m \rho_R^+ \d x  + \int_{\T} \langle \eps \partial^m \mathbb{P}^\perp G_R \cdot \zeta_{i-} \rangle \partial_i \partial^m \rho_R^- \d x \Big] \Big\} \\
		\no & + 4 \| \div_x \partial^m u_R \|^2_{L^2_x} + \frac{C}{\eps^2} \big( \| \partial^m \mathbb{P}^\perp G_R \|^2_{L^2_{x,v}} + \| \nabla_x \partial^m \mathbb{P}^\perp G_R \|^2_{L^2_{x,v}} \big) + C \eps^2 \| \mathcal{P}_{\mathfrak{B}} \partial^m H_R \|^2_{L^2_{x,v}} \\
		\no & + C \sum_{m' \leq m} \big( \| \partial^{m-m'} u_0 \|^2_{L^\infty_x} + \| \partial^{m-m'} B_0 \|^2_{L^\infty_x} \big) ( \| \partial^{m'} \mathbb{P}^\perp G_R \|^2_{L^2_{x,v}} + \| \partial^{m'} u_R \|^2_{L^2_x} + \| \partial^{m'} B_R \|^2_{L^2_x} ) \,.
	\end{align}

For any fixed integer $N \geq 1$ and $|m| \leq N$, the Sobolev embedding inequality $H^2_x (\T) \hookrightarrow L^\infty_x (\T)$ implies
	\begin{equation}\label{MM-31}
		\| \partial^m u_0 \|^2_{L^\infty_x} + \| \partial^m B_0 \|^2_{L^\infty_x} \leq C ( \| u_0 \|^2_{H^{N+2}_x} + \| B_0 \|^2_{H^{N+2}_x} ) \leq C \mathcal{E}^{\IN}_{0,N+2} \leq C \lambda_0 (N+2),
	\end{equation}
where the last two inequalities is guaranteed by Lemma \ref{lemm:bnd-NSMF}. Therefore, by adding \eqref{MM-19} and \eqref{MM-23} to \eqref{MM-30} (up to a constant multiple of $\frac{1}{4}$), and summing up for $|m| \leq N$, the inequality \eqref{MM-31} ensures us to get that
	\begin{align}\label{MM-32}
		\no \tfrac{1}{2} \| & \nabla_x u_R \|^2_{H^N_x} + \tfrac{1}{4} \| \nabla_x \theta_R \|^2_{H^N_x} + \tfrac{1}{8} \| \nabla_x \rho_R^+ \|^2_{H^N_x} + \tfrac{1}{8} \| \nabla_x \rho_R^- \|^2_{H^N_x} + \tfrac{1}{4} \| \div_x E_R \|^2_{H^N_x} \\[3pt]
		\leq & \eps \tfrac{\d }{\d t} \widetilde{\mathcal{A}}_N (t) + C \eps^2 \big( \tfrac{1}{2} \| \nabla_x u_R \|^2_{H^N_x} + \tfrac{1}{4} \| \nabla_x \theta_R \|^2_{H^N_x} + \tfrac{1}{8} \| \nabla_x \rho_R^+ \|^2_{H^N_x} + \tfrac{1}{8} \| \nabla_x \rho_R^- \|^2_{H^N_x} \big) \\[3pt]
		\no & \tfrac{C}{\eps^2} ( 1 + \eps^2 \lambda_0 (N+2) ) \| \mathbb{P}^\perp G_R \|^2_{H^{N+1}_x L^2_v} + C \eps^2 \| \mathcal{P}_{\mathfrak{B}} H_R \|^2_{H^N_x L^2_v} + C \mathcal{E}^{\IN}_{0,N+2} ( \| u_R \|^2_{H^N_x} + \| B_R \|^2_{H^N_x} ) \,,
	\end{align}
	where the quantity $ \widetilde{\mathcal{A}}_N (t) $ is defined in \eqref{Tilde-AN}.

By taking $f = u_R$, $\theta_R$ and $\rho_R^\pm$ in the following Poincar\'e inequality:
	\begin{equation}\label{MM-33}
		\| f \|^2_{L^2_x} \leq C \| \nabla_x f \|^2_{L^2_x} + C \Big( \int_{\T} f \d x \Big)^2 \,, \ \forall \, f(x) \in H^1_x \,,
	\end{equation}
and by taking $ \eps \in ( 0 , \min \{ 1, \tfrac{1}{\sqrt{2C}} \} ] $, we are led by the relation \eqref{MM-32} to \eqref{Mic-Mac-Inq}. The proof of Lemma \ref{Lm-Mic-Mac-decomp} is therefore finished.
\end{proof}

\subsection{Estimates on \texorpdfstring{$(E_R, B_R)$}{(E\_R, B\_R)}: Proof of Lemma \ref{Lm-Mxw-Dec}}
\label{Sec:Lmm-MD}

This subsection is devoted to deriving the estimates of $(E_R, B_R)$ based on the equations \eqref{Mxw-GR}, \eqref{Mxw-4} and \eqref{Mxw-5}.

\begin{proof}[Proof of Lemma \ref{Lm-Mxw-Dec}]
	By acting the operator $\partial^m$ on \eqref{Mxw-5} for $|m| \leq N-1$, taking the $L^2_x$-inner product by multiplying by $\partial_t \partial^m B_R$ and integrating by parts over $x \in \T$, it follows that
	\begin{equation}\label{DI-1}
			\tfrac{1}{2} \tfrac{\d}{\d t} \big( \| \partial_t \partial^m B_R \|^2_{L^2_x} + \| \nabla_x \partial^m B_R \|^2_{L^2_x} \big) + \sigma \| \partial_t \partial^m B_R \|^2_{L^2_x}
			= \int_{\T} \nabla_x \times \mathcal{K} ( \mathbb{P}^\perp G_R ) \cdot \partial_t \partial^m B_R \d x \,.
	\end{equation}
	If one replaces the multiplied vector $\partial_t \partial^m B_R$ by $\partial^m B_R$ in the above procedure, there holds
	\begin{align}\label{DI-2}
		\no & \tfrac{1}{2} \tfrac{\d}{\d t} \big( \| \partial_t \partial^m B_R + \partial^m B_R \|^2_{L^2_x} - \| \partial_t \partial^m B_R \|^2_{L^2_x} + (\sigma - 1) \| \partial^m B_R \|^2_{L^2_x} \big) \\
		\no & \quad - \| \partial_t \partial^m B_R \|^2_{L^2_x} + \| \nabla_x \partial^m B_R \|^2_{L^2_x} \\
		\no & = \int_{\T} \nabla_x \times \mathcal{K} ( \mathbb{P}^\perp G_R ) \cdot \partial^m B_R \d x \\
		& = \int_{\T} \mathcal{K} ( \mathbb{P}^\perp G_R ) \cdot ( \nabla_x \times \partial^m B_R  ) \d x \,.
	\end{align}

On the other hand, for $|m| \leq N-1$, from applying $\partial^m$ on the first two equations of \eqref{Mxw-GR}, taking $L^2$-inner product with $\partial^m E_R$ and $\partial^m B_R$, respectively, and combining with \eqref{Mxw-4}, it is derived that
	\begin{align}\label{DI-3}
		\no \tfrac{1}{2} \tfrac{\d}{\d t} & \big( \| \partial^m E_R \|^2_{L^2_x} + \| \partial^m B_R \|^2_{L^2_x} \big) = - \tfrac{1}{\eps} \int_{\T} \partial^m \langle G_R \cdot \mathcal{T}_1  v \rangle \cdot \partial^m E_R \d x \\
		= & - \sigma \| \partial^m E_R \|^2_{L^2_x} + \tfrac{1}{2} \sigma \int_{\T} \nabla_x \partial^m ( \rho_R^+ - \rho_R^- ) \cdot \partial^m E_R \d x - \int_{\T} \partial^m \mathcal{K} ( \mathbb{P}^\perp G_R ) \cdot \partial^m E_R \d x \,.
	\end{align}
Note that $\div_x E_R = \rho_R^+ - \rho_R^-$ implies $\int_{\T} \nabla_x \partial^m ( \rho_R^+ - \rho_R^- ) \cdot \partial^m E_R \d x = - \| \div_x \partial^m E_R \|^2_{L^2_x}$, then, it follows from the equality \eqref{DI-3} that
	\begin{align}\label{DI-4}
		\no \tfrac{1}{2} \tfrac{\d}{\d t} ( \| \partial^m E_R \|^2_{L^2_{x}} + \| \partial^m B_R \|^2_{L^2_{x}}) + \sigma \| \partial^m E_R \|^2_{L^2_{x}} + & \tfrac{1}{2} \sigma \| \div_x \partial^m E_R \|^2_{L^2_{x}} \\
		= & \int_{\T} \partial^m \mathcal{K} ( \mathbb{P}^\perp G_R ) \cdot \partial^m E_R \d x \,.
	\end{align}
Let $\delta = \tfrac{1}{2} \min\{ 1 , \sigma \} \in (0 , \tfrac{1}{2}] $. By multiplying by $\delta$ in \eqref{DI-2} and then adding it and \eqref{DI-4} to the equality \eqref{DI-1}, we obtain that
	\begin{align}\label{DI-5}
		\no & \tfrac{1}{2} \tfrac{\d}{\d t} \Big( \| E_R \|^2_{H^{N-1}_x} + ( 1 - \delta + \sigma \delta ) \| B_R \|^2_{H^{N-1}_x} + \| \nabla_x B_R \|^2_{H^{N-1}_x} \\
		\no & \qquad \qquad \qquad \qquad + (1 - \delta) \| \partial_t B_R \|^2_{H^{N-1}_x} + \delta \| \partial_t B_R + B_R \|^2_{H^{N-1}_x} \Big) \\
		\no & \qquad + \sigma \| E_R \|^2_{H^{N-1}_x} + \delta \| \nabla_x B_R \|^2_{H^{N-1}_x} + (\sigma - \delta) \| \partial_t B_R \|^2_{H^{N-1}_x} + \tfrac{1}{2} \sigma \| \div_x E_R \|^2_{H^{N-1}_x} \\
		& = \sum_{|m| \leq N-1} \int_{\T} \nabla_x \times \partial^m \mathcal{K} ( \mathbb{P}^\perp G_R ) \cdot \partial_t \partial^m B_R \d x \\
		\no & \quad + \sum_{|m| \leq N-1} \int_{\T}  \partial^m \mathcal{K} ( \mathbb{P}^\perp G_R ) \cdot ( \delta \nabla_x \times \partial^m B_R - \partial^m E_R ) \d x \\
		\no & = I_{N-1} + I\!I_{N-1} \,.
	\end{align}

For the term $I_{N-1}$ in \eqref{DI-5}, by recalling the definition of $\mathcal{K} ( \mathbb{P}^\perp G_R )$ in \eqref{Mxw-4} and using the Sobolev embedding inequality $H^2_x (\T) \hookrightarrow L^\infty_x (\T)$, a simple calculation implies
	\begin{align}\label{DI-6}
		\no I_{N-1} \leq & - \eps \sum_{|m| \leq N-1} \int_{\T} \nabla_x \times \partial_t \langle \partial^m \mathbb{P}^\perp G_R \cdot \mathcal{T}_1 \widetilde{\Phi} (v) \rangle \cdot \partial_t \partial^m B_R \d x \\
		& + C ( 1 + \| B_0 \|_{H^{N+2}_x}  ) \| \mathbb{P}^\perp G_R \|_{H^{N+1}_x L^2_v} \| \partial_t B_R \|_{H^{N-1}_x} + C \eps \| H_R \|_{H^N_x L^2_v} \| \partial_t B_R \|_{H^{N-1}_x} \\
		\no & + C \big( \| u_0 \|_{H^{N+1}_x} + \| B_0 \|_{H^{N+1}_x} \big) \big( \| u_R \|_{H^{N-1}_x} + \| B_R \|_{H^{N-1}_x} \big)  \| \partial_t B_R \|_{H^{N-1}_x} \,.
	\end{align}
Here the first term of the right-hand side of \eqref{DI-6} still needs to be handled. First, according to the equation \eqref{Mxw-5}, there holds
	\begin{align}\label{DI-7}
		\no - \eps & \sum_{|m| \leq N-1} \int_{\T} \nabla_x \times \partial_t \langle \partial^m \mathbb{P}^\perp G_R \cdot \mathcal{T}_1 \widetilde{\Phi} (v) \rangle \cdot \partial_t \partial^m B_R \d x \\
		& = - \eps \tfrac{\d}{\d t} \sum_{|m| \leq N-1} \int_{\T} \nabla_x \times \langle \partial^m \mathbb{P}^\perp G_R \cdot \mathcal{T}_1 \widetilde{\Phi} (v) \rangle \cdot \partial_t \partial^m B_R \d x \\
		\no & \quad + \eps \sum_{|m| \leq N-1} \int_{\T} \nabla_x \times \langle \partial^m \mathbb{P}^\perp G_R \cdot \mathcal{T}_1 \widetilde{\Phi} (v) \rangle \big[ \Delta_x \partial^m B_R - \sigma \partial_t \partial^m B_R + \nabla_x \times \partial^m \mathcal{K} (\mathbb{P}^\perp G_R) \big] \d x \,.
	\end{align}
	A straightforward computation implies
	\begin{align}\label{DI-8}
		\no \eps \sum_{|m| \leq N-1} & \int_{\T} \nabla_x \times  \langle \partial^m \mathbb{P}^\perp G_R \cdot \mathcal{T}_1 \widetilde{\Phi} (v) \rangle \cdot \Delta_x \partial^m B_R \d x \\
		= & - \eps \sum_{|m| \leq N-1} \int_{\T} \nabla_x ( \nabla_x \times \langle \partial^m \mathbb{P}^\perp G_R \cdot \mathcal{T}_1 \widetilde{\Phi} (v) \rangle ) : \nabla_x \partial^m B_R \d x \\
		\no \leq & C \eps \| \mathbb{P}^\perp G_R \|_{H^{N+1}_x L^2_v} \| \nabla_x B_R \|_{H^{N-1}_x} \,,
	\end{align}
	and similarly,
	\begin{equation}\label{DI-9}
		\begin{aligned}
			\eps \sum_{|m| \leq N-1} \int_{\T} \nabla_x \times  \langle \partial^m \mathbb{P}^\perp G_R \cdot \mathcal{T}_1 \widetilde{\Phi} (v) \rangle \cdot & ( - \sigma \partial_t \partial^m B_R ) \d x
		\leq C \eps \| \mathbb{P}^\perp G_R \|_{H^N_x L^2_v} \| \partial_t B_R \|_{H^{N-1}_x} \,.
		\end{aligned}
	\end{equation}

Furthermore, together with the definition of $\mathcal{K} ( \mathbb{P}^\perp G_R )$ in \eqref{Mxw-4}, a similar argument as that in \eqref{DI-6} yields
	\begin{align}\label{DI-10}
		\no \eps & \sum_{|m| \leq N-1} \int_{\T} \nabla_x \times \langle \partial^m \mathbb{P}^\perp G_R \cdot \mathcal{T}_1 \widetilde{\Phi} (v) \rangle \cdot \nabla_x \times \partial^m \mathcal{K} ( \mathbb{P}^\perp G_R ) \d x \\
		\leq & - \tfrac{\eps^2}{2} \tfrac{\d}{\d t} \|  \nabla_x \times \langle \mathbb{P}^\perp G_R \cdot \mathcal{T}_1 \widetilde{\Phi} (v) \rangle \|^2_{H^{N-1}_x} + C \eps ( 1 + \| B_0 \|_{H^{N+2}_x} ) \| \mathbb{P}^\perp G_R \|^2_{H^{N+1}_x L^2_v} \\
		\no & + C \eps^4 \| H_R \|^2_{H^N_x L^2_v} + C \eps ( \| u_0 \|^2_{H^{N+2}_x} + \| B_0 \|^2_{H^{N+2}_x} ) ( \| B_R \|^2_{H^N_x} + \| u_R \|^2_{H^N_x} ) \,.
	\end{align}
	By plugging the inequalities \eqref{DI-8}-\eqref{DI-10} into \eqref{DI-7}, one then has
	\begin{align}\label{DI-11}
		\no - \eps & \sum_{|m| \leq N-1} \int_{\T} \nabla_x \times \partial_t \langle \partial^m \mathbb{P}^\perp G_R \cdot \mathcal{T}_1 \widetilde{\Phi} (v) \rangle \cdot \partial_t \partial^m B_R \d x \\
		\no \leq & - \tfrac{\eps^2}{2} \tfrac{\d}{\d t} \|  \nabla_x \times \langle \mathbb{P}^\perp G_R \cdot \mathcal{T}_1 \widetilde{\Phi} (v) \rangle \|^2_{H^{N-1}_x} \\
		& + C \eps ( 1 + \lambda_0^\frac{1}{2} (N+2) ) \| \mathbb{P}^\perp G_R \|^2_{H^{N+1}_x L^2_v} + C \eps^4 \| H_R \|^2_{H^N_x L^2_v} \\
		\no & + C \eps ( \| \partial_t B_R \|^2_{H^{N-1}_x} + \| \nabla_x B_R \|^2_{H^{N-1}_x} ) + C \eps \lambda_0 (N+2) ( \| B_R \|^2_{H^N_x} + \| u_R \|^2_{H^N_x} ) \,,
	\end{align}
	where the inequality \eqref{Bnd-NSFM} in Lemma \ref{lemm:bnd-NSMF} has been used. Therefore, substituting \eqref{DI-11} into \eqref{DI-6} reduces to
	\begin{equation}\label{DI-12}
		\begin{aligned}
			I_{N-1} \leq & - \tfrac{\eps^2}{2} \tfrac{\d}{\d t} \| \nabla_x \times \langle \mathbb{P}^\perp G_R \cdot \mathcal{T}_1 \widetilde{\Phi} (v) \rangle \|^2_{H^{N-1}_x} + \tfrac{C}{\eps^2} \| \mathbb{P}^\perp G_R \|^2_{H^{N+1}_x L^2_v} + C \eps^4 \| H_R \|^2_{H^N_x L^2_v} \\
			& + C ( \eps + \xi + \sqrt{\mathcal{E}^{\IN}_{0,N+2}} ) ( \| \partial_t B_R \|^2_{H^{N-1}_x} + \| \nabla_x B_R \|^2_{H^{N-1}_x} ) + C \sqrt{\mathcal{E}^{\IN}_{0,N+2}} \| u_R \|^2_{H^N_x} \,,
		\end{aligned}
	\end{equation}
	where $\xi > 0$ is sufficiently small to be determined, and the inequality \eqref{Poincare-BR} in Lemma \ref{Lm-Integral-Bnd} has been used.

It remains to estimate the term $I\!I_{N-1}$ in \eqref{DI-5}. By the equation \eqref{Mxw-5}, the H\"older inequality, the Sobolev embedding theory and the Young inequality, we can infer that
	\begin{align}\label{DI-13}
		\no I\!I_{N-1} \leq & - \eps \sum_{|m| \leq N-1} \int_{\T} \partial_t \langle \partial^m \mathbb{P}^\perp G_R \cdot \mathcal{T} \widetilde{\Phi} (v) \rangle \cdot (\delta \nabla_x \times \partial^m B_R - \partial^m E_R ) \d x \\
		& + \tfrac{C}{\eps^2} \| \mathbb{P}^\perp G_R \|^2_{H^{N+1}_x L^2_v} + C(\xi) \eps^2 \| H_R \|^2_{H^N_x L^2_v} + C \sqrt{\mathcal{E}^{\IN}_{0,N+2}} \| u_R \|^2_{H^{N-1}_x} \\
		\no & + ( \eps^2 + \xi + C \sqrt{\mathcal{E}^{\IN}_{0,N+2}} ) ( \| \nabla_x B_R \|^2_{H^{N-1}_x} + \| E_R \|^2_{H^{N-1}_x} ) \,,
	\end{align}
	where the constant $\xi > 0$ is sufficiently small to be determined. For the first term of the right-hand side of \eqref{DI-13}, one has
	\begin{align}\label{DI-14}
		\no - \eps \sum_{|m| \leq N-1} & \int_{\T} \partial_t \langle \partial^m \mathbb{P}^\perp G_R \cdot \mathcal{T} \widetilde{\Phi} (v) \rangle \cdot (\delta \nabla_x \times \partial^m B_R - \partial^m E_R ) \d x \\
		= & \eps \tfrac{\d}{\d t} \sum_{|m| \leq N-1} \int_{\T}  \langle \partial^m \mathbb{P}^\perp G_R \cdot \mathcal{T} \widetilde{\Phi} (v) \rangle \cdot (\delta \nabla_x \times \partial^m B_R - \partial^m E_R ) \d x \\
		\no & + \eps \sum_{|m| \leq N-1} \int_{\T} \langle \partial^m \mathbb{P}^\perp G_R \cdot \mathcal{T} \widetilde{\Phi} (v) \rangle \cdot ( \delta \nabla_x \times \partial_t \partial^m B_R - \partial_t \partial^m E_R ) \d x \,,
	\end{align}
	where the last term in \eqref{DI-14} can be bounded by
	\begin{align}\label{DI-15}
		\no \eps \sum_{|m| \leq N-1} & \int_{\T} \langle \partial^m \mathbb{P}^\perp G_R \cdot \mathcal{T} \widetilde{\Phi} (v) \rangle \cdot ( \delta \nabla_x \times \partial_t \partial^m B_R - \partial_t \partial^m E_R ) \d x \\
		\no = & - \eps \delta \sum_{|m| \leq N-1} \int_{\T} \nabla_x \times [ \nabla_x \times \langle \partial^m \mathbb{P}^\perp G_R \cdot \mathcal{T} \widetilde{\Phi} (v) \rangle ] \cdot \partial^m E_R \d x \\
		\no & - \eps \sum_{|m| \leq N-1} \int_{\T} \langle \partial^m \mathbb{P}^\perp G_R \cdot \mathcal{T} \widetilde{\Phi} (v) \rangle \cdot \nabla_x \times \partial^m B_R \d x \\
		& + \sum_{|m| \leq N-1} \int_{\T} \langle \partial^m \mathbb{P}^\perp G_R \cdot \mathcal{T} \widetilde{\Phi} (v) \rangle \cdot \langle \partial^m \mathbb{P}^\perp G_R \cdot \mathcal{T} v \rangle \d x \\
		\no \leq & C \eps \| \mathbb{P}^\perp G_R \|_{H^{N+1}_x L^2_v} ( \| \nabla_x B_R \|_{H^{N-1}_x} + \| E_R \|_{H^{N-1}_x} ) + C \| \mathbb{P}^\perp G_R \|^2_{H^{N+1}_x L^2_v} \,,
	\end{align}
where the relation $\int_{\T} A \cdot \nabla_x \times ( \nabla_x \times B  ) \d x = \int_{\T} \nabla_x \times ( \nabla_x \times A ) \cdot B \d x$ and the equations \eqref{Mxw-GR} have been used. Collecting the inequalities \eqref{DI-13}-\eqref{DI-15} thereby implies
		\begin{align}\label{DI-16}
			I\!I_{N-1} \leq & - \eps \tfrac{\d}{\d t} \sum_{|m| \leq N-1} \int_{\T}  \langle \partial^m \mathbb{P}^\perp G_R \cdot \mathcal{T} \widetilde{\Phi} (v) \rangle \cdot (\delta \nabla_x \times \partial^m B_R - \partial^m E_R ) \d x \no \\
			& + \tfrac{C}{\eps^2} \| \mathbb{P}^\perp G_R \|^2_{H^{N+1}_x L^2_v} + C(\xi) \eps^2 \| H_R \|^2_{H^N_x L^2_v} + C \sqrt{\mathcal{E}^{\IN}_{0,N+2}} \| u_R \|^2_{H^{N-1}_x} \\ \no
			& + ( \eps^2 + \xi + C \sqrt{\mathcal{E}^{\IN}_{0,N+2}} ) ( \| \nabla_x B_R \|^2_{H^{N-1}_x} + \| E_R \|^2_{H^{N-1}_x} ),
		\end{align}
with an undetermined $\xi > 0$. Together with \eqref{DI-5}, \eqref{DI-12} and \eqref{DI-16}, one knows
	\begin{align}\label{DI-17}
		\no \tfrac{1}{2} \tfrac{\d}{\d t} & \Big( \| E_R \|^2_{H^{N-1}_x} + ( 1 - \delta + \sigma \delta  )  \| B_R \|^2_{H^{N-1}_x}  + \| \nabla_x B_R \|^2_{H^{N-1}_x} + (1 - \delta ) \| \partial_t B_R \|^2_{H^{N-1}_x} \\
		\no & \quad + \delta \| \partial_t B_R + B_R \|^2_{H^{N-1}_x} + \eps^2 \| \nabla_x \times \langle \mathbb{P}^\perp G_R \cdot \mathcal{T}_1 \widetilde{\Phi} (v) \rangle \|^2_{H^{N-1}_x} \Big) \\
		\no & + \sigma \| E_R \|^2_{H^{N-1}_x} + \delta \| \nabla_x B_R \|^2_{H^{N-1}_x} + ( \sigma - \delta ) \| \partial_t B_R \|^2_{H^{N-1}_x} + \tfrac{1}{2} \sigma \| \div_x  E_R \|^2_{H^{N-1}_x} \\
		\leq & \eps \tfrac{\d}{\d t} \sum_{|m| \leq N-1} \int_{\T}  \langle \partial^m \mathbb{P}^\perp G_R \cdot \mathcal{T} \widetilde{\Phi} (v) \rangle \cdot (\partial^m E_R - \delta \nabla_x \times \partial^m B_R ) \d x \\
		\no & + \tfrac{C}{\eps^2} \| \mathbb{P}^\perp G_R \|^2_{H^{N+1}_x L^2_v} + C(\xi) \eps^2 \| H_R \|^2_{H^N_x L^2_v} + C \sqrt{\mathcal{E}^{\IN}_{0,N+2}} \| u_R \|^2_{H^{N}_x} \\
		\no & + C ( \eps + \xi + C \sqrt{\mathcal{E}^{\IN}_{0,N+2}} ) ( \| \partial_t B_R \|^2_{H^{N-1}_x} + \| \nabla_x B_R \|^2_{H^{N-1}_x} + \| E_R \|^2_{H^{N-1}_x} )  \,.
	\end{align}
As a consequence, by picking $\mathcal{E}^{\IN}_{0,N+2}$, $\eps\,, \xi \ll 1$, the inequality \eqref{DI-17} concludes \eqref{Decay-Inq}. The proof of Lemma \ref{Lm-Mxw-Dec} is therefore completed.
\end{proof}

\subsection{Kinetic Dissipation of \texorpdfstring{$G_R$}{G\_R}: Proof of Lemma \ref{Lm-Unif-Spatial-Bnd}}
\label{Sec:Lmm-USB}

In this subsection, we aims at proving Lemma \ref{Lm-Unif-Spatial-Bnd}, namely, at deriving the kinetic dissipation contribution $\tfrac{1}{\eps^2} \sum_{|m| \leq N+1} \| \partial^m \mathbb{P}^\perp G_R \|^2_{L^2_{x,v} (\nu)}$ with singularity of order $\tfrac{1}{\eps^2}$, and finishing the spatial derivative estimates \eqref{Unif-Spatial-Bnd}. Remark firstly that the kinetic dissipation is one of the important structures which can absorb the singularity. In order to deal with the singularity, there are two key cancellation relations, as follows, (see also \eqref{Sp-9} and \eqref{Sp-11} below)
\begin{equation*}
	\begin{aligned}
		\langle [ \mathcal{T} (v \times \partial^{m - m'} B_R ) \cdot \nabla_v \partial^{m'} G_0 ] \cdot \partial^m \mathbb{P} G_R \rangle = \langle [ \mathcal{T} (v \times \partial^{m'} B_0) \cdot \nabla_v \partial^{m-m'} \mathbb{P} G_R ] \cdot \partial^m \mathbb{P} G_R \rangle = 0 \,.
	\end{aligned}
\end{equation*}

\begin{proof}[Proof of Lemma \ref{Lm-Unif-Spatial-Bnd}]
Take $\partial^m $ with $|m| \leq N+1$ on the first equation of the remainder system \eqref{Remd-Equ-GR}, and then we get
	\begin{align}\label{Sp-5}
		\no \partial_t \partial^m & G_R + \tfrac{1}{\eps} v \cdot \nabla_x \partial^m G_R - \tfrac{1}{\eps} ( \partial^m E_R \cdot v ) \mathcal{T}_1 + \tfrac{1}{\eps} \mathcal{T} ( v \times B_0 ) \cdot \nabla_v \partial^m G_R + \tfrac{1}{\eps^2} \mathbb{L} \partial^m G_R \\
		= & - \tfrac{1}{\eps} \sum_{0 \neq m' \leq m} C_{m}^{m'} \mathcal{T} (v \times \partial^{m'} B_0) \cdot \nabla_v \partial^{m-m'} G_R - \tfrac{1}{\eps} \partial^m [ \mathcal{T} (v \times B_R) \cdot \nabla_v G_0 ] + \partial^m H_R \,.
	\end{align}
	Taking $L^2_{x,v}$-inner product with $\partial^m G_R$ in \eqref{Sp-5} reduces to
	\begin{align}\label{Sp-6}
		\no \tfrac{1}{2} \tfrac{\d}{\d t} & \| \partial^m G_R \|^2_{L^2_{x,v}} - \int_{\T} \partial^m E_R \cdot \tfrac{1}{\eps} \langle \partial^m G_R \cdot \mathcal{T}_1 v \rangle \d x + \tfrac{1}{\eps^2} \int_{\T} \langle \mathbb{L} \partial^m G_R \cdot \partial^m G_R \rangle \d x \\
		= & - \tfrac{1}{\eps} \sum_{0 \neq m' \leq m} C_m^{m'} \int_{\T} \langle [ \mathcal{T} ( v \times \partial^{m'} B_0 ) \cdot \nabla_v \partial^{m-m'} G_R ] \cdot \partial^m G_R \rangle \d x \\
		\no & - \tfrac{1}{\eps} \int_{\T} \langle \partial^m [ \mathcal{T} (v \times B_R) \cdot \nabla_v G_0 ] \cdot \partial^m G_R \rangle \d x + \int_{\T} \langle \partial^m H_R \cdot \partial^m G_R \rangle \d x \,.
	\end{align}
	Recalling the Maxwell equations of \eqref{Remd-Equ-GR}, one has
	\begin{align}\label{Sp-7}
		\no - \int_{\T} \partial^m E_R \cdot \tfrac{1}{\eps} \langle \partial^m G_R \cdot \mathcal{T}_1 v \rangle \d x = & \int_{\T} \partial^m E_R \cdot ( \partial_t \partial^m E_R - \nabla_x \times \partial^m B_R ) \d x \\
		= & \tfrac{1}{2} \tfrac{\d}{\d t} \| \partial^m E_R \|^2_{L^2_x} - \int_{\T} ( \nabla_x \times \partial^m E_R) \cdot \partial^m B_R \d x \\
		\no = & \tfrac{1}{2} \tfrac{\d}{\d t} \| \partial^m E_R \|^2_{L^2_x} + \int_{\T} \partial_t \partial^m B_R \cdot \partial^m B_R \d x \\
		\no = & \tfrac{1}{2} \tfrac{\d}{\d t} \big( \| \partial^m E_R \|^2_{L^2_x} + \| \partial^m B_R \|^2_{L^2_x} \big) \,,
	\end{align}
where the second equality is implied by $\int_{\T} A \cdot ( \nabla_x \times B ) \d  x = \int_{\T} (\nabla_x \times A) \cdot B \d x$.

By the coercive property of the projection operator $\mathbb{L}$ in Lemma \ref{Lm-L-Propty}, there is a constant $\Lambda > 0$, such that
	\begin{equation}\label{Sp-8}
		\tfrac{1}{\eps^2} \int_{\T} \langle \mathbb{L} \partial^m G_R \cdot \partial^m G_R \rangle \d  x \geq \tfrac{\Lambda}{\eps^2} \| \partial^m \mathbb{P}^\perp G_R \|^2_{L^2_{x,v} (\nu)} \,.
	\end{equation}
	Since $G_0 = (g_0^+, g_0^-)^\top$ with $g_0^\pm = \rho_0^\pm + u_0 \cdot v + \theta_0 ( \tfrac{|v|^2}{2} - \tfrac{3}{2} )$ and $\mathbb{P} G_R$ is defined in \eqref{Projc-P}, a simple calculation implies
		\begin{multline}\label{Sp-9}
			\langle [ \mathcal{T} (v \times \partial^{m - m'} B_R ) \cdot \nabla_v \partial^{m'} G_0 ] \cdot \partial^m \mathbb{P} G_R \rangle \\
			= \left\langle \left[ (v \times \partial^{m-m'} B_R) \cdot \left( \begin{array}{c} u_0 \\ - u_0 \end{array} \right) \right] \cdot \left( \begin{array}{c} \partial^m \rho_R^+ + \partial^m u_R \cdot v + \partial^m \theta_R ( \tfrac{|v|^2}{2} - \tfrac{3}{2}) \\ \partial^m \rho_R^- + \partial^m u_R \cdot v + \partial^m \theta_R ( \tfrac{|v|^2}{2} - \tfrac{3}{2}) \end{array} \right) \right\rangle = 0 \,.
		\end{multline}
In fact, this is one of the key cancellations to deal with the singularity.

Splitting $G_R = \mathbb{P} G_R + \mathbb{P}^\perp G_R $, using the cancellation \eqref{Sp-9} and the Sobolev embedding theory lead us to
	\begin{align}\label{Sp-10}
		\no - \tfrac{1}{\eps} \int_{\T} & \langle \partial^m [\mathcal{T} (v \times B_R) \cdot \nabla_v G_0 ] \cdot \partial^m G_R \rangle \d x \\
		\no = & - \tfrac{1}{ \eps } \sum_{m' \leq m} C_m^{m'} \int_{\T} \langle [ \mathcal{T} (v \times \partial^{m-m'} B_R) \cdot \nabla_v \partial^{m'} G_0 ] \cdot  \mathbb{P}^\perp \partial^m G_R  \rangle \d x \\
		\leq & \tfrac{C}{\eps} \sum_{m' \leq m}  \int_{\T} |\partial^{m-m'} B_R| \, |\partial^{m'} u_0 | \Big( \int_{\R} |v|^2 M \d v \Big)^\frac{1}{2} \langle | \partial^m \mathbb{P}^\perp G_R | ^2 \rangle^\frac{1}{2} \d x \\
		\no \leq & \tfrac{C}{\eps} \| u_0 \|_{H^{N+2}_x} \| B_R \|_{H^{N+1}_x} \| \partial^m \mathbb{P}^\perp G_R \|_{L^2_{x,v}} \,.
	\end{align}

Another key cancellation relation concerning the singularity is similar as in \eqref{Sp-9}, i.e.,
	\begin{equation}\label{Sp-11}
		\langle [ \mathcal{T} (v \times \partial^{m'} B_0) \cdot \nabla_v \partial^{m-m'} \mathbb{P} G_R ] \cdot \partial^m \mathbb{P} G_R \rangle = 0 \,.
	\end{equation}
Then, combining the relation $G_R = \mathbb{P} G_R + \mathbb{P}^\perp G_R $ and \eqref{Sp-11} together yields that
	\begin{align}\label{Sp-12}
		\no - \tfrac{1}{\eps} \sum_{0 \neq m' \leq m} & C_m^{m'} \int_{\T} \langle [ \mathcal{T} ( v \times \partial^{m'} B_0 ) \cdot \nabla_v \partial^{m-m'} G_R ] \cdot \partial^m G_R \rangle \d x \\
		\no = & - \tfrac{1}{\eps} \sum_{0 \neq m' \leq m} C_m^{m'} \int_{\T} \langle [ \mathcal{T} ( v \times \partial^{m'} B_0 ) \cdot \nabla_v \partial^{m-m'} \mathbb{P} G_R ] \cdot \partial^m \mathbb{P}^\perp G_R  \rangle \d x \\
		& - \tfrac{1}{\eps} \sum_{0 \neq m' \leq m} C_m^{m'} \int_{\T} \langle [ \mathcal{T} ( v \times \partial^{m'} B_0 ) \cdot \nabla_v \partial^{m-m'} \mathbb{P}^\perp G_R ] \cdot  \partial^m \mathbb{P} G_R \rangle \d x \\
		\no & - \tfrac{1}{\eps} \sum_{0 \neq m' \leq m} C_m^{m'} \int_{\T} \langle [ \mathcal{T} ( v \times \partial^{m'} B_0 ) \cdot \nabla_v \partial^{m-m'} \mathbb{P}^\perp G_R ] \cdot  \partial^m \mathbb{P}^\perp G_R \rangle \d x \\
		\no := & I_1 + I_2 + I_3 \,.
	\end{align}

By the definition of $\mathbb{P} G_R$ in \eqref{Projc-P}, we deduce that for $|m| \leq N+1$,
	\begin{align}\label{Sp-13}
		\no I_1 = & - \tfrac{1}{\eps} \sum_{0 \neq m' \leq m} C_m^{m'} \int_{\T} \left\langle \left[  ( v \times \partial^{m'} B_0 ) \cdot  \partial^{m-m'} \left( \begin{array}{c} u_R \\ - u_R \end{array} \right) \right] \cdot \partial^m \mathbb{P}^\perp G_R  \right\rangle \d x \\
		\leq & \tfrac{C}{\eps} \| B_0 \|_{H^{N+2}_x} \| u_R \|_{H^{N+1}_x} \| \partial^m \mathbb{P}^\perp G_R \|_{L^2_{x,v}} \leq \tfrac{C}{\eps} \sqrt{\mathcal{E}^{\IN}_{0,N+2}} \| \partial^m \mathbb{P}^\perp G_R \|_{L^2_{x,v}} \| \mathbb{P} G_R \|_{H^{N+1}_x L^2_v} \,,
	\end{align}
where Lemma \ref{lemm:bnd-NSMF} and the inequality \eqref{Sp-1} have been used.
Similarly, we can derive from the integration by parts over $v \in \R^3$ and the symmetry of the integral on the $v$-variable, that
	\begin{equation}\label{Sp-14}
		\begin{aligned}
			I_2 = & - \tfrac{1}{\eps} \sum_{0 \neq m' \leq m} C_m^{m'} \int_{\T} \langle [ \mathcal{T} ( v \times \partial^{m'} B_0 ) \cdot \nabla_v \partial^m \mathbb{P} G_R ] \cdot \partial^{m-m'} \mathbb{P}^\perp G_R  \rangle \d x \\
			\leq & \tfrac{C}{\eps} \sqrt{\mathcal{E}^{\IN}_{0,N+2}} \| \partial^m \mathbb{P} G_R \|_{L^2_{x,v}} \| \mathbb{P}^\perp G_R \|_{H^{N+1}_x L^2_v} \,.
		\end{aligned}
	\end{equation}
For the term $I_3$, the H\"older inequality and the Sobolev embedding $H^2_x (\T) \hookrightarrow L^\infty_x(\T)$ imply
	\begin{equation}\label{Sp-15}
		\begin{aligned}
			I_3 \leq & \tfrac{C}{\eps} \| B_0 \|_{H^{N+2}_x} \| \partial^m \mathbb{P}^\perp G_R \|_{L^2_{x,v}} \| \nabla_v \mathbb{P}^\perp G_R \|_{H^N_x L^2_v} \\
			\leq & \tfrac{C}{\eps} \sqrt{\mathcal{E}^{\IN}_{0,N+2}} \| \partial^m \mathbb{P}^\perp G_R \|_{L^2_{x,v}} \| \nabla_v \mathbb{P}^\perp G_R \|_{H^N_x L^2_v} \,.
		\end{aligned}
	\end{equation}
Combining the relations \eqref{Sp-12}-\eqref{Sp-15} enables us to get
		\begin{align}\label{Sp-16}
			- \frac{1}{\eps} \sum_{0 \neq m' \leq m} & C_m^{m'} \int_{\T} \langle [ \mathcal{T} ( v \times \partial^{m'} B_0 ) \cdot \nabla_v \partial^{m-m'} \mathbb{P}^\perp G_R ] \cdot  \partial^m \mathbb{P}^\perp G_R \rangle \d x \\ \no
			\leq & \frac{C}{\eps} \sqrt{\mathcal{E}^{\IN}_{0,N+2}} \| \partial^m \mathbb{P}^\perp G_R \|_{L^2_{x,v}} \Big( \| \mathbb{P} G_R \|_{H^{N+1}_x L^2_v} + \| \nabla_v \mathbb{P}^\perp G_R \|_{H^N_x L^2_v} \Big) \,.
		\end{align}
From substituting the inequalities \eqref{Sp-7}, \eqref{Sp-8} and \eqref{Sp-10} into \eqref{Sp-6}, it follows,
	\begin{equation}\label{Sp-17}
		\begin{aligned}
			\tfrac{1}{2} \tfrac{\d}{\d t} & \Big( \| \partial^m G_R \|^2_{L^2_{x,v}} + \| \partial^m E_R \|^2_{L^2_x} + \| \partial^m B_R \|^2_{L^2_x} \Big) + \tfrac{\Lambda}{\eps^2} \| \partial^m \mathbb{P}^\perp G_R \|^2_{L^2_{x,v}} \\
			\leq & \tfrac{C}{\eps} \| u_0 \|_{H^{N+2}_x} \| B_R \|_{H^{N+1}_x} \| \partial^m \mathbb{P}^\perp G_R \|_{L^2_{x,v}} + \int_{\T} \langle \partial^m H_R \cdot \partial^m G_R \rangle \d x \\
			& + \tfrac{C}{\eps} \sqrt{\mathcal{E}^{\IN}_{0,N+2}} \| \partial^m \mathbb{P}^\perp G_R \|_{L^2_{x,v}} \Big( \| \mathbb{P} G_R \|_{H^{N+1}_x L^2_v} + \| \nabla_v \mathbb{P}^\perp G_R \|_{H^N_x L^2_v} \Big) \,.
		\end{aligned}
	\end{equation}

Next one deals with the term $\frac{C}{\eps} \| u_0 \|_{H^{N+2}_x} \| B_R \|_{H^{N+1}_x} \| \partial^m \mathbb{P}^\perp G_R \|_{L^2_{x,v}}$ in the right-hand side of \eqref{Sp-17}. From the conservation laws \eqref{Conservtn-g0} of the leading order term $(g_0^\pm , E_0 , B_0)$, one has
$$\int_{\T} \theta_0 \d x = 0, \text{ and } \int_{\T} ( 2 u_0 + E_0 \times B_0  ) \d  x = 0.$$
The Poincar\'e inequality implies
	\begin{align}\label{Sp-18}
		\no \| u_0 \|_{L^2_x} \leq & C \| \nabla_x u_0 \|_{L^2_x} + C \Big| \int_{\T} u_0 \d x \Big| \leq C \| \nabla_x u_0 \|_{L^2_x} + C \Big| \int_{\T} E_0 \times B_0 \d  x \Big| \\
		\leq & C \| \nabla_x u_0 \|_{L^2_x} + C \| E_0 \|_{L^2_x} \| B_0 \|_{L^2_x} \,,
	\end{align}
which reduces to
	\begin{equation}\label{Sp-19}
		\begin{aligned}
			\| u_0 \|_{H^{N+2}_x} \lesssim \| \nabla_x u_0 \|_{H^{N+1}_x} + \| \nabla_x u_0 \|_{L^2_x} +  \| E_0 \|_{L^2_x} \| B_0 \|_{L^2_x} \lesssim \Big( 1 + \sqrt{\mathcal{E}^{\IN}_{0,N+1}} \Big) \mathcal{D}^\frac{1}{2}_{0,N+1} (t) \,,
		\end{aligned}
	\end{equation}
where the last inequality is implied by Lemma \ref{lemm:bnd-NSMF}.

By summing up for $|m| \leq N+1$ in \eqref{Sp-17}, and combined the inequality \eqref{Sp-19} and the Young inequality, one has
	\begin{align}\label{Sp-20}
		\no \tfrac{\d}{\d t} \big( & \| G_R \|^2_{H^{N+1}_x L^2_v} + \| E_R \|^2_{H^{N+1}_x} + \| B_R \|^2_{H^{N+1}_x} \big) + \tfrac{\Lambda}{\eps^2} \| \mathbb{P}^\perp G_R \|^2_{H^{N+1}_x L^2_v (\nu)} \\
		\leq & C \Big( 1 + {\mathcal{E}^{\IN}_{0,N+1}} \Big) \mathcal{D}_{0,N+1} (t) \| B_R \|^2_{H^{N+1}_x} + \sum_{|m| \leq N+1} \int_{\T} \langle \partial^m H_R \cdot \partial^m G_R \rangle \d x \\
		\no & + C {\mathcal{E}^{\IN}_{0,N+2}}  \big( \| \mathbb{P} G_R \|^2_{H^{N+1}_x L^2_v} + \| \nabla_v \mathbb{P}^\perp G_R \|_{H^N_x L^2_v} \big) \,.
	\end{align}
Observe that the inequality \eqref{Poincare-BR} in Lemma \ref{Lm-Integral-Bnd} implies for $N \geq 1$,
	\begin{equation*}
		\| B_R \|^2_{H^N_x} \leq  C_2 \| \nabla_x B_R \|^2_{L^2_x} + \| \nabla_x B_R \|^2_{H^{N-1}_x} \leq ( 1 + C_2) \| \nabla_x B_R \|^2_{H^{N-1}_x}.
	\end{equation*}
Recall the inequality \eqref{Sp-1} and the initial condition \eqref{Sp-3}, i.e., $\mathcal{E}^{\IN}_{0,N+2} \leq \lambda_R (N+2)$, then we have
	\begin{align}\label{Sp-21}
			C_3 \mathcal{E}^{\IN}_{0,N+2} \| B_R \|^2_{H^N_x} \leq C_3 ( 1 + C_2 ) \lambda_R(N+2) \| \nabla_x B_R \|^2_{H^{N-1}_x} \,, \\ \no
			C_3 \sqrt{\mathcal{E}^{\IN}_{0,N+2}} \| u_R \|^2_{H^N_x} \leq \tfrac{C_3}{C_0} \sqrt{\lambda_R (N+2)} \| \mathbb{P} G_R \|^2_{H^{N+1}_x L^2_v} \,.
	\end{align}
If one chooses $\lambda_R (N+2) \in ( 0 , \lambda_1 (N+2) ] $, depending only on $N$, $\mu$, $\sigma$, $\kappa$, such that
$$C_3 ( 1 + C_2 ) \lambda_R(N+2) \leq \frac{\delta}{2}, \text{ and } \tfrac{C_3}{C_0} \sqrt{\lambda_R (N+2)} \leq \tfrac{C_0 C_3}{2 C_1},$$
then the inequality \eqref{Sp-21} reduces to
	\begin{equation}\label{Sp-22}
			C_3 \mathcal{E}^{\IN}_{0,N+2} \| B_R \|^2_{H^N_x} \leq \frac{\delta}{2} \| \nabla_x B_R \|^2_{H^{N-1}_x} \,, \ C_3 \sqrt{\mathcal{E}^{\IN}_{0,N+2}} \| u_R \|^2_{H^N_x} \leq \tfrac{C_0 C_3}{2 C_1} \| \mathbb{P} G_R \|^2_{H^{N+1}_x L^2_v} \,.
	\end{equation}

Therefore, from combining with inequalities \eqref{Sp-1}-\eqref{Sp-2} and \eqref{Sp-22}, plugging the inequality \eqref{Integral-Inqu} into \eqref{Mic-Mac-Inq}, multiplying by $\frac{C_3}{C_1}$ in \eqref{Mic-Mac-Inq} and adding them to the inequality \eqref{Decay-Inq}, we can derive that
	\begin{align}\label{Sp-23}
		\no \tfrac{\d}{\d t} \Big[ & \| E_R \|^2_{N^{N-1}_x} + (1-\delta + \sigma \delta) \| B_R \|^2_{N^{N-1}_x} + \| \nabla_x B_R \|^2_{N^{N-1}_x} + (1-\delta) \| \partial_t B_R \|^2_{N^{N-1}_x} \\
		\no & \quad + \delta \| \partial_t B_R + B_R \|^2_{N^{N-1}_x} + \eps^2 \| \nabla_x \times \langle \mathbb{P}^\perp G_R \cdot \mathcal{T}_1 \widetilde{\Phi}(v) \rangle \|^2_{H^{N-1}_x} \Big] + \sigma \| E_R \|^2_{N^{N-1}_x} \\
		\no & + \tfrac{\delta}{2} \| \nabla_x B_R \|^2_{N^{N-1}_x} + (\sigma - \delta) \| \partial_t B_R \|^2_{N^{N-1}_x} + \tfrac{C_0 C_3}{2 C_1} \| \mathbb{P} G_R \|^2_{H^{N+1}_x L^2_v} + \tfrac{C_3}{C_1} \| \div_x E_R \|^2_{H^N_x} \\
		\no \leq & C_3 \eps \tfrac{\d}{\d t} \mathcal{A}_N (t) + \tfrac{2 C_3}{\eps^2} \| \mathbb{P}^\perp G_R \|^2_{H^{N+1}_x L^2_v} + C_3 \eps^2 \| H_R \|^2_{H^N_x L^2_v} + C_3 \eps^2 \| \mathcal{P}_{\mathfrak{B}} H_R \|^2_{H^N_x L^2_v}  \\
		& + C_2 C_3 \big( \mathcal{D}_{0,2} (t) + \mathcal{D}_{1,2} (t) \big) + C_2 C_3 \eps^2 ( \| E_R \|^4_{L^2_x} + \| B_R \|^2_{L^2_x} \| \nabla_x B_R \|^2_{L^2_x} ) \,,
	\end{align}
	where the scalar functional $\mathcal{A}_N (t)$ is defined as in \eqref{Sp-4}.

Multiplying by $\frac{\Lambda + 2 C_3}{\Lambda}$ in \eqref{Sp-20} and by $C_2 C_3$ in \eqref{bnd:Spec-energy-1} in Lemma \ref{lemm:bnd-linearMaxwl} for the case $M = N$, respectively, and adding them to \eqref{Sp-23}, one has
	\begin{align}\label{Sp-24}
		\no \tfrac{\d}{\d t} \Big[ & \| E_R \|^2_{H^{N-1}_x} + (1-\delta + \sigma \delta) \| B_R \|^2_{H^{N-1}_x} + \| \nabla_x B_R \|^2_{H^{N-1}_x} + (1-\delta) \| \partial_t B_R \|^2_{H^{N-1}_x} \\
		\no & \quad + \delta \| \partial_t B_R +  B_R \|^2_{H^{N-1}_x} + \eps^2 \| \nabla_x \times \langle \mathbb{P}^\perp G_R \cdot \mathcal{T}_1 \widetilde{\Phi}(v) \rangle \|^2_{H^{N-1}_x} - C_3 \eps \tfrac{\d}{\d t} \mathcal{A}_N (t) \\
		\no & + \tfrac{\Lambda + 2 C_3}{\Lambda} \big( \| G_R \|^2_{H^{N+1}_x L^2_v} + \| E_R \|^2_{H^{N+1}_x} + \| B_R \|^2_{H^{N+1}_x} \big) + C_2 C_3 \big( \mathcal{E}_{1,N} (t) + \widetilde{C}_N \mathcal{E}_{0,N+2} (t) \big) \Big] \\
		\no & + \sigma \| E_R \|^2_{H^{N-1}_x} + \tfrac{\delta}{2} \| \nabla_x B_R \|^2_{H^{N-1}_x} + (\sigma - \delta) \| \partial_t B_R \|^2_{H^{N-1}_x} + \mathcal{D}_{0,N+2} (t) \\
		\no & + \tfrac{C_0 C_3}{4 C_1} \| \mathbb{P} G_R \|^2_{H^{N+1}_x L^2_v} + \tfrac{C_3}{C_1} \| \div_x E_R \|^2_{H^{N}_x} + \tfrac{\Lambda}{\eps^2} \| \mathbb{P}^\perp G_R \|^2_{H^{N+1}_x L^2_v (\nu)} \\
		\no \leq & C_3 \eps^2 ( \| H_R \|^2_{H^N_x L^2_v} + \| \mathcal{P}_{\mathfrak{B}} H_R \|^2_{H^N_x L^2_v} ) + \tfrac{\Lambda + 2 C_3}{\Lambda} \sum_{|m|\leq N+1} \int_{\T} \langle \partial^m H_R \cdot \partial^m G_R  \rangle \d x \\
		\no & + \tfrac{C(\Lambda + 2 C_3)}{\Lambda} \Big[ ( 1 + \lambda_0(N+2) ) \| B_R \|^2_{H^{N+1}_x} \mathcal{D}_{0,N+2}(t) + \mathcal{E}^{\IN}_{0,N+2} \| \nabla_v \mathbb{P}^\perp G_R \|^2_{H^N_x L^2_v} \Big] \\
		& + C_2 C_3 \eps^2 \big( \| E_R \|^4_{L^2_x} + \| B_R \|^2_{L^2_x} \| \nabla_x B_R \|^2_{L^2_x} \big) \,.
	\end{align}
We take now $C_4 = \tfrac{\Lambda + 2 C_3}{\Lambda} > 1$, $C_5= C_2 C_3 > 0$, $C_6 = \min \{ 1 , \tfrac{\delta}{2} , \sigma - \delta, \tfrac{C_0 C_3}{4 C_1} , \tfrac{C_3}{C_1} , \Lambda \} > 0$ and $C_7 = \max \{ C_3 , \tfrac{C (\Lambda + 2 C_3)}{\Lambda} ( 1 + \lambda_0 (N+2) ) , \tfrac{\Lambda + 2 C_3}{\Lambda}, C_2 C_3 \} > 0$, then the inequality \eqref{Sp-24} concludes \eqref{Unif-Spatial-Bnd}, and hence the proof of Lemma \ref{Lm-Unif-Spatial-Bnd} is completed.
\end{proof}

\subsection{Mixed Derivative Estimates: Proof of Lemmas \ref{lemm:remainder-apriori}-\ref{lemm:remainder-apriori2}}
\label{Sec:Lmm-RA12}

In this subsection, in order to control the quantities related to $H_R$ and the term $\mathcal{E}_{0, N+2}^{\IN} \sum_{|m| \leq N} \| \nabla_v \partial^m \mathbb{P}^\perp G_R \|^2_{L^2_{x,v}}$, one will derive the mixed derivative estimates of $G_R$ in weighted spaces. First, we focus on $\partial^m_\beta G_R$ with $|m| + |\beta| \leq N$, corresponding to Lemma \ref{lemm:remainder-apriori}. Second, we consider $\partial^m_\beta G_R$ with $|m| + |\beta| \leq N+1$ and $\beta \neq 0$, corresponding to Lemma \ref{lemm:remainder-apriori2}.

\begin{proof}[Proof of Lemma \ref{lemm:remainder-apriori}]
From applying the mixed multi-derivative operator $\p^m_{\beta}$ with $|m| + |\beta| \leq N$ to the $G_R$-equation in \eqref{Remd-Equ-GR} and combining with the decomposition $G_R = \mathbb{P}G_R  + \mathbb{P}^\perp G_R$, we can get
\begin{align}
	[\p_t & + \tfrac{1}{\eps} v \cdot \nabla_x + \tfrac{1}{\eps} \mathcal{T} (v \times B_0) \cdot \nabla_v] \p^m_{\beta} \mathbb{P}^\perp G_R + \tfrac{1}{\eps^2} \p^m_{\beta} \mathbb{L} \mathbb{P}^\perp G_R
	\\\no
	= & - \big[ \p_t + \tfrac{1}{\eps} v \cdot \nabla_x + \tfrac{1}{\eps} \mathcal{T} (v \times B_0) \cdot \nabla_v \big] \p^m_{\beta} \mathbb{P}G_R
	\\\no &
	- \tfrac{1}{\eps} \sum_{|\beta_1|=1} C_{\beta}^{\beta_1} \p_{\beta_1} v \cdot \nabla_x \p^m_{\beta- \beta_1}G_R
	- \tfrac{1}{\eps} \p^m_{\beta} [\mathcal{T} (v \times B_0) \cdot \nabla_v G_0] + \tfrac{1}{\eps} \p^m E_R \cdot \p_{\beta} v \mathcal{T}_1
	\\\no &
	- \tfrac{1}{\eps} \sum_{ |\beta_1|=1, m_1\le m} C_m^{m_1} C_{\beta}^{\beta_1} \mathcal{T} (\p_{\beta_1} v \times \p^{m_1} B_0) \cdot \nabla_v \p^{m-m_1}_{\beta- \beta_1} G_R + \p^m_{\beta} H_R.
\end{align}
Taking the $L^2_{x,v}$-inner product with $\w^{2l} \p^m_{\beta} \mathbb{P}^\perp G_R $, one has
\begin{align}\label{eq:weakform-re}
	\no \tfrac{1}{2} \tfrac{\d}{\d t} & \| \w^l \p^m_{\beta} \mathbb{P}^\perp G_R \|^2_{L^2_{x,v}} + \tfrac{1}{\eps^2} \int_{\mathbb{T}^3} \langle \w^{2l} \p^m_{\beta} \mathbb{L} \mathbb{P}^\perp G_R \cdot \p^m_{\beta} \mathbb{P}^\perp G_R \rangle \d x \\
	\no = & \tfrac{1}{\eps} \int_{\mathbb{T}^3} \p^m E_R \langle \p_{\beta} v \mathcal{T}_1 \cdot \w^{2l} \p^m_{\beta} \mathbb{P}^\perp G_R \rangle \d x + \tfrac{1}{\eps} \int_{\mathbb{T}^3} \langle \mathcal{T} ( v \times B_0 ) \cdot  \nabla_v \w^{l} \p^m_{\beta} \mathbb{P}^\perp G_R \cdot \w^{l} \p^m_{\beta} \mathbb{P}^\perp G_R \rangle \d x \\
	\no & - \iint_{\mathbb{T}^3 \times \mathbb{R}^3} \big( \p_t + \tfrac{1}{\eps} v \cdot \nabla_x + \tfrac{1}{\eps} \mathcal{T} ( v \times B_0 ) \cdot \nabla_v \big) \p^m_{\beta} \mathbb{P} G_R \cdot \w^{2l} \p^m_{\beta} \mathbb{P}^\perp G_R M \d v \d x \\
	& - \tfrac{1}{\eps} \sum_{|\beta_1|=1} C_{\beta}^{\beta_1} \int_{\mathbb{T}^3} \langle \p_{\beta_1} v \cdot \nabla_x \p^m_{\beta- \beta_1} G_R \cdot \w^{2l} \p^m_{\beta} \mathbb{P}^\perp G_R \rangle \d x \\
	\no & - \tfrac{1}{\eps} \int_{\mathbb{T}^3} \langle \p^m_{\beta} [ \mathcal{T} (v \times B_0) \cdot \nabla_v G_0 ] \cdot \nabla_x \p^m_{\beta- \beta_1} G_R \cdot \w^{2l} \p^m_{\beta} \mathbb{P}^\perp G_R \rangle \d x \\
	\no & - \tfrac{1}{\eps} \sum_{\stackrel{|\beta_1|=1}{m_1\le m}} C_m^{m_1} C_{\beta}^{\beta_1} \int_{\mathbb{T}^3} \langle \mathcal{T} (\p_{\beta_1} v \times \p^{m_1} B_0) \cdot \nabla_v \p^{m-m_1}_{\beta- \beta_1} G_R \cdot \nabla_x \p^m_{\beta- \beta_1}G_R \cdot \w^{2l} \p^m_{\beta} \mathbb{P}^\perp G_R \rangle \d x \\
	\no & + \int_{\mathbb{T}^3} \langle \p^m_{\beta} H_R \cdot \w^{2l} \p^m_{\beta} \mathbb{P}^\perp G_R \rangle \d x \\
	\no := & I_1 + I_2 + \cdots + I_7 \,.
\end{align}
By the inequality \eqref{Coercive-L} in Lemma \ref{Lm-L-Q}, one obtains
\begin{align}
	\tfrac{1}{\eps^2} \int_{\mathbb{T}^3} \langle \w^{2l} \p^m_{\beta} \mathbb{L} \mathbb{P}^\perp G_R \cdot \p^m_{\beta} \mathbb{P}^\perp G_R \rangle \d x
	\ge \tfrac{1}{2 \eps^2} \|\w^l \p^m_{\beta} \mathbb{P}^\perp G_R\|_{L^2_{x,v}(\nu)}
	- \tfrac{C}{\eps^2} \|\w^l \p^m \mathbb{P}^\perp G_R\|_{L^2_{x,v}(\nu)}.
\end{align}
Now we are left to estimate the terms $I_i$ with $1\le i \le 6$ in the right-hand side of the equality {\eqref{eq:weakform-re}}, while estimating the last term $I_7$ is postponed behind.

For the term $I_1$, we observe that $I_1=0$ for the case $|\beta|\ge 2$. This means we only need to consider the cases $\beta =0$ and $|\beta|=1$. In the case of $\beta=0$, we have, for $|m|+|\beta| \le N$, that
\begin{align}
	I_1 = & \tfrac{1}{\eps} \int_{\mathbb{T}^3} \p^m E_R \int_{\mathbb{R}^3} v \mathcal{T}_1 \cdot \w^{2l} \p^m \mathbb{P}^\perp G_R M\d v \d x \\
	\no = & - \tfrac{1}{\eps} \sum_{|e|=1} \int_{\mathbb{T}^3} \p^{m-e} E_R \int_{\mathbb{R}^3} v \mathcal{T}_1 \cdot \w^{2l} \p^{m+e} \mathbb{P}^\perp G_R M\d v \d x \,.
\end{align}
Then the H\"older inequality yields
\begin{align}
	I_1 \le & \tfrac{C}{\eps} \sum_{|e|=1} \int_{\mathbb{T}^3} |\p^{m-e} E_R| \left\{ \int_{\mathbb{R}^3} |v|^2 \w^{4l} M\d v \right\}^{\frac{1}{2}} \left\{ \int_{\mathbb{R}^3} |\p^{m+e} \mathbb{P}^\perp G_R|^2 M\d v \right\}^{\frac{1}{2}} \d x \\
	\no \le & \tfrac{C}{\eps} \|E_R\|_{H^{N-1}_x} \| \mathbb{P}^\perp G_R\|_{H^{N+1}_x L^2_v} \,.
\end{align}
In case $|\beta|=1$, there holds $|m| \le N-1$ (since $|m|+|\beta| \le N$). It can be deduced that
\begin{align}
	I_1 \le & \tfrac{C}{\eps} \int_{\mathbb{T}^3} |\p^m E_R| \left\{ \int_{\mathbb{R}^3} \w^{2l} M \d v \right\}^{\frac{1}{2}} \left\{ \int_{\mathbb{R}^3} \w^{2l} |\p^m_{\beta} \mathbb{P}^\perp G_R|^2 M\d v \right\}^{\frac{1}{2}} \d x \\
	\no \le & \tfrac{C}{\eps} \|E_R\|_{H^{N-1}_x} \|\w^{l} \p^m_{\beta} \mathbb{P}^\perp G_R\|_{L^2_{x,v}} \,.
\end{align}
In summary, the term $I_1$ can be bounded by
\begin{align}\label{eq:term-I1}
	I_1 \le \tfrac{C}{\eps} \|E_R\|_{H^{N-1}_x} (\|\w^{l} \p^m_{\beta} \mathbb{P}^\perp G_R\|_{L^2_{x,v}} + \| \mathbb{P}^\perp G_R\|_{H^{N+1}_x L^2_v}) \,.
\end{align}

Notice that $I_2=0$ for $l=0$. It suffices to estimate $I_2$ for the case $l >0$. Observe that $|v| \le \w(v)$ and $|\nabla_v \w^l(v)| = |\w^{l-1}(v) \frac{v}{\w(v)} | \le \w^{l-1}(v)$. From the Sobolev embedding inequality $H^2_x(\mathbb{T}^3) \hookrightarrow L^\infty_x(\mathbb{T}^3)$ and the result \eqref{Bnd-NSFM} in Lemma \ref{lemm:bnd-NSMF}, it is derived that
\begin{align}\label{eq:term-I2}
	I_2 \le & \tfrac{C}{\eps} \iint_{\mathbb{T}^3 \times \mathbb{R}^3} |B_0| |v| \w^{l-1} |\p^m_{\beta}\mathbb{P}^\perp G_R |^2 w^l M\d v \d x \\
	\no \le & C \eps \|B_0\|_{H^2_x} \cdot \tfrac{1}{\eps^2} \|\w^l \p^m_{\beta}\mathbb{P}^\perp G_R\|^2_{L^2_{x,v}} \le C \eps \lambda_0(2) \cdot \tfrac{1}{\eps^2} \|\w^l \p^m_{\beta}\mathbb{P}^\perp G_R\|^2_{L^2_{x,v}} \,.
\end{align}

We next turn to the estimate on $I_3$. Recalling the definition $\mathbb{P}G_R$ in \eqref{Projc-GR}, i.e., $\mathbb{P} G_R = \rho_R^+ (1,0)^\top + \rho_R^- (0,1)^\top + u_R \cdot v (1,1)^\top + \theta_R (\tfrac{|v|^2}{2} - \tfrac{3}{2}) (1,1)^\top$, one deduces that
\begin{align}\label{eq:term-I3-temp}
	\no I_3 \le & C (\| \p_t \p^m \rho_R^+ \|_{L^2_x} + \| \p_t \p^m \rho_R^- \|_{L^2_x} + \| \p_t \p^m u_R \|_{L^2_x} + \| \p_t \p^m \theta_R \|_{L^2_x}) \| \w^l \p^m_{\beta} \mathbb{P}^\perp G_R \|_{L^2_{x,v}} \\
	& + \tfrac{C}{\eps} ( 1 + \| B_0 \|_{H^2_x} ) \| \mathbb{P} G_R \|_{H^{N+1}_x L^2_v} \| \w^l \p^m_{\beta} \mathbb{P}^\perp G_R \|_{L^2_{x,v}} \,,
\end{align}
where the Sobolev embedding inequality $H^2_x(\mathbb{T}^3) \hookrightarrow L^\infty_x(\mathbb{T}^3)$ and $ \| \p^m \rho_R^+\|_{L^2_x} + \|\p^m \rho_R^-\|_{L^2_x} + \|\p^m u_R\|_{L^2_x} + \|\p^m \theta_R\|_{L^2_x} \sim \|\p^m \mathbb{P}G_R\|_{L^2_{x,v}}$ have been used. By the relation \eqref{MM-10} and the H\"older inequality, it follows that for $|m|\le N$,
\begin{align}\label{eq:term-I3-temp2}
	\|\p_t \p^m & \rho_R^+\|_{L^2_x} + \|\p_t \p^m \rho_R^-\|_{L^2_x} + \|\p_t \p^m u_R\|_{L^2_x} + \|\p_t \p^m \theta_R\|_{L^2_x}
	\\\no
	\le & \tfrac{C}{\eps} \|\nabla_x \p^m \mathbb{P} G_R\|_{L^2_{x,v}}
	+ \tfrac{C}{\eps} (1+\|B_0\|_{H^{N+2}_x}) \| \mathbb{P}^\perp G_R\|_{H^{N+1}_x L^2_v}
	+ C \|\p^m H_R\|_{L^2_{x,v}}.
\end{align}
Then, based on \eqref{eq:term-I3-temp}-\eqref{eq:term-I3-temp2} and \eqref{Bnd-NSFM}, $I_3$ can be bounded by
\begin{align}\label{eq:term-I3}
	I_3 \lesssim & \tfrac{1}{\eps} \| \mathbb{P}G_R\|_{H^{N+1}_x L^2_v} \|\w^l \p^m_{\beta}\mathbb{P}^\perp G_R\|_{L^2_{x,v}} \\
	\no & + \tfrac{1}{\eps} \| \mathbb{P}^\perp G_R\|_{H^{N+1}_x L^2_v} \|\w^l \p^m_{\beta}\mathbb{P}^\perp G_R\|_{L^2_{x,v}} + \|\p^m H_R\|_{L^2_{x,v}} \|\w^l \p^m_{\beta}\mathbb{P}^\perp G_R\|_{L^2_{x,v}} \,.
\end{align}

For the term $I_4$, one has
\begin{align}\label{eq:term-I4}
	I_4 = & \tfrac{1}{\eps} \sum_{|\beta_1|=1} C_{\beta}^{\beta_1} \int_{\mathbb{T}^3} \langle \p_{\beta_1} v \cdot \nabla_x \p^m_{\beta- \beta_1} (\mathbb{P} G_R + \mathbb{P}^\perp G_R) \cdot \w^{2l} \p^m_{\beta}\mathbb{P}^\perp G_R \rangle \d x
	\\
	\no \le & \tfrac{C}{\eps} \big( \| \mathbb{P} G_R\|_{H^{N+1}_x L^2_v} + \| \mathbb{P}^\perp G_R\|_{H^N_{x,v}(\w^{2l})} \big) \|\w^l \p^m_{\beta}\mathbb{P}^\perp G_R\|_{L^2_{x,v}} \,.
\end{align}

We now estimate the term $I_5$. Recall \eqref{g0-0}, i.e., $G_0 = \rho_0^+ (1,0)^\top + \rho_0^- (0,1)^\top + u_0 \cdot v (1,1)^\top +\theta_0 (\tfrac{|v|^2}{2}-\tfrac{3}{2}) (1,1)^\top $, where $- \theta_0 = \rho_0 = \tfrac{\rho_0^+ + \rho_0^-}{2}$ and $n_0 = \rho_0^+ - \rho_0^-$. From the Sobolev embedding inequality, \eqref{Bnd-NSFM} and \eqref{Poincare-BR} in Lemma \ref{Lm-Integral-Bnd}, we can infer that
\begin{align}\label{eq:term-I5}
	I_5 = & \tfrac{1}{\eps} \int_{\mathbb{T}^3} \langle \p^m_{\beta} [ ( v \times B_R ) \cdot u_0 \mathcal{T}_1 ] \cdot \w^{2l} \p^m_{\beta}\mathbb{P}^\perp G_R \rangle \d x \\
	\no \le & \tfrac{C}{\eps} \|u_0\|_{H^{N+2}_x} \|B_R\|_{H^N_x} \|\w^l \p^m_{\beta}\mathbb{P}^\perp G_R\|_{L^2_{x,v}} \le \tfrac{C}{\eps} \lambda_0^{\frac{1}{2}}(N+2) \|\nabla_x B_R\|_{H^{N-1}_x} \|\w^l \p^m_{\beta}\mathbb{P}^\perp G_R\|_{L^2_{x,v}} \,.
\end{align}
By the similar argument as that in estimating $I_4$ and $I_5$, $I_6$ can be bounded by
\begin{align}\label{eq:term-I6}
	I_6 \le \tfrac{C}{\eps} \lambda_0^{\frac{1}{2}}(N+2) \big( \| \mathbb{P} G_R\|_{H^{N+1}_x L^2_v} + \| \mathbb{P}^\perp G_R\|_{H^N_{x,v} (\w^{2l}) } \big)  \|\w^l \p^m_{\beta}\mathbb{P}^\perp G_R\|_{L^2_{x,v}} \,.
\end{align}

Combining the above estimates on $I_i$ ($1\le i \le 6$) as in \eqref{eq:term-I1}, \eqref{eq:term-I2}, \eqref{eq:term-I3}, \eqref{eq:term-I4}, \eqref{eq:term-I5}, \eqref{eq:term-I6}, and noticing the fact $\|\cdot\|_{L^2_{x,v}} \le \|\cdot\|_{L^2_{x,v}(\nu)}$, it follows that
\begin{align}
	\no \tfrac{1}{2} \tfrac{\d}{\d t} & \|\w^l \p^m_{\beta}\mathbb{P}^\perp G_R\|^2_{L^2_{x,v}} + \tfrac{1}{2 \eps^2} \|\w^l \p^m_{\beta}\mathbb{P}^\perp G_R\|^2_{L^2_{x,v}(\nu)} \\
	\no \le & \tfrac{C}{\eps} \| \p^m \mathbb{P}^\perp G_R \|^2_{L^2_{x,v}(\nu)} + \tfrac{C}{\eps} \|E_R\|_{H^{N-1}_x} \| \mathbb{P}^\perp G_R \|_{H^{N+1}_x L^2_v (\nu)} + C \|\p^m H_R\|_{L^2_{x,v}} \|\w^l \p^m_{\beta}\mathbb{P}^\perp G_R\|_{L^2_{x,v}(\nu)} \\
	\no & + \tfrac{C}{\eps} ( \| \mathbb{P} G_R\|_{H^{N+1}_x L^2_v} + \| \mathbb{P}^\perp G_R \|_{H^{N+1}_x L^2_v(\nu)}) \|\w^l \p^m_{\beta}\mathbb{P}^\perp G_R\|_{L^2_{x,v}(\nu)} + C \eps \tfrac{1}{\eps^2} \| \mathbb{P}^\perp G_R\|^2_{H^N_{x,v} (\w^{2l} \nu) } \\
	\no & + \tfrac{C}{\eps} \|\nabla_x B_R\|_{H^{N-1}_x} \|\w^l \p^m_{\beta}\mathbb{P}^\perp G_R\|_{L^2_{x,v}(\nu)} + \int_{\mathbb{T}^3} \langle \p^m_{\beta} H_R \cdot \w^{2l} \p^m_{\beta}\mathbb{P}^\perp G_R \rangle \d x \,.
\end{align}
By summing up for $|m|+|\beta|\le N$ and using the H\"older inequality, the previous inequality concludes \eqref{esm:remainder-mix}. The proof of Lemma \ref{lemm:remainder-apriori} is therefore completed.
\end{proof}

\begin{proof}[Proof of Lemma \ref{lemm:remainder-apriori2}]

The process is analogous to the proof of Lemma \ref{lemm:remainder-apriori}, with some necessary modifications to match the case $|m|+|\beta| \le N+1, \beta \neq 0$. Note that in this case, the index $m$ satisfies $|m| \le N$, which ensures it is still possible to control the higher-order derivatives with respect to the spatial variables.

We start from the formulation \eqref{eq:weakform-re}. For term $I_1$, it suffices to deal with the case $|\beta| =1$, namely,
\begin{align}
	I_1 = & \tfrac{1}{\eps} \sum_{|e|=1} \int_{\mathbb{T}^3} \p^{m-e} E_R \int_{\mathbb{R}^3} \p_{e} (\p_{\beta} v \mathcal{T}_1 \cdot \w^{2l} M) \p^{m+e}_{\beta-e} \mathbb{P}^\perp G_R \d v \d x \\
	\no \le & \tfrac{C}{\eps} \sum_{|e|=1} \int_{\mathbb{T}^3} |\p^{m-e} E_R| \left\{ \int_{\mathbb{R}^3} \w^{2l} (\w^{l-1} + \w^{l+1})^2 M\d v \right\}^{\frac{1}{2}} \left\{ \int_{\mathbb{R}^3} | \p^{m+e}_{\beta-e} \mathbb{P}^\perp G_R|^2 M\d v \right\}^{\frac{1}{2}} \d x \\
	\no \le & \tfrac{C}{\eps} \|E_R\|_{H^{N-1}_x} \| \mathbb{P}^\perp G_R\|_{H^{N+1}_x L^2_v} \,,
\end{align}
where $|\nabla_v (\w^{2l} M)| \le C(\w^{2l-1} + \w^{2l+1}) M $ has been employed.

On the other hand, the estimation on $I_3$ in \eqref{eq:term-I3} can be modified as
\begin{multline}\label{eq:term-I3-higher}
	I_3 \le \tfrac{C}{\eps} (\| \mathbb{P} G_R\|_{H^{N+1}_x L^2_v} + \| \mathbb{P}^\perp G_R\|_{H^{N+1}_x L^2_v}) \|\w^l \p^m_{\beta}\mathbb{P}^\perp G_R\|_{L^2_{x,v}(\nu)} \\
	+ C \|\p^m H_R\|_{L^2_{x,v}} \|\w^l \p^m_{\beta}\mathbb{P}^\perp G_R\|_{L^2_{x,v}(\nu)}.
\end{multline}
Furthermore, a similar argument as in proving \eqref{eq:term-I4}-\eqref{eq:term-I6} enables us to get
\begin{align}
	I_4 + I_5 + I_6 
	\le \tfrac{C}{\eps} \big( \|\nabla_x B_R\|_{H^{N-1}_x} + \| \mathbb{P} G_R\|_{H^{N+1}_x L^2_v} + \| \mathbb{P}^\perp G_R\|_{\widetilde{\mathbf{H}}^{N+1}_{x,v}(\w^{2l} \nu)}\big) \|\w^l \p^m_{\beta}\mathbb{P}^\perp G_R\|_{L^2_{x,v}(\nu)} \,.
\end{align}
From combining all above inequalities with the estimate on $I_2$ \eqref{eq:term-I2}, and summing up for $|m|+|\beta| \le N+1$ with $ \beta \neq 0$, we can conclude that the estimate \eqref{esm:remainder-mix2} holds. The proof of Lemma \ref{lemm:remainder-apriori2} is thus completed.
\end{proof}

\subsection{Estimates on \texorpdfstring{$H_R$}{H\_R}: Proof of Lemmas \ref{Lm-HR-Square}-\ref{Lm-HR-Inner-Product}}
\label{Sec:HR}

In this subsection, the goal is to derive the bounds of the quantities associated to the higher order term $H_R$ in \eqref{Remnd-H}, namely, to prove Lemmas \ref{Lm-HR-Square}-\ref{Lm-HR-Inner-Product}.

\begin{proof}[Proof of Lemma \ref{Lm-HR-Square}]
	Recalling the definition of $H_R$ in \eqref{Remnd-H} (the subscript $\eps$ is already dropped), one first decomposes $H_R$ as
	\begin{equation}\label{HR-2}
		H_R = \tfrac{1}{\eps} \Gamma_0 G_R + H_R^{(1)} + H_R^{(2)} + H_R^{(3)} \,,
	\end{equation}
	where $\Gamma_0 G_R$ is given in \eqref{Remd-Gamma0}, the symbols $H_R^{(i)}$ with $i=1,2,3$ are expressed as follows:
	\begin{align}\label{HR-3}
		\no H_R^{(1)} = & \eps \left( \begin{array}{c}
			\mathcal{Q} ({g}_2^+ , g_R^+ + g_R^-) + \Q (g_{R, \eps}^+, {g}_2^+ + {g}_2^- ) \\
			\mathcal{Q} ({g}_2^- , g_R^+ + g_R^-) + \Q (g_{R, \eps}^-, {g}_2^+ + {g}_2^- )
		\end{array} \right) + \left( \begin{array}{c}
		\mathcal{Q} (g_R^+ , g_R^+ + g_R^-) \\ \mathcal{Q} (g_R^- , g_R^+ + g_R^-) \end{array} \right) \\
		& + \left(
		\begin{array}{c}
			\mathcal{Q} ({g}_1^+ , g_R^+ + g_R^-) + \mathcal{Q} (g_R^+ , {g}_1^+ + {g}_1^-) \\
			\mathcal{Q} ({g}_1^- , g_R^+ + g_R^-) + \mathcal{Q} (g_R^- , {g}_1^+ + {g}_1^-)
		\end{array} \right) \,,
	\end{align}
	\begin{align}\label{HR-4}
		\no H_R^{(2)} = & - \eps \mathcal{T} {E}_1 \cdot \nabla_v G_{R} + \eps \mathcal{T} {E}_1 \cdot v G_{R} - \eps \mathcal{T} E_{R} \cdot \nabla_v G_{R} + \eps \mathcal{T} E_{R} \cdot v G_{R} \\
		\no & -  \mathcal{T} E_0 \cdot \nabla_v G_{R} +  \mathcal{T} E_0 \cdot v G_{R} - \mathcal{T} (v \times {B}_1) \cdot \nabla_v G_{R} - \mathcal{T} (v \times B_{R}) \cdot \nabla_v G_{R} \\
		& - \mathcal{T} E_{R} \cdot \nabla_v ( G_0 + \eps G_1 + \eps^2 G_2 ) + \mathcal{T} E_{R} \cdot v ( G_0 +  \eps G_1 + \eps^2 G_2 )  \\
		\no & - \mathcal{T} (v \times B_{R}) \cdot \nabla_v ( G_1 + \eps G_2 ) \,,
	\end{align}
	and
	\begin{equation}\label{HR-5}
		H_R^{(3)} = (\mathcal{R}^+, \mathcal{R}^-)^\top \,,
	\end{equation}
	where $\mathcal{R}^\pm$ are defined in \eqref{Remd-Terms}.

	Observe that for $|m| \leq N$,
	\begin{align*}
		|\partial^m \Gamma_0 G_R|^2 \leq & C \Big( |\partial^m \mathcal{Q} (g_0^+ , g_R^+ + g_R^-)|^2 + |\partial^m \mathcal{Q}( g_R^+, g_0^+ + g_0^- )|^2 \\
		& + |\partial^m \mathcal{Q} (g_0^- , g_R^+ + g_R^-)|^2 + |\partial^m \mathcal{Q}( g_R^-, g_0^+ + g_0^- )|^2 \Big) \,.
	\end{align*}
	Then \eqref{Bnd-Q-2} in Lemma \ref{Lm-L-Q} implies
	\begin{align}\label{HR-6}
		\no \| \partial^m \Gamma_0 G_R \|^2_{L^2_v}
		\no \leq & C \sum_{m_1 + m_2 = m} \Big( \| \partial^{m_1} g_0^+ \|^2_{L^2_v }+ \| \partial^{m_1} g_0^- \|^2_{L^2_v} + \| \w^{2 \gamma} \partial^{m_1} g_0^+ \|^2_{L^2_v} + \| \w^{2 \gamma} \partial^{m_1} g_0^- \|^2_{L^2_v} \Big) \\
		& \qquad \times \Big( \| \partial^{m_2} g_R^+ \|^2_{L^2_v }+ \| \partial^{m_2} g_R^- \|^2_{L^2_v} + \| \w^{2 \gamma} \partial^{m_2} g_R^+ \|^2_{L^2_v} + \| \w^{2 \gamma} \partial^{m_2} g_R^- \|^2_{L^2_v} \Big) \\
		\no \leq & C \sum_{m_1 + m_2 = m} \Big( \| \partial^{m_1} G_0 \|^2_{L^2_v } + \| \w^{2 \gamma} \partial^{m_1} G_0 \|^2_{L^2_v} \Big) \Big( \| \partial^{m_2} G_R \|^2_{L^2_v }+ \| \w^{2 \gamma} \partial^{m_2} G_R \|^2_{L^2_v} \Big) \,.
	\end{align}

Recalling \eqref{g0-0}, i.e., $G_0 = \rho_0^+ (1,0)^\top + \rho_0^- (0,1)^\top  u_0 \cdot v (1,1)^\top + \theta_0 (\tfrac{|v|^2}{2} - \tfrac{3}{2}) (1,1)^\top$, where $- \theta_0 = \rho_0 = \tfrac{\rho_0^+ + \rho_0^-}{2}$, $n_0 = \rho_0^+ - \rho_0^-$, one deduces that
	\begin{equation}\label{HR-7}
		\begin{aligned}
			\| \partial^{m_1} G_0 \|^2_{L^2_v } + \| \w^{2 \gamma} \partial^{m_1} G_0 \|^2_{L^2_v} \leq  C \big( \| \partial^{m_1} u_0 \|^2_{L^\infty_x} + \| \partial^{m_1} \theta_0 \|^2_{L^\infty_x} + \| \partial^{m_1} n_0 \|^2_{L^\infty_x} \big) \,.
		\end{aligned}
	\end{equation}
	Then the inequalities \eqref{HR-6} and \eqref{HR-7} reduce to
	\begin{equation}\label{HR-8}
		\begin{aligned}
			\| \partial^m \Gamma_0 G_R \|^2_{L^2_{x,v}} \leq C \sum_{m_2 \leq m} \| \partial^{m-m_2} ( u_0, \theta_0, n_0  ) \|^2_{L^\infty_x} \| \w^{2 \gamma} \partial^{m_2} G_R \|^2_{L^2_{x,v}} \,.
		\end{aligned}
	\end{equation}

By splitting $G_R = \mathbb{P} G_R + \mathbb{P}^\perp G_R$ and combining the inequality \eqref{Weight-Hydro-Control}, one knows
	\begin{equation}\label{HR-9}
		\begin{aligned}
			\| \w^{2 \gamma} \partial^{m_2} G_R \|^2_{L^2_{x,v}} \leq & C \big( \| \w^{2 \gamma} \partial^{m_2} \mathbb{P} G_R \|^2_{L^2_{x,v}} + \| \w^{2 \gamma} \partial^{m_2} \mathbb{P}^\perp G_R \|^2_{L^2_{x,v}} \big) \\
			\leq & C \big( \| \partial^{m_2} \mathbb{P} G_R \|^2_{L^2_{x,v}} + \| \w^{2 \gamma} \partial^{m_2} \mathbb{P}^\perp G_R \|^2_{L^2_{x,v}} \big) \,.
		\end{aligned}
	\end{equation}
	Together with \eqref{HR-8} and \eqref{HR-9}, there holds
	\begin{equation}\label{HR-10}
		\begin{aligned}
			\| \Gamma_0  G_R \|^2_{H^N_x L^2_v} \lesssim & \| ( u_0, \theta_0, n_0 ) \|^2_{H^{N+2}_x} \big( \| \mathbb{P} G_R \|^2_{H^N_x L^2_v} + \| \w^{2 \gamma} \partial^{m} \mathbb{P}^\perp G_R \|^2_{H^N_x L^2_v (\w^{4 \gamma} \nu)} \big) \,,
		\end{aligned}
	\end{equation}
	where we have used the Sobolev embedding $H^2_x (\T) \hookrightarrow L^\infty_x (\T)$ and Lemma \ref{lemm:bnd-NSMF}.

By the similar scheme of \eqref{HR-6}, we get
	\begin{align}\label{HR-11}
		\no \| \partial^m H_R^{(1)} \|^2_{L^2_v} \leq & C \eps^2 \sum_{m_1 \leq m} \| \w^{2 \gamma} \partial^{m-m'} G_2 \|^2_{L^2_v} \| \w^{2 \gamma} \partial^{m_1} G_R \|^2_{L^2_v} \\
		\no & + C \sum_{m_1 \leq m} \| \w^{2 \gamma} \partial^{m-m'} G_1 \|^2_{L^2_v} \| \w^{2 \gamma} \partial^{m_1} G_R \|^2_{L^2_v} \\
		& + C \sum_{m_1 \leq m} \| \w^{2 \gamma} \partial^{m-m_1} G_R \|^2_{L^2_v} \| \w^{2 \gamma} \partial^{m_1} G_R \|^2_{L^2_v} : = I\!I\!I_1 + I\!I\!I_2 + I\!I\!I_3 \,.
	\end{align}

Recalling the definition of ${g}_1^\pm$ in \eqref{g1-g2}, it is derived from the  Sobolev embedding $H^2_x(\T) \hookrightarrow L^\infty_x (\T)$ and Lemma \ref{lemm:bnd-NSMF}-\ref{lemm:bnd-linearMaxwl} that
	\begin{align}\label{HR-12}
		\no \| G_1 \|^2_{H^N_x L^2_v} \lesssim & \| ( {u}_1 , {\theta}_1, {n}_1) \|^2_{H^{N+2}_x} + \| ( u_0 , \theta_0 , n_0 , E_0) \|^2_{H^{N+3}_x}   + \| (u_0 , \theta_0, n_0 , B_0) \|^4_{H^{N+3}_x} \\
		\no \lesssim & ( 1 + \mathcal{E}^{\IN}_{0,N+2} ) \mathcal{D}_{0,N+2}(t) + ( \mathcal{E}^{\IN}_{0,N+4} + \mathcal{E}^{\IN}_{1,N+2} ) + \mathcal{E}^{\IN}_{0,N+3} ( 1 + \mathcal{E}^{\IN}_{0,N+3} ) \\
		\leq & C(\mathcal{E}^{\IN}_{0,N+4} , \mathcal{E}^{\IN}_{1,N+2}) < \infty \,,
	\end{align}
	where the fact $\mathcal{D}_{0,N+2} (t) \leq C \mathcal{E}_{0,N+3} (t)$ is utilized. Similarly, there holds
	\begin{align}\label{HR-13}
		\| G_2 \|^2_{H^N_x L^2_v} \leq C (  \mathcal{E}^{\IN}_{0,N+5} , \mathcal{E}^{\IN}_{1,N+3}) < \infty \,.
	\end{align}
	Consequently, the inequalities \eqref{HR-9}, \eqref{HR-11}, \eqref{HR-12} and \eqref{HR-13} reduce to
	\begin{equation}\label{HR-14}
		\int_{\T} ( I\!I\!I_1 + I\!I\!I_2 ) \d  x \leq C (  \mathcal{E}^{\IN}_{0,N+5} , \mathcal{E}^{\IN}_{1,N+3}) \big( \| \mathbb{P} G_R \|^2_{H^N_x L^2_v} + \| \mathbb{P}^\perp G_R \|^2_{H^N_x L^2_v (\w^{4 \gamma})} \big) \,.
	\end{equation}

Next we consider the integral $ \int_{\T}  I\!I\!I_3 \d  x $. By the Sobolev embedding theory and the H\"older inequality, it is easily derived that
	\begin{align}\label{HR-17}
		\no \int_{\T}  I\!I\!I_3 \d  x & \leq C \| G_R \|^2_{H^N_x L^2_v (\w^{4 \gamma})} \| \w^{2 \gamma} \partial^{m} G_R \|^2_{H^N_x L^2_v (\w^{4 \gamma} \nu)} \\
		& \leq C \big( \| \mathbb{P} G_R \|^2_{H^N_x L^2_v} + \| \mathbb{P}^\perp G_R \|^2_{H^N_x L^2_v (\w^{4 \gamma})} \big) \big( \| \mathbb{P} G_R \|^2_{H^N_x L^2_v} + \| \mathbb{P}^\perp G_R \|^2_{H^N_x L^2_v (\w^{4 \gamma} \nu)} \big) \,,
	\end{align}
	where we have used the inequality \eqref{Weight-Hydro-Control} in the last line. Consequently, combining with \eqref{HR-11}, \eqref{HR-14} and \eqref{HR-17}, we know
	\begin{align}\label{HR-18}
		\no \| H_R^{(1)} \|^2_{H^N_x L^2_v} \leq & C (  \mathcal{E}^{\IN}_{0,N+5} , \mathcal{E}^{\IN}_{1,N+3}) \big( \| \mathbb{P} G_R \|^2_{H^N_x L^2_v} + \| \mathbb{P}^\perp G_R \|^2_{H^N_x L^2_v ( \w^{4 \gamma} \nu)} \big) \\
		& + C \big( \| \mathbb{P} G_R \|^2_{H^N_x L^2_v} + \| \mathbb{P}^\perp G_R \|^2_{H^N_x L^2_v ( \w^{4 \gamma} \nu)} \big)^2 \,.
	\end{align}

	We then estimate the quantity $ \| H_R^{(2)} \|^2_{H^N_x L^2_v}$. By the similar argument as that in \eqref{HR-17}, one has
	\begin{align}\label{HR-20}
		\sum_{|m| \leq N} & \int_{\T} \int_{\R^3} [ - \eps \partial^m ( \mathcal{T} {E}_1 \cdot \nabla_v G_R ) ]^2 M \d v \d x \\
		\no \leq & C \eps^2 \sum_{|m| \leq N} \sum_{m_1 \leq m} \int_{\T} |\partial^{m-m_1} {E}_1|^2 \| \partial^{m_1} \nabla_v G_R \|^2_{L^2_v} \d x \leq C \eps^2 \| {E}_1 \|^2_{H^N_x} \| \nabla_v G_R \|^2_{H^N_x L^2_v} \,.
	\end{align}
	Similarly, there holds
	\begin{align}\label{HR-21}
		\sum_{|m| \leq N} \int_{\T} \int_{\R^3} [ \eps \partial^m ( \mathcal{T} {E}_1 \cdot v G_R ) ]^2 M \d v \d x \leq C \eps^2 \| {E}_1 \|^2_{H^N_x} \| G_R \|^2_{H^N_xx L^2_v (\w^2)}  \,,
	\end{align}
	\begin{equation}\label{HR-22}
		\begin{aligned}
			\sum_{|m| \leq N} \int_{\T} \int_{\R^3} & | \eps \partial^m ( - \mathcal{T} E_R \cdot \nabla_v G_R + \mathcal{T} E_R \cdot v G_R ) |^2 M \d v \d x \\
			\leq & C  \eps^2 \| E_R \|^2_{H^N_x} \big( \| \nabla_v G_R \|^2_{H^N_x L^2_v} + \| \w \partial^m G_R \|^2_{H^N_x L^2_v (\w^2)} \big) \,,
		\end{aligned}
	\end{equation}
	\begin{equation}\label{HR-23}
		\begin{aligned}
			\sum_{|m| \leq N} \int_{\T} \int_{\R^3} & | \partial^m ( - \mathcal{T} E_0 \cdot \nabla_v G_R + \mathcal{T} E_0 \cdot v G_R ) |^2 M \d v \d x \\
			\leq & C \| E_0 \|^2_{H^N_x} \big( \| \nabla_v G_R \|^2_{H^N_x L^2_v} + \| G_R \|^2_{H^N_x L^2_v (\w^2 )}  \big) \,,
		\end{aligned}
	\end{equation}
	and
		\begin{align}\label{HR-24}
			\no \sum_{|m| \leq N} & \int_{\T} \langle \big| \partial^m [ - \mathcal{T} E_R \cdot \nabla_v ( G_0 + \eps G_1 + \eps^2 G_2 ) + \mathcal{T} E_R \cdot v ( G_0 + \eps G_1 + \eps^2 G_2 ) ] \big|^2 \rangle \d x \\
			\leq & C \| E_R \|^2_{H^N_x} \big(  \| \nabla_v ( G_0 + \eps G_1 + \eps^2 G_2 ) \|^2_{H^N_x L^2_v} + \| G_0 + \eps G_1 + \eps^2 G_2 \|^2_{H^N_x L^2_v (\w^2)} \big) \,.
		\end{align}

Recalling the definitions of $G_k$ ($k = 0,1,2$), (see \eqref{g0-0} and \eqref{g1-g2},) it is derived from the Sobolev embedding theory, the H\"older inequality and Lemmas \ref{lemm:bnd-NSMF}-\ref{lemm:bnd-linearMaxwl} that
	\begin{equation}\label{HR-25}
		\begin{aligned}
			\| \nabla_v ( G_0 + \eps G_1 + \eps^2 G_2 ) \|^2_{H^N_x L^2_v} + \| G_0 + \eps G_1 + \eps^2 G_2 \|^2_{H^N_x L^2_v (\w^2)} \\
			\leq C ( 1 + \mathcal{E}^{\IN}_{0,N+3} ) \big( \mathcal{D}_{0,N+2} (t) + \mathcal{D}_{1,N+2} (t) \big)\,,
		\end{aligned}
	\end{equation}
	which means, by \eqref{HR-24}, that
		\begin{align}\label{HR-26}
			\no \sum_{|m| \leq N} \int_{\T} \langle \big| \partial^m [ - \mathcal{T} E_R \cdot \nabla_v ( G_0 + \eps G_1 + & \eps^2 G_2 ) + \mathcal{T} E_R \cdot v ( G_0 + \eps G_1 + \eps^2 G_2 ) ] \big|^2 \rangle \d x \\
			\leq & C \| E_R \|^2_{H^N_x} \big( \mathcal{D}_{0,N+2} (t) + \mathcal{D}_{1,N+2} (t) \big),
		\end{align}
	for $N \geq 4$.

Furthermore, if $N \geq 4$, from the Sobolev embedding and the H\"older inequality, it follows that
	\begin{equation}\label{HR-27}
		\begin{aligned}
			& \sum_{|m| \leq N} \int_{\T} \int_{\R^3} |\partial^m [ - \mathcal{T} (v \times B_R) \cdot \nabla_v G_R ]|^2 M \d v \d x \leq C \| B_R \|^2_{H^N_x} \| \nabla_v G_R \|^2_{H^N_x L^2_v (\w^2)} \,.
		\end{aligned}
	\end{equation}
	and additionally, by Lemma \ref{lemm:bnd-linearMaxwl}, that
	\begin{align}\label{HR-28}
		\no & \sum_{|m| \leq N} \int_{\T} \int_{\R^3} | \partial^m [ - \mathcal{T} (v \times {B}_1) \cdot \nabla_v G_R ] |^2 M \d v \d x \\
		\leq & C \| {B}_1 \|^2_{H^N_x} \| \nabla_v G_R \|^2_{H^N_x L^2_v (\w^2)} \leq C( \mathcal{E}^{\IN}_{0,N+2} + \mathcal{E}^{\IN}_{1,N+2} ) \| \nabla_v G_R \|^2_{H^N_x L^2_v (\w^2)} \,,
	\end{align}
	and
	\begin{align}\label{HR-29}
		\no & \sum_{|m| \leq N} \int_{\T} \int_{\R^3} |\partial^m [ - \mathcal{T} ( v \times B_R ) \cdot \nabla_v ( G_1 + \eps  G_2 ) ]|^2 M \d v \d x \\
		\leq & C \| B_R \|^2_{H^N_x} \| \nabla_v ( G_1 + \eps G_2 ) \|^2_{H^N_x L^2_v (\w^2)} \leq C(\mathcal{E}^{\IN}_{0,N+5} + \mathcal{E}^{\IN}_{1,N+3}) \| \nabla_x B_R \|^2_{H^{N-1}_x} \,,
	\end{align}
	where the last inequality is implied by the inequalities \eqref{Poincare-BR}, \eqref{HR-12} and \eqref{HR-13}.

	In summary, together with \eqref{Weight-Hydro-Control} and Lemmas \ref{lemm:bnd-NSMF}-\ref{lemm:bnd-linearMaxwl}, it is derived from the inequalities \eqref{HR-20}-\eqref{HR-23} and \eqref{HR-26}-\eqref{HR-29} that
	\begin{align}\label{HR-30}
		\no \| H_R^{(2)} \|^2_{H^N_x L^2_v} \leq & C \| E_R \|^2_{H^N_x} \big( \mathcal{D}_{0,2} (t) + \mathcal{D}_{1,2} (t) \big) \\
		\no & + C \| ( E_R , B_R ) \|^2_{H^N_x} \big( \| \nabla_v G_R \|^2_{H^N_x L^2_v} + \| G_R \|^2_{H^N_x L^2_v(\w^2 \nu)} \big) \\
		& + C(\mathcal{E}^{\IN}_{0,N+5} , \mathcal{E}^{\IN}_{1,N+3}) \big[ \| \nabla_x B_R \|^2_{H^{N-1}_x} +  \| \nabla_v G_R \|^2_{H^N_x L^2_v} + \| G_R \|^2_{H^N_x L^2_v(\w^2 \nu)} \big] \,.
	\end{align}

	Finally, we are ready to estimate the quantity $ \| H_R^{(3)} \|^2_{H^N_x L^2_v}$. Recall that $ H_R^{(3)} = (  \mathcal{R}^+ , \mathcal{R}^- )^\top $ is defined by \eqref{Remd-Terms}. From Lemma \ref{lemm:bnd-NSMF}, \ref{lemm:bnd-linearMaxwl}, the Sobolev embedding theory and the H\"older inequality, it follows that
	\begin{equation}\label{HR-31}
		\| H_R^{(3)} \|^2_{H^N_x L^2_v} \leq C(\mathcal{E}^{\IN}_{0,N+5} , \mathcal{E}^{\IN}_{1,N+3}) \big( \mathcal{D}_{0,N+5} (t) + \mathcal{D}_{1,N+3} (t) \big) \,.
	\end{equation}
	Thus, combining together with the relations \eqref{HR-10}, \eqref{HR-18}, \eqref{HR-30} and \eqref{HR-31}, one concludes \eqref{HR-Square}. The proof of Lemma \ref{Lm-HR-Square} is therefore completed.
\end{proof}

\begin{proof}[Proof of Lemma \ref{Lm-HR-Inner-Product}]
	For notational convenience, denote by $\mathbb{P}^\perp G_R = ( p_R^+ , p_R^- )^\top$. For $|m| + |\beta| \leq N+1$ and $\beta \neq 0$, from \eqref{Bnd-Q-1} in Lemma \ref{Lm-L-Q}, there holds
	\begin{align}\label{HRI-1}
		\no \int_{\T} & \langle \tfrac{1}{ \eps } \partial^m_\beta \Gamma_0 G_R \cdot \w^{2l} \partial^m_\beta \mathbb{P}^\perp G_R \rangle \d x \\
		\no = & \sum_{\tau = \pm} \int_{\T} \tfrac{1}{\eps} \langle ( \partial^m_\beta \mathcal{Q} (g_0^\tau , g_R^+ + g_R^-) + \partial^m_\beta \mathcal{Q} (g_R^\tau , g_0^+ + g_0^-) ) \cdot \w^{2l} \partial^m_\beta p_R^\tau \rangle \d x \\
		\no \leq & \tfrac{C}{\eps} \sum_{\substack{m_1 + m_2 = m \\ \beta_1 + \beta_2 = \beta }} \sup_{x \in \T} \| \w^l \partial^{m_1}_{\beta_1} G_0 \|_{L^2_v (\nu)} \int_{\T} \| \w^l \partial^{m_2}_{\beta_2} G_R \|_{L^2_v(\nu)} \| \w^l \partial^m_\beta \mathbb{P}^\perp G_R \|_{L^2_v (\nu)} \d x \\
		\leq & \tfrac{C}{\eps} \sum_{\substack{ m_2 \leq m \\ \beta_2 \leq \beta }} \sqrt{\mathcal{E}^{\IN}_{0,N+5}} \big( \| \partial^{m_2} \mathbb{P} G_R \|_{L^2_{x,v}} + \| \w^l \partial^{m_2}_{\beta_2} \mathbb{P}^\perp G_R \|_{L^2_{x,v}(\nu)} \big) \| \w^l \partial^{m}_{\beta} \mathbb{P}^\perp G_R \|_{L^2_{x,v}(\nu)} \\
		\no \leq & \tfrac{C}{\eps} \| \mathbb{P} G_R \|_{H^N_x L^2_v} \| \mathbb{P}^\perp G_R \|_{\widetilde{\mathbf{H}}^{N+1}_{x,v}(\w^{2l} \nu)} + C \eps \Big( \tfrac{1}{\eps^2} \| \mathbb{P}^\perp G_R \|_{H^N_{x,v} (\w^{2l} \nu)} + \tfrac{1}{\eps^2} \| \mathbb{P}^\perp G_R \|_{\widetilde{\mathbf{H}}^{N+1}_{x,v}(\w^{2l} \nu)} \Big) \,,
	\end{align}
	where the inequality \eqref{Weight-Hydro-Control} and Lemma \ref{lemm:bnd-NSMF} are also used. Similar calculations as that in \eqref{HRI-1} imply that, for $|m| + |\beta| \leq N+1$ and $\beta \neq 0$,
	\begin{align}\label{HRI-2}
		\no \int_{\T} & \langle \partial^m_\beta H_R^{(1)} \cdot \w^{2l} \partial^m_\beta \mathbb{P}^\perp G_R \rangle \d x \\
		\no \leq & C \eps^2 \Big( 1 + \| \mathbb{P} G_R \|_{H^N_x L^2_v} + \| \mathbb{P}^\perp G_R \|_{H^N_{x,v} (\w^{2l})L^2_{x,v}} + \| \mathbb{P}^\perp G_R \|_{\widetilde{\mathbf{H}}^{N+1}_{x,v}(\w^{2l} )} \Big) \\
		& \times \Big( \| \partial^m \mathbb{P} G_R \|^2_{H^N_x L^2_v} + \tfrac{1}{\eps^2} \| \mathbb{P}^\perp G_R \|^2_{H^N_{x,v} (\w^{2l} \nu)} + \tfrac{1}{\eps^2} \| \mathbb{P}^\perp G_R \|^2_{\widetilde{\mathbf{H}}^{N+1}_{x,v}(\w^{2l} \nu)} \Big) \,.
	\end{align}

We next turn to consider the quantity $ \sum\limits_{\substack{|m|+|\beta| \leq N+1 \\ \beta \neq 0}} \int_{\T} \langle \partial^m_\beta H_R^{(2)} \cdot \w^{2l} \partial^m_\beta \mathbb{P}^\perp G_R \rangle \d x $. First, one has
	\begin{align}\label{HRI-3}
		\no - \eps & \langle \partial^m_\beta ( \mathcal{T} E_R \cdot \nabla_v G_R ) \cdot \w^{2l} \partial^m_\beta \mathbb{P}^\perp G_R \rangle \\
		\no = & - \eps \langle \partial^m ( \mathcal{T} E_R \cdot \nabla_v \partial_\beta \mathbb{P} G_R ) \cdot \w^{2l} \partial^m_\beta \mathbb{P}^\perp G_R \rangle - \eps \langle  ( \mathcal{T} E_R \cdot \nabla_v \partial^m_\beta \mathbb{P}^\perp G_R ) \cdot \w^{2l} \partial^m_\beta \mathbb{P}^\perp G_R \rangle \\
		& - \eps \sum_{0 \neq m_1 \leq m} \langle  ( \mathcal{T} \partial^{m_1} E_R \cdot \nabla_v \partial^{m-m_1}_\beta \mathbb{P}^\perp G_R ) \cdot \w^{2l} \partial^m_\beta \mathbb{P}^\perp G_R \rangle : = A_1 + A_2 + A_3 \,.
	\end{align}
For the term $A_1$, it is derived from the Sobolev embedding theory, the inequality \eqref{Weight-Hydro-Control} and the definition of $\mathbb{P}$ in \eqref{Projc-P} that, for $|m| + |\beta| \leq N+1$ and $\beta \neq 0$,
	\begin{align}\label{HRI-4}
		\no \int_{\T} A_1 \d x \leq & C \eps \sum_{m_1 \leq m} \int_{\T} |\partial^{m-m_1} E_R | \, |\partial_{m_1} \theta_R| \, \| \w^l \partial^m_\beta \mathbb{P}^\perp G_R \|_{L^2_v} \d x \\
		\leq & C \eps \| E_R \|_{H^N_x} \| \mathbb{P} G_R \|_{H^N_x L^2_v} \| \mathbb{P}^\perp G_R \|^2_{\widetilde{\mathbf{H}}^{N+1}_{x,v}(\w^{2l} )} \,.
	\end{align}
As for the term $A_2$, a simple calculation implies
	\begin{align}\label{HRI-5}
		\no \int_{\T} A_2 \d x = & - \frac{\eps}{2} \int_{\T} \langle E_R \cdot \nabla_v [ \mathcal{T} \partial^m_\beta \mathbb{P}^\perp G_R \cdot \partial^m_\beta \mathbb{P}^\perp G_R ] \w^{2l} (v) \rangle \d x \\
		\leq & C \eps \| E_R \|_{L^\infty_x} \int_{\T} \left\langle |v| \w^{2l}(v) \big| 1 - \tfrac{2l}{\w^2(v)} \big| \, |\partial^m_\beta \mathbb{P}^\perp G_R|^2 \right\rangle \d x \\
		\no \leq & C \eps \| E_R \|_{H^N_x} \| \w^l \partial^m_\beta \mathbb{P}^\perp G_R \|^2_{L^2_{x,v}(\nu)} \,,
	\end{align}
	where $ \sup_{v \in \R^3} \big| 1 - \tfrac{2l}{\w^2(v)} \big| \leq 2 l + 1 $ and the Sobolev embedding $H^2_x(\T) \hookrightarrow L^\infty_x(\T)$ have been used. For the term $A_3$, using again the Sobolev embedding theory implies
	\begin{equation}\label{HRI-6}
			\int_{\T} A_3 \d x \leq C \eps \| E_R \|_{H^N_x} \| \mathbb{P}^\perp G_R \|^2_{\widetilde{\mathbf{H}}^{N+1}_{x,v}(\w^{2l} )} \,.
	\end{equation}
Collecting the inequalities \eqref{HRI-3}-\eqref{HRI-6} reduces to
	\begin{multline}\label{HRI-7}
			\int_{\T} - \eps \langle \partial^m_\beta ( \mathcal{T} E_R \cdot \nabla_v G_R ) \cdot \w^{2l} \partial^m_\beta \mathbb{P}^\perp G_R \rangle \d x \\
			\leq C \eps^2 \| E_R \|_{H^N_x} \Big( \| \mathbb{P} G_R \|^2_{H^N_x L^2_v} + \tfrac{1}{\eps^2} \| \mathbb{P}^\perp G_R \|^2_{\widetilde{\mathbf{H}}^{N+1}_{x,v}(\w^{2l} )} \Big) \,.
	\end{multline}

Concerning the term $- \mathcal{T} (v \times B_R) \cdot \nabla_v G_R$ in \eqref{HR-4}, one deduces that for $|m| + |\beta| \leq N+1$ with $\beta \neq 0$,
	\begin{align}\label{HRI-8}
		\no \int_{\T} & \langle \partial^m_\beta ( \eps \mathcal{T} E_R \cdot v G_R ) \cdot \w^{2l} \partial^m_\beta \mathbb{P}^\perp G_R \rangle \d  x \\
		= & \sum_{m_1 \leq m} \int_{\T} \langle ( \eps \mathcal{T} \partial^{m-m_1} E_R \cdot v \partial^{m_1}_\beta ( \mathbb{P} G_R + \mathbb{P}^\perp G_R ) ) \cdot \w^{2l} \partial^m_\beta \mathbb{P}^\perp G_R  ) \rangle \d x \\
		\no & + \sum_{m_1 \leq m, |\beta_1| = 1} \int_{\T} \langle ( \eps \mathcal{T} \partial^{m-m_1} E_R \cdot \partial_{\beta_1}v \partial^{m_1}_{\beta - \beta_1} ( \mathbb{P} G_R + \mathbb{P}^\perp G_R ) ) \cdot \w^{2l} \partial^m_\beta \mathbb{P}^\perp G_R  ) \rangle \d x \\
		\no \leq & C \eps^2 \| E_R \|_{H^N_x} \Big( \tfrac{1}{\eps^2} \| \mathbb{P}^\perp G_R \|^2_{\widetilde{\mathbf{H}}^{N+1}_{x,v}(\w^{2l} \nu)} + \| \mathbb{P} G_R \|^2_{H^N_x L^2_v} + \tfrac{1}{\eps^2} \| \mathbb{P}^\perp G_R \|^2_{H^N_{x,v} (\w^{2l} \nu)} \Big) \,.
	\end{align}
By the similar arguments as in \eqref{HRI-7} and \eqref{HRI-8}, it is implied that for $|m|+|\beta| \leq N+1$ with $\beta \neq 0$,
	\begin{align}\label{HRI-9}
		\int_{\T} & \langle \partial^m ( - \eps \mathcal{T} {E}_1 \cdot \nabla_v G_R + \eps \mathcal{T} {E}_1 \cdot v G_R  ) \cdot \w^{2l} \partial^m_\beta \mathbb{P}^\perp G_R \rangle \d  x \\
		\no \leq & C \eps^2 \sqrt{\mathcal{E}^{\IN}_{0,N+5}} \Big( \tfrac{1}{\eps^2} \| \mathbb{P}^\perp G_R \|^2_{\widetilde{\mathbf{H}}^{N+1}_{x,v}(\w^{2l} \nu)} + \| \mathbb{P} G_R \|^2_{H^N_x L^2_v} + \tfrac{1}{\eps^2} \| \mathbb{P}^\perp G_R \|^2_{H^N_{x,v} (\w^{2l} \nu)} \Big) \,,
	\end{align}
	and
	\begin{align}\label{HRI-10}
		\int_{\T} & \langle \partial^m ( - \eps \mathcal{T} E_0 \cdot \nabla_v G_R + \eps \mathcal{T} E_0 \cdot v G_R  ) \cdot \w^{2l} \partial^m_\beta \mathbb{P}^\perp G_R \rangle \d  x \\
		\no \leq & C \eps^2 \sqrt{\mathcal{E}^{\IN}_{0,N}} \Big( \tfrac{1}{\eps^2} \| \mathbb{P}^\perp G_R \|^2_{\widetilde{\mathbf{H}}^{N+1}_{x,v}(\w^{2l} \nu)} + \| \mathbb{P} G_R \|^2_{H^N_x L^2_v} + \tfrac{1}{\eps^2} \| \mathbb{P}^\perp G_R \|^2_{H^N_{x,v} (\w^{2l} \nu)} \Big) \,,
	\end{align}
	where Lemma \ref{lemm:bnd-linearMaxwl} and Lemma \ref{lemm:bnd-NSMF} have also been used.

	For the term $- \mathcal{T} (v \times B_R) \cdot \nabla_v G_R$ in \eqref{HR-4}, it is derived that for $|m|+|\beta| \leq N+1$ with $\beta \neq 0$,
	\begin{align}\label{HRI-11}
		\no \int_{\T} & \langle \partial^m_\beta ( - \mathcal{T} (v \times B_R) \cdot \nabla_v G_R ) \cdot \w^{2l} \partial^m_\beta \mathbb{P}^\perp G_R \rangle \d x \\
		\no = & \int_{\T} \langle \partial^m ( - \mathcal{T} (v \times B_R) \cdot \nabla_v \partial_\beta \mathbb{P} G_R ) \cdot \w^{2l} \partial^m_\beta \mathbb{P}^\perp G_R \rangle \d x \\
		\no & + \sum_{|\beta_1|=1} \int_{\T} \langle \partial^m ( - \mathcal{T} (\partial_{\beta_1} v \times B_R) \cdot \nabla_v \partial_{\beta - \beta_1} \mathbb{P} G_R ) \cdot \w^{2l} \partial^m_\beta \mathbb{P}^\perp G_R \rangle \d x \\
		\no & + \int_{\T} \langle \partial^m ( - \mathcal{T} (v \times B_R) \cdot \nabla_v \partial_\beta \mathbb{P}^\perp G_R ) \cdot \w^{2l} \partial^m_\beta \mathbb{P}^\perp G_R \rangle \d x \\
		\no & + \sum_{|\beta_1|=1} \int_{\T} \langle \partial^m ( - \mathcal{T} (\partial_{\beta_1} v \times B_R) \cdot \nabla_v \partial_{\beta - \beta_1} \mathbb{P}^\perp G_R ) \cdot \w^{2l} \partial^m_\beta \mathbb{P}^\perp G_R \rangle \d x \\
		\leq & C \eps^2 \| B_R \|_{H^N_x} \Big( \| \mathbb{P} G_R \|^2_{H^N_x L^2_v} + \tfrac{1}{\eps^2} \| \mathbb{P}^\perp G_R \|^2_{\widetilde{\mathbf{H}}^{N+1}_{x,v}(\w^{2l} \nu)} \Big) \,,
	\end{align}
	where we have used the cancellation
	\begin{align*}
		\langle  ( - & \mathcal{T} (v \times B_R) \cdot \nabla_v \partial^m_\beta \mathbb{P}^\perp G_R ) \cdot \w^{2l} \partial^m_\beta \mathbb{P}^\perp G_R \rangle \\
		= & \tfrac{1}{2} \langle [ ( v \times B_R ) \cdot ( 2l \w^{2l - 2}(v) - 1 ) v ] \mathcal{T} \partial^m_\beta \mathbb{P}^\perp G_R \cdot \partial^m_\beta \mathbb{P}^\perp G_R \rangle = 0 \,.
	\end{align*}
Similarly, one has
	\begin{align}\label{HRI-12}
		\no \sum_{\substack{|m|+|\beta| \leq N+1 , \beta \neq 0}} & \int_{\T} \langle \partial^m_\beta ( - \mathcal{T} (v \times {B}_1 ) \cdot \nabla_v G_R ) \cdot \w^{2l} \partial^m_\beta \mathbb{P}^\perp G_R \rangle \d x \\
		\leq & C \eps^2 \sqrt{\mathcal{E}^{\IN}_{0,N+2}} \Big( \| \mathbb{P} G_R \|^2_{H^N_x L^2_v} + \tfrac{1}{\eps^2} \| \mathbb{P}^\perp G_R \|^2_{\widetilde{\mathbf{H}}^{N+1}_{x,v}(\w^{2l} \nu)} \Big) \,,
	\end{align}
	and
	\begin{align}\label{HRI-13}
		\no \sum_{\substack{|m|+|\beta| \leq N+1 , \beta \neq 0}} & \int_{\T} \langle \partial^m_\beta ( - \mathcal{T} (v \times B_R ) \cdot \nabla_v ( G_1 + \eps G_2) ) \cdot \w^{2l} \partial^m_\beta \mathbb{P}^\perp G_R \rangle \d x \\
		\leq & C \eps^2  \Big( \| \nabla_x B_R \|^2_{H^{N-1}_x} + \tfrac{1}{\eps^2} \| \mathbb{P}^\perp G_R \|^2_{\widetilde{\mathbf{H}}^{N+1}_{x,v}(\w^{2l} \nu)} \Big) \,.
	\end{align}

	Recalling the definitions of $G_k$ ($k = 0,1,2$) in \eqref{g0-0} and \eqref{g1-g2} again, Lemmas \ref{lemm:bnd-NSMF}-\ref{lemm:bnd-linearMaxwl} imply that
		\begin{align}\label{HRI-14}
			\sum_{\substack{|m|+|\beta| \leq N+1 , \beta \neq 0}} & \int_{\T} \langle \partial^m_\beta [ - \mathcal{T} E_R \cdot \nabla_v (G_0 + \eps G_1 + \eps^2 G_2) \no \\
			& \qquad \qquad + \mathcal{T} E_R \cdot v (G_0 + \eps G_1 + \eps^2 G_2) ] \cdot \w^{2l} \partial^m_\beta \mathbb{P}^\perp G_R \rangle \d x \\ \no
			\leq & C \eps^2 \| E_R \|_{H^N_x} \Big( \mathcal{D}_{0,N+5} (t) +  \mathcal{D}_{1,N+3} (t) + \tfrac{1}{\eps^2} \| \mathbb{P}^\perp G_R \|^2_{\widetilde{\mathbf{H}}^{N+1}_{x,v}(\w^{2l} \nu)} \Big)  \,.
		\end{align}
which, together with the inequalities \eqref{HRI-7}-\eqref{HRI-14}, yields
		\begin{align}\label{HRI-15}
			\no \sum_{\substack{|m|+|\beta| \leq N+1 , \beta \neq 0}} & \int_{\T} \langle \partial^m_\beta H_R^{(2)} \cdot \w^{2l} \partial^m_\beta \mathbb{P}^\perp G_R \rangle \d x \\
			\no \leq & C \eps^2 ( 1 + \| ( E_R, B_R ) \|_{H^N_x} ) \Big( \mathcal{D}_{0,N+5} (t) +  \mathcal{D}_{1,N+3} (t)  + \| \mathbb{P} G_R \|^2_{H^N_x L^2_v} + \| \nabla_x B_R \|^2_{H^{N-1}_x} \\
			& + \tfrac{1}{\eps^2} \| \mathbb{P}^\perp G_R \|^2_{H^N_{x,v} (\w^{2 l} \nu)} + \tfrac{1}{\eps^2} \| \mathbb{P}^\perp G_R \|^2_{\widetilde{\mathbf{H}}^{N+1}_{x,v}(\w^{2l} \nu)} \Big) \,.
		\end{align}

As for the last term $H_R^{(3)}$ in \eqref{HR-5}, one easily derives from Lemma \ref{lemm:bnd-NSMF} and \ref{lemm:bnd-linearMaxwl} that
	\begin{equation}\label{HRI-16}
		\begin{aligned}
			\sum_{\substack{|m|+|\beta| \leq N+1 , \beta \neq 0}} & \int_{\T} \langle \partial^m_\beta H_R^{(3)} \cdot \w^{2l} \partial^m_\beta \mathbb{P}^\perp G_R \rangle \d x \\
			\leq & C \eps^2 \Big( \mathcal{D}_{0,N+5} (t) +  \mathcal{D}_{1,N+3} (t) + \tfrac{1}{\eps^2} \| \mathbb{P}^\perp G_R \|^2_{\widetilde{\mathbf{H}}^{N+1}_{x,v}(\w^{2l} \nu)} \Big) \,.
		\end{aligned}
	\end{equation}
	Consequently, combining the inequalities \eqref{HRI-2}, \eqref{HRI-15} and \eqref{HRI-16} reduces to
	\begin{align}\label{HRI-17}
		\Bigg| \sum_{\substack{|m| + |\beta| \leq N+1 , \beta \neq 0 }} \int_{\T} \langle H_R \cdot \w^{2l} \partial^m_\beta \mathbb{P}^\perp G_R \rangle \d x \Bigg| \leq \textrm{RHS of \eqref{HR-Inner-Product}} \,.
	\end{align}
Here ``RHS'' stands for ``right-hand side''. By the similar arguments as that of \eqref{HRI-17}, one can also deduce that
	\begin{align}\label{HRI-18}
		\Bigg| \sum_{|m| + |\beta| \leq N} \int_{\T} \langle H_R \cdot \w^{2l} \partial^m_\beta \mathbb{P}^\perp G_R \rangle \d x \Bigg| \leq \textrm{RHS of \eqref{HR-Inner-Product}} \,.
	\end{align}
	Observe the fact
	\begin{equation}\label{HRI-19}
		\langle \Gamma_0 G_R \cdot \mathbb{P} G_R \rangle = \langle H_R^{(1)} \cdot \mathbb{P} G_R \rangle = 0 \,,
	\end{equation}
	where $\Gamma_0 G_R$ and $H_R^{(1)}$ are given in \eqref{Remd-Gamma0} and \eqref{HR-3}, respectively. This enables us to obtain, by a similar derivation as that of \eqref{HRI-17}, that
	\begin{align}\label{HRI-20}
		\Bigg| \sum_{|m| \leq N+1} \int_{\T} \langle \partial^m H_R \cdot \partial^m G_R \rangle \d x \Bigg| \leq \textrm{RHS of \eqref{HR-Inner-Product}} \,.
	\end{align}
	Therefore, combining the inequalities \eqref{HRI-17}, \eqref{HRI-18} and \eqref{HRI-20} together concludes \eqref{HR-Inner-Product}. The proof of Lemma \ref{Lm-HR-Inner-Product} is thus finished.
\end{proof}

\appendix

\section{Estimates on \texorpdfstring{$(\rho_1, u_1, \theta_1, E_1, B_1)$}{(rho\_1, u\_1, theta\_1, E\_1, B\_1)}: Proof of Lemma \ref{lemm:bnd-linearMaxwl} }
\label{Appendix-A}

This section aims at deriving the energy bounds of the linear Maxwell equations \eqref{Od:eps(10)}-\eqref{Spec-Chos*} with initial data \eqref{IC-LM}, namely, at proving Lemma \ref{lemm:bnd-linearMaxwl}.

\begin{proof}[Proof of Lemma \ref{lemm:bnd-linearMaxwl}]

The proof will be divided into two steps.

\noindent\textbf{Step 1. Estimates on $(\rho_1, u_1, \theta_1)$.}

We first prove the bound \eqref{bnd:Spec-u} of $ u_1$. Recalling the relations \eqref{Spec-Chos*}, one sees that
\begin{equation*}
	\begin{aligned}
		u_1 (t,x) = \nabla_x \phi (t,x) \,, \ \Delta_x \phi = \partial_t \theta_0 \,, \int_{\T} \phi \d x = 0 \,.
	\end{aligned}
\end{equation*}
The standard elliptic theory implies $\|\phi(t,\cdot)\|_{H_x^{M+2}} \lesssim \|\p_t \theta_0(t,\cdot) \|_{H_x^M}$, provided that $\p_t \theta_0(t,\cdot) \in {H_x^M}$. Hence, it follows from the third equation of the NSMF system \eqref{Limit-Equ} that, for $M \ge 1$,
\begin{multline}
	\| u_1 \|_{H^{M+1}_x} = \|\nabla_x \phi\|_{H^{M+1}_x}(t) \lesssim \| \p_t \theta_0(t,\cdot) \|_{H_x^M} \\
	= \| \kappa \Delta_x \theta_0 - u_0 \cdot \nabla_x \theta_0\|_{H^M_x}
	\lesssim \|\nabla_x \theta_0\|_{H^{M+1}_x} + \|u_0\|_{H^M_x} \|\nabla_x \theta_0\|_{H^M_x} \,.
\end{multline}
Then, from the definitions of $\mathcal{E}_{0,M+1} (t)$ and $\mathcal{D}_{0,M+1} (t) $ in \eqref{E0s} and \eqref{D0s} with $s = M+1$, the inequality \eqref{bnd:Spec-u} holds.

We next derive the inequality \eqref{bnd:Spec-rho}. By \eqref{Spec-Chos*}, $\rho_1$ and $\theta_1$ satisfy
\begin{equation*}
	\begin{aligned}
		\rho_1 = \theta_1 \,, \ \Delta_x \rho_1 = \tfrac{1}{6} \Delta_x |u_0|^2 - \tfrac{1}{2} \div_x ( u_0 \cdot \nabla_x u_0 - \tfrac{1}{2} j_0 \times B_0) \,, \ \int_{\T} \rho_1 \d x = 0 \,.
	\end{aligned}
\end{equation*}
Then, together with the last third equation of the NSMF system \eqref{Limit-Equ}, the elliptic theory yields
\begin{align}
	\| \rho_1 \|_{H^{M+1}_x} = \| \theta_1 \|_{H^{M+1}_x}
	\lesssim & \| \tfrac{1}{6} \Delta_x |u_0|^2 - \tfrac{1}{2} \div_x ( u_0 \cdot \nabla_x u_0 - \tfrac{1}{2} j_0 \times B_0 ) \|_{H_x^{M-1}} \\
	\no \lesssim & \big( \| \nabla_x u_0 \|_{H_x^{M+1}} + \| \nabla_x n_0 \|_{H_x^{M+1}} + \| E_0 \|_{H_x^{M+1}} \big)
	\big( \| u_0 \|_{H_x^{M+1}} + \| B_0 \|_{H_x^{M+1}} \big) \\
	\no & + \big( \| \nabla_x u_0 \|_{H_x^{M+1}} + \| n_0 \|_{H_x^{M+1}} \big) \big( \| u_0 \|^2_{H_x^{M+1}} + \| B_0 \|^2_{H_x^{M+1}} \big) \,.
\end{align}
Thus, the inequality \eqref{bnd:Spec-rho} follows from the definition of $\mathcal{E}_{0,M+1} (t) $ and $\mathcal{D}_{0,M+1} (t) $ in \eqref{E0s} and \eqref{D0s}, respectively.

\noindent\textbf{Step 2. Estimates on $(E_1, B_1)$.}

It remains to derive the inequalities \eqref{bnd:Spec-energy-1}-\eqref{bnd:Spec-energy-2}. The key point is to seek some essential decay structures. In \eqref{Od:eps(10)}, denote by $j_1 = u_2^+ - u_2^-$.

From applying $\div_x$ to the first equation of \eqref{Od:eps(10)} and using $\div_x E_1 = n_1$, it follows that
\begin{align}\label{eq:Spec-Chos-n}
	\p_t n_1 - \tfrac{1}{2}\sigma \Delta_x n_1 + \sigma n_1
	= - \div_x j_1.
\end{align}
Moreover, from taking the time derivative on the $ B_1$-equation in \eqref{Od:eps(10)},
\begin{align*}
	\p_{tt} B_1 + \nabla_x \times (\nabla_x \times B_1) - \sigma \nabla_x \times E_1
	= \nabla_x \times j_1 \,.
\end{align*}
Noticing that $\div_x B_1 = 0$ yields $\nabla_x \times	(\nabla_x B_1) = - \Delta_x B_1$, one has
\begin{align}\label{eq:Spec-Chos-B}
	\p_{tt} B_1 - \Delta_x B_1 + \sigma \p_t B_1
	= \nabla_x \times j_1 \,,
\end{align}
which possesses a decay effect $\sigma \p_t B_1$.

Then, we turn to focus on the system \eqref{Od:eps(10)}, \eqref{eq:Spec-Chos-n} and \eqref{eq:Spec-Chos-B}. For all $|m| \le M$, applying $\p^m$ to the first and second equations of \eqref{Od:eps(10)}, and then taking $L^2_x$-inner product with $\p^m E_1$ and $\p^m B_1$, respectively, one infers that
\begin{align*}
	\tfrac{1}{2} \tfrac{\d}{\d t} \big( \|\p^m E_1 \|^2_{L^2_x} + \|\p^m B_1 \|^2_{L^2_x} \big) + \sigma \|\p^m E_1 \|^2_{L^2_x} + \tfrac{1}{2} \sigma \| \p^m n_1 \|^2_{L^2_x}
	= - \int_{\T} \p^m j_1 \cdot \p^m E_1 \d x \,,
\end{align*}
where the cancellation relation $\int_{\mathbb{T}^3} (\nabla_x \times \p^m  E_1 ) \cdot \p^m B_1 \d x = \int_{\mathbb{T}^3} ( \nabla_x \times \p^m B_1 ) \cdot \p^m E_1 \d x $ and the third equation of \eqref{Od:eps(10)}, i.e., $\div_x E_1 = n_1$ have been used. Then, summing up for $|m| \le M$ implies
\begin{align}\label{eq:Spec-temp1}
	\tfrac{1}{2} \tfrac{\d}{\d t} \big( \| E_1 \|^2_{H^M_x} + \| B_1 \|^2_{H^M_x} \big) + \sigma \| E_1 \|^2_{H^M_x} + \tfrac{1}{2} \sigma \| n_1 \|^2_{H^M_x} = - \sum_{|m| \le M} \int_{\mathbb{T}^3} \p^m j_1 \cdot \p^m E_1 \d x \,.
\end{align}
Performing the similar procedure as above to the equation \eqref{eq:Spec-Chos-n}, one has
\begin{align}\label{eq:Spec-temp2}
	\tfrac{1}{2} \tfrac{\d}{\d t} \| n_1 \|^2_{H^M_x} + \sigma \| n_1 \|^2_{H^M_x} + \tfrac{1}{2} \sigma \| \nabla_x n_1 \|^2_{H^M_x} = \sum_{|m| \le M} \int_{\mathbb{T}^3} \p^m j_1 \cdot \nabla_x \p^m n_1 \d x \,.
\end{align}
For the equation \eqref{eq:Spec-Chos-B}, from taking the operator $\p^m$, multiplying by $\p_t \p^m B_1$, integrating by parts over $x \in \mathbb{T}^3$ and summing up for $|m| \leq M$, it follows that
\begin{align}\label{eq:Spec-temp3}
	\no \tfrac{1}{2} \tfrac{\d}{\d t} \big( \| \p_t B_1 \|^2_{H^M_x} + \| \nabla_x B_1 \|^2_{H^M_x} \big) + & \sigma \| \p_t \p^m B_1 \|^2_{H^M_x} \\
	& = \sum_{|m| \le M} \int_{\mathbb{T}^3} ( \nabla_x \times \p^m j_1 ) \cdot \p_t \p^m B_1 \d x \,.
\end{align}
Furthermore, if one multiplies by $\p^m B_1$ instead of $\p_t \p^m B_1$ in the arguments of \eqref{eq:Spec-temp3}, there holds
\begin{align}\label{eq:Spec-temp4}
	\tfrac{1}{2} \tfrac{\d}{\d t} \big( \| \p_t B_1 + B_1 \|^2_{H^M_x} - \| \p_t B_1 \|^2_{H^M_x} - ( 1 - \sigma ) \| B_1 \|^2_{H^M_x} \big) - \| \p_t B_1 \|^2_{H^M_x}+ \| \nabla_x B_1 \|^2_{H^M_x} \\
	\no = \sum_{ |m| \le M} \int_{\mathbb{T}^3} \p^m j_1 \cdot ( \nabla_x \times \p^m B_1) \d x \,,
\end{align}
where the equality $( \nabla_x \times \p^m j_1) \cdot \p^m B_1 = \div_x ( \p^m j_1 \times \p^m B_1 ) + \p^m j_1 \cdot ( \nabla_x \times \p^m B_1)$ has been utilized.

Choosing a positive constant $\delta = \tfrac{1}{2} \min \{ 1,\sigma \} \in ( 0, \tfrac{1}{2} ] $, multiplying equality \eqref{eq:Spec-temp4} by $\delta$, and adding the resultant to the relations \eqref{eq:Spec-temp1}, \eqref{eq:Spec-temp2} and \eqref{eq:Spec-temp3}, one has
\begin{align}\label{bnd:I-M}
	\no \tfrac{1}{2} \tfrac{\d}{\d t} & \big( \| E_1 \|^2_{H^M_x} + \| n_1 \|^2_{H^M_x} + ( 1 - \delta + \delta \sigma ) \| B_1 \|^2_{H^M_x}
	+ \| \nabla_x B_1 \|^2_{H^M_x} \\
	\no & \hspace*{6em} + ( 1 - \delta ) \| \p_t B_1 \|^2_{H^M_x} + \delta \| \p_t B_1 + B_1 \|^2_{H^M_x} \big) \\
	& + ( \sigma - \delta ) \| \p_t B_1 \|^2_{H^M_x} + \delta \| \nabla_x B_1 \|^2_{H^M_x} + \sigma \| E_1 \|^2_{H^M_x} + \tfrac{1}{2} \sigma \| \nabla_x n_1 \|^2_{H^M_x} + \tfrac{3}{2} \sigma \| n_1 \|^2_{H^M_x} \\
	\no = & \sum_{|m| \le M} \int_{\mathbb{T}^3} [ \p^m j_1 \cdot ( \nabla_x \p^m n_1 - \p^m n_1 + \delta \nabla_x \times \p^m B_1 ) + ( \nabla_x \times \p^m j_1 ) \cdot \p_t \p^m B_1 ] \d x
	\\ \no
	\le & \tfrac{\sigma}{2} \| n_1 \|^2_{H^M_x} + \tfrac{\sigma}{4} \| \nabla_x n_1 \|^2_{H^M_x} + \tfrac{\delta}{2} \| \nabla_x B_1 \|^2_{H^M_x} + \tfrac{ \sigma - \delta}{2} \| \p_t B_1 \|^2_{H^M_x} + C \| j_1 \|^2_{H^M_x} + C \| \nabla_x \times j_1 \|^2_{H^M_x} \,,
\end{align}
where we have used the H\"older inequality in the last line.

Recalling the expression of $j_1 = u_2^+ - u_2^-$ in \eqref{Od:eps(10)} and the fact $\nabla_x \times (\nabla_x n_1) = 0$, it follows from the Sobolev embedding theory that
\begin{align}
	\no \| j_1 \|^2_{H^M_x} + \| \nabla_x \times j_1 \|^2_{H^M_x}
	\leq & \| u_1 \|^2_{H^{M+1}_x} + C \| B_0 \|^2_{H^{M+1}_x} \| u_1 \|^2_{H^{M+1}_x} + C \sum \| \Gamma_0^- \|^2_{H^{M+1}_x} \\
	\no & + C ( \| u_0 \|^2_{H^{M+1}_x} + \| \theta_0 \|^2_{H^{M+1}_x} ) ( \| n_1 \|^2_{H^M_x} + \| \nabla_x n_1 \|^2_{H^M_x} + \| \nabla_x B_1 \|^2_{H^M_x} ) \,,
\end{align}
where the Poincar\'e inequality $\| B_1 \|_{L^2_x} \le C \| \nabla_x B_1 \|_{L^2_x}$ has been used. Recalling the definition of $\Gamma_0^-$ in \eqref{Od:eps(7)}, one can deduce from the Sobolev embedding theory that
\begin{align}
	\no \sum \|\Gamma_0^-\|^2_{H^{M+1}_x}
	\lesssim & \| \nabla_x u_0 \|^2_{H^{M+2}_x} + \| \nabla_x \theta_0 \|^2_{H^{M+2}_x} + \| E_0 \|^2_{H^{M+2}_x} \\
	\no & \quad + \big( \| \nabla_x n_0 \|^2_{H^{M+2}_x} + \| \nabla_x \theta_0 \|^2_{H^{M+2}_x} + \| E_0 \|^2_{H^{M+2}_x} \big) \\
	\no & \qquad \times \big( \| n_0 \|^2_{H^{M+2}_x} + \| \theta_0 \|^2_{H^{M+2}_x} + \| u_0 \|^2_{H^{M+2}_x} + \| B_0 \|^2_{H^{M+2}_x} \big) \\
	\no & \quad + ( \| n_0 \|^2_{H^{M+2}_x} + \| \nabla_x u_0 \|^2_{H^{M+2}_x}) \cdot ( \| \theta_0 \|^2_{H^{M+2}_x} + \| u_0 \|^2_{H^{M+2}_x} + \| B_0 \|^2_{H^{M+2}_x} ) \\
	\no \lesssim & \mathcal{D}_{0,M+2} (t) + \mathcal{E}_{0,M+2} (t) \big( 1 + \mathcal{E}_{0,M+2} (t) \big) \mathcal{D}_{0,M+2} (t) \,,
\end{align}
where the Ohm's law $j_0 = n_0 u_0 + \sigma (-\tfrac{1}{2} \nabla_x n_0 + E_0 + u_0 \times B_0) $ in \eqref{Limit-Equ} has been used. Therefore,
\begin{align}\label{bnd:j-1-wd}
	\| j_1 \|^2_{H^M_x} + \| \nabla_x \times j_1 \|^2_{H^M_x}
	\lesssim & \big( 1 + \mathcal{E}_{0,M+2} (t) \big) \| u_1 \|^2_{H^{M+1}_x} + \big( 1 + \mathcal{E}_{0,M+2} (t) \big)^2 \mathcal{D}_{0,M+2} (t) \\
	\no & + \mathcal{E}_{0,M+2} (t) \big( \| n_1 \|^2_{H^M_x} + \| \nabla_x n_1 \|^2_{H^M_x} + \| \nabla_x B_1 \|^2_{H^{M}_x} \big) \,.
\end{align}
Plugging the inequality \eqref{bnd:j-1-wd} into \eqref{bnd:I-M} implies
\begin{align}\label{bnd:I-M-left}
	\no \tfrac{1}{2} \tfrac{\d}{\d t} & \big( \| E_1 \|^2_{H^M_x} + \| n_1 \|^2_{H^M_x} + (1 - \delta + \delta \sigma ) \| B_1 \|^2_{H^M_x}
	+ \| \nabla_x B_1 \|^2_{H^M_x} \\
	\no & \hspace*{6em} + ( 1 - \delta ) \| \p_t B_1 \|^2_{H^M_x} + \delta \| \p_t B_1 + B_1 \|^2_{H^M_x} \big) \\
	& + ( \sigma - \delta ) \| \p_t B_1 \|^2_{H^M_x} + \delta \| \nabla_x B_1 \|^2_{H^M_x} + \sigma \| E_1 \|^2_{H^M_x} + \tfrac{1}{2} \sigma \| \nabla_x n_1 \|^2_{H^M_x} + \tfrac{3}{2} \sigma \| n_1 \|^2_{H^M_x} \\
	\no \le & C \big( 1 + \mathcal{E}_{0,M+2} (t) \big) \| u_1 \|^2_{H^{M+1}_x} + C \big( 1 + \mathcal{E}_{0,M+2} (t) \big)^2 \mathcal{D}_{0,M+2} (t) \\
	\no & + C \mathcal{E}_{0,M+2} (t) \big( \| n_1 \|^2_{H^M_x} + \| \nabla_x n_1 \|^2_{H^M_x} + \| \nabla_x B_1 \|^2_{H^M_x} \big) \,.
\end{align}
By the monotonicity of $\mathcal{E}_{0,s}(t)$ and $\mathcal{D}_{0,s}(t)$ with respect to the index $s \ge 0$, the inequality \eqref{bnd:Spec-u} reduces to
\begin{align}
	\| u_1 \|^2_{H^{M+1}_x} (t) \le C \big( 1 + \mathcal{E}_{0, M+2} (t) \big) \mathcal{D}_{0,M+2} (t) \,,
\end{align}
which, combined with the fact that $\mathcal{E}_{0,M+2}(t) \le C \mathcal{E}_{0,M+2}^{\IN} \le C \lambda_0 (M+2)$ due to Lemma \ref{lemm:bnd-NSMF}, enables us to get
\begin{align}
	\no \tfrac{1}{2} \tfrac{\d}{\d t} &
	\big( \| E_1 \|^2_{H^M_x} + \| n_1 \|^2_{H^M_x} + ( 1 - \delta + \delta \sigma ) \| B_1 \|^2_{H^M_x} + \| \nabla_x B_1 \|^2_{H^M_x} \\
	\no & \hspace*{6em} + ( 1 - \delta ) \| \p_t B_1 \|^2_{H^M_x} + \delta \| \p_t B_1 + B_1 \|^2_{H^M_x} \big) \\
	\no & + ( \sigma - \delta ) \| \p_t B_1 \|^2_{H^M_x} + \delta \| \nabla_x B_1 \|^2_{H^M_x} + \sigma \| E_1 \|^2_{H^M_x} + \tfrac{1}{2} \sigma \| \nabla_x n_1 \|^2_{H^M_x} + \tfrac{3}{2} \sigma \| n_1 \|^2_{H^M_x} \\
	\no \le & C ( 1 + \mathcal{E}_{0,M+2}^{\IN} )^2 \mathcal{D}_{0,M+2} (t) + C \mathcal{E}_{0,M+2}^{\IN} \big( \| n_1 \|^2_{H^M_x} + \| \nabla_x n_1 \|^2_{H^M_x} + \| \nabla_x B_1 \|^2_{H^M_x} \big) \,.
\end{align}
As a consequence, by the definitions of $\mathcal{E}_{1,M}(t)$ and $\mathcal{D}_{1,M}(t)$ in \eqref{E1M} and \eqref{D1M}, respectively, it follows that
\begin{align}\label{esm:I-M}
	\tfrac{1}{2} \tfrac{\d}{\d t} \mathcal{E}_{1,M} (t) + 2 \mathcal{D}_{1,M} (t) \le C ( 1 + \mathcal{E}_{0,M+2}^{\IN} )^2 \mathcal{D}_{0,M+2} (t) + C \mathcal{E}_{0,M+2}^{\IN} \mathcal{D}_{1,M} (t) \,.
\end{align}

Applying the energy inequality \eqref{esm:E-0s} with $s=M+2$ in Lemma \ref{lemm:bnd-NSMF} yields
\begin{align}
	\tfrac{\d}{\d t}[ \mathcal{E}_{1, M} (t) + \widetilde{C}_M \mathcal{E}_{0, M+2} (t) ]  + \left[ 2 \mathcal{D}_{1, M} (t) + \mathcal{D}_{0, M+2} (t) \right] \le C \mathcal{E}_{0,M+2}^{\IN} \mathcal{D}_{1,M} (t) \,,
\end{align}
with $\widetilde {C}_M= 1+ C ( 1 + \lambda_0 (M+2))^2 \ge 1$. Therefore, by choosing a small positive constant $\lambda_1(M+2) \in ( 0, \lambda_0(M+2)] $ such that if $\mathcal{E}_{0,M+2}^{\IN} \leq \lambda_1(M+2)$, then $C \mathcal{E}_{0,M+2}^{\IN} \le 1$, the above inequality will conclude \eqref{bnd:Spec-energy-1}. The bound \eqref{bnd:Spec-energy-2} also follows immediately. The proof of Lemma \ref{lemm:bnd-linearMaxwl} is completed.
\end{proof}


\section*{Acknowledgments}

The author N. Jiang is supported by grants from the National NSFC under contract Nos. 11971360 and 11731008, and the Strategic Priority Research Program of Chinese Academy of Sciences, Grant No. XDA25010404. Y.-L. Luo is supported by the Starting Research Fund from South China University of Technology (Double First-Class Construction Project) under contract No. D6211300. T.-F. Zhang is supported by grants from the National Natural Science Foundation of China under contract No. 11701534, and No. 11871203.

\bigskip


\begin{thebibliography}{99}

\bibitem{AIM-ARMA-2015} D. Ars\'enio, S. Ibrahim and N. Masmoudi,
A derivation of the magnetohydrodynamic system from Navier-Stokes-Maxwell systems. \emph{Arch. Ration. Mech. Anal.} \textbf{216} (2015), no. 3, 767-812.

\bibitem{Arsenio-SRM} D. Ars\'enio, L. Saint-Raymond,
Compactness in kinetic transport equations and hypoellipticity. \emph{J. Funct. Anal.} \textbf{261} (2011), no. 10, 3044-3098.

\bibitem{Arsenio-SaintRaymond} D. Ars\'enio and L. Saint-Raymond, \emph{From the Vlasov-Maxwell-Boltzmann system to incompressible viscous electro-magneto-hydrodynamics.} Vol. 1. EMS Monographs in Mathematics. European Mathematical Society (EMS), Zürich, 2019.

\bibitem {BGL-CPAM1993} C. Bardos, F. Golse, and C. D. Levermore,
{Fluid Dynamic Limits of Kinetic Equations. II. Convergence Proofs
for the Boltzmann Equation,} \emph{Commun. Pure and Appl. Math.} \textbf{46}
(1993), 667-753.

\bibitem{Biskamp} D. Biskamp,
\emph{Nonlinear magnetohydrodynamics.} Cambridge Monographs on Plasma Physics, \textbf{1}. Cambridge University Press, Cambridge, 1993.

\bibitem{Davidson} P. A. Davidson,
\emph{An introduction to magnetohydrodynamics. }
Cambridge Texts in Applied Mathematics. Cambridge University Press, Cambridge, 2001.


\bibitem{DL-CPAM1989} R. DiPerna and P.-L. Lions,
Global weak solutions of Vlasov-Maxwell systems. \emph{Comm. Pure Appl. Math.} \textbf{42} (1989), no. 6, 729-757.

\bibitem{DL-Annals1989} R. DiPerna and P.-L. Lions,
On the Cauchy problem for Boltzmann equations: global existence and weak stability. \emph{Ann. of Math.} (2) \textbf{130} (1989), no. 2, 321-366.

\bibitem{DLYZ-KRM2013} R. J. Duan, S. Q. Liu, T. Yang and H. J. Zhao,
Stability of the nonrelativistic Vlasov-Maxwell-Boltzmann system for angular non-cutoff potentials. \emph{Kinet. Relat. Models} \textbf{6} (2013), no. 1, 159-204.

\bibitem{DLYZ-CMP2017} R-J. Duan, Y-J Lei, T. Yang and H.J. Zhao,
The Vlasov-Maxwell-Boltzmann system near Maxwellians in the whole space with very soft potentials. \emph{Comm. Math. Phys.} \textbf{351} (2017), no. 1, 95-153.

\bibitem{FLLZ-2018} Y. Z. Fan,  Y. J. Lei, S.Q. Liu and H.J. Zhao,
The non-cutoff Vlasov-Maxwell-Boltzmann system with weak angular singularity. \emph{Sci. China Math.} \textbf{61} (2018), no. 1, 111-136.

\bibitem{GIM2014} P. Germain, S. Ibrahim, and N. Masmoudi,
Well-posedness of the Navier-Stokes-Maxwell equations. \emph{Proc. Roy. Soc. Edinburgh Sect. A} \textbf{144} (2014), no. 1, 71-86.

\bibitem{Glassey-1996} R. T. Glassey, \emph{The Cauchy problem in kinetic theory.} Society for Industrial and Applied Mathematics (SIAM), Philadelphia, 1996.

\bibitem{G-SRM-Invent2004} F. Golse and L. Saint-Raymond,
{The Navier-Stokes limit of the Boltzmann equation for bounded
collision kernels.} \emph{Invent. Math.} \textbf{155} (2004), no. 1, 81--161.

\bibitem{G-SRM2009} F. Golse and L. Saint-Raymond,
{The Incompressible Navier-Stokes Limit of the Boltzmann
Equation for Hard Cutoff Potentials.}  \emph{J. Math. Pures Appl. (9)} \textbf{
91} (2009), no. 5, 508--552.

\bibitem{Guo-Inventions2003} Y. Guo,
The Vlasov-Maxwell-Boltzmann system near Maxwellians. \emph{Invent. Math.} \textbf{153} (2003), no. 3, 593-630.

\bibitem{Guo-CPAM-2006} Y. Guo, Boltzmann diffusive limit beyond the Navier-Stokes approximation. \emph{Comm. Pure Appl. Math.} \textbf{59} (2006), no. 5, 626-687.

\bibitem{Ibrahim-Keraani-2011-SIMA} S. Ibrahim  and S. Keraani,
Global small solutions for the Navier-Stokes-Maxwell system. \emph{SIAM J. Math. Anal.} \textbf{43} (2011), no. 5, 2275-2295.

\bibitem{Ibrahim-Yoneda-2012-JMAA}S. Ibrahim and T. Yoneda,
Local solvability and loss of smoothness of the Navier-Stokes-Maxwell equations with large initial data. \emph{J. Math. Anal. Appl.} \textbf{396} (2012), no. 2, 555-561.

\bibitem{Jang-ARMA2008} J. Jang,
Vlasov-Maxwell-Boltzmann diffusive limit. \emph{Arch. Ration. Mech. Anal.} \textbf{194} (2009), no. 2, 531-584.

\bibitem{Jang-Masmoudi-SIMA2012} J. Jang and N. Masmoudi,
Derivation of Ohm's law from the kinetic equations.  \emph{SIAM J. Math. Anal.} \textbf{44} (2012), no. 5, 3649-3669.

\bibitem{JL-CMS-2018} N. Jiang and Y.-L. Luo,
Global classical solutions to the two-fluid incompressible Navier-Stokes-Maxwell system with Ohm's law. \emph{Commun. Math. Sci.}, \textbf{16} (2018), no. 2, 561-578.

\bibitem{JL-2022annPDE} N. Jiang and Y.-L. Luo,
From Vlasov-Maxwell-Boltzmann system to two-fluid incompressible Navier-Stokes-Fourier-Maxwell system with Ohm's law: convergence for classical solutions. {\em Ann. PDE} {\bf 8} (2022), no. 1, Paper No. 4, 126 pp.

\bibitem{JM-CPAM2017} N. Jiang and N. Masmoudi,
Boundary layers and incompressible Navier-Stokes-Fourier limit of the Boltzmann equation in bounded domain I. \emph{Comm. Pure Appl. Math.} \textbf{70} (2017), no. 1, 90-171.

\bibitem{JXZ-IUM-2018} N. Jiang, C.-J. Xu and H. J. Zhao. Incompressible Navier-Stokes-Fourier limit from the Boltzmann equation: classical solutions. \emph{Indiana Univ. Math. J.}, \textbf{67} (2018), no. 5, 1817-1855.

\bibitem{LS-2010-KRM} C. D. Levermore and W. Sun. Compactness of the gain parts of the linearized Boltzmann operator with weakly cutoff kernels. \emph{Kinet. Relat. Models}, \textbf{3} (2010), no. 2, 335-351.

\bibitem{Lions-Kyoto1994} P.-L. Lions,
Compactness in Boltzmann's equation via Fourier integral operators and applications. I, II. \emph{J. Math. Kyoto Univ.} \textbf{34} (1994), no. 2, 391-427, 429-461.

\bibitem{Lions-Kyoto1994-2} P.-L. Lions,
Compactness in Boltzmann's equation via Fourier integral operators and applications. III. \emph{J. Math. Kyoto Univ.} \textbf{34} (1994), no. 3, 539-584.

\bibitem{Lions-1996} P.-L. Lions,
\emph{Mathematical topics in fluid mechanics. Vol. 1. Incompressible models.} Oxford Lecture Series in Mathematics and its Applications, 3. Oxford Science Publications. The Clarendon Press, Oxford University Press, New York, 1996.

\bibitem{Masmoudi-JMPA2010} N. Masmoudi,
Global well posedness for the Maxwell-Navier-Stokes system in 2D. \emph{J. Math. Pures Appl. (9)} \textbf{93} (2010), no. 6, 559-571.

\bibitem{Masmoudi-SRM-CPAM2003} N. Masmoudi and L. Saint-Raymond,
From the Boltzmann equation to the Stokes-Fourier system in a bounded domain.
\emph{Comm. Pure Appl. Math.} \textbf{56} (2003), no. 9, 1263-1293.

\end{thebibliography}

\end{document}